\documentclass[12pt]{article}
\usepackage{amsmath,amssymb}
\usepackage[mathscr]{euscript}
\usepackage[OT1,T1]{fontenc}
\usepackage[all]{xy}
\usepackage{hyperref}
\usepackage{cancel}
\usepackage[normalem]{ulem}

\usepackage{xcolor}

\newtheorem{prop}{Proposition}[section]
\newtheorem{rema}{Remark}[section]
\newtheorem{defi}{D{\'e}finition}[section]
\newtheorem{lemm}{Lemma}[section]
\newtheorem{theo}{Theorem}[section]
\newtheorem{coro}{Corollary}[section]

\newcommand{\iN}{\hbox{ {\leaders\hrule\hskip.2cm}{\vrule height .22cm} }}

\newcommand{\R}{\mathbb{R}}
\newcommand{\C}{\mathbb{C}}
\newcommand{\Z}{\mathbb{Z}}

\newcommand{\N}{\mathbb{N}}

\newcommand{\dR}{\hbox{d}}
\newcommand{\goG}{\mathfrak{G}}
\newcommand{\goL}{\mathfrak{L}}
\newcommand{\goP}{\mathfrak{P}}
\newcommand{\gog}{\mathfrak{g}}
\newcommand{\gos}{\mathfrak{s}}
\newcommand{\gol}{\mathfrak{l}}
\newcommand{\gop}{\mathfrak{p}}
\newcommand{\gou}{\mathfrak{u}}

\newcommand{\celg}{{\underline{\mathfrak{g}}}}
\newcommand{\cell}{{\underline{\,\mathfrak{l}\,}}}
\newcommand{\cels}{{\underline{\mathfrak{s}}}}
\newcommand{\celu}{{\underline{\mathfrak{u}}}}

\newcommand{\celp}{{\underline{\mathfrak{p}}}}

\newcommand{\celv}{\underline{V}}
\newcommand{\celw}{\underline{W}}

\title{Gauge and Gravity theories on a dynamical principal bundle}
\author{Fr{\'e}d{\'e}ric \textsc{H{\'e}lein}\footnote{Institut de Math{\'e}matiques de Jussieu,
UMR CNRS 7586 Universit{\'e} Denis Diderot --- Paris 7,
UFR de Math{\'e}matiques,  B{\^a}timent Sophie Germain
75205 Paris Cedex 13, France, \textsf{helein@math.univ-paris-diderot.fr}}
}

\begin{document}

\maketitle
\textbf{Abstract} --- \emph{In this paper we present original variational formulations of Yang-Mills, Einstein's gravitation and Kaluza-Klein theories, where, in the spirit of General Relativity, the principal bundle structure over the space-time is not fixed a priori but is dynamical. In the Yang-Mills case only a topological fibration is given a priori. In the gravity and the Kaluza-Klein theories no fibration is assumed: any critical point of the action functional defines a foliation of the manifold and the leaves make up the space-time. The latter is naturally equipped with a pseudo-Riemannian metric and, under some hypotheses, this foliation is actually a fibration. In all cases the apparition of a (at least local) principal bundle structure and a connection follows from the dynamics. Moreover the metric and the connection thus constructed are solutions of the Yang-Mills, the Einstein-Cartan or the Yang-Mills-Einstein equations, depending on the model. A crucial point is that we face the difficulty that some Lagrange multiplier fields (which are responsible for the foliation, the principal bundle structure and the connection) create unwanted terms in the equations. This difficulty is overcome by the observation that, if the structure group is compact, these terms miraculously cancel.}
\pagebreak


\section{Introduction}
A large part of theoretical Physics is based on the principle of gauge symmetry, which itself amounts to postulate the existence of principal bundles over the space-time, at a more or less formal level. However there is no fundamental rationale for explaining this postulate. This is in contrast to General Relativity, the fundamental principle of which is the equivalence principle, which results in the covariance of the theory with respect to diffeomorphisms and which do not postulate the existence of a structure as particular as that of a principal bundle. This lack of justification of the principal bundle structure is particularly evident in Kaluza-Klein theories, aiming to combine General Relativity with gauge theories: the most common hypothesis to explain the symmetry breaking at the origin of gauge fields goes back to O. Klein, it consists in assuming that the fibers of the total space above the space-time are tiny and is not completely satisfactory.  Moreover although in General Relativity the principal bundle structure may appear as non essential for pure gravity, it becomes necessary for a correct description of the fermions on a curved space-time, through the introduction of the Spin bundle.

In this paper we present alternative theories in which the principal bundle structure is not given \emph{a priori} but derives from a solution of the equations of dynamics. These theories sit on a manifold which is a candidate to get a principal bundle structure. This bundle structure will be constructed out of a dynamical field which is a 1-form with coefficients in the Lie algebra of the structure group, which could also be interpreted as a connection form on a trivial vector bundle on the manifold. Auxiliary fields are introduced in order to force integrability conditions allowing to construct a foliation which, under certain assumptions, will form a principal bundle over a quotient space, equipped with an equivariant connection. The quotient space can then be identified with a space-time manifold and the constructed fields can then be shown to be the solutions of some gauge theoretical system of equations (such as, e.g., the Yang--Mills equations) over this space-time.

However the auxiliary fields, which play the role of Lagrange multiplier for imposing non holonomic constraints, could possibly spoil the theory since they create unwanted sources in the r.h.s. of the dynamical equations. A crucial step in the study of the Euler--Lagrange equations is to prove that, under some hypotheses, these sources actually vanish. The main hypothesis in order to achieve this cancellation is to assume that the structure group is compact and simply connected. The cancellation phenomenon is then a consequence of the fact that, after a suitable trivialization of the bundle and a gauge transform, the unwanted sources are simultaneously constant on each fibers and equal to the integral of an exact form of maximal degree over the fiber, which thus cancels thanks to Stokes theorem.

For instance, in the case of the Yang--Mills theory, the Lagrangian which will be used is invariant by the group of diffeomorphisms which preserve the fibers of a submersion. This large symmetry group reduces to the standard gauge group acting on a principal bundle on classical solutions. Similarly the Lagrangian of the 4-dimensional Gravitation theory which follows is invariant by diffeomorphisms of a manifold of dimension 10 (i.e. the dimension of the Poincaré group).
Combining properties of both approaches leads to unify the gravity and the Yang--Mills fields in the spirit of Kaluza--Klein theories but without the need to assume \emph{a priori} a fibration and the equivariance of the fields along the fibration.

These various models follows the same main lines: given some Lie algebra $\gog$ of finite dimension $\hbox{dim}\gog=r$, they involve three \emph{dynamical} objects:
\begin{enumerate}
 \item a manifold $\mathcal{F}$, of dimension $N\geq r$;
 \item a 1-form $\theta^\gog$ on $\mathcal{F}$ with coefficients
 in $\gog$ and of rank $r$ everywhere;
\item an $(N-2)$-form $\pi_\gog$ on $\mathcal{F}$ with
coefficients in the dual space $\gog^*$;
\end{enumerate}
The main, naked term in the action functional is
\begin{equation}\label{Actiongenerique}
\mathcal{A}[\mathcal{F},\theta,\pi]:=
 \int_\mathcal{F}\left\langle\pi_\gog \wedge
\left(\dR\theta^\gog+\frac{1}{2}[\theta^\gog\wedge \theta^\gog]\right) \right\rangle
=
 \int_\mathcal{F}\left\langle\pi_\gog \wedge
\Theta^\gog\right\rangle,
\end{equation}
where $\langle\cdot,\cdot\rangle$ denotes the duality pairing between $\gog^*$ and $\gog$ and $\Theta^\gog:= \dR\theta^\gog+\frac{1}{2}[\theta^\gog\wedge \theta^\gog]$.

We note that the critical points of the action (\ref{Actiongenerique}) satisfy the Euler--Lagrange equations
\begin{equation}\label{systemLie3}
  \left\{\begin{array}{ccl}
         \dR\theta^\gog+\frac{1}{2}[\theta^\gog\wedge \theta^\gog] & = & 0 \\
         \dR\pi_\gog + \hbox{ad}_{\theta^\gog}^*\wedge \pi_\gog & = & 0
        \end{array}\right.
\end{equation}
The first equation (obtained by using $\pi_\gog$ as a Lagrange multiplier) is the Maurer--Cartan one. Assume that $\hbox{dim}\mathcal{F} = N=r=\hbox{dim}\gog$ and that the rank of $\theta^\gog$ is maximal everywhere. This allows, by integrating $\theta^\gog$, to construct a diffeomorphism from any neighbourhood of a point in $\mathcal{F}$ to a neighbourhood of the identity in the Lie group $\goG$, the Lie algebra of which is $\gog$. Hence $\mathcal{F}$ has locally the same structure as $\goG$. This corresponds to a local version of the \emph{Cartan--Lie theorem} asserting that a finite dimensional Lie algebra can be integrated to produce a corresponding Lie group.

Variants of this mechanism, obtained by imposing some constraints on the fields $\pi_\gog$, lead to less rigid conditions on $\theta^\gog$ and thus to identify, at least locally, $\mathcal{F}$ with a principal bundle. Indeed if we further add some extra terms in the action, then critical points $(\theta^\gog,\pi_\gog)$ correspond to solutions of gauge theoretical problems (e.g. Maxwell, Yang--Mills, Einstein--Cartan) on $\mathcal{X}$.

\subsection{Principal bundle structure starting from a submersion}
Assume now that $\mathcal{F}$ is a manifold of dimension $N = r+n$, where $n>0$, set $\gos:= \R^n$ and let $(\mathcal{X},\textbf{g})$ be a pseudo Riemannian manifold of dimension $n$. Assume that there is a submersion $\mathcal{F}\xrightarrow{P}\mathcal{X}$. We suppose that there is exists a 1-form $\underline{\beta}^\gos$ on $\mathcal{X}$ with coefficients in $\gos$, the components of which are an \emph{orthonormal} coframe on $(\mathcal{X},\textbf{g})$ and we denote by $\beta^\gos$ the pull-back by $\mathcal{F}\xrightarrow{P}\mathcal{X}$ of $\underline{\beta}^\gos$.

Consider dynamical fields which are pairs $(\theta^\gog,\pi_\gog)$, where $\theta^\gog$ is a 1-form on $\mathcal{F}$ with coefficients in the Lie algebra $\gog$ (with components $\theta^i$ in a basis) and $\pi_\gog$ is a $(N-2)$-form on $\mathcal{F}$ with coefficients in the Lie  $\gog^*$ (with components $\pi_i$). We also assume that the rank of $(\beta^\gos,\theta^\gog)$ is $N$ everywhere, so that its components $(\beta^a,\theta^i)_{1\leq a \leq n< i \leq N}$ in a basis of $\gos\oplus \gog$ provides us with a coframe on $\mathcal{F}$. This defines a volume $N$-form $\beta^{(n)}\wedge \theta^{(r)}$ on $\mathcal{F}$, where $\beta^{(n)}$ and $\theta^{(r)}$ are the exterior products of the components of, respectively, $\beta^\gos$ and $\theta^\gog$.
We then look at pairs $(\theta^\gog,\pi_\gog)$ which are critical points of $\mathscr{A}$ given by (\ref{Actiongenerique}) under the constraint that for all $1\leq a,b\leq n$ and $n<i\leq N$, the coefficient $\pi{_i}^{ab}$ such that $\beta^a\wedge \beta^b\wedge \pi_i = \pi{_i}^{ab} \beta^{(n)}\wedge \theta^{(r)}$ vanishes. Then the Euler--Lagrange equations are
\begin{equation}\label{systemLie3YM}
  \left\{\begin{array}{ccl}
         \dR\theta^\gog+\frac{1}{2}[\theta^\gog\wedge \theta^\gog] & = & \frac{1}{2}\Theta{^\gog}_{ab}\beta^a\wedge\beta^b \\
         \dR\pi_\gog + \hbox{ad}_{\theta^\gog}^*\wedge \pi_\gog & = & 0
        \end{array}\right.
\end{equation}
Here the first equation means that, if we decompose $\dR\theta^\gog+\frac{1}{2}[\theta^\gog\wedge\theta^\gog]$ by using the coframe $(\beta^a,\theta^i)_{1\leq a \leq n< i \leq N}$, the coefficients of $\beta^a\wedge\theta^j$ and $\theta^i\wedge \theta^j$ vanish.
This relation allows to identify locally each fiber of the submersion $\mathcal{F}\xrightarrow{P}\mathcal{X}$ with an open subset of $\goG$ and hence to endow $\mathcal{F}$ with a \emph{local} structure of principal bundle with structure group $\goG$ and base manifold $\mathcal{X}$. Moreover $\theta^\gog$ defines a connection on this bundle.

Now assume that, instead of assuming the constraints $\pi{_i}^{ab} = 0$ as previously, we add to the functional $\mathscr{A}$ in (\ref{Actiongenerique}) the integral $\int_\mathcal{F}\frac{1}{4}|\pi{_\gog}^{\gos\gos}|^2\beta^{(n)}\wedge \theta^{(r)}$, where $\pi{_\gog}^{\gos\gos}$ is the tensor, the components of which are $(\pi{_i}^{ab})_{i,a,b}$ and $|\pi{_\gog}^{\gos\gos}|^2$ is its norm  computed by using the pseudo Riemannian metric on $\mathcal{X}$ and an $\hbox{Ad}_\goG$-invariant metric on $\gog$. Then the Euler--Lagrange equations imply that the components of $\pi{_\gog}^{\gos\gos}$ correspond to (minus) the Hodge dual of $\dR\theta^\gog+\frac{1}{2}[\theta^\gog\wedge\theta^\gog]$. Moreover instead of the second equation in (\ref{systemLie3YM}) we get
\[
 \dR\pi_i + \hbox{ad}_{\theta^\gog}^*\wedge \pi_i = \frac{1}{2}|\pi{_\gog}^{\gos\gos}|^2\beta^{(n)}\wedge \theta^{(r-1)}_i
\]
where $\theta^{(r-1)}_i:= \frac{1}{(r-1)!} \epsilon_{ii_2\cdots i_r}\theta^{i_1}\wedge \cdots\wedge \theta^{i_r}$.
It turns out that on can deduce from this equation that the connection is a solution of a Yang--Mills equation with \emph{a priori} non vanishing sources which come from components of $\pi_\gog$ which are different from $\pi{_\gog}^{\gos\gos}$.

However a second mechanism comes into play and leads, under some general hypotheses (in particular that the group $\goG$ is compact), to the conclusion that these sources actually vanish, so that actually we obtain a solution of the Yang--Mills equation in vacuum. Thanks to this \emph{cancellation} phenomenon we obtain the following results, proved in Section \ref{sectionGauge}.\\

\noindent\textbf{Theorem \ref{statementTheoYMgeneral}} ---
\emph{Let $\gog$ be a Lie algebra of dimension $r$. Let $(\mathcal{X}, \mathbf{g})$ be a connected pseudo Riemannian manifold of dimension $n$, $\mathcal{F}$ a smooth manifold of dimension $N=n+r$ such that there exists a smooth submersion $\mathcal{F}\xrightarrow{P} \mathcal{X}$
with connected fibers. Let $(\beta^a)_{1\leq a\leq n}$ be the pull-back image by $P$ of a given orthonormal moving coframe on $(\mathcal{X}, \mathbf{g})$.}

\emph{Let $\theta^\gog$ be a 1-form on $\mathcal{F}$ with coefficient in $\gog$ of maximal rank everywhere and $\pi_\gog$ an $(N-2)$-form on $\mathcal{F}$ with coefficient in $\gog^*$. Assume that $(\theta^\gog,\pi_\gog)$ is a $\mathscr{C}^2$ critical point of}
\[
 \int_\mathcal{F}\left\langle\pi_\gog \wedge
\left(\dR\theta^\gog+\frac{1}{2}[\theta^\gog\wedge \theta^\gog]\right) \right\rangle
+ \frac{1}{4}|\pi{_\gog}^{\gos\gos}|^2\beta^{(n)}\wedge \theta^{(r)}
\] 
\emph{Assume that either, (i)
$\gog = u(1)$ and at least one fiber $P^{-1}(\{x\})$ is compact or, (ii) $\gog$ is the Lie algebra of a \textbf{compact}, \textbf{simply connected} Lie group $\widehat{\goG}$.}

\emph{Then $\theta^\gog$ endows $\mathcal{F}$ with a  principal bundle structure with a structure group $\goG$, which is either $U(1)$ in Case (i), or a quotient of $\widehat{\goG}$ by a finite subgroup in Case (ii). Moreover it defines a connection on $\mathcal{F}\xrightarrow{P} \mathcal{X}$ which is either  a solution of the Maxwell equation on $(\mathcal{X},\textbf{g})$ in Case (i), or a solution of the Yang--Mills equation in Case (ii).}\\

\noindent
Examples of compact simply connected groups are the groups $SU(k)$, for $k\geq 2$. However $U(1)$ is not simply connected.

\subsection{Principal bundle structure starting from nothing}

It is possible to dispense with the assumption that there exists a submersion from $\mathcal{F}$ to a lower dimensional manifold $\mathcal{X}$. For that purpose we assume that we replace $\beta^\gos$ (which was previously given \emph{a priori}) by $\theta^\gos$, which is now a dynamical fields. This amounts to embedd the Lie algebra $\gog$ in the larger one $\gou:= \gos\oplus \gog$, such that $[\gog,\gog]\subset \gog$ and $[\gog,\gos]\subset \gos$  and to consider $(\theta^\gou,\pi_\gou) = (\theta^\gos + \theta^\gog, \pi_\gos + \pi_\gog)$ as dynamical fields, with coefficients in, respectively, $\gou$ and $\gou^*$. We then assume that $\theta^\gou$ has a maximal rank, so that its components $(\theta^I)_{1\leq I\leq N} = (\theta^a)_{1\leq a\leq n}\cup (\theta^i)_{n<i\leq N}$ provides us with a coframe on $\mathcal{F}$. We also impose the constraint $\pi{_\gou}^{\gos\gos} = 0$, where $\pi{_\gou}^{\gos\gos}$ is the tensor with components $(\pi{_I}^{ab})_{1\leq a,b\leq n;1\leq I\leq N}$ which are
defined by $\theta^a\wedge \theta^b\wedge \pi_I = \pi{_I}^{ab}\theta^{(N)}$, where $\theta^a$ and $\theta^b$ are components of $\theta^\gos$ and $\theta^{(N)}$ is the exterior product of all components of $\theta^\gou$.
Under these assumptions a critical point of $\mathscr{A}$ satisfies the Euler--Lagrange equations
\begin{equation}\label{systemLie3Fibre}
 \left\{
 \begin{array}{ccl}
  \dR\theta^\gos + \frac{1}{2}[\theta^{\gos}\wedge \theta^{\gos}]^\gos + [\theta^{\gog}\wedge \theta^{\gos}]
  & = & \frac{1}{2}\Theta{^\gos}_{ab}\theta^a\wedge \theta^b \\
 \dR\theta^{\gog} + \frac{1}{2}[\theta^{\gos}\wedge \theta^{\gos}]^\gog + \frac{1}{2}[\theta^{\gog}\wedge \theta^{\gog}]
 & = & \frac{1}{2}\Theta{^{\gog}}_{ab}\theta^a\wedge \theta^b \\
 \dR\pi_\gou + \hbox{ad}^*_{\theta^\gou}
 \wedge \pi_\gou & = & \Psi{_\gou}^i\theta^{(N-1)}_i
 \end{array}\right.
\end{equation}
where the $\Psi{_\gou}^i$'s are coefficients in $\gou^*$, the components of which are $\Psi{_J}^i = \Theta{^L}_{JK}\pi{_L}^{iK}$.
By considering the $r$-dimensional submanifolds $\textsf{f}$ which are solutions of the exterior differential system $\theta^\gos|_\textsf{f} = 0$ we obtain a foliation of $\mathcal{F}$. This leads to endow a neighbourhood of any point of $\mathcal{F}$ with a \emph{local} principal bundle structure with structure group $\goG$ over some quotient manifold $\mathcal{X}$ of dimension $n$ (the space of leaves) and to construct a pseudo Riemannian metric and a $\gog$-value connection 1-form on $\mathcal{X}$.

\subsubsection{Kaluza--Klein theory}
Assume that the subspace $\gos\subset \gou$ is in the center of $\gou$ (it leads to simplifications in the two first equations in (\ref{systemLie3Fibre})) and fix a metric $\textsf{h}$ on $\gou$ which is invariant by the adjoint action of $\goG$ and such that $\gos\perp\gog$ . We further append to the dynamical fields $(\theta^\gou,\pi_\gou)$ a 1-form $\varphi^\gol$ with coefficients in the Lie algebra $\gol:= so(\gou,\textsf{h})$ and we add
the Palatini Lagrangian
$\int_\mathcal{F}\frac{1}{2}\theta^{(N-2)}_{IJ}\wedge \Phi^{IJ}$
to the action $\int_\mathcal{F}\left\langle\pi_\gou \wedge \Theta^\gou \right\rangle$. (Here the $\Phi^{IJ}$'s are the components of $\Phi^\gol:= \dR\varphi^\gol + \frac{1}{2}[\varphi^\gol\wedge \varphi^\gol]$ and, still, $\Theta^\gou:= \dR\theta^\gou+\frac{1}{2}[\theta^\gou\wedge \theta^\gou]$.) Then critical points of $\int_\mathcal{F}\left\langle\pi_\gou \wedge \Theta^\gou \right\rangle + \frac{1}{2}\theta^{(N-2)}_{IJ}\wedge \Phi^{IJ}$ under the constraint $\pi{_\gou}^{\gos\gos} = 0$ satisfy the system
\begin{equation}\label{systemLie3KK}
 \left\{
 \begin{array}{ccl}
  \dR\theta^\gos
  & = & \frac{1}{2}\Theta{^\gos}_{ab}\theta^a\wedge \theta^b \\
 \dR\theta^{\gog} + \frac{1}{2}[\theta^{\gog}\wedge \theta^{\gog}]
 & = & \frac{1}{2}\Theta{^{\gog}}_{ab}\theta^a\wedge \theta^b \\
 \dR\pi_\gou + \hbox{ad}^*_{\theta^\gou}
 \wedge \pi_\gou  & = & \Psi{_\gou}^i\theta^{(N-1)}_i + \frac{1}{2}\theta^{(N-3)}_{IJ\gou}\wedge \Phi^{IJ} \\
 \dR \theta^I + \varphi{^I}_J\wedge \theta^J & = & 0
 \end{array}\right.
\end{equation}
where $\Psi{_I}^i:= \Theta{^J}_{Ia}\pi{_J}^{ia}$.
By using the two first equations, if $\goG$ is compact and under mild topological hypotheses, we can construct a principal bundle structure on $\mathcal{F}$ and a pseudo Riemannian metric and a $\gog$-valued connection 1-form on the quotient manifold $\mathcal{X}$. Thanks to the last equation
we can identify $\varphi^{so}$ with the Levi-Civita connection on $\mathcal{F}$ with the metric $(\theta^\gou)^*\textsf{h}$. Thus $\frac{1}{2}\theta^{(N-3)}_{ab\gou}\wedge \Phi^{ab}$ can be interpreted as the Einstein tensor on $(\mathcal{F},(\theta^\gou)^*\textsf{h})$.
Hence the third equation means that $(\mathcal{F},(\theta^\gou)^*\textsf{h})$ is a solution of the Einstein equation with a complicated source equal to $\dR\pi_\gou + \hbox{ad}^*_{\theta^\gou} \wedge \pi_\gou  - \Psi{_\gou}^i\theta^{(N-1)}_i$.

By analyzing the latter equation (in a local trivialization) we deduce that the metric $(\theta^\gos)^*\textsf{h}$ (not $(\theta^\gou)^*\textsf{h}$ !) and the connection on $\mathcal{X}$ are solutions of an Einstein--Yang--Mills system of equations. Here again a subtle \emph{cancellation} mechanism comes into play which allow to let the sources of this system vanish.
We can hence realize the Kaluza--Klein programme without assuming any fibration \emph{a priori}, under some generic topological hypotheses. The following result is proved in Section \ref{SectionKK}. \\

\noindent\textbf{Theorem \ref{theorem0}} ---
\emph{Assume that $\widehat{\goG}$ is a \textbf{simply connected} Lie group of dimension $r$. Let $\gou = \gos \oplus \gog$ (where $\gos:= \R^n$) and let $\emph{\textsf{h}}$ a metric on $\gou$ such that $\gos\perp\gog$ and which is invariant by the adjoint action of $\widehat{\goG}$. Let $\mathcal{Y}$ be a connected oriented manifold of dimension $N=n+r$. Let  $\theta^\gou= \theta^\gos + \theta^\gog,$ be a 1-form on $\mathcal{Y}$ with coefficients in $\gou$ of rank $N$ everywhere, $\pi_\gos = \pi_\gou + \pi_\gog$ be a $(N-2)$-form on $\mathcal{Y}$ with coefficients in $\gou^*$ and $\varphi^\gol$ be a 1-form on $\mathcal{Y}$ with coefficients in $\gol = so(\gou,\emph{\textsf{h}})$. Assume that $(\theta^\gou,\pi_\gou,\varphi^\gol)$ is a critical point of class $\mathscr{C}^2$ of
\[
 \int_\mathcal{Y}\left\langle\pi_\gou \wedge
\left(\dR\theta^\gou+\frac{1}{2}[\theta^\gou\wedge \theta^\gou]\right) \right\rangle + \frac{1}{2}\theta^{(N-2)}_{IJ}\wedge \Phi^{IJ}
\]
under the constraint that $\pi_\gou\wedge \theta^a\wedge \theta^b = 0$, for any components $\theta^a$ and $\theta^b$ of $\theta^\gos$.}

\emph{Then $\mathcal{Y}$ is foliated by submanifolds $\textsf{\emph{f}}$ of codimension $n$ which are diffeomorphic to a Lie group $\goG$ which is a quotient of $\widehat{\goG}$ by a finite subgroup and on which $\widehat{\goG}$ acts. }

\emph{If furthermore $\goG$ is \textbf{compact}, then the leaves are the fibers of a 
principal bundle $\mathcal{Y}\xrightarrow{P}\mathcal{X}$ over an $n$-dimensional manifold $\mathcal{X}$ with structure group $\mathfrak{G}$.
Moreover $\theta^\gou$ encodes a pseudo Riemannian metric $\emph{\textbf{g}}$ on $\mathcal{X}$ and a $\gog$-valued connection 1-form $\emph{\textbf{A}}^\gog$ on $\mathcal{X}$, which are solutions of the Einstein--Yang--Mills system
  \[
   \left\{
 \begin{array}{ccl}
  \mathbf{R}(\mathbf{g}){^\gos}_\gos
  - \frac{1}{2}\mathbf{R}\delta{^\gos}_\gos
  + \Lambda\delta{^\gos}_\gos & = &
\frac{1}{2}\emph{\textbf{F}}{_\celg}^{\gos\cels}\emph{\textbf{F}}{^\celg}_{\gos\cels}
-\frac{1}{4} |\emph{\textbf{F}}|^2\delta{^\gos}_\gos\\
\nabla^{T\mathcal{X},\emph{\textbf{A}}}_\cels \emph{\textbf{F}}_{\gog}{^{\gos\cels}} & = & 0
 \end{array}\right.
\]
with some
cosmological constant equal to $\Lambda$.}\\

\noindent
A special case is for $\widehat{\goG} = \R$. Then, if one leaf is compact we obtain a principal bundle with structure group $U(1)$ and a solution of the Einstein--Maxwell system. However if $\widehat{\goG} = SU(k)$, then $\goG$ is necessarily compact and all conclusions of the theorem are satisfied.

\subsubsection{Gravitation on the principal bundle of frames}\label{paragraphIntroGrav}
In the two previous situations the group $\goG$ played the role of a structure group for a Yang--Mills gauge theory. For gravity theories we replace $\goG$ by a 'Lorentz' group, i.e. a group $\goL:= SO(\gos,\textsf{b})$ of isometries of some fixed Euclidean or Minkowski space $(\gos,\textsf{b})$ of dimension $n$, or its spin group $Spin(\gos,\textsf{b})$. We also introduce the 'Poincaré' group $\goP:= \goL\ltimes \gos$ and we denote by $\gol$ and $\gop = \gol \oplus \gos$ the Lie algebras of, respectively, $\goL$ and $\goP$. Then on a given manifold $\mathcal{P}$ of dimension $N:= n + \frac{n(n-1)}{2} = \hbox{dim}\goP$ we consider a pair of fields $(\varphi^\gop,\pi_\gop)$, where $\varphi^\gop$ is a 1-form of rank $N$ on $\mathcal{P}$ with coefficients in $\gop$ and $\pi_\gop$ is a $(N-2
)$-form with coefficients in $\gop^*$. Since $\varphi^\gop$ has a maximal rank, its components provide us with a coframe on $\mathcal{P}$ and by the splitting $\gop = \gol\oplus \gos$ we have the decompositions $\varphi^\gop = \varphi^\gol+\varphi^\gos$ and $\pi_\gop = \pi_\gol + \pi_\gos$.

As previously we consider the action functional $\mathscr{A}[\varphi^\gop,\pi_\gop] = \int_\mathcal{P}\left\langle\pi_\gop \wedge
\left(\dR\varphi^\gop+\frac{1}{2}[\varphi^\gop\wedge \varphi^\gop]\right) \right\rangle$ and let us first impose to $(\varphi^\gop,\pi_\gop)$ to satisfy the constraint
$\pi{_\gop}^{\gos\gos} = 0$, meaning
that, for any components $\varphi^a,\varphi^b$ of $\varphi^\gos$, $\varphi^a\wedge \varphi^b\wedge \pi_\gop = 0$. Then a critical point of $\mathscr{A}$ under these constraints satisfies exactly the system (\ref{systemLie3Fibre}), by replacing $\theta^\gos$, $\theta^\gog$, $\Theta^\gos$ and $\Theta^\gog$ by, respectively, $\varphi^\gos$, $\varphi^\gol$, $\Phi^\gos$ and $\Phi^\gol$. This
allows to locally identify $\mathfrak{P}$ with a principal bundle with structure group $\goL$ and a base manifold of dimension $n$. The fields $\varphi^\gos$ and $\varphi^\gol$ also define respectively a metric and a metric preserving connection on the quotient manifold. Hence we obtain a local structure of \emph{Cartan geometry} (see \S \ref{IntroCartanGeo} below).

Now we further add to the action $\mathscr{A}$ the 'Palatini' term $\int_\mathcal{P} \frac{1}{2}\varphi^{(N-2)}_{ab}\wedge \Phi^{ab}$, where the $\Phi^{ab}$'s are the components of $\Phi^\gop:= \dR\varphi^\gop + \frac{1}{2}[\varphi^\gop\wedge \varphi^\gop]$. Then (using the assumption $[\gog,\gos]\subset \gos$) the critical points of the total action $\int_\mathcal{P}\left\langle\pi_\gop \wedge
\Phi^\gop \right\rangle + \frac{1}{2}\varphi^{(N-2)}_{ab}\wedge \Phi^{ab}$ are solutions of
\begin{equation}\label{systemLie3Grav}
 \left\{
 \begin{array}{ccl}
  \dR\varphi^\gos + \frac{1}{2}[\varphi^{\gos}\wedge \varphi^{\gos}]^\gos + [\varphi^{\gog}\wedge \varphi^{\gos}]
  & = & \frac{1}{2}\Phi{^\gos}_{ab}\varphi^a\wedge \varphi^b \\
 \dR\varphi^{\gog} + \frac{1}{2}[\varphi^{\gos}\wedge \varphi^{\gos}]^\gog + \frac{1}{2}[\varphi^{\gog}\wedge \varphi^{\gog}]
 & = & \frac{1}{2}\Phi{^{\gog}}_{ab}\varphi^a\wedge \varphi^b \\
 \dR\pi_\gou + \hbox{ad}^*_{\varphi^\gou}
 \wedge \pi_\gou & = & \Psi{_\gou}^i\theta^{(N-1)}_i
 - \frac{1}{2}\Psi \theta^{(N-1)}_\gou
 \end{array}\right.
\end{equation}
where $\Psi:= \Psi{_a}^a = \Phi{^I}_{ab}\kappa{_I}^{ab}$. This leads to define a local Cartan geometry. Moreover the metric and the connexion on the local quotient manifold $\mathcal{X}$ are solutions of an \emph{Einstein--Cartan} system of equations.
As in the previous situations some sources (coming from the complicated structure of the third equation in (\ref{systemLie3Grav})) may  appear \emph{a priori} in these Einstein--Cartan equations (involving the Einstein tensor and the torsion). They may however vanish thanks to the cancellation phenomenon and under some assumptions.

In the following the total action $\int_\mathcal{P} \left\langle\pi_\gop \wedge
\Phi^\gop \right\rangle + \frac{1}{2}\varphi^{(N-2)}_{ab}\wedge \Phi^{ab}$ is replaced by the equivalent one $\int_\mathcal{P}\left\langle\pi_\gop \wedge
\Phi^\gop \right\rangle$ provided that, instead of the constraint $\varphi^a\wedge\varphi^b \wedge \pi_\gog = 0$, we impose that $\varphi^a\wedge\varphi^b \wedge \pi_\gos = 0$ and $\varphi^a\wedge\varphi^b \wedge \pi_\gol = \kappa{_\gol}^{ab}\varphi^{(N)}$, where the $\kappa{_\gol}^{ab}$'s are the components of a tensor $\kappa{_\gol}^{\gos\gos} \in \gol^*\otimes \gos\wedge \gos$ which encodes the canonical identification of $\gol = so(\gos,\textsf{b})$ with $\gos\wedge\gos$. This approach leads to the following, which is proved in Section \ref{sectionGravity}:\\

\noindent
\textbf{Theorem \ref{theoBigOne}} --- \emph{Let $\widehat{\goP}$ be a Lie group of dimension $N$ and $\widehat{\goL}\subset \widehat{\goP}$ a simply connected Lie subgroup of dimension $r$. Assume that their respective Lie algebras $\gop$ and $\gol$ are \emph{unimodular} and that there exists a vector subspace $\gos\subset \gop$ which is stable by $\hbox{\emph{Ad}}_\goL$ and such that $\gop = \gos\oplus \gol$ (i.e. $\widehat{\goP}/\widehat{\goL}$ is \emph{reductive}).
Let $\kappa{_\gol}^{\gos\gos}$ be a tensor in $\gop^*\otimes \gos\wedge \gos$ which is invariant by the adjoint action of $\goL$.}

\emph{Let $\varphi^\gop$ be a 1-form with coefficients in $\gop$ on $\mathcal{P}$ of rank $N$ everywhere and $\pi_\gop$ be a $(N-2)$-form  with coefficients in $\gop^*$ on $\mathcal{P}$. Assume that $(\pi_\gop,\varphi^\gop)$ is a smooth critical point of
\[
 \int_\mathcal{P}\left\langle\pi_\gop \wedge
\left(\emph{\dR}\varphi^\gop + \frac{1}{2}[\varphi^\gop\wedge\varphi^\gop]\right) \right\rangle
\]
under the constraint that $\varphi^a\wedge\varphi^b \wedge \pi_\gos = 0$ and $\varphi^a\wedge\varphi^b \wedge \pi_\gol = \kappa{_\gol}^{ab}\varphi^{(N)}$, where the $\kappa{_\gol}^{ab}$ are the components of $\kappa{_\gol}^{\gos\gos}$ in a basis of $\gos$.}

\emph{Then $\mathcal{F}$ is foliated by smooth leaves of dimension $r:= \hbox{\emph{dim}}\gol$ covered by $\widehat{\goL}$ and, in a sufficiently small open subset of $\mathcal{F}$, we can identify the set of leaves with a quotient manifold $\mathcal{X}$ of dimension $n:= \hbox{\emph{dim}}\gos$.
Moreover $\varphi^\gop$ encodes a local principal structure on the leaves and a metric and a connection on the local quotient manifold which are solutions of a \emph{generalized} Einstein--Cartan system of equations, the sources of which are total divergences \emph{on each leaf}.}\\

In the case where $\widehat{\goP}$ and $\widehat{\goL}$ are respectively the spin  Poincaré group and the spin  group and if $\kappa{_\gol}^{\gos\gos}$ encodes the canonical identification of $so(\gos,\textsf{b})$ with $\gos\wedge\gos$, then the generalized Einstein--Cartan system of equations coincides with the standard one, with sources which are total divergences.

More can be said under the additional hypothesis that the foliation is actually a fibration: the quotient manifold $\mathcal{X}$ (which represents the space-time) has then a manifold structure and the critical point produces a solution of an Einstein--Cartan system on $\mathcal{X}$ in presence of a stress-energy tensor and an angular momentum tensor. Lastly if we assume further that $\goL$ is \textbf{compact} (which is not the case if $\goL$ is the Lorentz group!) or that the fields $\pi{_\gou}$ decay at infinity, we can then conclude that the sources of the Einstein--Cartan system actually \emph{vanish}.

\subsection{Cartan geometries}\label{IntroCartanGeo}
As alluded in \S \ref{paragraphIntroGrav} a pair $(\varphi^\gop,\pi_\gop)$ which is a critical point of $\int_\mathcal{P}\left\langle\pi_\gop \wedge
\left(\dR\varphi^\gop+\frac{1}{2}[\varphi^\gop\wedge \varphi^\gop]\right) \right\rangle$ under the constraints $\varphi^a\wedge \varphi^b\wedge \pi_\gop = 0$ defines locally a structure of  \emph{Cartan geometry} on $\mathcal{P}$.

The relevance of Cartan geometry for General Relativity has been highlighted for instance in \cite{neemanregge,Wise2009}. It is based on the fact that, in the moving frame approach on General Relativity, the $so(1,3)$-valued spin connection form $\omega^\gol$ and the $\R^4$-valued soldering form $\theta^\gos$ should be understood as the two components of a single 1-form with coefficients in the Lie algebra $so(1,3)\ltimes \R^4$ of the Poincar{\'e} group (as in \cite{macdowellmansouri}). However the right geometric interpretation requires to consider all these forms as defined on the principal bundle of orthonormal frames over the space-time $\mathcal{X}$ and to understand $\textbf{A}^\gop = \omega^\gol+\theta^\gos$ as the expression of a \emph{Cartan connection} $\varphi^\gop$ in a particular trivialization of the frame bundle.

In a few words each Cartan geometry can be seen as a deformation of a rigid geometric model, called a \emph{Klein geometry}, which can be defined as a homogeneous space $\goP/\goL$, where $\goP$ is a Lie Group and $\goL$ a Lie sub-group of $\goP$. The space $\goP/\goL$ has the canonical principal bundle structure $\goL \longrightarrow \goP \longrightarrow \goP/\goL$ and $\goP$ is canonically endowed with the (left invariant) Maurer--Cartan 1-form $\eta^\gop$ with coefficients in the Lie algebra $\gop$ of $\goP$ (if $\goP$ is a matrix group, $\eta^\gop_g = g^{-1}\hbox{d}g$). A Cartan geometry is described by a principal fiber bundle $\goL \longrightarrow \mathcal{P} \longrightarrow \mathcal{X}$ and the Maurer--Cartan form $\eta^\gop$ is there replaced by a 1-form $\varphi^\gop$ defined on $\mathcal{P}$ with coefficients in the Lie algebra $\gop$ which has a maximal rank and is normalized and equivariant under the action of $\goL$ on $\mathcal{P}$. The form $\varphi^\gop$ is called a \emph{Cartan connection} and is a concept different from the well-known so-called Ehresmann connection.
The value at a point of the curvature 2-form $\dR\varphi^\gop + \frac{1}{2}[\varphi^\gop\wedge\varphi^\gop]$ measures the obstruction for $(\goL \longrightarrow \mathcal{P} \longrightarrow \mathcal{X},\varphi^\gop)$ to coincide at first order at this point with the model $(\goL \longrightarrow \goP \longrightarrow \goP/\goL,\eta^\gop)$.

The most natural situation is when $\goP = SO(n)\ltimes \R^n$ is the group of affine Euclidean isometries of the Euclidean space of dimension $n$ and $\goL=SO(n)$. Then $\goP/\goL$ is just the Euclidean space of dimension $n$ and the corresponding Cartan geometry is just another way to look at the standard Riemannian geometry. Replacing $SO(n)$ by the Lorentz group $SO(1,n-1)$ then leads to the pseudo Riemannian geometry, the framework for General Relativity. Another interesting application to General Relativity is that, by replacing the Minkowski space as a model by the de Sitter space ($\simeq SO(1,n)/SO(1,n-1)$) or the anti-de Sitter space ($\simeq SO(2,n-1)/SO(1,n-1)$), we get the Einstein equations with a positive (respectively negative) cosmological constant, as seen by S.W. MacDowell and F. Mansouri \cite{macdowellmansouri} (see \S \ref{soussectioncosmo}).
More comments on Cartan geometry are presented in \S \ref{paragraphLeviCivita} in this paper and is e.g. expounded in details in  \cite{sharpe}.
Recent accounts of its relation with General Relativity can found in \cite{Wise2009,catren}.

\subsection{A crucial point: the cancellation of the sources}\label{introcancel}

One can notice in the examples expounded in this paper that the field $\pi_\gou$ is not connected \emph{a priori} with any physically observable quantity. Indeed this field plays the role of a Lagrange multiplier for forcing the foliation and the equivariance property along the fibers. However $\pi_\gog$ has also the effect to create unwanted sources in the Euler--Lagrange equations (at least if we want to recover the standard equations of Physics or of Geometry).
A crucial step is to ensure that, under some reasonable hypotheses, these sources vanish. Here a subtle mechanism comes into play to cancel these sources, based on the facts that, on the one hand, the average of these sources on each fiber vanishes because it is the integral of a closed form and, on the other hand, these sources are constant on each fiber. However, in order to observe this cancellation, a local trivialization based on a gauge transformation is required, which requires a  delicate computation. An alternative approach have been developped by J. Pierard de Maujouy in \cite{pierard1}.

Although this mechanism works perfectly if the fibers are compact, we meet some difficulties for using it when $\goG$ is not compact: we then need to assume that the field $\pi_\gou$ and its first derivatives decay at infinity in each fiber for being able to exploit it. This is the reason why, in Theorem \ref{theoBigOne} we cannot conclude in full generality that the sources (the stress-energy and the relativistic angular momentum tensors) vanish if $\goL$ is not compact.

\subsection{Further comments}

\subsubsection{Origin of the variational formulations}
The various constructions in this paper do not come out of the blue, but have been derived first in the two papers \cite{helein14,heleinvey15} motivated by natural questions in the framework of \emph{multisymplectic geometry}. This framework generalizes the symplectic geometry in the sense that it provides a geometrical description of the Hamiltonian structure of solutions of problems in the Calculus of Variations in \emph{several variables} without depending on the choice of a particular system of coordinates (such as, for instance, a time coordinates for evolution problems). The Yang--Mills and the gravitation formulations were obtained, first, by lifting in an equivariant way the standard Lagrangian formulation of these theories on the principal bundle (see \S \ref{relationshipYM} and \S \ref{paragraphLiftingdebase}) and, second, by performing a Legendre transformation (in the multisymplectic context) by taking into account the equivariance of the connexion.
The extra field $\pi_\gog$ appears then naturally as the (multi)momentum variable conjugate to the gauge field and the constraints on $\pi{_\gou}^{\gos\gos}$ are consequences of the equivariance of the connexion (and thus reflects the gauge invariance of the initial problem). Hence the action in  (\ref{Actiongenerique}) may be viewed as the analogue for gauge theories of the integral $\int p_\mu dq^\mu$ in Mechanics. It is important to notice that the interpretation of $\pi_\gog$ as a (multi)momentum variable was a reliable indication of its relevance and importance.

The Kaluza--Klein formulation was constructed afterwards in \cite{helein2020} by combining ingredients from both theories.

\subsubsection{Perspectives}

\noindent
\textbf{Kaluza--Klein theory} ---
The Kaluza--Klein theory has a long history, starting from the work of T. Kaluza \cite{kaluza} in 1921 and O. Klein \cite{klein} in 1926, for the structure group $\R$ or $U(1)$. Some inconsistency were observed and fixed through the introduction of an additional fields (\emph{radion} or \emph{dilaton}) independently by P. Jordan \cite{jordan} in 1947 and Y. Thiry \cite{thiry} in 1948. The addition of this field may be avoided by renouncing to impose the Einstein equation on the total space of the bundle and instead by looking for the critical points of the Einstein--Hilbert action on the fiber bundle under the equivariance constraint. By following this alternative option the theory was extended to Yang--Mills fields by R. Kerner \cite{kerner} in 1968, leading to the Einstein--Yang--Mills system. Our theory is connected with the latter approach.

The most commonly used explanation for the fact that the universe we observe is 4-dimensional is basically due to Klein and relies on the hypothesis that the extra dimension is tiny and hence impossible to observe at our scale (this is reinforced at the quantum level by Heisenberg's uncertainty principle).
Our formulation does not need this assumption.\\

\noindent
\textbf{Gravity theory} ---
A physical motivation behind our gravity theory in Section \ref{sectionGravity} is to build a framework for relativistic theories which is not restricted to the set of events in space-time, but which also includes all possible frames of reference at each events. This idea was proposed for quantum field theory by F. Lur{\c c}at in 1964 \cite{lurcat}. Later on it was implemented for gravity theories by M. Toller \cite{toller} and, independently, by Y. Ne'eman and T. Regge \cite{neemanregge} in 1978. The latter work (which used ideas related to Cartan geometry) was motivated by supergravity theories and was followed by a series of papers \cite{ddfr1980,df1980,df1982}. These papers
proposed variational formulations for producing dynamically principal bundle structures (called there \emph{group manifolds}) and solutions of the Einstein--Cartan system of equations.
However they differ from our approach since their action functionals involve an integral over an $n$-dimensional section of the principal bundle (where $n$ is the dimension of space-time) and, as Ne'eman and Regge noted in \cite{neemanregge}, \S 5, no way to `\emph{extend the integration to the entire group space}' was known at that time. Under the hypothesis that the cancellation phenomenon holds (see below and  \S \ref{introcancel}) our result Theorem \ref{theoBigOne} answers positively to this question.

Our method is based on the introduction of Lagrange multiplier fields $\pi_\gop$ and most of the results in this paper involve the cancellation phenomenon (see \S \ref{introcancel}) in order to remove these fields from physical observation. However this cancellation phenomenon might not take place in gravitational theories because the Lorentz group is not compact. If so this would lead to modify the physics thus modelled, by adding new matter fields. The question of analyzing such possibilities and their possible physical relevance is quite difficult, due to the complexity of the equations. This is why we endeavored to derive the complete equations \ref{newugly} and its consequence (\ref{conservationdesdeuxtenseurs}) in a geometrical language.

Our point of view shares similarities with the interesting recent work by S. Gielen an D.K.  Wise \cite{gielen-wise}. Here the fundamental geometrical framework is the bundle of unit time-like vectors on the space-time manifold. A variational formulation of gravity is also proposed. The Authors remark also that the latter fields may also create unwanted sources to the equations.

The models proposed here do not include fermions, i.e. Dirac fields. It is however an essential question to incorporate them in, e.g., a gravity theory. It is also natural in our framework by choosing the structure group for the principal bundle to be the spin group. This question is addressed in \cite{pierard2}.

Lastly this paper addresses only \emph{classical} solutions of our models and shows that they do not differ from standard classical solutions under mild assumptions.
However it is possible that their quantification leads to different physical phenomena.


\subsection{Content of this paper}

Many results presented here were partially proved or sketched in  \cite{helein14,heleinvey15,helein2020,helein22}. However we have endeavored to simplify the computations of the Euler--Lagrange equations which were relatively tedious and to give more precise informations about these equations and their structure through the introduction of a general framework. In this process we developped a more general approach, leading to some generalizations and improvements. In particular we present the first complete and rigorous proof of the existence of a fibration in our Yang--Mills and Kaluza--Klein models.

Section \ref{Section2} is mainly pedagogical and is devoted to recall the relationship between the standard geometry of connections and metric viewed on the manifold and its lift to a principal bundle. We also discuss Cartan geometry and about the Palatini functional.

Section \ref{paragraphConventions} expounds notations and conventions which are used afterwards. Some useful technical lemmas are also stated and proven.

Section \ref{sectionGauge} is devoted to the pure Yang--Mills theory.
For pedagogical reasons we start by proving first Theorem \ref{statementTheoYMgeneral} for Maxwell fields, i.e. for $\goG = U(1)$, on the flat Minkowski space. This result is new and its proof allows to understand the cancellation phenomenon in a simple context (although some arguments are different from the case where $\goG$ is compact simply connected).
We prove afterwards Theorem \ref{statementTheoYMgeneral} on Yang--Mills fields. This result generalizes the one in \cite{helein14} since it allows more general hypotheses, for, in \cite{helein14}, we made the assumption that the 1-form $\theta^\gog$ is \emph{normalized}.

Section \ref{SectionKK} is devoted to the proof of
Theorem \ref{theorem0} on Kaluza--Klein models. This result was proved in \cite{helein2020}. Here we reproduce most of the computations of this paper in a, hopefully, more transparent and direct language and derive the complete system of equations, including some of these which were hidden in \cite{helein2020}, and  complete proof. Moreover we incorporate a cosmological term in the action. 

Section \ref{sectionGravity} contains the proof of
Theorem \ref{theoBigOne}
on gravitation, a result which extends to a larger class of groups $(\goP,\goL)$ the result in \cite{heleinvey15} which was specialized to the case where $\goP = SO(1,3)\ltimes \R^4$ is the Poincar{\'e} group and $\goL = SO(1,3)$ is the Lorentz group (or their spin covers $Spin(1,3)\ltimes \R^4$
and $Spin(1,3)$). We give applications of these results to the case where $\goP$ is $SO(1,n)$, $SO(1,n-1)\ltimes \R^n$ or $SO(2,n-1)$ and $\goL=SO(n-1)$. For $n=4$, we also show that one can deform the standard gravity by introducing the Barbero--Immirzi parameter, through different choices of the tensor $\kappa{_\gop}^{\gos\gos}$.

\section{Generalities on connections}\label{Section2}

\subsection{Connections in gauge theories and Ehresmann connections}\label{sectiongeneraliteYM}
Assume that $\mathcal{X}$ is an $n$-dimensional manifold and that $\goG$ is a finite dimensional Lie group. Let's denote by
$\gog$ its Lie algebra.
\subsubsection{In the physics literature}
A gauge field on a manifold $\mathcal{X}$ is described by a
1-form $\textbf{A}^\gog$ on $\mathcal{X}$ with coefficient in $\gog$,
i.e. $\textbf{A}^\gog\in \gog\otimes\Omega^1(\mathcal{X})$. Note that this means implicitely that the associated principal bundle is trivial. Using local coordinates
$x^\mu$ on $\mathcal{X}$, one can decompose $\textbf{A}^\gog = \textbf{A}{^\gog}_\mu dx^\mu$ (where the summation over $\mu$ is assumed), and each $\textbf{A}{^\gog}_\mu$ is a
$\gog$-valued function on $\mathcal{X}$. Its curvature is~:
\[
 \textbf{F}{^\gog}:= \dR \textbf{A}{^\gog} + \frac{1}{2}[\textbf{A}{^\gog}\wedge \textbf{A}{^\gog}] = \frac{1}{2}\left(\frac{\partial \textbf{A}{^\gog}_\nu}{\partial x^\mu}
 - \frac{\partial \textbf{A}{^\gog}_\mu}{\partial x^\nu} + [\textbf{A}{^\gog}_\mu,\textbf{A}{^\gog}_\nu]\right)dx^\mu\wedge dx^\nu.
\]
But since in all physically relevant cases $\goG$ can be represented by matrices, we can also
write $\textbf{F}{^\gog} = \dR \textbf{A}{^\gog} +\textbf{A}{^\gog}\wedge \textbf{A}{^\gog}$.

Let us fix some Riemannian or pseudo-Riemannian metric $\textsf{g}$ on $\mathcal{X}$
and an $\hbox{ad}_\goG$-invariant metric $\textsf{k}$ on $\gog$.
Then the Yang--Mills action is defined on the set of $\gog$-valued forms on $\mathcal{X}$
by
\[
 \mathcal{YM}[\textbf{A}{^\gog}]:= -\frac{1}{4}\int_\mathcal{X} |\textbf{F}{^\gog}|^2\dR\hbox{vol},
\]
where $\dR\hbox{vol}$ is the Riemannian volume form on $\mathcal{X}$ and $|\textbf{F}{^\gog}|$ is the Hilbert--Schmidt norm of $\textbf{F}{^\gog}$ computed using
$\textsf{g}$ and $\textsf{k}$. It is well-kown that $\mathcal{YM}$ is invariant by gauge
transformations:
\[
 \left\{\begin{array}{ccl}
         \textbf{A}{^\gog} & \longmapsto & g^{-1}\dR g +g^{-1}\textbf{A}{^\gog}g \\
         \textbf{F}{^\gog} & \longmapsto & g^{-1}\textbf{F}{^\gog}g,
        \end{array}
\right.
\]
for any smooth map $g$ from $\mathcal{X}$ to $\goG$.
Actually $\gog$-valued forms correspond to
connections on a principal bundle over $\mathcal{X}$ as described below.

\subsubsection{Geometric viewpoint: principal bundles}
A way to represent connections consists with working in an associated principal principal
bundle over $\mathcal{X}$ with structure group $\goG$:
\[
\goG\longrightarrow \mathcal{F} \xrightarrow{\hbox{ }P\hbox{ }} \mathcal{X}
\]
Here, if $r:= \hbox{dim}\goG$, $\mathcal{F}$ is an $(n+r)$-dimensional
manifold equipped with a submersion $P:\mathcal{F}\longrightarrow \mathcal{X}$, such that, for any $\textsf{x}\in \mathcal{X}$, the fiber
$\mathcal{F}_\textsf{x}:= P^{-1}(\{\textsf{x}\})$ is diffeomorphic to $\goG$ and there exists a right action of $\goG$ on $\mathcal{F}$
\[
 \begin{array}{ccc}
  \mathcal{F}\times \goG & \longmapsto & \mathcal{F}\\
  (\textsf{z},g)& \longmapsto & \textsf{z}\cdot g
 \end{array}
\]
such that the $\goG$-orbit of any point $\textsf{z}\in\mathcal{F}$ coincides with the fiber
$\mathcal{F}_{P(\textsf{z})}$ containing $\textsf{z}$.
We hence get a representation of $\gog$
in the space of tangent vector fields $\mathcal{X}(\mathcal{F})$
\[
 \begin{array}{ccc}
  \mathcal{F}\times \gog & \longmapsto & T\mathcal{F}\\
  (\textsf{z},\xi^\gog)& \longmapsto & (\textsf{z},\textsf{z}\cdot \xi^\gog),
 \end{array}
\]
where $\textsf{z}\cdot \xi^\gog:= \frac{d(\textsf{z}\cdot e^{\varepsilon \xi^\gog})}{d\varepsilon}(0)$, which induces a vector space isomorphism $T_\textsf{z}(\mathcal{F}_{P(\textsf{z})}) \simeq \gog$. As a consequence of these
definitions, for any $\textsf{z}\in \mathcal{F}$, the kernel of $\dR P_\textsf{z}$
in $T_\textsf{z}\mathcal{F}$ coincides with the \emph{vertical subspace}
$V_\textsf{z}:= \textsf{z}\cdot \gog:= \{\textsf{z}\cdot \xi^\gog;\; \xi^\gog\in \gog\}$.

\subsubsection{\emph{General} Ehresmann connections}\label{paragraphEhresmann}

A general \emph{Ehresmann connection} (as defined in \cite{Ehresmann50}) is
a distribution $(H_\textsf{z})_{\textsf{z}\in \mathcal{F}}$ of subspaces of $T\mathcal{F}$ such that, $\forall \textsf{z}\in \mathcal{F}$,
$H_\textsf{z}\oplus V_\textsf{z}= T_\textsf{z}\mathcal{F}$. We call each subspace $H_\textsf{z}$ a \emph{horizontal subspace}.
It can be completely defined by a 1-form $\theta^\gog\in \Omega^1(\mathcal{F})\otimes \gog$, with coefficients in $\gog$, such that
\[
\forall \textsf{z}\in \mathcal{F},\quad
 \hbox{Ker}\theta^\gog_\textsf{z} = H_\textsf{z}.
\]
This form is not unique. However if we impose a \textbf{normalization} condition
\begin{equation}\label{normalization0}
\forall \textsf{z}\in \mathcal{F},\forall \xi^\gog\in \gog,\quad
 \theta^\gog_\textsf{z}(\textsf{z}\cdot \xi^\gog) = \xi^\gog,
\end{equation}
then $\theta^\gog$ is uniquely defined.

Note that, for a general Ehresmann connection, the dependence of $H_\textsf{z}$ in $\textsf{z}$, where $\textsf{z}$ runs
in a fiber $\mathcal{F}_\textsf{x}$, may be completely arbitrary. Hence
this notion is more general than the standard connection used in Physics.
Indeed it turns out that the standard connections in Physics
and in Mathematics satisfy the further \textbf{equivariance condition}
$\left(e^{t(\cdot\xi^\gog)}\right)^*\theta^\gog = \hbox{Ad}_{e^{-t\xi^\gog}}\theta^\gog$, $\forall t$, which implies
\begin{equation}\label{equivariance0}
 L_{\textsf{z}\cdot\xi^\gog}\theta^\gog + [\xi^\gog,\theta^\gog]=0.
\end{equation}
A key observation is that, if (\ref{normalization0}) is
satisfied, then
$L_{\textsf{z}\cdot\xi^\gog}\theta^\gog  + [\xi^\gog,\theta^\gog]= \dR(\textsf{z}\cdot\xi^\gog\iN \theta^\gog)
+ \textsf{z}\cdot\xi^\gog\iN \dR\theta^\gog + [\xi^\gog,\theta^\gog] = 0
+ \textsf{z}\cdot\xi^\gog\iN(\dR\theta^\gog + \frac{1}{2}[\theta^\gog\wedge \theta^\gog])$.
Hence
\[
\left\{\begin{array}{ccc}
    \textsf{z}\cdot \xi^\gog\iN \theta^\gog & = & \xi^\gog \\
    L_{\textsf{z}\cdot\xi^\gog}\theta^\gog + [\xi^\gog,\theta^\gog] & = & 0
    \end{array}\right.
\quad \Longleftrightarrow \quad
\left\{\begin{array}{lcc}
    \textsf{z}\cdot \xi^\gog\iN\ \theta^\gog & = & \xi^\gog \\
    \textsf{z}\cdot\xi^\gog\iN(\dR\theta^\gog + \frac{1}{2}[\theta^\gog\wedge\theta^\gog])
    & = & 0
    \end{array}\right.
\]
Beware that in most references the term 'Ehresmann connection' is used for meaning '\emph{normalized equivariant} Ehresmann connection'.

\subsubsection{Relationship between both points of view}\label{relationshipYM}
Consider a section $\sigma$ of $\mathcal{F}$ over some open subset of $\mathcal{X}$.
For avoiding clumsiness we assume that $\sigma$ is defined globally on $\mathcal{X}$,
i.e.  $\sigma:\mathcal{X}\longrightarrow \mathcal{F}$.
Then, for any $\theta^\gog\in \Omega^1(\mathcal{F})\otimes \gog$ which is \emph{normalized} and \emph{equivariant},
$\textbf{A}{^\gog} =\sigma^*\theta^\gog\in \Omega^1(\mathcal{X})\otimes \gog$ represents a standard connection.
Moreover if $\tilde{\sigma}:\mathcal{X}\longrightarrow \mathcal{F}$
is another section, then $\tilde{\textbf{A}}{^\gog}:= \tilde{\sigma}^*\theta^\gog$ is another connection and
$\textbf{A}{^\gog}$ and $\tilde{\textbf{A}}{^\gog}$ are related by a gauge transformation.

Actually any section $\sigma:\mathcal{X}\longrightarrow \mathcal{F}$ gives us a
diffeomorphism
\[
 \begin{array}{ccc}
  \mathcal{X}\times \goG & \longrightarrow & \mathcal{F}\\
  (\textsf{x},g) & \longmapsto & \sigma(\textsf{x})\cdot g
 \end{array}
\]
the inverse of which provides us with a local chart
\[
 \begin{array}{ccl}
  \mathcal{F} & \longrightarrow & \mathcal{X}\times \goG\\
  \textsf{z} & \longmapsto & (\textsf{x},g)\hbox{ s.t. } \textsf{z} = \sigma(\textsf{x})\cdot g
 \end{array}
\]
In these coordinates the normalization condition (\ref{normalization0}) reads:
$\exists \textbf{A}{^\gog}\in\Omega^1(\mathcal{F})\otimes \gog$ s.t.
\begin{equation}\label{normalization1}
 \theta^\gog = g^{-1}\dR g + g^{-1}\textbf{A}{^\gog}g\quad\hbox{and}\quad
 (\textsf{z}\cdot \xi^\gog)\iN \textbf{A}{^\gog} =0,\quad \forall \xi^\gog\in \gog
\end{equation}
and, if so, the equivariance condition (\ref{equivariance0}) reads
\begin{equation}\label{equivariance1}
 L_{z\cdot \xi^\gog}\textbf{A}{^\gog} =0,\quad \forall \xi^\gog\in \gog.
\end{equation}
Note that (\ref{normalization1}) means that $\textbf{A}{^\gog}$ has the decomposition $\textbf{A}{^\gog} = \textbf{A}{^\gog}_\mu dx^\mu$,
where each $\textbf{A}{^\gog}_\mu$ is a function on $\mathcal{F}$ (i.e. depending on the coordinates $x$ and $g$),
whereas the equivariance condition (\ref{equivariance1}) then means that actually the functions
$\textbf{A}{^\gog}_\mu$ depend only on $\textsf{x}$.

\subsection{Gravity and Cartan connections}
\subsubsection{Levi-Civita connections in orthonormal moving frames}\label{paragraphLeviCivita}

Let $\mathcal{X}$ be a manifold of dimension $n$ and $\gos$ a vector space of the same dimension $n$. Assume we are given $e^\gos\in \gos\otimes \Omega^1(\mathcal{X})$, an $\gos$-valued 1-form of rank $n$ everywhere. It provides us with a \emph{solder form}, i.e., at any point $\textsf{x}\in \mathcal{X}$, an isomorphism $T_\textsf{x}\mathcal{X}\longrightarrow \gos$. By choosing a basis $(E_1,\cdots,E_n)$ of $\gos$ we decompose $e^\gos$ as $e^\gos = e^aE_a$. Then the components $(e^1,\cdots,e^n)$ form a coframe on $\mathcal{X}$. We will thus call \emph{coframe} or \emph{soldering form} any $e^\gos\in \gos\otimes \Omega^1(\mathcal{X})$ of maximal rank (see Definition \ref{definitioncofram31}).
By the same token we define the dual frame $(e_1,\cdots,e_n)$. Then any connection $\nabla$ on $T\mathcal{X}$ can be characterized by an $\hbox{End}(\gos)$-valued 1-form $\gamma^{gl(\gos)}$, the components in a basis $(E_1,\cdots ,E_n)$ of which are $\left(\gamma{^a}_b\right)_{1\leq a,b\leq n}$ so that $\nabla$ is given by  $\nabla_Xe_a = \gamma{^b}_a(X)e_b$, for any smooth vector field $X$. We define $\gamma{^a}_{bc}:= \gamma{^a}_b(e_c)$, so that
\begin{equation}\label{exportKK}
 \gamma{^a}_b = \gamma{^a}_{bc}e^c
\end{equation}
This connection is \emph{torsion free} iff $\dR e^a + \gamma{^a}_b\wedge e^b = 0$.

If furthermore $\gos$ is endowed with a non degenerate bilinear form $\textsf{b}$, then $\mathcal{X}$ is endowed with the
pseudo-Riemannian metric $\textbf{g}:= (e^\gos)^*\textsf{b} =
\textsf{b}_{ab}e^a\otimes e^b$, where $\textsf{b}_{ab}:= \textsf{b}(E_a,E_b)$. A connection $\nabla$, which is defined by $\gamma^{gl(\gos)}$, \emph{respects the metric} $\textbf{g}$ iff the coefficients of $\gamma^{gl(\gos)}$ are in  $so(\gos,\textsf{b})$, i.e. $\gamma^{ab}:= \gamma{^a}_{b'}\textsf{b}^{b'b}$ is skewsymmetric.

The Levi-Civita connection $\nabla^{T\mathcal{X}}$ on $T\mathcal{X}$ on $(\mathcal{X},\textbf{g})$ is the unique
connection which is torsion free and respects the metric.

\subsubsection{The Palatini formulation of gravity}\label{paragraphPalatinidebase}
The previous framework allows us to set the so-called 'Palatini' (also called 'Trautman' in \cite{neemanregge}) formulation of gravity theories and its $n$-dimensional generalizations as follows. Suppose we are given some model $n$-dimensional space $(\gos,\textsf{b})$ as in \S \ref{paragraphLeviCivita} and an $n$-dimensional \emph{oriented} manifold $\mathcal{X}$. Let $\gol:= so(\gos,\textsf{b})$. Consider the set of pairs $\left( e^{\gos},\gamma^\gol\right)$, where $e^{\gos} \in \gos\otimes \Omega^1(\mathcal{X})$ is a solder form on $\mathcal{X}$ and $\gamma^\gol\in \gol\otimes \Omega^1(\mathcal{X})$ is a 1-form with coefficients in $so(\gos,\textsf{b})$. Using a decomposition of $\gos$ in a basis as in \S \ref{paragraphLeviCivita}, the Palatini action is then given by
\[
 \mathscr{A}_P\left( e^{\gos},\gamma^\gol\right) =
 \int_\mathcal{X} \frac{1}{2}e^{(N-2)}_{ab}\wedge\left(
 \dR \gamma{^a}_d + \gamma{^a}_c\wedge \gamma{^c}_d\right) \textsf{b}^{db}
\]
where $e^{(n-2)}_{a_1a_2} =
\frac{1}{(n-2)!} \epsilon_{a_1\cdots a_n}
e^{a_3}\wedge \cdots \wedge e^{a_n}$ and $\epsilon_{a_1\cdots a_n}$ is the completely antisymmetric
tensor such that $\epsilon_{1\cdots n} = 1$.
Actually the  expression $\gamma{^a}_c\wedge \gamma{^c}_{b}$ is nothing but a component of $\frac{1}{2}[\gamma^\gol\wedge \gamma^\gol]$, where $[\cdot,\cdot]$ is the Lie bracket of $so(\gos,\textsf{b})$. By setting $\Gamma^\gol:= \dR \gamma^\gol + \frac{1}{2}[\gamma^\gol\wedge \gamma^\gol]$ and $\Gamma^{ab}:= \Gamma{^a}_{b'}\textsf{b}^{b'b}$, the Palatini action reads $\mathscr{A}_P\left( e^{\gos},\gamma^\gol\right) =
 \int_\mathcal{X} \frac{1}{2}e^{(N-2)}_{ab}\wedge\Gamma^{ab}$.

It is well-known that critical points of $\mathscr{A}_P$ correspond to solutions of the Einstein equations in the vacuum on $\mathcal{X}$: to $e^\gos\simeq \left(e^a\right)_{1\leq a\leq n}$ and $\gamma^\gol$, it corresponds a pseudo metric $\textbf{g}:= (e^\gos)^*\textsf{b}$ and a connection $\nabla$
on $T\mathcal{X}$.
The vanishing of the first variation $\mathscr{A}_P$ with respect to variations of $\gamma^\gol$ implies that $\nabla$ is the Levi-Civita connection. The vanishing of the first variation with respect to $e^\gos$ reads as the Einstein equation.

This description requires the existence of a moving frame on $\mathcal{X}$, which is possible only \emph{locally} in general, for topological reasons. This can be fixed by, e.g., replacing the $\R^n$-valued form $e^\gos$ by a form with values in some vector bundle $V\mathcal{X}$ of rank $n$ equipped with a pseudo metric and which is topologically equivalent to $T\mathcal{X}$ and $\gamma^\gol$ by a 1-form with coefficients in the bundle $so(V\mathcal{X})$. An alternative way to fix this point would be to work on the frame bundle.

\subsubsection{Lifting on the frame bundle}\label{paragraphLiftingdebase}
As in \S \ref{sectiongeneraliteYM} one can associate to any connection $\nabla^{T\mathcal{X}}$ on $T\mathcal{X}$ a normalized and equivariant Ehresmann connection $\nabla^{\mathcal{P}}$ on a principal bundle $\pi:\mathcal{P}\longrightarrow \mathcal{X}$ associated to $T\mathcal{X}$. The simplest choice for $\mathcal{P}$ is the bundle $F(T\mathcal{X})$ defined as follows, which can be identified with the following subset of $\gos\otimes T^*\mathcal{X}$:
\[
F(T\mathcal{X}) :=
(\gos\otimes T^*\mathcal{X})_{iso} :=
 \{(\textsf{x}, A^\gos)\in
 \gos\otimes T^*\mathcal{X} \,;\,
 A^\gos: T_\textsf{x}\mathcal{X} \longrightarrow \gos\hbox{ is an isomorphism}\} 
\]
The group $GL(\gos)$ of linear automorphisms of $\gos$ acts on the right on $F(T\mathcal{X})$ through $(g,A^\gos)\longmapsto A^\gos \cdot g:=  g^{-1}A^\gos$.

The canonical fibration map $\pi:\gos\otimes T^*\mathcal{X} \longrightarrow \mathcal{X}$, $(\textsf{x}, A^\gos)\longmapsto \textsf{x}$, defines a \emph{canonical} $\gos$-valued 1-form $\varphi^\gos$ on $\gos\otimes T^*\mathcal{X}$, given by $\varphi^\gos:= \pi^* A^\gos$. Its restriction on $(\gos\otimes T^*\mathcal{X})_{iso}$ (which we still denote by $\varphi^\gos$) is the canonical soldering form on $(\gos\otimes T^*\mathcal{X})_{iso} = F(T\mathcal{X})$.

Now consider a (possibly local) section $\alpha:\mathcal{X}\longrightarrow F(T\mathcal{X})$. It allows us to trivialize $F(T\mathcal{X})$, i.e. to construct a diffeomorphism
\[
 \begin{array}{ccc}
  \mathcal{X}\times GL(n,\R) & \longrightarrow & F(T\mathcal{X}) \\
  (\textsf{x},g) & \longmapsto & (\textsf{x},\alpha_\textsf{x}\cdot g)
 \end{array}
\]
Then any connection $\nabla^{T\mathcal{X}}$ on $T\mathcal{X}$ is defined by a 1-form $\gamma^{gl(\gos)}$ with coefficients in $gl(\gos)$ by setting that, for any smooth tangent vector fields $X,Y$ on $\mathcal{X}$,
\[
 \langle\alpha,\nabla^{T\mathcal{X}}_XY\rangle = L_X \langle\alpha, Y\rangle + \gamma^{gl(\gos)}(X) \langle\alpha, Y\rangle
\]
where $\langle\cdot,\cdot \rangle$ is the $\gos$-valued pairing between $\gos\otimes T^*\mathcal{X}$ and $T\mathcal{X}$. 

Assume now that we are given a basis $(E_1,\cdots, E_n)$ of $\gos$. Then, to any $(\textsf{x}, A^\gos)\in (\gos\otimes T^*\mathcal{X})_{iso}$ it corresponds a unique frame in $T_\textsf{x}\mathcal{X}$ which is the inverse image of $(E_1,\cdots, E_n)$ by $A^\gos$. By applying this in particular for $A^\gos = (\alpha_\textsf{x})^*\varphi^\gos$ we get a moving frame $(e_1,\cdots, e_n)$ on $\mathcal{X}$ and hence the matrix representation $(\gamma{^a}_b)_{a,b}$ of $\gamma^{gl(\gos)}$ in this basis. Then the previous relation translates as $\nabla^{T\mathcal{X}}_XY^a = L_XY^a + \gamma{^a}_b(X)Y^b$, where $X = X^ae_a$ and $Y = Y^be_b$.

Moreover, still by using the trivialization, we can define a 1-form $\varphi^{gl(\gos)}$ on $F(T\mathcal{X})$ with coefficients in $gl(\gos)$ by
\[
 \varphi^{gl(\gos)}:= g^{-1}\gamma^{gl(\gos)} g + g^{-1}\dR g\in gl(\gos)\otimes \Omega^1(F(T\mathcal{X}))
\]
This 1-form is obviously normalized and equivariant and its restriction to the image of $\alpha$ coincides with $\gamma^{gl(\gos)}$. As in \ref{paragraphEhresmann} $\varphi^{gl(\gos)}$ defines at each point $\textsf{z}\in F(T\mathcal{X})$ a horizontal subspace $H_\textsf{z}:= \hbox{Ker}\varphi^{gl(\gos)}_\textsf{z}\in T_\textsf{z}F(T\mathcal{X})$ and hence
a normalized and equivariant Ehresmann connection $\nabla^{F(T\mathcal{X})}$ on $F(T\mathcal{X})$. On the other hand it is clear that the restriction of $\varphi^\gos_\textsf{z}$ on the horizontal space $H_\textsf{z}$ is an isomorphism.
Hence the rank of $\varphi^\gos+\varphi^{gl(\gos)}\in (\gos\oplus gl(\gos))\otimes \Omega^1((\gos\otimes T^*\mathcal{X})_{iso})$ is maximal everywhere, which means that it provides us with a coframe on $(\gos\otimes T^*\mathcal{X})_{iso}$.

Assume furthermore that $\mathcal{X}$ is pseudo Riemannian and, for simplicity, is oriented and that the connection $\nabla^{T\mathcal{X}}$ respects the metric. Then we can reduce $F(T\mathcal{X})$ to the bundle of orthonormal frames $SO(T\mathcal{X})$ and
replace $gl(\gos)$ by $\gol:= so(\gos,\textsf{b})$. Hence
$\varphi^\gol = \varphi^{gl(\gos)}$ has coefficients in $\gol$.
We remark that $\varphi^\gos+\varphi^\gol$ encodes exactly the pair $(e^\gos,\gamma^\gol)$ which are the dynamical fields in the Palatini formulation of gravity. The 1-form $\varphi^\gos+\varphi^\gol$ is a particular case of a \emph{Cartan connection}. In the case where the bundle $F(T\mathcal{X})$ admits a two-sheeted spin cover $Spin(T\mathcal{X})$ we can extend these definitions by considering the pull-back images of $\varphi^\gos$ and $\varphi^\gol$ by the cover map $Spin(T\mathcal{X}) \longrightarrow SO(T\mathcal{X})$.

\subsubsection{Cartan connections and Cartan geometries}

Cartan connexions were defined by Ehresmann in \cite{Ehresmann50}. A comprehensive presentation of Cartan geometries and of their relationship with gravity theories can be found in \cite{Wise2009} and a full treatise in \cite{sharpe}.

Cartan geometries can be seen as smooth deformations of \emph{Klein geometries} which, themselves, are a way to understand and generalize Euclidean spaces or the Minkowski space as symmetric spaces. Within Klein geometry the relevant properties of a space are encoded in the group of symmetry $\goP$ (like \emph{Poincaré}) acting on the space on the right. Moreover the subgroups of $\goP$ which leave a given point invariant can be identified with a subgroup $\goL$ (like \emph{Lorentz}) of $\goP$. As a consequence the space can be identified with the coset $\goP/\goL$.
All that defines a principal right bundle $\goL\longrightarrow\goP \longrightarrow \goP/\goL$ over $\goP/\goL$ with structure group $\goL$. The infinitesimal structure of this geometry is encoded by the canonical Maurer--Cartan 1-form $g^{-1}\dR g$ on $\goP$.

A Cartan geometry is described by principal right bundle $\goL\longrightarrow\mathcal{P}\longrightarrow \mathcal{X}$ over a manifold $\mathcal{X}$ of dimension equal to $\hbox{dim}(\goP/\goL)$, which is endowed with a \emph{Cartan connection} which can be seen as a deformation of the canonical Maurer--Cartan 1-form $g^{-1}\dR g$ on $\goP$:
let $\gop$ and $\gol$ be the Lie algebras of, respectively, $\goP$ and $\goL$.
A \emph{Cartan connection} on $\goL\longrightarrow\mathcal{P} \longrightarrow \mathcal{X}$ is a 1-form $\varphi^\gop\in \gop\otimes \Omega^1(\mathcal{P})$ with \emph{maximal rank everywhere} (i.e. a coframe on $\mathcal{P}$), which is \emph{equivariant} with respect to the right action of $\gol$ on $\mathcal{P}$, i.e. such that,
\begin{equation}\label{Cartangeoequiv}
 \forall \textsf{z}\in \mathcal{P},\forall \xi^\gol\in \gol,\quad L_{\textsf{z}\cdot\xi^\gol} \varphi^\gop + [\xi^\gol,\varphi^\gop] = 0
\end{equation}
and which is \emph{normalized} with respect to the right action of $\gol$ on $\mathcal{P}$, i.e. such that
\begin{equation}\label{Cartangeonorma}
 \forall \textsf{z}\in \mathcal{P},\forall \xi^\gol\in \gol,\quad \varphi^\gop_\textsf{z}(\textsf{z}\cdot \xi^\gol) = \xi^\gol,
\end{equation}
Note that the latter relation implies in particular that the restriction of $\varphi^\gop$ to a vertical subspace $T_\textsf{z}\mathcal{P}_x$ takes values in $\gol\subset \gop$. We note also that, as in \S \ref{paragraphEhresmann}, conditions (\ref{Cartangeoequiv}) and (\ref{Cartangeonorma}) are equivalent to the conditions
\[
\forall \textsf{z}\in \mathcal{P},
 \forall \xi^\gol\in \gol,\quad
 \textsf{z}\cdot \xi^\gol \iN \varphi^\gop = \xi^\gol\quad \hbox{and}\quad
 \textsf{z}\cdot \xi^\gol\iN \left(\dR \varphi^\gop + \frac{1}{2} [\varphi^\gop\wedge \varphi^\gop]\right) =0
\]
A \emph{Cartan geometry} is a principal bundle $\goL\longrightarrow\mathcal{P}\longrightarrow\mathcal{X}$ endowed with a \emph{Cartan connection} $\varphi^\gop$. Its curvature 2-form $\dR \varphi^\gop + \frac{1}{2} [\varphi^\gop\wedge \varphi^\gop]$ is an obstruction for $\mathcal{X}$ to be locally identified with $\goP/\goL$.

We consider here \emph{reductive} Cartan geometries: a Cartan geometry $(\goL\longrightarrow\mathcal{P}\longrightarrow\mathcal{X},\nabla)$ modeled on $\goL\longrightarrow\goP \longrightarrow \goP/\goL$ is \emph{reductive} if there exists a vector space decomposition
\[
 \gop = \gol \oplus \gos
\]
which is \emph{invariant by the adjoint action} of $\goL$ on $\gop$.

The example in \S \ref{paragraphLiftingdebase} corresponds to a reductive Cartan geometry with $\goP = SO(\gos,\textsf{b})\ltimes \gos$, $\goL = SO(\gos,\textsf{b})$. In this case a Cartan connection $\varphi^\gop$
on $\goL\longrightarrow\mathcal{P}\longrightarrow\mathcal{X}$ describes a pseudo Riemannian structure and a metric preserving connection on $\mathcal{X}$. Through the $\hbox{Ad}_\goL$-invariant splitting $\gol\oplus \gos$, a Cartan connection $\varphi^\gop$ can be decomposed as $\varphi^\gop = \varphi^\gos + \varphi^\gol$. We recover hence the description of \S \ref{paragraphLiftingdebase}. The standard General Relativity theory corresponds to the case where where $(\gos,\textsf{b})$ is the 4-dimensional Minkowski space, $\goP = SO(1,3)\ltimes \gos$ and $\goL = SO(1,3)$.

Note that if a pseudo Riemannian manifold $\mathcal{X}$ admits a spin structure, we can replace the bundle $SO(T\mathcal{X})$ by its 2-sheeted covering $Spin(T\mathcal{X})$ its suffices to define its Cartan connection as the pull-back of $\varphi^\gop\in so(T\mathcal{X})\otimes \Omega^1(\mathcal{P})$ by the covering map $Spin(T\mathcal{X}) \longrightarrow SO(T\mathcal{X})$.

\subsubsection{Generalized Palatini models}\label{paragraph1.2.5}
We may generalize the Palatini model in \S \ref{paragraphPalatinidebase} by replacing the Klein model $SO(1,3)\longrightarrow SO(1,3)\ltimes \gos\longrightarrow \gos$ of a Minkowski space by a reductive Klein model $\goL\longrightarrow\goP \longrightarrow \goP/\goL$.
For instance keeping $\goL = SO(1,3)$ but replacing $SO(1,3)\ltimes \gos$ by $SO(1,4)$ or $SO(2,3)$ leads to other gravitation theories with a non vanishing cosmological constant (see \S \ref{soussectioncosmo}).

The extra ingredient is a constant tensor $\kappa{_\gol}^{\gos\gos}\in \gol^*\otimes\gos\wedge\gos$ which is invariant by the adjoint action of $\goL$, i.e. such that, $\forall g\in \goL$, $\hbox{Ad}^*_g\otimes \hbox{Ad}_g \otimes \hbox{Ad}_g(\kappa{_\gol}^{\gos\gos}) = \kappa{_\gol}^{\gos\gos}$. We set
\[
 \mathscr{A}_P(\theta^\gos,\varphi^\gol) =
 \int_\mathcal{X}\frac{1}{2}\kappa{_\cell}^{\cels\cels}\theta^{(n-2)}_{\cels\cels}\wedge \Phi^\cell,
\]
where $\Phi^\gol:= \dR\varphi^\gol + \frac{1}{2}[\varphi^\gol\wedge \varphi^\gol]\in \gol\otimes \Omega^2(\mathcal{X})$ is the curvature 2-form of $\varphi^\gol$ and we use the conventions of \S \ref{paragraphIndices} for $\kappa{_\cell}^{\cels\cels}\theta^{(n-2)}_{\cels\cels}\wedge \Phi^\cell$: it means that if $\left(E_a\right)_{1\leq a\leq n}$ is a basis of $\gos$ and if $\left(\textbf{t}_i\right)_{1\leq i\leq r}$ is a basis of $\gol$, if we let $\kappa{_i}^{ab}$ be the coefficients such that $\kappa{_\gol}^{\gos\gos}:= \kappa{_i}^{ab}\textbf{t}^i\otimes E_a\otimes E_b$ and if $\theta^\gos = \theta^aE_a$ and $\Phi^\gol = \Phi^i\textbf{t}_i$,
then
\begin{equation}\label{avantProposition2point1}
  \frac{1}{2}\kappa{_\cell}^{\cels\cels}\theta^{(n-2)}_{\cels\cels}\wedge \Phi^\cell
:= \frac{1}{2}\kappa{_i}^{ab}\theta^{(n-2)}_{ab}\wedge \Phi^i.
\end{equation}
Here it is worth to introduce a specific basis of $\gol$ in the case where $\goL = SO(\gos,\textsf{b})$, through the following, the proof of which is straightforward.
\begin{prop}[basis of $\gol = so(\gos,\textbf{b})$]
Let $(\gos,\textbf{b})$ be a vector space endowed with a symmetric non degenerate bilinear form $\emph{\textbf{b}}$. Let $\left(E_a\right)_{1\leq a\leq n}$ be a basis of $\gos$. Then there exists a unique basis of $\gol:= so(\gos, \emph{\textbf{b}})$, which we denote by $\left(\emph{\textbf{t}}^{ab}\right)_{1\leq a<b\leq n}$, such that:
for any $\xi^\gol\in so(\gos,\emph{\textbf{b}})$, if $(\xi{^a}_b)_{1\leq a,b,\leq n}$ is the matrix of $\xi^\gol$ in $\left(E_a\right)_{1\leq a\leq n}$, i.e. such that $\xi^\gol(E_a) = \xi{^b}_aE_b$, and if we let $\xi^{ab}:= \xi{^a}_{b'}\emph{\textsf{b}}^{b'b}$, then $\xi^\gol = \sum_{1\leq a<b\leq n}\xi^{ab}\emph{\textbf{t}}_{ab}$.
Moreover since $\xi^{ab}+\xi^{ba} =0$ (because $\xi^\gol\in \gol$), by defining $\emph{\textbf{t}}_{ba}:= - \emph{\textbf{t}}_{ab}$ for $1\leq b\leq a \leq n$, we can write
\begin{equation}\label{conventionIndicesso}
  \xi^\gol = \frac{1}{2}\xi^{ab}\emph{\textbf{t}}_{ab}
 := \frac{1}{2}\sum_{1\leq a,b\leq n}\xi^{ab}\emph{\textbf{t}}_{ab}
\end{equation}
\end{prop}
Thus the set of indices $\{i\in\N\,|\,1\leq i\leq r\}$ in (\ref{avantProposition2point1}) is actually the set of ordered pairs $\{ab=[a,b]\in \N^2\,|\,1\leq a < b\leq n\}$.
Back to (\ref{avantProposition2point1}), by choosing
\begin{equation}\label{kappasimple}
\kappa{_\gol}^{\gos\gos} :=
\frac{1}{2}\textbf{t}^{ab}\otimes \textbf{t}_a\wedge \textbf{t}_b,
\hbox{ i.e. }
\kappa{_{[c,d]}}^{ab} =
\delta^{ab}_{cd}:= \delta^a_c\delta^b_d -\delta^a_d\delta^b_c
\end{equation}
we recover the standard Palatini action
\[
\int_\mathcal{X}\frac{1}{2}\kappa{_\cell}^{\cels\cels}\theta^{(n-2)}_{\cels\cels}\wedge \Phi^\cell =
\int_\mathcal{X} \frac{1}{4}\kappa{_{[c,d]}}^{ab}\theta^{(N-2)}_{ab}\wedge \Phi^{cd}
= \int_\mathcal{X} \frac{1}{2}\theta^{(N-2)}_{ab}\wedge\Phi^{ab}
\]

\subsection{Towards variational formulations on the principal bundle}\label{paragraphTowards}
The basic ideas behind the variational theories expounded in this paper is to find a variational formulation of Yang--Mills equations or of gravitation sitting on the principal bundle. In the case of Yang--Mills theories, a simple way to proceed is based on the fact that, roughly speaking, if the structure group $\goG$ is compact, for any $\theta^\gog \in \gog\otimes\Omega^1(\mathcal{F})$ which is normalized and equivariant we can write
$\theta^\gog = g^{-1}\dR g + g^{-1}\textbf{A}^\gog g$ in a trivialization, where $\textbf{A}^\gog = \textbf{A}^\gog_\mu(x)\dR x^\mu$ and hence
\[
 \int_\mathcal{X}\left|\textbf{F}^\gog\right|^2
 \hbox{dvol}_\mathcal{X}
 = \frac{1}{\hbox{vol}(\goG)}
 \int_{\mathcal{F}}\left|\Theta^\gog\right|^2
 \hbox{dvol}_\mathcal{F}
\]
where $\textbf{F}^\gog = \dR \textbf{A}^\gog + \frac{1}{2}[\textbf{A}^\gog\wedge\textbf{A}^\gog]$ and $\Theta^\gog:= \dR \theta^\gog + \frac{1}{2}[\theta^\gog\wedge\theta^\gog]$.

This tells us that we may replace the standard Yang--Mills action by
$\int_{\mathcal{F}}\left|\Theta^\gog\right|^2
\hbox{dvol}_\mathcal{F}$ \emph{provided that we assume the constraint that $\theta^\gog$ is equivariant and normalised}.
The delicate point is to impose these constraints: $\textsf{z}\cdot \xi^\gog \iN \theta^\gog = \xi^\gog$ and $\textsf{z}\cdot \xi^\gog \iN (\dR \theta^\gog + \frac{1}{2}[\theta^\gog\wedge \theta^\gog]) = 0$. This is more or less what is done in the action functional in Section \ref{sectionGauge} through the introduction of auxiliary fields which play the role of Lagrange multipliers.

\section{Notations, conventions and some useful results}\label{paragraphConventions}
Through the paper the interior product of a vector with an exterior differential form is denoted by $\iN$.
Some gothic letters have been chosen in relation to their possible physical meaning:
\[
 \begin{array}{ll}
\gos & \hbox{like \emph{space}} \\
  \goG\hbox{ and }\gog & \hbox{for a structure \emph{group} (e.g. }SU(m)\hbox{ or the \emph{Lorentz} group) and its Lie algebra}
  \\
  \goL\hbox{ and }\gol & \hbox{like the \emph{Lorentz} group and its Lie algebra}  \\
  \goP\hbox{ and }\gop & \hbox{like the \emph{Poincaré} group and its Lie algebra}
 \end{array}
\]
Underlined letters $\cels,\celg,\cell,\celp,\celu$ refer to pairs of repeated indices, i.e. duality pairings, see \S \ref{paragraphIndices} below.

\subsection{Linear representations and tensor products of representations}\label{representations}
In the following $\goG$ is a finite dimensional Lie group of dimension $r$ and  $\gog$ its Lie algebra.\\

\noindent
\textbf{Dual representations} ---
Let $V$ be a finite dimensional vector space and let $V^*$ be its dual space.
Let $\hbox{R}:\goG\longrightarrow GL(V)$, $g\longmapsto \hbox{R}_g$,
be a linear representation of $\goG$. We define its dual representation
$\hbox{R}^*:\goG\longrightarrow GL(V^*)$ by:
$\forall g\in \goG$,
\[
\forall \alpha \in V^*,
\forall u\in V,\quad
 \langle \hbox{R}_g^*\alpha,u\rangle :=
 \langle \alpha,\hbox{R}_{g^{-1}}u\rangle
\]
where $\langle \cdot,\cdot \rangle$ denotes the duality pairing.
Similarly given a linear representation $\rho:\gog\longrightarrow gl(V)$ of
$\gog$, we define its dual representation
$\rho^*:\gog\longrightarrow gl(V^*)$ by:
$\forall \xi\in \gog$,
\[
\forall \alpha \in V^*,
\forall u\in V,\quad
 \langle \rho^*(\xi)\alpha,u\rangle :=
 - \langle \alpha,\rho(\xi) u\rangle
\]
These definitions give rise to the relations
\begin{equation}\label{Adjointpairing}
\forall \alpha \in V^*,
\forall u\in V,\quad
 \langle \hbox{R}_g^*\alpha,\hbox{R}_gu\rangle = \langle\alpha,u\rangle
\end{equation}
and
\begin{equation}\label{adjointpairing}
\forall \alpha \in V^*,
\forall u\in V,\quad
 \langle \rho^*(\xi)\alpha,u\rangle
 + \langle \alpha,\rho(\xi) u\rangle = 0
\end{equation}

\noindent
\textbf{Adjoint and coadjoint representations} --- The \emph{adjoint} representation of $\goG$ maps any
$g\in \goG$ to $\hbox{Ad}_g\in GL(\gog)$ defined by: $\forall \zeta\in \gog$,
$\hbox{Ad}_g\zeta:= \frac{d}{dt}\left( ge^{t\zeta}g^{-1}\right)|_{t=0}\in \gog$.
If we assume that $\goG$ is a matrix group (which is always the case in our context) then $\hbox{Ad}_g\zeta = g\zeta g^{-1}$.
Following the previous definitions its dual representation\footnote{The definition given here for the \emph{adjoint} representation of $\goG$
on $\gog^*$ coincides with the standard definition of the so-called
\emph{coadjoint representation}, denoted by most Authors by $\hbox{Ad}^*$. Beware that the sign convention is opposite to the definition used by the author in \cite{helein14} and \cite{heleinvey15}.} is the \emph{co-adjoint} representation
$\hbox{Ad}^*:\goG\longrightarrow GL(\gog^*)$, defined by:
$\forall g \in \goG$,
\[
\forall \lambda\in \gog^*,
\forall \zeta\in \gog,\quad
\langle \hbox{Ad}^*_g\lambda,\zeta\rangle  := \langle \lambda,\hbox{Ad}_{g^{-1}}\zeta\rangle,
\quad
\forall \xi\in \gog.
\]
The adjoint representation of $\gog$ maps any
$\xi\in \gog$ to $\hbox{ad}_\xi\in gl(\gog)$ defined by: $\forall \zeta\in \gog$,
$\hbox{ad}_\xi\zeta:= \frac{d}{dt}\left( \hbox{Ad}_{e^{t\xi}}\zeta\right)|_{t=0}= [\xi,\zeta]\in \gog$.
Its dual representation is
$\hbox{ad}^*:\gog\longrightarrow gl(\gog^*)$, defined by:
$\forall \xi \in \gog$,
\[
\forall \lambda\in \gog^*,
\forall \zeta\in \gog,\quad
\langle \hbox{ad}^*_\xi\lambda,\zeta\rangle :=
- \langle \lambda,\hbox{ad}_\xi\zeta\rangle
\]
As a consequence of (\ref{Adjointpairing}) and
(\ref{adjointpairing}) these representations satisfy the relations
$\langle \hbox{Ad}_g^*\lambda,\hbox{Ad}_g\zeta\rangle = \langle\lambda,\zeta\rangle$
and $\langle \hbox{ad}_\xi^*\lambda,\zeta\rangle
 + \langle \lambda,\hbox{ad}_\xi\zeta\rangle = 0$,
 $\forall \lambda\in \gog^*$,
$\forall \zeta\in \gog$. \\

\noindent
\textbf{Use of bases} --- Let $(\textbf{t}_i)_{1\leq i\leq r}$ be a basis of $\gog$ and let $(\textbf{t}^i)_{1\leq i\leq r}$ be its dual basis. Then the Lie algebra structure is encoded in the structure coefficients $\textbf{c}^k_{ij}$, i.e. such
that $[\textbf{t}_i,\textbf{t}_j] = \textbf{t}_k\textbf{c}^k_{ij}$.
Then
\begin{equation}\label{expressionadjointaction}
  \hbox{ad}_{\textbf{t}_i}\textbf{t}_j
 = \textbf{c}^k_{ij}\textbf{t}_k
\end{equation}
and, for the coadjoint representation
$\hbox{ad}^*:\gog\longrightarrow gl(\gog^*)$,
\begin{equation}\label{expressioncoadjointaction}
  \hbox{ad}_{\textbf{t}_i}^*\textbf{t}^j
 = - \textbf{c}^j_{ik}\textbf{t}^k
\end{equation}

\noindent
\textbf{Tensor products} --- Given a finite family of representations of $\goG$,
$\hbox{R}_i:\goG\longrightarrow GL(V_i)$, for $1\leq i\leq a$, we define their tensor product
$\hbox{R}:= \hbox{R}_1\otimes\cdots\otimes\hbox{R}_a$ to be the map
$\hbox{R}:\goG \longrightarrow GL\left(V_1\otimes \cdots\otimes V_a\right)$ such that:
$\forall g\in \goG$, $\forall (u_1,\cdots,u_a)\in V_1\times \cdots\times V_a$
\[
 \hbox{R}_g(u_1\otimes \cdots\otimes u_a)
 = (\hbox{R}_1)_g(u_1)\otimes \cdots\otimes (\hbox{R}_a)_g(u_a)
\]
Given a finite family of representations of $\gog$,
$\rho_i:\gog\longrightarrow gl(V_i)$, for $1\leq i\leq a$, we define their tensor product
$\rho:= (\rho_1\otimes 1\otimes \cdots \otimes 1) + \cdots
+ (1\otimes \cdots \otimes 1\otimes \rho_a)$
to be the map $\rho:\gog\longrightarrow
gl\left(V_1\otimes \cdots\otimes V_a\right) $ such that: $\forall \xi\in \gog$, $\forall (u_1,\cdots,u_a)\in V_1\times \cdots\times V_a$
\begin{equation}\label{adjointontensorproduct}
 \rho(\xi)(u_1\otimes \cdots\otimes u_a) =
 (\rho_1(\xi)u_1)\otimes u_2\otimes \cdots\otimes u_a + \cdots + u_1\otimes \cdots\otimes (\rho_a(\xi)u_a)
\end{equation}

\subsection{Intrinsic indices and some standard tensors}
The proofs of our results rely on expressions involving tensors with many indices. In order to limit the proliferation of indices we adopt the following conventions.
\begin{enumerate}
 \item Given a vector space $V$, $x^V$ represents a vector in $V$.
 \item If $V_1,\cdots ,V_a$ are vector spaces, $x^{V_1\cdots V_a}$ represents a tensor in $V_1\otimes\cdots \otimes V_a$.
 \item If $V^*$ is the dual space of $V$ we may denote by $\ell_V$ (instead of $\ell^{V^*}$) an element of $V^*$. 
 \item We use this convention for any tensor: any index $V$ will refer to $V$ or to its dual space, according to its position, upper or lower, respectively.
 \end{enumerate}
Using this convention, if $\gog$ is a Lie algebra with basis $(\mathbf{t}_1,\cdots,\mathbf{t}_r)$ and dual basis $(\mathbf{t}^1,\cdots,\mathbf{t}^r)$, if we denote $\mathbf{c}^k_{ij}$ its structure coefficients in this basis, we define 
\begin{equation}\label{cgggref}
 \mathbf{c}{^\gog}_{\gog\gog}:=
 \mathbf{c}^k_{ij}\mathbf{t}_k\otimes \mathbf{t}^j\otimes \mathbf{t}^j
 \quad \in \gog\otimes \gog^*\otimes \gog^*
\end{equation}
If $(\gos,\textsf{b})$ is an $n$-dimensional Euclidean or Minkowski vector space with basis $(E_1,\cdots,E_n)$, we denote
\begin{equation}\label{deltass}
 \delta{^\gos}_\gos:= \delta^a_b E_a\otimes E^b \in \gos\otimes \gos^*
 \quad \hbox{and}\quad
 \delta{_\gos}^\gos:= \delta^b_a E^a\otimes E_b\in \gos^*\otimes \gos
\end{equation}
where $\delta^a_b$ is the Kronecker symbol.

If $(\mathcal{M},g)$ is a (pseudo-)Riemannian manifold of the same dimension as $\gos$ and $(e^1,\cdots ,e^n)$ is a (possibly locally defined) orthonormal moving frame on $\mathcal{M}$, we set $e^\gos = e^1E_1+\cdots + e^nE_n$.

A connection 1-form in this frame reads $\omega^\gol$, or $\omega{^\gos}_\gos = \omega{^a}_b \mathbf{t}_a\otimes \mathbf{t}^b$, through the identification of  $\gol = so(\gos,\textsf{b})$, the Lie algebra of isometries of $(\gos,\textsf{b})$, with a subspace of $\gos\otimes \gos^*$.

Its curvature 2-form reads $\Omega^\gol = \hbox{d}\omega^\gol + \frac{1}{2}[\omega^\gol\wedge \omega^\gol]$ or
$\Omega{^\gos}_\gos = \hbox{d}\omega{^\gos}_\gos + \omega{^\gos}_a\wedge \omega{^a}_\gos$ ($=\hbox{d}\omega{^\gos}_\gos + \omega{^\gos}_\cels\wedge \omega{^\cels}_\gos$, according to the conventions in the next paragraph). It can be represented by
$\Omega^{\gos\gos} := \Omega{^\gos}_a \textsf{b}^{a\gos}$ ($=\Omega{^\gos}_\cels \textsf{b}^{\cels\gos}$).
The decomposition of $\Omega^{\gos\gos}$ in the basis $\left(e^a\wedge e^b\right)_{1\leq a<b\leq n}$ involves the coefficients of the Riemann curvature tensor $\mathbf{R}{^{\gos\gos}}_{\gos\gos}$:
\[
 \Omega^{\gos\gos}
 = \frac{1}{2}\sum_{a,b=1}^n \mathbf{R}{^{\gos\gos}}_{ab}e^a\wedge e^b
 \left( = \frac{1}{2} \mathbf{R}{^{\gos\gos}}_{\cels\cels}e^\cels\wedge e^\cels\right)
\]
from which we define the Ricci tensor and the scalar curvature:
\begin{equation}\label{notationRicci}
 \mathbf{R}{^{\gos}}_{\gos}:= \mathbf{R}{^{\gos a}}_{\gos a} = \mathbf{R}{^{\gos\cels}}_{\gos\cels}
 \quad \hbox{and}\quad 
 \mathbf{R}:= \mathbf{R}{^{a}}_{a} = \mathbf{R}{^{\cels}}_{\cels}
\end{equation}
and the Einstein tensor
\begin{equation}\label{notationEinstein}
 \mathbf{E}{_\gos}^\gos:= \mathbf{R}{_\gos}^\gos - \frac{1}{2}\mathbf{R}\delta{_\gos}^\gos
\end{equation}

\subsection{Contractions of tensors and intrinsic indices}\label{paragraphIndices}
Using the previous conventions, in order to help to identify which pair of indices are summed when summations on repeated indices occur, we introduce the following conventions (recall that the summation over pairs of repeated indices corresponds to a duality product). 

For any integer $a\in \N^*$, let
$[\![1,a]\!]:= [1,a]\cap \N$. Let $a,b\in \N^*$ and let $(V_1,\cdots,V_a)$ and
$(W_1,\cdots,W_b)$ be two lists of vector spaces (possibly with repetition).
Let $c\in \N^*$ such that $c\leq \hbox{min}(a,b)$ and let
$\sigma:[\![1,c]\!] \longrightarrow [\![1,a]\!]$ and
$\tau:[\![1,c]\!] \longrightarrow [\![1,b]\!]$ be
two \emph{one-to-one} maps.
Assume that, $\forall i\in [\![1,c]\!]$,
$V_{\sigma(i)}$ is in duality with
$W_{\tau(i)}$.
We then define the \textbf{contracted tensor product} to be the bilinear map
\[
\textsf{C}_{\sigma,\tau}:
\left(V_1\otimes\cdots\otimes V_a\right)
\times
\left(W_1\otimes\cdots\otimes W_b\right)
\longrightarrow
Z_1\otimes\cdots\otimes Z_{a+b-2c}
\]
where $Z_1, \cdots , Z_{a+b-2c}$ is the list of vector spaces obtained after removing all vector spaces $V_{\sigma(i)}$ and $W_{\tau(i)}$ for $i\in [\![1,c]\!]$, from the list $(V_1,\cdots,V_a, W_1,\cdots,W_b)$.
For $S\in V_1\otimes\cdots\otimes V_a$ and $T\in W_1\otimes\cdots\otimes W_b$,
$\textsf{C}_{\sigma,\tau}(S, T)$ is the tensor in $Z_1\otimes\cdots\otimes Z_{a+b-2c}$ obtained by contracting, in the tensor product
$S\otimes T$, all pairs of indices associated to the positions $(\sigma(i),\tau(i))$, for $i\in [\![1,c]\!]$.

A precise definition is given at the end of this paragraph. However it may be more illuminating to start by illustrating this definition through examples.\\

\noindent
\textbf{A list of examples}\\

In the following $(\textbf{v}_i)_i$ is a basis of $V$ and $(\textbf{v}^i)_i$ is its dual basis.
\begin{enumerate}
 \item we denote the duality product between $x^V$ and $\ell_V$ by
 \[
  \ell_{\celv} x^{\celv}
  := \textsf{C}_{\sigma,\tau}(\ell_V,x^V)
  := \ell_1x^1+\cdots +\ell_kx^k\quad \in \R
 \]
where we use the underlined out letter
${\underline{V}}$ repeated twice to indicate the duality pairing, i.e. the summation over repeated indices. \emph{Here $\sigma$ and $\tau$ are such that $a=b=c=1$ and $(\sigma(1), \tau(1)) = (1,1)$.}
\item however if two indices $V$ are repeated but not underlined, then it means that we  consider their tensor product. Hence
\[
 \ell_V x^V :=\ell_V \otimes x^V
 = \ell_ix^j\ \textbf{v}^i\otimes \textbf{v}_j\quad \in V^*\otimes V
\]
Beware it is not commutative!
\end{enumerate}
These rules are then extended to tensors as follows:
\begin{enumerate}
 \item[(iii)] Suppose that $V$ and $W$ are two different vector spaces. Consider for example two tensors
 $S{^V}_{WW} = S{^i}_{AB}\textbf{v}_i\otimes \textbf{w}^A\otimes\textbf{w}^B\in V\otimes W^*\otimes W^*$ and $T{_V}^{WW} = T{_j}^{CD}\textbf{v}^j\otimes \textbf{w}_C\otimes\textbf{w}_D\in V^*\otimes W\otimes W$, then
 \[
 \begin{array}{ccccc}
   S{^{\celv}}_{WW}\
  T{_{\celv}}^{WW}
  & := &  S{^i}_{AB}\ T{_i}^{CD}
  \textbf{w}^A\otimes\textbf{w}^B\otimes \textbf{w}_C\otimes\textbf{w}_D & \in & W^*\otimes W^*\otimes W\otimes W \\
   S{^V}_{\celw W}\
  T{_V}^{\celw W}
  & := & S{^i}_{AB}\ T{_j}^{AD}
  \textbf{v}_i\otimes\textbf{w}^B\otimes \textbf{v}^j\otimes\textbf{w}_D & \in & V\otimes W^*\otimes V^*\otimes W \\
   S{^V}_{\celw W}\
  T{_V}^{W\celw}
  & := & S{^i}_{AB}\ T{_j}^{CA}
  \textbf{v}_i\otimes\textbf{w}^B\otimes \textbf{v}^j\otimes\textbf{w}_C & \in & V\otimes W^*\otimes V^*\otimes W \\
   S{^{\celv}}_{\celw W}\
  T{_{\celv}}^{W\celw}
  & := & S{^i}_{AB}\ T{_i}^{CA}
  \textbf{w}^B\otimes\textbf{w}_C & \in & W^*\otimes W
 \end{array}
 \]
\emph{Here $a=b=3$ and the expression on the left hand side is equal to $\emph{\textsf{C}}_{\sigma,\tau}(S,T)$, where: on the first line, $c=1$ and $(\sigma(1), \tau(1)) = (1,1)$, on the second line, $c=1$ and $(\sigma(1), \tau(1)) = (2,2)$, on the third line, $c=1$ and $(\sigma(1),\tau(1)) = (2,3)$ and, on the last line, $c=2$ and $(\sigma(1), \tau(1)) = (1,1)$ and $(\sigma(2), \tau(2)) = (2,3)$.}
\item[(iv)] if the same vector space occurs several times in each factor and several pairings occur between these factors \emph{by respecting the order}, we also struck the indices corresponding to these factors. For instance, for  $S{^V}_{WW}$ and $T{_V}^{WW}$ as before,
 \[
 \begin{array}{ccccc}
   S{^V}_{\celw\celw}\
  T{_V}^{\celw\celw}
  & :=  & S{^i}_{AB}\ T{_j}^{AB}
  \textbf{v}_i\otimes\textbf{v}^j & \in & V\otimes  V^*  \\
   S{^{\celv}}_{\celw\celw}\
  T{_{\celv}}^{\celw\celw}
  & := &  S{^i}_{AB}\ T{_i}^{AB} & \in & \R
 \end{array}
 \]
\emph{Here $a=b=3$ and the expression on the left hand side is equal to $\emph{\textsf{C}}_{\sigma,\tau}(S,T)$, where: on the first line, $c=2$, $(\sigma(1), \tau(1)) = (2,2)$ and $(\sigma(2),\tau(2)) = (3,3)$, on the second line, $c=3$ and $(\sigma(1),\tau(1)) = (1,1)$, $(\sigma(2),\tau(2)) = (2,2)$ and $(\sigma(3),\tau(3)) = (3,3)$.}
 \item[(v)] in case of ambiguity, e.g. if the same vector space occurs several times in each factor and several pairings occur between these factors but the pairing between these factors does not respect the order, we label the factors by integers in order to indicate the right couplings. For instance, for  $S{^V}_{WW}$ and $T{_V}^{WW}$ as before,
 \[
 \begin{array}{ccccc}
   S{^V}_{\celw_1\celw_2}\
  T{_V}^{\celw_2\celw_1}
  & :=  & S{^i}_{AB}\ T{_j}^{BA}
  \textbf{v}_i\otimes\textbf{v}^j & \in & V\otimes  V^* \\
   S{^{\celv}}_{\celw_1\celw_2}\
  T{_{\celv}}^{\celw_2\celw_1}
  & := &  S{^i}_{AB}\ T{_i}^{BA} & \in & \R
 \end{array}
 \]
\emph{Here $a=b=3$ and the expression on the left hand side is equal to $\emph{\textsf{C}}_{\sigma,\tau}(S,T)$, where: on the first line, $c=2$, $(\sigma(1),\tau(1)) = (2,3)$ and $(\sigma(2),\tau(2)) = (3,2)$, on the second line, $c=3$,
  $(\sigma(1),\tau(1)) = (1,1)$, $(\sigma(2),\tau(2)) = (2,3)$ and $(\sigma(3),\tau(3)) = (3,2)$.}

\item[(vi)] Lastly by using the definition of $\mathbf{c}{^\gog}_{\gog\gog}$ given by (\ref{cgggref}),  (\ref{expressionadjointaction}) translates as
\begin{equation}
\forall \xi^\gog,\zeta^\gog\in \gog,\quad  \hbox{ad}_{\xi^\gog} \zeta^\gog = \mathbf{c}{^\gog}_{\celg\celg}\xi^\celg \zeta^\celg
\end{equation}
and (\ref{expressioncoadjointaction}) as
\begin{equation}
 \forall \xi^\gog\in \gog,\forall \ell_\gog\in \gog^*,\quad \hbox{ad}_{\xi^\gog} \ell_\gog = - \mathbf{c}{^{\celg_1}}_{\celg_2\gog}\xi^{\celg_2} \ell_{\celg_1}
\end{equation}
\end{enumerate}
Note that all these conventions are independant of the choice of the bases of the vector spaces.

We extend this operation to any pair of differential forms with coefficients in tensor products: for $a,b\in \N^*$, $V_1,\cdots,V_a$,
$W_1,\cdots,W_b$, $c\in \N^*$ and
$(\sigma,\tau):[\![1,c]\!] \longrightarrow [\![1,a]\!]\times [\![1,b]\!]$ as previously, for $p,q\in \N$, we define the \emph{contracted wedge product} to be the unique bilinear map
\[
 V_1\otimes\cdots\otimes V_a
\otimes \Omega^p(\mathcal{N})
\times
W_1\otimes\cdots\otimes W_b
\otimes \Omega^q(\mathcal{N})
\xrightarrow{\textsf{C}_{\sigma,\tau}(\cdot \wedge \cdot)}
Z_1\otimes\cdots\otimes Z_{a+b-2c}
\otimes \Omega^{p+q}(\mathcal{N})
\]
such that, $\forall S\in V_1\otimes\cdots\otimes V_a$, $\forall T\in W_1\otimes\cdots\otimes W_b$, $\forall \alpha\in \Omega^p(\mathcal{N})$ and $\forall \beta\in \Omega^q(\mathcal{N})$
\[
 \textsf{C}_{\sigma,\tau}(S\otimes \alpha\wedge T\otimes \beta) =
 \textsf{C}_{\sigma,\tau}(S,T)\  \alpha\wedge \beta.
\]
\begin{defi}
Let $a,b,c\in \N^*$, $(V_1,\cdots, V_a)$ and $(W_1, \cdots, W_b)$, two lists of vector
spaces, $\sigma:[\![1,c]\!] \longrightarrow [\![1,a]\!]$ and
$\tau:[\![1,c]\!] \longrightarrow [\![1,b]\!]$ as in \S \ref{paragraphIndices}. Let $i_1,\cdots,i_c\in [\![1,a]\!]$ such that $i_1<\cdots<i_c$ and
$\{i_1,\cdots,i_c\}:=\sigma([\![1,c]\!])$ and,
similarly, $j_1,\cdots,j_c\in [\![1,b]\!]$ such that $j_1<\cdots<j_c$ and
$\{j_1,\cdots,j_c\}:=\tau([\![1,c]\!])$.
Let $(v_1,\cdots,v_a,w_1,\cdots,w_b)\in
V_1\times \cdots\times V_a\times
W_1\times \cdots\times W_b$, $p,q\in \N$ and
$\alpha\in \Omega^p(\mathcal{N})$ and
$\beta\in \Omega^q(\mathcal{N})$. We then set
\begin{equation}\label{definitionsuperwedge}
\begin{array}{l}
 \displaystyle
 \textsf{\emph{C}}_{\sigma,\tau}\left(
 \ v_1\otimes\cdots\otimes v_a\otimes\alpha\
 \wedge
 \ w_1\otimes \cdots\otimes w_b\otimes\beta\ \right):= \\
 \displaystyle \quad
 \prod_{k=1}^c\langle v_{\sigma(k)},w_{\tau(k)}\rangle
 \left(\bigotimes_{i\in [\![1,a]\!]\setminus \sigma([\![1,c]\!])}v_i\right)
 \otimes \left(\bigotimes_{j\in [\![1,b]\!]\setminus \tau([\![1,c]\!])}w_j\right) \otimes \alpha \wedge \beta
\end{array}
\end{equation}
Then the \textbf{contracted wedge product} is the unique extension of
\[
\begin{array}{r}
\textsf{\emph{C}}_{\sigma,\tau}(\cdot \wedge \cdot ):
V_1\otimes \cdots\otimes V_a\otimes \Omega^p(\mathcal{N})\times
W_1\otimes \cdots\otimes W_b\otimes \Omega^q(\mathcal{N})
\longrightarrow \quad\quad \\
\displaystyle
 \left(\bigotimes_{i\in [\![1,a]\!]\setminus \sigma([\![1,c]\!])}V_i\right)
 \otimes \left(\bigotimes_{j\in [\![1,b]\!]\setminus \tau([\![1,c]\!])}W_i\right) \otimes
 \Omega^{p+q}(\mathcal{N})
\end{array}
\]
which is bilinear.
\end{defi}

\subsection{Vector and tensor valued forms and coframes on a manifold}\label{paragraphCoframes}
\begin{defi}\label{definitioncofram31}
Let $\mathcal{N}$ be a manifold of dimension $N$ and $V,V_1,\cdots, V_a$ be vector spaces and $p\in \N$.
Let $\mathcal{O}\subset \mathcal{N}$ an open subset.
\begin{enumerate}
 \item A \emph{vector space valued $p$-form} $e^V$ on $\mathcal{O}$ is an element of $V\otimes \Omega^p(\mathcal{O})$, i.e. a $p$-form on $\mathcal{O}$ with coefficients in $V$.
 \item If $p=1$, $\emph{\hbox{dim}}\mathcal{N} = \emph{\hbox{dim}} V$ and $e^V\in V\otimes \Omega^1(\mathcal{O})$ has a maximal rank everywhere, then $e^V$ is a \textbf{coframe} on $\mathcal{O}$.
 \item If $V = V_1\oplus V_2$ and $e^V\in V\otimes \Omega^p(\mathcal{O})$, then $e^{V_1}\in V_1\otimes \Omega^1(\mathcal{O})$ and $e^{V_2}\in V_2\otimes \Omega^1(\mathcal{O})$ are the projections of $e^V$ to, respectively, $V_1$ and $V_2$, through the splitting $V = V_1\oplus V_2$, so that $e^V = e^{V_1} + e^{V_2}$.
\end{enumerate}
\end{defi}
Note that we will also meet tensor valued $p$-forms on $\mathcal{O}$, i.e. elements of $V_1\otimes
  \cdots \otimes V_a\otimes  \Omega^p(\mathcal{O})$.
Most of the time we will not specify the domain $\mathcal{O}\subset \mathcal{N}$.

Consider a vector valued 1-form $e^V$ on $\mathcal{N}$ and let's choose a basis $\left(\textbf{v}_A\right)_{1\leq A\leq m}$ of $V$, with dual basis $\left(\textbf{v}^A\right)_{1\leq A\leq m}$. Then $e^V$ decomposes in this basis as $e^V:= e^A\textbf{v}_A$ for some collection $\left(e^A\right)_{1\leq A\leq m}$ of 1-forms on $\mathcal{N}$. Obviously $e^V$ is a coframe iff $\left(e^A\right)_{1\leq A\leq m}$ is a basis of $T^*\mathcal{N}$.
For any $k\in \N$ we set
\[
\begin{array}{r}
  \underbrace{V\wedge \cdots \wedge V}_{k}
 := \{ S^{A_1\cdots A_k}\textbf{v}_{A_1}\otimes \cdots \otimes \textbf{v}_{A_k}
 \in V\otimes \cdots \otimes V\; | \quad\quad \\
S^{\sigma(A_1)\cdots \sigma(A_k)}
= (-1)^{|\sigma|}S^{A_1\cdots A_k}, \forall \sigma \in \mathfrak{S}(k)\}
\end{array}
\]
We define
\begin{equation}\label{deffifi}
 \begin{array}{ccccccc}
e^{VV} & := & e^V\wedge e^V & = & \textbf{v}_A\otimes \textbf{v}_Be^A\wedge e^B & \in & V\wedge V\otimes \Omega^2(\mathcal{N}) \\
e^{VVV} & := & e^V\wedge e^V \wedge e^V & = & \textbf{v}_A\otimes \textbf{v}_B\otimes \textbf{v}_Ce^A\wedge e^B\wedge e^C & \in &
V\wedge V\wedge V\otimes \Omega^3(\mathcal{N})
 \end{array}
\end{equation}
and so on.
We also set
$e^{(m)}:= e^1\wedge \cdots \wedge e^m\in \Omega^m(\mathcal{N})$ and
\[
 \begin{array}{ccll}
      e^{(m-1)}_A & := &
 \frac{1}{(m-1)!}\epsilon_{AA_2\cdots A_m}e^{A_2}\wedge \cdots
 \wedge e^{A_m} & \in \Omega^{m-1}(\mathcal{N}) \\
 e^{(m-2)}_{AB} & := &
 \frac{1}{(m-2)!}\epsilon_{ABA_3\cdots A_m}e^{A_3}\wedge \cdots
 \wedge e^{A_m} & \in \Omega^{m-2}(\mathcal{N})
 \\
 e^{(m-3)}_{ABC} & :=  &
\frac{1}{(m-3)!}\epsilon_{ABCA_4\cdots A_m}e^{A_4}\wedge \cdots
\wedge e^{A_m} & \in \Omega^{m-3}(\mathcal{N})
\end{array}
\]
and we define
\begin{equation}\label{defantififi}
  \begin{array}{ccccc}
e_V^{(m-1)} & := & \textbf{v}^Ae^{(m-1)}_A & \in & V^*\otimes \Omega^{m-1}(\mathcal{F})  \\
e_{VV}^{(m-2)} & := & \textbf{v}^A\otimes\textbf{v}^Be^{(m-2)}_{AB}& \in & V^*\wedge V^*\otimes\Omega^{m-2}(\mathcal{F}) \\
e_{VVV}^{(m-3)} & := &  \textbf{v}^A\otimes\textbf{v}^B\otimes\textbf{v}^Ce^{(m-3)}_{ABC}& \in & V^*\wedge V^*\wedge V^*\otimes\Omega^{m-3}(\mathcal{F})
\end{array}
\end{equation}
In the following we assume that $e^V$ is a \textbf{coframe}. Then
\begin{itemize}
 \item any 1-form $\alpha\in \Omega^1(\mathcal{N})$ can be decomposed as $\alpha = \alpha_Ae^A$ and we associate to it the $V^*$-valued function $\alpha_V:= \alpha_A\textbf{v}^A\in V^*\otimes \mathscr{C}^\infty(\mathcal{N})$.
 \item any 2-form $\beta\in \Omega^2(\mathcal{N})$ can be decomposed as $\beta = \frac{1}{2}\beta_{AB}e^A\wedge e^B$, with $\beta_{AB} + \beta_{BA} = 0$, we associate to it the $V^*\wedge V^*$-valued function $\beta_{VV}:= \beta_{AB}\textbf{v}^A\otimes \textbf{v}^B\in V^*\wedge V^*\otimes \mathscr{C}^\infty(\mathcal{N})$ .
 \item the generalization of these conventions to forms of arbitrary degress is straightforward.
\end{itemize}
Hence the following isomorphisms, which are independant of the choice of basis:
\[
 \begin{array}{ccccccc}
  \Omega^1(\mathcal{N}) & \ni & \alpha = \alpha_Ae^A & \longmapsto &
  \alpha_V:= \textbf{v}^A\alpha_A & \in & V^*\otimes\mathscr{C}^\infty(\mathcal{N}) \\
  \Omega^2(\mathcal{N})& \ni &\beta = \frac{1}{2!}\beta_{AB}e^A\wedge e^B & \longmapsto &\beta_{VV}:= \textbf{v}^A \otimes\textbf{v}^B\beta_{AB} & \in & V^*\otimes V^*\otimes \mathscr{C}^\infty(\mathcal{N})
 \end{array}
\]
Then by using the convention of \S \ref{paragraphIndices}
we have
\begin{equation}\label{ConventionV}
  \alpha = \alpha_{\celv}e^{\celv}\hbox{ for }\alpha \in \Omega^1(\mathcal{N}) \hbox{ , }\quad \beta
 = \frac{1}{2!}\beta_{\celv\celv}\ e^{\celv\celv}\hbox{ for }\beta\in \Omega^2(\mathcal{N})
\end{equation}
and so on.

Forms of degree $N-p$ and for small values of $p$ (e.g. $p=1$, 2 or 3) also decompose in the bases
\[
\left(e^{(N-1)}_A\right)_{1\leq A\leq N},\quad \left(e^{(N-2)}_{AB}\right)_{1\leq A<B\leq N}\quad \hbox{ and } \quad \left(e^{(N-3)}_{ABC}\right)_{1\leq A<B<C\leq N}
\]
of respectively $\Omega^{N-1}(\mathcal{N})$,  $\Omega^{N-2}(\mathcal{N})$ and $\Omega^{N-3}(\mathcal{N})$.
This allows us to decompose any form in $\Omega^{N-p}(\mathcal{N})$ for $p=1,2$ or 3 and leads to the following isomorphisms, which depends of the choice of the basis of $V^*$ only through the $m$-form $\textbf{v}^1 \wedge \cdots \wedge \textbf{v}^m$:
\[
 \begin{array}{ccccccc}
  \Omega^{N-1}(\mathcal{N}) & \ni & \alpha = \alpha^Ae^{(N-1)}_A & \longmapsto &
  \alpha^V:= \textbf{v}_A\alpha^A & \in & V\otimes\mathscr{C}^\infty(\mathcal{N})\\
  \Omega^{N-2}(\mathcal{N}) & \ni & \beta = \frac{1}{2!}\alpha^{AB}e^{(N-2)}_{AB} & \longmapsto &
  \beta^{VV}:= \textbf{v}_A\otimes\textbf{v}_B\beta^{AB} & \in & V\otimes V\otimes\mathscr{C}^\infty(\mathcal{N})\\
  \Omega^{N-3}(\mathcal{N}) & \ni & \gamma = \frac{1}{3!}\gamma^{ABC}e^{(N-3)}_{ABC} & \longmapsto &
  \gamma^{VVV}:= \textbf{v}_A\otimes\textbf{v}_B\otimes\textbf{v}_C\gamma^{ABC} & \in & V\otimes V\otimes V\otimes\mathscr{C}^\infty(\mathcal{N})
 \end{array}
\]
We hence can write

\begin{equation}\label{ConventionV*}
 \alpha = \alpha^{\celv}e_{\celv}^{(N-1)},
 \quad \beta = \frac{1}{2}\beta^{\celv\celv}e_{\celv\celv}^{(N-2)}
 \ \hbox{ and }\
 \gamma = \frac{1}{3!}\gamma^{\celv\celv\celv}e_{\celv\celv\celv}^{(N-3)}
\end{equation}
Note that, if we let
$(e_1,\cdots,e_N)$
be the moving frame on $\mathcal{N}$ which is dual to $(e^1,\cdots,e^N)$, then
$e^{(N-1)}_A:= e_A\iN e^{(N)}$, $e^{(N-2)}_{AB}:= e_B\iN e^{(N-1)}_A$ and  $e^{(N-3)}_{ABC}:= e_C\iN e^{(N-2)}_{AB}$.

\subsection{Connections}
\label{appendixpseudoconnections}

Let $\gog$ be a Lie algebra and consider a $\gog$-valued-form $\omega^\gog\in \gog\otimes \Omega^1 (\mathcal{N})$ defined on a smooth manifold $\mathcal{N}$. Let $V$ be a vector space representation of $\gog$ and denote by
$\rho: \gog\longrightarrow gl(V)$ the associated morphism. On the trivial vector bundle $\mathcal{N}\times V$ we define the \emph{connection} associated to $\omega^\gog$, to be the first order differential operator
\[
 \dR^{\omega^\gog}: V\otimes \mathscr{C}^\infty(\mathcal{N})
 \longrightarrow V\otimes \Omega^{1}(\mathcal{N})
\]
defined by
\[
 \forall f^\gog \in \gog\otimes \mathscr{C}^\infty(\mathcal{N}),\quad
 \dR^{\omega^\gog} f^\gog:= \dR f^\gog + (\rho\omega^\gog)f^\gog
\]
and we extend this operator as $\dR^{\omega^\gog}: V\otimes \Omega^p(\mathcal{N})
 \longrightarrow V\otimes \Omega^{p+1}(\mathcal{N})$
by
\[
 \forall\alpha^V \in V\otimes \Omega^p(\mathcal{N}),\quad
 \dR^{\omega^\gog}\alpha^V:= \dR\alpha^V + (\rho\omega^\gog)\wedge \alpha^V
\]
where, if $(\textbf{t}_i)_{1\leq i\leq r}$ is a basis of $\gog$ and $\omega^\gog = \omega^i\textbf{t}_i$,
$(\rho\omega^\gog)\wedge \alpha^\gog
:= \omega^i\wedge (\rho\textbf{t}_i) \alpha^\gog$.
We define the \emph{curvature 2-form} of $\dR^{\omega^\gog}$ to be
\[
\Omega^\gog:= \dR \omega^\gog +
\frac{1}{2}[\omega^\gog\wedge \omega^\gog]
\quad \in \gog\otimes \Omega^2(\mathcal{N})
\]
it satisfies the property that $\dR^{\omega^\gog}\circ \dR^{\omega^\gog} = (\rho\Omega^\gog) \wedge$.
Most of the representations used within this paper are the adjoint $\hbox{ad}_{\omega^\gog}$ and the coadjoint ones $\hbox{ad}^*_{\omega^\gog}$ and their tensor products. Recall that, if $\textbf{c}^k_{ij}$ are the structure coefficients of $\gog$, so that $[\textbf{t}_i,\textbf{t}_j] = \textbf{c}^k_{ij}\textbf{t}_k$ (see Section \ref{representations}), then
\[
 \forall\alpha^\gog \in \gog\otimes \Omega^p(\mathcal{N}),\quad
 \hbox{ad}_{\omega^\gog}\wedge \alpha^\gog =
 [\omega^\gog \wedge \alpha^\gog]
 = \textbf{c}^k_{ij}\  \omega^i\wedge\alpha^j\ \textbf{t}_k
\]
and if $(\textbf{t}^i)_{1\leq i\leq r}$ is the dual basis of $\gog^*$,
\[
 \forall\alpha_\gog \in \gog^*\otimes \Omega^p(\mathcal{N}),\quad
 \hbox{ad}^*_{\omega^\gog}\wedge \alpha_\gog = - \textbf{c}^k_{ij}\  \omega^i\wedge\alpha_k\ \textbf{t}^j
\]
As a consequence of these definitions and of (\ref{adjointontensorproduct}), if $\rho_1:\gog\longrightarrow gl(V_1)$, $\cdots$, $\rho_k:\gog\longrightarrow gl(V_k)$ are vector space representations and if $V = V_1\otimes \cdots \otimes V_k$, then
$\forall\alpha^V \in V\otimes \Omega^p(\mathcal{N})$,
\begin{equation}\label{domegatensoriel}
 \dR^{\omega^\gog}\alpha^V:= \dR\alpha^V + \left(\rho_1\omega^\gog\otimes 1\otimes \cdots \otimes 1\right)\wedge
 \alpha^V
 + \cdots +
 \left(1\otimes 1\otimes \cdots \otimes \rho_k\omega^\gog\right)\wedge
 \alpha^V
\end{equation}
Through a decomposition of $\alpha^V$ by using bases of the spaces $V_1,\cdots, V_k$ and by denoting by $(\rho_\ell\omega^\gog)^{i_\ell}_{j_\ell}$ the matrix coefficients of $\rho_\ell(\omega^\gog)$ in each basis, the latter relation reads
\[
 \dR^{\omega^\gog}\alpha^{i_1\cdots i_k}
 = \dR\alpha^{i_1\cdots i_k}
 + (\rho_1\omega^\gog)^{i_1}_{j_\ell}\alpha^{j_1i_2\cdots i_k}
 + \cdots +
 (\rho_k\omega^\gog)^{i_k}_{j_k}\alpha^{i_1\cdots i_{k-1} j_k}
\]
Most of the time, in order to lighten the notations we will write $\dR^{\omega^\gog} = \dR^\omega$, if there is no ambiguity.

\subsection{Some useful results}
\subsubsection{Exterior differential algebra}
\begin{lemm}\label{lemmefifi}
Let $V$ be a vector space of dimension $N$. Let $e^V\in V\otimes \Omega^1(\mathcal{N})$ be a vector valued 1-form over a manifold $\mathcal{N}$ and let $e^A$, $e_A^{(N-1)}$, $e_{AB}^{(N-2)}$ and $e_{ABC}^{(N-3)}$ as in (\ref{defantififi}). Then
 \begin{equation}\label{fifi}
\left\{\begin{array}{cclc}
 e^A\wedge e^{(N-1)}_{A'} & = & \delta^A_{A'}e^{(N)}
 & \hbox{(a)} \\
 e^A\wedge e^{(N-2)}_{A'B'} & = & \delta^A_{B'}e^{(N-1)}_{A'}
 - \delta^A_{A'}e^{(N-1)}_{B'}& \hbox{(b)} \\
 e^A\wedge e^{(N-3)}_{A'B'C'} & = & \delta^A_{C'}e^{(N-2)}_{A'B'}
 + \delta^A_{B'}e^{(N-2)}_{C'A'} + \delta^A_{A'}e^{(N-2)}_{B'C'}& \hbox{(c)}\\
    e^A\wedge e^B\wedge e^{(N-2)}_{A'B'} & = & \delta^{AB}_{A'B'}e^{(N)} & \hbox{(d)} \\
    e^A\wedge e^B\wedge e^{(N-3)}_{A'B'C'} & = &
    \delta^{AB}_{B'C'}e^{(N-1)}_{A'}
    +\delta^{AB}_{C'A'}e^{(N-1)}_{B'}
    + \delta^{AB}_{A'B'}e^{(N-1)}_{C'}
    & \hbox{(e)}
\end{array}
\right.
\end{equation}
where $\delta^{AB}_{CD}:= \delta^A_C\delta^B_D - \delta^A_D\delta^B_C$.
Moreover
\begin{equation}\label{deAB}
 \left\{\begin{array}{ccl}
         \emph{\dR} e^{(N-1)}_A & = & \emph{\dR} e^B\wedge e^{(N-2)}_{AB} \\
         \emph{\dR} e^{(N-2)}_{AB} & = & \emph{\dR} e^C\wedge e^{(N-3)}_{ABC} \\
\emph{\dR} e^{(N-3)}_{ABC} & = & \emph{\dR} e^D\wedge e^{(N-3)}_{ABCD}
        \end{array}\right.
\end{equation}
\end{lemm}
\emph{Proof} ---
Relation (\ref{fifi}) is a consequence of the following elementary results. Let $(\textbf{v}_A)_{1\leq A\leq N}$ be a basis of $V$. We denote by $(\textbf{v}^A)_{1\leq A\leq N}$ the basis of $V^*$ which is dual to $(\textbf{v}_A)_{1\leq A\leq N}$. Set $\textbf{v}^{(N)}:= \textbf{v}^1\wedge \cdots \wedge \textbf{v}^N
 = \frac{1}{N!}\epsilon_{A_1\cdots A_N}\textbf{v}^{A_1}\wedge \cdots \wedge \textbf{v}^{A_N}\in \Lambda^NV^*
$ and
\[
 \begin{array}{ccll}
      \textbf{v}^{(N-1)}_A & = &
 \frac{1}{(N-1)!}\epsilon_{AA_2\cdots A_N}\textbf{v}^{A_2}\wedge \cdots
 \wedge \textbf{v}^{A_N} & \in \Lambda^{N-1}V^* \\
 \textbf{v}^{(N-2)}_{AB} & =  &
 \frac{1}{(N-2)!}\epsilon_{ABA_3\cdots A_N}\textbf{v}^{A_3}\wedge \cdots
 \wedge \textbf{v}^{A_N} & \in \Lambda^{N-2}V^*
 \\
 \textbf{v}^{(N-3)}_{ABC} & = &
 \frac{1}{(N-3)!}\epsilon_{ABCA_4\cdots A_N}\textbf{v}^{A_4}\wedge \cdots
 \wedge \textbf{v}^{A_N} & \in \Lambda^{N-3}V^*
         \end{array}
\]
A key observation is that
\begin{equation}\label{trmoinsunA}
 \textbf{v}^{(N-1)}_A:= \textbf{v}_A\iN \textbf{v}^{(N)},\quad
 \textbf{v}^{(N-2)}_{AB}:=\textbf{v}_B\iN \textbf{v}^{(N-1)}_A,\quad
 \textbf{v}^{(N-3)}_{ABC}:= \textbf{v}_C\iN \textbf{v}^{(N-2)}_{AB},
\end{equation}
from which we can easily deduce the following
\begin{equation}\label{thetatheta}
 \begin{array}{cclc}
 \textbf{v}^A\wedge \textbf{v}^{(N-1)}_{A'} & = & \delta^A_{A'}\textbf{v}^{(N)} & \hbox{(a)} \\
 \textbf{v}^A\wedge \textbf{v}^{(N-2)}_{A'B'} & = & \delta^A_{B'}\textbf{v}^{(N-1)}_{A'}
 - \delta^A_{A'}\textbf{v}^{(N-1)}_{B'} & \hbox{(b)} \\
 \textbf{v}^A\wedge \textbf{v}^{(N-3)}_{A'B'C'} & = & \delta^A_{C'}\textbf{v}^{(N-2)}_{A'B'}
 + \delta^A_{B'}\textbf{v}^{(N-2)}_{C'A'} + \delta^A_{A'}\textbf{v}^{(N-2)}_{B'C'} & \hbox{(c)} \\
    \textbf{v}^A\wedge \textbf{v}^B\wedge \textbf{v}^{(N-2)}_{A'B'} & = & \delta^{AB}_{A'B'}\textbf{v}^{(N)} & \hbox{(d)} \\
    \textbf{v}^A\wedge \textbf{v}^B\wedge \textbf{v}^{(N-3)}_{A'B'C'} & = &
    \delta^{AB}_{B'C'}\textbf{v}^{(N-1)}_{A'}
    +\delta^{AB}_{C'A'}\textbf{v}^{(N-1)}_{B'}
    + \delta^{AB}_{A'B'}\textbf{v}^{(N-1)}_{C'}
    & \hbox{(e)}
\end{array}
\end{equation}
where $\delta^{AB}_{A'B'}:= \delta^A_{A'}\delta^B_{B'} - \delta^A_{B'}\delta^B_{A'}$.
To prove (a) it suffices to
developp the relation
$0 = \textbf{v}_{A'}\iN 0
=\textbf{v}_{A'}\iN (\textbf{v}^A\wedge \textbf{v}^{(N)})$ and to use the graded Leibniz rule for the interior product.
Computing the interior product by $\textbf{v}_{B'}$
to both sides of (a) leads to (b) and
computing the interior product by $\textbf{v}_{C'}$
to both sides of (b) leads to (c). Then (d) follows from (a) and (b) and (e) follows from (b) and (c).
Lastly (\ref{fifi}) follows from (\ref{thetatheta}) by taking the pull-back by $e^V$, since $e^A = (e^V)^*\textbf{v}^A$, $e^{(N-1)}_A =  (e^V)^*\textbf{v}^{(N-1)}_A$, $e^{(N-2)}_{AB} =  (e^V)^*\textbf{v}^{(N-2)}_{AB}$ and $e^{(N-3)}_{ABC} = (e^V)^*\textbf{v}^{(N-3)}_{ABC}$.

Relations (\ref{deAB}) are easy consequence of the graded Leibniz rule for the exterior derivative.  
\hfill $\square$

\begin{lemm}\label{lemmadvarphiA}
 Let $e^V\in V\otimes \Omega^1(\mathcal{N})$ be a smooth frame over a manifold $\mathcal{N}$ and let $1\leq p\leq m-1$. Then
\begin{equation}\label{deVV}
 \emph{\dR} e^{(m-p)}_{V\cdots V}
  = \emph{\dR} e^{\underline{V}}\wedge e^{(m-p-1)}_{V\cdots V\underline{V}}
\end{equation}
(for instance $\emph{\dR} e^{(m-1)}_V = \emph{\dR} e^{\underline{V}}\wedge e^{(m-2)}_{V\underline{V}}$ and
$\emph{\dR} e^{(m-2)}_{VV} = \emph{\dR} e^{\underline{V}}\wedge e^{(m-3)}_{VV\underline{V}}$).
\end{lemm}
\emph{Proof} --- This relation is a translation of (\ref{deAB}). \hfill $\square$

\begin{lemm}\label{lemmeLeibniz}
Let $\gog$ be a Lie algebra and $\omega^\gog\in \gog\otimes \Omega^1(\mathcal{N})$.
Then $\emph{\dR}^\omega$ satisfies  the graded \emph{Leibniz rule} with respect to the contracted wedge product, which means the following.

Let $a,b\in \N^*$ and let $(V_1,\cdots,V_a)$,
$(W_1,\cdots,W_b)$ be two lists of vector spaces which are all linear \textbf{representations} of $\gog$. Let $c\in \N^*$ and $\sigma:[\![1,c]\!] \longrightarrow [\![1,a]\!]$ and
$\tau:[\![1,c]\!] \longrightarrow [\![1,b]\!]$ be two one-to-one maps. Let $p,q\in \N$. Then $\forall \beta\in  V_1\otimes\cdots\otimes V_a
\otimes \Omega^p(\mathcal{N})$, $\forall \gamma\in W_1\otimes\cdots\otimes W_b
\otimes \Omega^q(\mathcal{N})$,
\begin{equation}\label{twistedLeibniz}
  \emph{\dR}^\omega \emph{\textsf{C}}_{\sigma,\tau}(\beta \wedge \gamma)
  = \emph{\textsf{C}}_{\sigma,\tau}(\emph{\dR}^\omega \beta \wedge \gamma)
  + (-1)^p\emph{\textsf{C}}_{\sigma,\tau}(\beta \wedge \emph{\dR}^\omega \gamma)
\end{equation}
\end{lemm}
\emph{Proof} --- It is a consequence of the Leibniz rule for the exterior differential $\dR$ and of elementary properties of representations (\ref{adjointpairing}) and (\ref{adjointontensorproduct}).\hfill $\square$\\

For example let $\omega^\gog\in \gog\otimes \Omega^1(\mathcal{N})$, let $V$ be a vector space representation of $\gog$ and consider any
$\beta{^\gog}_{VV}\in
\gog\otimes V^* \otimes V^*\otimes \Omega^p(\mathcal{N})$ and any
$\gamma{^{VV}}_{\gog}
\in  V\otimes V\otimes\gog^*\otimes \Omega^q(\mathcal{N})$. Then
\[
 \begin{array}{ccc}
 \dR^\omega\left(\beta{^\gog}_{V\underline{V}}
 \wedge \gamma{^{V\underline{V}}}_{\gog} \right)
 & = & \left(\dR^\omega\beta{^\gog}_{V\underline{V}}
  \right)\wedge \gamma{^{V\underline{V}}}_{\gog}
 + (-1)^p \beta{^\gog}_{V\underline{V}}
 \wedge \left(\dR^\omega\gamma{^{V\underline{V}}}_{\gog} \right) \\
 \dR^\omega\left(\beta{^\gog}_{\underline{V}\underline{V}}
 \wedge \gamma{^{\underline{V}\underline{V}}}_{\gog} \right)
 & = & \left(\dR^\omega\beta{^\gog}_{\underline{V}\underline{V}}
  \right)\wedge \gamma{^{\underline{V}\underline{V}}}_{\gog}
 + (-1)^p \beta{^\gog}_{\underline{V}\underline{V}}
 \wedge \left(\dR^\omega\gamma{^{\underline{V}\underline{V}}}_{\gog} \right)
 \end{array}
\]

\begin{lemm}\label{dAegg}
 Let $\rho:\gog\longrightarrow gl(V)$ be a linear representation and assume it is \textbf{unimodular}, i.e. $\hbox{\emph{tr}}(\rho\xi) =
(\rho\xi)^B_{B} = 0$, $\forall \xi\in\gog$. Let
 $e^V\in V\otimes\Omega^1(\mathcal{N})$ and $\emph{\dR}^\omega:= \emph{\dR} + (\rho\omega)\wedge$. and let $1\leq p\leq m-1$. Then
\begin{equation}\label{nablaLeibniz}
 \emph{\dR}^\omega e^{(m-p)}_{V\cdots V}
  = \emph{\dR}^\omega e^{\underline{V}}\wedge e^{(m-p-1)}_{V\cdots V\underline{V}}
\end{equation}
(for instance $\emph{\dR}^\omega e^{(m-1)}_V = \emph{\dR} e^{\underline{V}}\wedge e^{(m-2)}_{V\underline{V}}$ and
$\emph{\dR}^\omega e^{(m-2)}_{VV} = \emph{\dR} e^{\underline{V}}\wedge e^{(m-3)}_{VV\underline{V}}$).
\end{lemm}
\emph{Proof} --- Let us prove, for instance, (\ref{nablaLeibniz}) for $p=2$. It amounts to prove $\dR^\omega  e^C\wedge e^{(m-3)}_{ABC} = \dR^\omega e^{(m-2)}_{AB}$. We use (\ref{fifi}) and (\ref{deAB}) in the follow computation
\[
 \begin{array}{ccl}
  \dR^\omega  e^C\wedge e^{(m-3)}_{ABC} & = & \left(\dR e^C + (\rho\omega)^C_{D}\wedge e^{D}\right)
  \wedge e^{(m-3)}_{ABC} \\
  & = & \dR e^C \wedge e^{(m-3)}_{ABC} + (\rho\omega)^C_{D}\wedge
  \left(\delta^{D}_A e^{(m-2)}_{BC}
  - \delta^{D}_B e^{(m-2)}_{AC} +\delta^{D}_C e^{(m-2)}_{AB} \right)\\
  & = & \dR e^{(m-2)}_{AB} + (\rho\omega)^C_A\wedge e^{(m-2)}_{BC}
  - (\rho\omega)^C_B\wedge  e^{(m-2)}_{AC}
  + (\rho\omega)^C_C\wedge e^{(m-2)}_{AB}
 \end{array}
\]
But since $\rho$ is unimodular, $(\rho\omega)^C_C = 0$ and thus, by permuting indices,
\[
 \dR^\omega e^C\wedge e^{(m-3)}_{ABC} = \dR e^{(m-2)}_{AB}
  - (\rho\omega)^{C}_A\wedge  e^{(m-2)}_{CB}
 - (\rho\omega)^{C}_B\wedge  e^{(m-2)}_{AC}
\]
which is the expression for $\dR^\omega e^{(m-2)}_{AB}$.\hfill $\square$

\subsubsection{Gauge transformations}
\begin{lemm}\label{lemme2point3}
 Let $g\in \mathscr{C}^\infty(\mathcal{N},\goG)$, $\emph{\hbox{R}}:\goG\longrightarrow GL(V)$ be a linear representation map of $\goG$ on a vector space $V$. Let $e^V, f^V\in V\otimes \Omega^1(\mathcal{N})$ such that
\begin{equation}\label{definitioneenfonctiondephi}
 e^V:= \emph{\hbox{R}}_gf^V
\end{equation}
then,
\begin{enumerate}
 \item by using Notation (\ref{deffifi})
\begin{equation}\label{eetf55}
 e^{VV} = \emph{\hbox{R}}_g\otimes \emph{\hbox{R}}_g\ f^{VV}
\end{equation}
\item if $\Phi^\gog, \Omega^\gog\in \gog\otimes\Omega^2(\mathcal{N})$ decompose as $\Phi^\gog = \frac{1}{2} \Phi{^\gog}_{\underline{V}\underline{V}}f^{\underline{V}\underline{V}}$ and
$\Omega^\gog = \frac{1}{2} \Omega{^\gog}_{\underline{V}\underline{V}}e^{\underline{V}\underline{V}}$, then
\begin{equation}\label{57}
  \hbox{\emph{Ad}}_g\Phi{^\gog}
  = \Omega{^\gog}
 \quad \Longleftrightarrow \quad
\hbox{\emph{Ad}}_g\otimes\emph{\hbox{R}}_g^*\otimes\emph{\hbox{R}}_g^*\left( \Phi{^\gog}_{VV}\right)
 = \Omega{^\gog}_{VV}
\end{equation}
\item if, furthermore, $\emph{\hbox{R}}:\goG\longrightarrow GL(V)$ is \textbf{unimodular}, then by using notations (\ref{defantififi}),
$e^{(m)}_V = f^{(m)}_V$ and
\begin{equation}\label{eetf56}
 e^{(m-1)}_{V} = \emph{\hbox{R}}_g^*f^{(m-1)}_{V},\quad
 e^{(m-2)}_{VV} = \emph{\hbox{R}}_g^*\otimes \emph{\hbox{R}}_g^* f^{(m-2)}_{VV},\quad
 e^{(m-3)}_{VVV} = \emph{\hbox{R}}_g^*\otimes\emph{\hbox{R}}_g^*\otimes\emph{\hbox{R}}_g^* f^{(m-3)}_{VVV}
\end{equation}
\item if $\emph{\hbox{R}}:\goG\longrightarrow GL(V)$ is \textbf{unimodular} and $\pi_\gog,p_\gog\in \Omega^{N-2}(\mathcal{N})\otimes \gog^*$ decompose as
$\pi_\gog = \frac{1}{2}\pi{_\gog}^{\celv\celv}f^{(N-2)}_{\celv\celv}$ and $p_\gog = \frac{1}{2}p{_\gog}^{\celv\celv}e^{(N-2)}_{\celv\celv}$, then
\begin{equation}\label{58}
 \hbox{\emph{Ad}}_g^*\pi_\gog
 = p_\gog
 \quad \Longleftrightarrow \quad
 \hbox{\emph{Ad}}_g^*\otimes\emph{\hbox{R}}_g\otimes\emph{\hbox{R}}_g\left( \pi{_\gog}^{VV}\right)
 = p{_\gog}^{VV}
\end{equation}
\end{enumerate}
\end{lemm}
\emph{Remark} --- If $f^V$ is a frame then $e^V$ is so and decompositions
$\Phi^\gog = \frac{1}{2} \Phi{^\gog}_{\celu\celu}f^{\celu\celu}$,
$\Omega^\gog = \frac{1}{2} \Omega{^\gog}_{\celu\celu}e^{\celu\celu}$,
$\pi_\gog = \frac{1}{2}\pi{_\gog}^{\celu\celu}f^{(N-2)}_{\celu\celu}$ and $p_\gog = \frac{1}{2}p{_\gog}^{\celu\celu}e^{(N-2)}_{\celu\celu}$ in (ii) and (iv)
are always possible.\\
\emph{Proof of Lemma \ref{lemme2point3}} --- The proofs of (i) and (iii) are straightforward. Assertion (ii) follows by using (\ref{Adjointpairing}) from
\[
\begin{array}{ccl}
 \hbox{Ad}_g\Phi{^\gog}
 & = & \frac{1}{2}\hbox{Ad}_g\left( \Phi{^\gog}_{\underline{V}\underline{V}}f^{\underline{V}\underline{V}}\right)
 = \frac{1}{2}\left(
 \hbox{Ad}_g\otimes\hbox{R}_g^*\otimes\hbox{R}_g^*\
 \Phi{^\gog}_{\underline{V}\underline{V}}\right)
 \left( \hbox{R}_g\otimes\hbox{R}_g\ f^{\underline{V}\underline{V}}\right) \\
 & = & \frac{1}{2}\left( \hbox{Ad}_g\otimes\hbox{R}_g^*\otimes\hbox{R}_g^*\ \Phi{^\gog}_{\underline{V}\underline{V}}\right)
 e^{\underline{V}\underline{V}}
\end{array}
\]
Assertion (iv) follows from $p_\gog = \frac{1}{2}p{_\gog}^{\underline{V}\underline{V}}e^{(N-2)}_{\underline{V}\underline{V}}$ and
\[
\begin{array}{ccl}
 \hbox{Ad}_g^*\pi_\gog & = & \frac{1}{2}\hbox{Ad}_g^*\left( \pi{_\gog}^{\underline{V}\underline{V}}f^{(N-2)}_{\underline{V}\underline{V}}\right)
 = \frac{1}{2}\left(\hbox{Ad}_g^*\otimes\hbox{R}_g\otimes\hbox{R}_g\ \pi{_\gog}^{\underline{V}\underline{V}}\right)\left(\hbox{R}_g^*\otimes\hbox{R}_g^*\ f^{(N-2)}_{\underline{V}\underline{V}}\right) \\
& = & \frac{1}{2}\left(\hbox{Ad}_g^*\otimes\hbox{R}_g\otimes\hbox{R}_g\ \pi{_\gog}^{\underline{V}\underline{V}}\right)e^{(N-2)}_{\underline{V}\underline{V}}
\end{array}
\]
\hfill $\square$

\begin{lemm}\label{mainlemma}
 Let $g\in \mathscr{C}^\infty(\mathcal{N},\goG)$ and $\theta^\gog,\omega^\gog\in \gog\otimes\Omega^1(\mathcal{N})$ such that
\begin{equation}\label{definitionAenfonctiondephi}
 \omega^\gog = \hbox{\emph{Ad}}_g\theta^\gog - \emph{\dR} g\,g^{-1}
\end{equation}
\begin{enumerate}
 \item then
 \begin{equation}
 \emph{\dR}\omega^\gog + \frac{1}{2}[\omega^\gog\wedge \omega^\gog]
 = \hbox{\emph{Ad}}_g\left(\emph{\dR}\theta^\gog + \frac{1}{2}[\theta^\gog\wedge \theta^\gog] \right)
\end{equation}
\item for any $\phi^\gog\in \gog\otimes \Omega^p(\mathcal{N})$,
\begin{equation}\label{novo25}
 \emph{\dR}^\omega\left(\hbox{\emph{Ad}}_g\phi^\gog\right) = \hbox{\emph{Ad}}_g\left(\emph{\dR}^\theta \phi^\gog \right)
\end{equation}
 \item for any $\pi_\gog\in \gog^*\otimes \Omega^p(\mathcal{N})$,
\begin{equation}\label{domegap=Adgdvarphivarpi}
  \emph{\dR}^\omega \left(\hbox{\emph{Ad}}_g^*\pi_\gog\right)
  = \hbox{\emph{Ad}}_g^*\left(\emph{\dR}^\theta\pi_\gog\right).
\end{equation}
\end{enumerate}
\end{lemm}
\emph{Proof} --- Result (i) is standard. The proof of (ii) is obtained as follows
\[
 \begin{array}{ccl}
  \dR^\omega\left(\hbox{Ad}_g\phi^\gog\right) & = &
  \dR \left( g\,\phi^\gog g^{-1}\right)
  + \hbox{ad}_{g\theta^\gog g^{-1}-\scriptsize{\dR} g\,g^{-1}}\wedge \left( g\,\phi^\gog g^{-1}\right) \\
  & = & \left[ \dR g\,g^{-1}\wedge g\,\phi^\gog g^{-1}\right] + g\,\dR \phi^\gog\,g^{-1}
  +g[\theta^\gog\wedge\phi^\gog]g^{-1}
  - \left[\dR g\,g^{-1}\wedge g\phi^\gog g^{-1}\right] \\
  & = & g\left(\dR \phi^\gog + [\theta^\gog\wedge \phi^\gog]\right)g^{-1}
 \end{array}
\]
We now deduce (iii).
Let $\pi_\gog\in \gog^*\otimes \Omega^p(\mathcal{N})$. Then for any $\phi^\gog\in \gog\otimes\Omega^q(\mathcal{N})$, by using (\ref{twistedLeibniz}) and (\ref{Adjointpairing}) and by applying (ii), i.e. (\ref{novo25}), to $\phi^\gog$, we obtain
\[
 \begin{array}{ccl}
  \left(\dR^\omega\hbox{Ad}_g^*\pi_\celg\right)\wedge \hbox{Ad}_g\phi^\celg
  & = & \dR^\omega\left(\hbox{Ad}_g^*\pi_\celg \wedge \hbox{Ad}_g\phi^\celg \right)
  - (-1)^p\hbox{Ad}_g^*\pi_\celg \wedge \left(\dR^\omega\hbox{Ad}_g\phi^\celg\right) \\
  & = & \dR^\omega\left(\pi_\celg \wedge \phi^\celg \right)
  - (-1)^p\hbox{Ad}_g^*\pi_\celg \wedge \hbox{Ad}_g\left(\dR^\theta\phi^\celg\right) \\
  & = & \dR^\theta\left(\pi_\celg  \wedge \phi^\celg \right)
  - (-1)^p\pi_\celg\wedge \left(\dR^\theta\phi^\celg \right)\\
  & = & \left(\dR^\theta\pi_\celg \right) \wedge \phi^\celg
  = \left(\hbox{Ad}_g^*\dR^\theta\pi_\celg \right) \wedge \hbox{Ad}_g\phi^\celg
 \end{array}
\]
Since this is true for any $\phi^\gog$, we deduce $\dR^\omega\hbox{Ad}_g^*\pi_\gog
= \hbox{Ad}^*_g\left(\dR^\theta\pi_\gog\right)$.\hfill $\square$

\begin{lemm}\label{lemma2point5}
Let $\omega^\gog,e^\gog \in \gog\otimes\Omega^1(\mathcal{N})$ and
$g\in \mathscr{C}^\infty(\mathcal{N},\goG)$ such that
\begin{equation}\label{relationaAdgg}
 e^\gog = \omega^\gog + \emph{\dR} g\ g^{-1}
\end{equation}
Then, by setting $\Omega^\gog:= \emph{\dR}\omega^\gog + \frac{1}{2}[\omega^\gog\wedge \omega^\gog]$,
 \begin{equation}\label{dAeagalFee}
  \emph{\dR}^\omega e^\gog = \Omega^\gog + \frac{1}{2}[e^\gog\wedge e^\gog],
 \end{equation}
\end{lemm}
\emph{Proof} --- This is a computation which uses $\dR\left( \dR g\ g^{-1}\right)
= \frac{1}{2}[\dR g\ g^{-1}\wedge \dR g\ g^{-1}]$
\[
 \begin{array}{ccl}
  \dR^\omega e^\gog & = & \hbox{d}e^\gog + [\omega^\gog\wedge e^\gog]
  = \dR\left(\dR g\ g^{-1} + \omega^\gog\right)  + [\omega^\gog\wedge (\dR g\ g^{-1} + \omega^\gog)] \\
  & = & \dR\left( \dR g\ g^{-1}\right) + \dR\omega^\gog + [\omega^\gog\wedge \dR g\ g^{-1}] + [\omega^\gog\wedge \omega^\gog] \\
  & = & \left(\frac{1}{2}[\dR g\ g^{-1}\wedge \dR g\ g^{-1}] + [\omega^\gog\wedge \dR g\ g^{-1}] + \frac{1}{2}[\omega^\gog\wedge \omega^\gog]\right) + \left(\dR\omega^\gog
  + \frac{1}{2}[\omega^\gog\wedge \omega^\gog]\right) \\
  & = & \frac{1}{2}[e^\gog\wedge e^\gog] +  \Omega^\gog
 \end{array}
\]
\hfill $\square$\\
\emph{Remark} --- Hypothesis (\ref{relationaAdgg}) occurs for instance if there exists some
$\theta^\gog\in \gog\otimes\Omega^1(\mathcal{N})$ such that
$\omega^\gog:= \hbox{Ad}_g\theta^\gog
 - \dR g\ g^{-1}$ and $e^\gog:= \hbox{Ad}_g\theta^\gog$.

\section{Gauge theories}\label{sectionGauge}
\subsection{General framework}
Assume we are given a vector space $\gos\simeq \R^n$,
endowed with a nondegenerate symmetric bilinear form $\textsf{b}$, and a smooth oriented pseudo Riemannian manifold $(\mathcal{X},\textbf{g})$
of dimension $n$,
such that, $\forall \textsf{x}\in \mathcal{X}$,
$(T_\textsf{x}\mathcal{X},\textsf{g}_\textsf{x})$ is isometric to $(\gos,\textsf{b})$.
In applications $(\gos,\textsf{b})$ will be either an Euclidean space (then $\mathcal{X}$ is Riemannian) or a Minkowski space and $(\mathcal{X},\textbf{g})$ a curved space or space-time. We fix a basis
$(E_1,\cdots, E_n)$ of $\gos$ and we set $\textsf{b}_{ab}:= \textsf{b}(E_a,E_b)$.

We are also given a compact (hence unimodular) Lie group $\goG$
of dimension $r$,
with Lie algebra $\gog$. We assume that
$\gog$ is endowed with a positive
$\hbox{Ad}_\goG$-invariant metric $\textsf{k}$, i.e.
such that $\textsf{k}(\hbox{Ad}_g\xi,\hbox{Ad}_g\zeta)
= \textsf{k}(\xi,\zeta)$, $\forall g\in \goG$,
$\forall \xi,\zeta\in \gog$. We let $(\textbf{t}_1,\cdots,\textbf{t}_r)$
be a basis of $\gog$ and $(\textbf{t}^1,\cdots,\textbf{t}^r)$ its dual basis of $\gog^*$. We set
\[
 N = n+r\quad \hbox{ and }\quad
 \gou:= \gos\oplus \gog
\]
A basis of $\gou$ is $(\textbf{u}_1,\cdots,\mathbf{u}_N
) = (E_1,\cdots, E_n,\textbf{t}_1,\cdots,\textbf{t}_r)$.

We are going to build a generalized gauge theory on $\mathcal{X}$ with group structure $\goG$, starting from a smooth submersion
$P:\mathcal{F}\longrightarrow\mathcal{X}$ with connected fibers over $\mathcal{X}$, where $\mathcal{P}$ is a smooth manifold of dimension $N$ (thus the dimension of the fibers is $r$).

The \emph{dynamical fields} of the problem are:
\begin{enumerate}
\item a $\gog$-valued 1-form $\theta^\gog$ on $\mathcal{F}$ such that, $\forall \textsf{x}\in \mathcal{X}$,
 the rank of the restriction of $\theta^\gog_\textsf{x}$ on the fiber
 $\mathcal{F}_\textsf{x}:= P^{-1}(\{\textsf{x}\})$ is
 equal to $r$ (thus $\theta^\gog$ induces
a connection on $\mathcal{F}$
in the general sense of Ehresmann);
\item a \emph{dual} $(N-2)$-form $\pi_\gog$ on $\mathcal{F}$ with
coefficients in $\gog^*$.
\end{enumerate}
We shall see that if $\theta^\gog$ is a classical solution of our dynamical equations, it will impose constraints on the geometry of $\mathcal{F}$. Hence the geometry of $\mathcal{F}$ is also a part of the dynamical variables, a similarity with General Relativity.
More precisely, assuming some generic hypotheses, any solution $(\theta^\gog,\pi_\gog)$ of the dynamical equations will define a $\goG$-principal bundle structure on  $\mathcal{F}$ and also a solution of the Yang--Mills system of equations on $\mathcal{X}$. One hypothesis will based on the following notion.
\begin{defi}\label{defigcomplete}
 Let $\gou$ be a vector space and $\gos,\gog\subset \gou$ be two vector subspaces such that $\gou = \gos \oplus \gog$. Let $\mathcal{Y}$ be a manifold of such that $\hbox{dim}\mathcal{N} = \hbox{dim}\gou$ and $\theta^\gou = \theta^\gos+\theta^\gog\in \gou\otimes \Omega^1(\mathcal{N})$ be a coframe.
 We say that $(\mathcal{Y},\theta^\gos,\theta^\gog)$ is \textbf{$\gog$-complete} if, for any continuous
map $v^\gog$ from $[0,1]$ to $\mathfrak{g}$ and for any
point $\emph{\textsf{y}}\in \mathcal{Y}$, there exists
an unique $\mathscr{C}^1$ map $\gamma:[0,1]\longrightarrow \mathcal{Y}$,
which is a solution of $(\gamma^*\theta^\gos)_t= 0$ and $(\gamma^*\theta^\gog)_t= v^\gog(t)\hbox{\emph{d}}t$, $\forall t\in [0,1]$, with the
initial condition $\gamma(0) = \emph{\textsf{y}}$.
\end{defi}

\subsubsection{Presentation of the model}\label{firstpresentationYM}

Working locally if necessary, we assume that there exists an
oriented \emph{orthonormal} coframe $\underline{\beta}^\gos$ on $\mathcal{X}$,
such that, in particular, $(\underline{\beta}^\gos)^*\textsf{b} = \textsf{g}$. We define the lifted forms $\beta^\gos:= P^*\underline{\beta}^\gos\in \Omega^1(\mathcal{F})$ which can be decomposed as $\beta^{\gos}:= \beta^aE_a$ and we let
$\beta^{(n)}:= \beta^1\wedge \cdots \wedge \beta^n\in \Omega^n(\mathcal{F})$.
Similarly $\theta^\gog\in \gog\otimes \Omega^1(\mathcal{F})$
decomposes as $\theta^\gog = \theta^i\textbf{t}_i$
and the $\gog^*$-valued $(N-2)$-form $\pi_\gog$ decomposes as $\pi{_\gog} = \pi{_i}\textbf{t}^i$. We set $\theta^{(r)}:= \theta^1\wedge \cdots \wedge \theta^r$.
We can consider $\beta^{\gos}$ and $\theta^\gog$ as the two components of the
1-form $f^\gou = f^A\textbf{u}_A =  \beta^{\gos} + \theta^\gog\in \gou\otimes \Omega^1(\mathcal{F})$
and we have $f^{(N)} = \beta^{(n)}\wedge \theta^{(r)}$.

The set of dynamical fields is
\begin{equation}\label{definitionSetEpourYM}
 \mathscr{E}:= \{ (\theta^\gog,\pi_\gog)\in \left(\gog\otimes\Omega^1(\mathcal{F})\right) \times \left(\gog^*\otimes\Omega^{N-2}(\mathcal{F})\right) \hbox{ of class }\mathscr{C}^2\, ;\, f^{(N)}\neq 0\}
\end{equation}
Observe that the condition $f^{(N)} = \beta^{(n)}\wedge \theta^{(r)}\neq 0$ ensures that $f^\gou = \beta^\gos+\theta^\gog$ is a coframe on $T^*\mathcal{F}$ and that $\hbox{rank}(\theta^\gog|_{\mathcal{F}_\textsf{x}}) = r$.
We denote by $(\frac{\partial}{\partial \beta^1},\cdots, \frac{\partial}{\partial \beta^n},\frac{\partial}{\partial \theta^1},
\cdots,\frac{\partial}{\partial \theta^r})$ its dual basis.
We also define $f^{(N-1)}_\gou$,  $f^{(N-2)}_{\gou\gou}$, $\beta^{(n-1)}_\gos$, $\beta^{(n-2)}_{\gos\gos}$, $\theta^{(r-1)}_\gog$
and $\theta^{(r-2)}_{\gog\gog}$ by following the rules in (\ref{defantififi}).
By applying the convention (\ref{ConventionV*})
we can decompose $\pi_\gog$ as
$\pi_\gog = \frac{1}{2}\pi{_\gog}^{\celu\celu}f_{\celu\celu}^{(n-2)}$.
By splitting $\pi{_\gog}^{\gou\gou} = \pi{_\gog}^{\gos\gos} + \pi{_\gog}^{\gos\gog} + \pi{_\gog}^{\gog\gos} + \pi{_\gog}^{\gog\gog}$, this gives also
\begin{equation}\label{decompositionofpiwithindices}
 \pi_\gog
 = \frac{1}{2}\pi{_\gog}^{\celu\celu}f_{\celu\celu}^{(n-2)}
 = \frac{1}{2}\pi{_\gog}^{\cels\cels}\beta_{\cels\cels}^{(n-2)} \wedge \theta^{(r)}
 - (-1)^n \pi{_\gog}^{\cels\celg}\beta_{\cels}^{(n-1)} \wedge \theta^{(r-1)}_\celg
 +\frac{1}{2}\pi{_\gog}^{\celg\celg}\beta^{(n)} \wedge \theta^{(r-2)}_{\celg\celg}
\end{equation}
The coefficient $\pi{_\gog}^{\gos\gos} = \pi{_i}^{ab}\textbf{t}^i\otimes E_a\otimes E_b
\in \gog^*\otimes \gos\otimes \gos$ which is also defined implicitely by
\begin{equation}\label{implicitpimumnu0}
 \pi{_\gog}\wedge \beta^\gos\wedge \beta^\gos = \pi{_\gog}^{\gos\gos}\beta^{(n)}\wedge \theta^{(r)}.
\end{equation}
plays a special role.
It defines the map
\begin{equation}\label{Qmap0}
 \begin{array}{cccc}
  \textsf{Q}{_\gog}^{\gos\gos}: &
  \mathscr{E} & \longrightarrow &
  \gog^*\otimes \gos\wedge\gos\otimes \mathscr{C}^\infty(\mathcal{F}) \\
  & (\theta^\gog,\pi_\gog) & \longmapsto &
  \pi{_\gog}^{\gos\gos}
 \end{array}
\end{equation}
We set
\[
 |\pi{_\gog}^{\gos\gos}|^2 =
 |\textsf{Q}{_\gog}^{\gos\gos}(\theta^\gog,\pi_\gog)|^2
 :=  \textsf{k}^{ij}\textsf{b}_{aa'} \textsf{b}_{bb'} \pi{_i}^{ab}\pi{_j}^{a'b'}
\]
or, by setting $\textsf{b}_{\gos\gos}:= \textsf{b}_{ab}E^a\otimes E^b\in \gos^*\otimes \gos^*$ and $\textsf{k}_{\gog\gog}:= \textsf{k}_{ij}\textbf{t}^i\otimes \textbf{t}^j\in \gog^*\otimes \gog^*$;
\begin{equation}\label{pigss}
  |\pi{_\gog}^{\gos\gos}|^2 :=
 \frac{1}{2}
 \pi{_\celg}^{\cels\cels}
 \ \pi{^\celg}_{\cels\cels}
 \quad \hbox{ where }
 \pi{^\gog}_{\gos\gos}
:=
\textsf{k}^{\gog\celg} \otimes\textsf{b}_{\gos\cels}\otimes
\textsf{b}_{\gos\cels}\left( \pi{_\celg}^{\cels\cels}\right)
\end{equation}
Lastly we define
\begin{equation}\label{actionYMnaivegeneral}
 \mathcal{A}[\theta^\gog,\pi{_\gog}]:=
 \int_\mathcal{F} \frac{1}{2}|\pi{_\gog}^{\gos\gos}|^2\beta^{(n)}\wedge \theta^{(r)}
 + \pi_\celg\wedge \left(\dR\theta{^\celg}+ \frac{1}{2}[\theta^\gog\wedge \theta^\gog]{^\celg}\right).
\end{equation}
We will prove the following result.
\begin{theo}\label{statementTheoYMgeneral}
Let $\gog$ a Lie algebra of dimension $r$.
 Let $\mathcal{F}$ and $\mathcal{X}$ be two smooth connected manifolds of dimensions $N:= n+r$ and $n$, respectively and $P:\mathcal{F}\longrightarrow \mathcal{X}$ be a smooth submersion with \textbf{connected} fibers.
 Consider the set $\mathscr{E}$ defined by (\ref{definitionSetEpourYM}).
 Assume that
 \begin{enumerate}
  \item either $\gog=u(1)\simeq \R$ and the fibers $\mathcal{F}_x:= P^{-1}(\{x\})$ are compact;
  \item or $\gog$ is the Lie algebra of a \textbf{compact, simply connected} Lie group $\widehat{\goG}$.
 \end{enumerate}
Let $(\theta^\gog,\pi_\gog)\in \mathscr{E}$ be a critical point of the functional (\ref{actionYMnaivegeneral}) and assume that $(\mathcal{F},\beta^\gos,\theta^\gog)$ is $\gog$-complete.
Then $(\theta^\gog,\pi_\gog)$
endows $P$ with a $\goG$-principal bundle structure, where $\goG$ is a compact connected Lie group. In case (i) $\goG = U(1)$, in case (ii) $\goG$ is a quotient of $\widehat{\goG}$ by a finite subgroup.

Moreover for any point in $\mathcal{X}$ there exist an open neighbourhood $\mathcal{O}$ of this point in $\mathcal{X}$ and a 
$\goG$-valued map $g$ defined on $\mathcal{O}$ such that, if $\mathbf{A}^\gog:= \hbox{\emph{Ad}}_g\theta^\gog - \hbox{\emph{d}}g\,g^{-1}$, $\mathbf{F}^\gog:= \hbox{\emph{d}}\mathbf{A}^\gog + \frac{1}{2}[\mathbf{A}^\gog\wedge \mathbf{A}^\gog]$ and $p{_\gog}^{\gog\gog}$ and $p{_\gog}^{\gog\gos}$ are the coefficients of $p_\gog:= \hbox{\emph{Ad}}^*_g\pi_\gog$ in the decomposition by using the coframe $e^\gog:= \hbox{\emph{Ad}}_g\theta^\gog$, then these fields are solutions of the system 
\begin{equation}\label{enonceYs}
\left\{
\begin{array}{rcl}
 \partial^{\gamma,\mathbf{A}}_{\cels} \mathbf{F}{_\gog}^{\gos\cels}
& = & 0 \quad \hbox{ (Yang--Mills)} \\
\partial^{\gamma,\mathbf{A}}_{\cels}  p{_\gog}^{\gog\cels} + \left(\partial_\celg p{_\gog}^{\gog\celg} + \frac{1}{2}\mathbf{c}{^\gog}_{\celg_1\celg_2}p{_\gog}^{\celg_1\celg_2}
\right)
& = & \frac{1}{2}|\mathbf{F}{_\gog}^{\gos\gos}|^2 \delta{_\gog}^\gog
- \frac{1}{2} 
 \mathbf{F}{_\gog}^{\cels_1\cels_2}\mathbf{F}{^{\gog}}_{\cels_1\cels_2}
 \end{array}\right.
\end{equation}
\end{theo}
Note that, in Case (i), where $\gog=u(1)$, the Yang--Mills system reduces to the Maxwell equations and the second equation in (\ref{enonceYs}) reduces to $\partial^{\gamma}_{\cels}p^\gos  = -\frac{1}{2}|\mathbf{F}^{\gos\gos}|^2$, where $p^\gos:= p{_\gog}^{\gog\gos}$.

The proofs of both cases follow similar key steps, although some arguments differ. As a warm up we first show Case (i) by assuming for simplicity that $\mathcal{X}$ is the flat Minkowski space $\gos$ of dimension 4, since it allows to get rid of unimportant details which can be fixed easily.
After introducing some extra notations,
we will then address Case (ii) in full generality. The crucial property that any compact Lie group is \textbf{unimodular} will used repeatedly.

\subsection{Study of the Maxwell case}
As announced we assume here that $\gog = u(1)=\R$ and $\mathcal{X} = \gos = \R^4$. Since the fibers of $\mathcal{F}\xrightarrow{P}\mathcal{X}$ are compact, connected and 1-dimensional they are all topologically equivalent to a circle. Hence
the manifold $\mathcal{F}$ is diffeomorphic to $\R^4\times S^1$. This allows us to choose global coordinates $(x^\mu,y) = (x^0,x^1,x^2,x^3,y)$, where $(x^0,x^1,x^2,x^3)\in \R^4$ and $y\in S^1\simeq \R/2\pi\Z$. We can then choose the coframe $(\beta^0,\beta^1,\beta^2,\beta^3)$ to be equal to $(\dR x^0,\dR x^1,\dR x^2,\dR x^3)$. We set
$\dR x^{(4)}:= \dR x^0\wedge \dR x^1\wedge \dR x^2\wedge \dR x^3$, $\dR x^{(3)}_\mu:= \frac{\partial}{\partial x^\mu}\iN \dR x^{(4)}$ and
$\dR x^{(3)}_{\mu\nu}:= \frac{\partial}{\partial x^\nu}\iN \dR x^{(4)}_\mu$.

Obviously we can identify $\gog^*\simeq \gog$ and the
metric $\textbf{h}$ with the standard
metric on $\R$. We can also dropp the index $\gog$ in $\theta^\gog$,  $\pi_\gog$ and $\textsf{Q}{_\gog}^{\gos\gos}$.
The set (\ref{definitionSetEpourYM}) reads here
\[
\mathscr{E}_{\textsf{Maxwell}}
 = \{ (\theta,\pi);\ \theta\in \Omega^1(\mathcal{F}),\pi\in \Omega^3(\mathcal{F}),
 \dR x^{(4)}\wedge \theta \neq 0\}
\]
Condition $\dR x^{(4)}\wedge \theta \neq 0$ means that, if we decompose
\[
\theta= \theta_0\dR x^0 + \theta_1\dR x^1 +
\theta_2\dR x^2 + \theta_3\dR x^3 + \theta_4\dR y,
\]
where the coefficients $\theta_A$ are functions of $(x^\mu,y)$,
then $\theta_4$ does not vanish.
Without loss of generality (since $\mathcal{F}$ is connected) we assume that
$\theta_4>0$.
The 3-form $\pi$ decomposes \emph{a priori} as
\[
 \pi = \frac{1}{2}\pi^{\mu\nu}\dR x^{(2)}_{\mu\nu}\wedge \theta
 - \pi^\mu \dR x^{(3)}_\mu
\]
which implies (see (\ref{thetatheta})) $\dR x^\mu\wedge \dR x^\nu\wedge\pi =
\pi^{\mu\nu}\dR x^{(4)}\wedge\theta$.
The quantity $|\pi^{\gos\gos}|^2$ reads
$|\pi^{\gos\gos}|^2 = \frac{1}{2}\textsf{b}_{\mu\mu'}\textsf{b}_{\nu\nu'}
 \pi^{\mu\nu}\pi^{\mu'\nu'}
 = \frac{1}{2}\pi^{\mu\nu}\pi_{\mu\nu}$,
where $\pi_{\mu\nu}:= \textsf{b}_{\mu\mu'}\textsf{b}_{\nu\nu'}\pi^{\mu'\nu'}$,
and the action is
\[
 \mathcal{A}[\theta,\pi]
 = \int_\mathcal{F}\frac{1}{2}
 |\pi^{\gos\gos}|^2\dR x^{(4)}\wedge \theta + \pi\wedge \dR\theta.
\]
The 'curvature' 2-form is simply $\Theta:= \dR\theta$, which we decompose as
\[
 \dR\theta = \Theta = \frac{1}{2}\Theta_{\mu\nu}\dR x^\mu\wedge \dR x^\nu
 + \Theta_\mu \dR x^\mu\wedge \theta.
\]
Hence by using (\ref{thetatheta}),
$\pi\wedge \dR\theta = \left(\frac{1}{2}\Theta_{\mu\nu}\pi^{\mu\nu} + \Theta_{\mu}\pi^\mu\right) \dR x^{(4)}\wedge \theta$.

\subsubsection{Study of the first variation}\label{paragraphMaxwellUl}

\noindent
\textbf{First variation with respect to $\pi$} ---
We write that the action is stationary with respect to variations $(\theta,\pi)\longmapsto (\theta,\pi+\varepsilon \delta\pi)$, for $\varepsilon$ small. This means that $\delta\theta = 0$ and the variations of $\pi$ are induced by the variations  $\delta\pi_{\mu\nu}$ and $\delta\pi^\mu$ of, respectively, $\pi_{\mu\nu}$ and $\pi^\mu$. We obtain straightforwardly (note that $\frac{1}{4}\pi_{\mu\nu}\pi^{\mu\nu}$
is quadratic in $\pi$, whereas $\Theta_{\mu\nu}\pi^{\mu\nu}$ is linear)
\begin{equation}\label{ELvarpi}
\left\{ \begin{array}{cccc}
         \pi_{\mu\nu} + \Theta_{\mu\nu} & = & 0 & \hbox{(a)} \\
         \Theta_{\mu} & = & 0 & \hbox{(b)}
        \end{array}\right.
\end{equation}
(equivalentely
$\frac{\partial \theta_\nu}{\partial x^\mu}
- \frac{\partial \theta_\mu}{\partial x^\nu} =
-\pi_{\mu\nu}$
and
$\frac{\partial \theta_4}{\partial x^\mu}
= \frac{\partial \theta_\mu}{\partial y}$).
Equation (b) means that $\frac{\partial}{\partial y}\iN \dR\theta = 0$ and has the following consequence: let $\mathcal{F}_{x_1}$ and $\mathcal{F}_{x_2}$
be two fibers over $x_1$ and $x_2\in \R^4$ respectively. Both are diffeomorphic to the circle $S^1$. Consider a path $\Gamma$ joining $x_1$ to $x_2$ in $\R^4$.
Its lift $\mathcal{S}:= P^{-1}(\Gamma)$ is a surface (having the topology of a cylinder)
the boundary of which is $\partial \mathcal{S} = \mathcal{F}_{x_2} - \mathcal{F}_{x_1}$
(choosing the orientation in an appropriate way). Thus
\begin{equation}\label{quantizationMaxwell}
  \int_{\mathcal{F}_{x_2}}\theta - \int_{\mathcal{F}_{x_1}}\theta
 = \int_{\partial\mathcal{S}}\theta = \int_{\mathcal{S}}\dR\theta = 0,
\end{equation}
where we have used $\dR\theta|_\mathcal{S} = 0$, because $\frac{\partial}{\partial y}$ is tangent to
$\mathcal{S}$ and $\frac{\partial}{\partial y}\iN \dR\theta = 0$. Since $\R^4$
is connected, this leads to a normalization of the fibers:
$\exists q\in (0,+\infty)$ such that
\[
 q = \int_{\mathcal{F}_x}\theta, \quad \forall x\in\R^4.
\]
Thus we can thus define a map
$f:\mathcal{F}\longrightarrow \R/q\Z$ such that
$\forall x\in \R^4$, $\dR f|_{\mathcal{F}_x} = \theta|_{\mathcal{F}_x}$, i.e. $\frac{\partial f}{\partial y} = \theta_4$, by
setting e.g.\footnote{One may as well
define $f$ by $f(x,y) = \int_0^y\theta_4(x,\sigma(x)+y')dy'$, where $\sigma:\R^4\longrightarrow \R/2\pi\Z$ is any section of $\mathcal{F}$.}
$f(x,y) = \int_0^y\theta_4(x,y')dy'$.
Then the map
\[
 \begin{array}{cccc}
  T: & \mathcal{F} & \longrightarrow & \R^4\times (\R/q\Z) \\
  & (x,y) & \longmapsto & (x,f(x,y))
 \end{array}
\]
is a diffeomorphism. We denote by $(x^\mu,s)$ coordinates on $\R^4\times (\R/q\Z)$.
Moreover
\[
 \theta= \left(\theta_\mu - \frac{\partial f}{\partial x^\mu}\right)\dR x^\mu + \dR f.
\]
and hence, by setting
\begin{equation}\label{AparTmoinsun}
 \textbf{A}_\mu:=
  \left(\theta_\mu - \frac{\partial f}{\partial x^\mu}\right)\circ T^{-1},
  \quad \hbox{ for }0\leq \mu\leq 3,
\end{equation}
and $\textbf{A}:= \textbf{A}_\mu \dR x^\mu$, we have
\begin{equation}\label{varphi=ds}
\theta= (\textbf{A}_\mu\circ T) \dR x^\mu
+ \dR f
=  T^*\left(\textbf{A} + \dR s\right).
\end{equation}
In particular $\textbf{A}+\dR s$ is \emph{normalized} (i.e. $\frac{\partial}{\partial s}\iN (\textbf{A}+\dR s) = 1$).

Moreover since
$T^*\left( T_*\frac{\partial}{\partial y}\iN \dR\textbf{A}\right)
= \frac{\partial}{\partial y}\iN T^*\dR\textbf{A}
= \frac{\partial}{\partial y}\iN
\dR\theta = 0$ by (\ref{ELvarpi}) and
$T_*\frac{\partial}{\partial y} =
(\theta_4\circ T^{-1})\frac{\partial}{\partial s}$, (\ref{ELvarpi}, b) translates as
\[
 \frac{\partial}{\partial s}\iN \dR\textbf{A} = 0
\]
Since we have obviously
$\frac{\partial}{\partial s}\iN \textbf{A} = 0$
we also get that
$L_{\frac{\partial}{\partial s}}\textbf{A}
 = \dR\left(\frac{\partial}{\partial s}\iN \textbf{A}\right)
 + \frac{\partial}{\partial s}\iN \dR\textbf{A}
 = 0$, i.e. $\frac{\partial \textbf{A}_\mu}{\partial s} = 0$,
$\forall \mu$, i.e. $\textbf{A}_\mu$ is a function of $x\in \R^4$ only.

Lastly we define $\textbf{F}:= \dR\textbf{A}$, so that $\Theta= \dR\theta = T^*\textbf{F}$
and we deduce from the previous results the
decomposition
\begin{equation}\label{FmunuMaxwell}
 \textbf{F} = \frac{1}{2}\textbf{F}_{\mu\nu}
 \dR x^\mu\wedge \dR x^\nu
\end{equation}
where the coefficients $\textbf{F}_{\mu\nu}$
are functions of $x\in \R^4$ only.
Equation (a) in (\ref{ELvarpi}) translates then as
$\textbf{F}_{\mu\nu} =
\frac{\partial \textbf{A}_\nu}{\partial x^\mu}
- \frac{\partial \textbf{A}_\mu}{\partial x^\nu} =
-p_{\mu\nu}$. \\

\noindent
\textbf{First variation with respect to $\theta$} ---
Here we write the condition for the action to be stationary with respect to variations $(\theta,\pi)\longmapsto (\theta+\varepsilon\delta \theta,\pi)$, for $\varepsilon$ small (hence $\delta\pi=0$).
We decompose
\[
 \delta\theta = \tau_\mu \dR x^\mu + \tau\theta
\]
which induces the variation $\delta(\dR x^{(4)}\wedge \theta) =
\tau \dR x^{(4)}\wedge \theta$.
Since $\delta(\dR x^\mu\wedge \dR x^\nu\wedge\pi) = 0$ and we must respect the constraint
$\dR x^\mu\wedge \dR x^\nu\wedge\pi =
\pi^{\mu\nu}\dR x^{(4)}\wedge\theta$,
this forces to have
\[
 0 = \delta\left(\pi^{\mu\nu}\dR x^{(4)}\wedge\theta\right)
 = \delta\pi^{\mu\nu}\dR x^{(4)}\wedge\theta
 + \pi^{\mu\nu}\dR x^{(4)}\wedge\delta\theta
 = \left(\delta\pi^{\mu\nu} + \tau \pi^{\mu\nu}\right)
 \dR x^{(4)}\wedge\theta
\]
Hence we must impose $\delta\pi^{\mu\nu} + \tau \pi^{\mu\nu} = 0$. The
induced variations on $|\pi^{\gos\gos}|^2$ is
$\delta |\pi^{\gos\gos}|^2 = -2\tau|\pi^{\gos\gos}|^2$. Hence
\[
 \delta\left(\frac{1}{2}|\pi^{\gos\gos}|^2\dR x^{(4)}\wedge \theta\right)
 = -\frac{1}{2}|\pi^{\gos\gos}|^2\lambda \dR x^{(4)}\wedge \theta
 = -\frac{1}{2}|\pi^{\gos\gos}|^2\delta\theta\wedge \dR x^{(4)}
\]
Moreover $\delta(\pi\wedge \dR\theta) = \dR(\delta\theta)\wedge \pi
= \dR(\delta\theta\wedge \pi) + \delta\theta\wedge \dR\pi$, hence
the vanishing of the first variation of $\mathcal{A}$ leads to
\[
0 = \int_\mathcal{F}\dR(\delta\theta\wedge \pi) +
 \delta\theta\wedge\left( \dR\pi -\frac{1}{2}|\pi^{\gos\gos}|^2\dR x^{(4)}\right),
\quad \forall \delta\theta
\]
i.e., if $\delta\theta$ has compact support,
\begin{equation}\label{ELvarphi}
 \dR\pi = \frac{1}{2}|\pi^{\gos\gos}|^2\dR x^{(4)}
\end{equation}
By using (\ref{varphi=ds}) we can write (see (\ref{thetatheta}))
\[
 \begin{array}{ccl}
  \pi & = & \frac{1}{2}\pi^{\mu\nu}\dR x^{(2)}_{\mu\nu}\wedge
  \left((\textbf{A}_\lambda\circ T) \dR x^\lambda+\dR f\right)  - \pi^\mu \dR x^{(3)}_\mu \\
 & = & \pi^{\mu\nu}(\textbf{A}_\nu\circ T) \dR x^{(3)}_\mu
  + \frac{1}{2}\pi^{\mu\nu}\dR x^{(2)}_{\mu\nu}\wedge \dR f
  - \pi^\mu \dR x^{(3)}_\mu \\
  & = & \frac{1}{2}\pi^{\mu\nu}\dR x^{(2)}_{\mu\nu}\wedge (T^*\dR s)
  + \left( \pi^{\mu\nu}(\textbf{A}_\nu\circ T) - \pi^\mu\right) \dR x^{(3)}_\mu
 \end{array}
\]
thus, by defining $p^{\mu\nu}$ and $p^\mu$ such that
$p^{\mu\nu}\circ T:= \pi^{\mu\nu}$, $p^\mu\circ T:= \pi^\mu
 - (\pi^{\mu\nu})(\textbf{A}_\nu\circ T)$
and $p:= \frac{1}{2}p^{\mu\nu}\dR x^{(2)}_{\mu\nu}\wedge \dR s - p^\mu \dR x^{(3)}_\mu$,
we obtain
\[
 \pi = T^*p
 = T^*\left(\frac{1}{2}p^{\mu\nu}\dR x^{(2)}_{\mu\nu}\wedge \dR s
- p^\mu \dR x^{(3)}_\mu \right).
\]
Then $\dR\pi = T^*\dR p$ with
\[
 \dR p = \dR\left(\frac{1}{2}p^{\mu\nu}\dR x^{(2)}_{\mu\nu}\wedge \dR s
 - p^\mu \dR x^{(3)}_\mu\right)
 = \frac{1}{2}\dR p^{\mu\nu} \wedge \dR x^{(2)}_{\mu\nu}\wedge \dR s
 - \dR p^\mu \wedge \dR x^{(3)}_\mu
\]
Thus setting $\dR p^{\mu\nu} = \partial_{\lambda} p^{\mu\nu}\dR x^\lambda
 + \partial_sp^{\mu\nu}\dR s$ and
$\dR p^\mu = \partial_{\lambda} p^\mu \dR x^\lambda
 + \partial_s p^\mu \dR s$, we get
\[
 \begin{array}{ccl}
  \dR p & = & \frac{1}{2}\left(\partial_{\lambda} p^{\mu\nu}\dR x^\lambda
  + \partial_s p^{\mu\nu} \dR s \right) \wedge \dR x^{(2)}_{\mu\nu}\wedge \dR s
  - \left(\partial_{\lambda} p^\mu \dR x^\lambda + \partial_s p^\mu \dR s\right)
  \wedge \dR x^{(3)}_\mu \\
  & = & \partial_\nu p^{\mu\nu}\dR x^{(3)}_\mu\wedge \dR s
  - \partial_\mu p^\mu \dR x^{(4)} - \partial_s p^\mu \dR s\wedge \dR x^{(3)}_\mu \\
  & = & \left(\partial_\nu p^{\mu\nu} + \partial_sp^\mu\right)\dR x^{(3)}_\mu\wedge \dR s
  - \partial_\mu p^\mu \dR x^{(4)}
 \end{array}
\]
We also note that
$p^{\mu\nu}\circ T = \pi^{\mu\nu}$
implies $T^*\left(\frac{1}{2}|p^{\gos\gos}|^2\dR x^{(4)}\right) = \frac{1}{2}|\pi^{\gos\gos}|^2\dR x^{(4)}$.
Hence (\ref{ELvarphi}) reads
$T^*\dR p = T^*\left(\frac{1}{2}|p^{\gos\gos}|^2\dR x^{(4)}\right)$, which is equivalent to
$\dR p = \frac{1}{2}|p^{\gos\gos}|^2\dR x^{(4)}$.
In view of the previous computations, this
is equivalent to the system
\begin{equation}\label{ELvarphi2}
\left\{ \begin{array}{cccc}
         \partial_{\nu} p^{\mu\nu}  & = & -\partial_sp^\mu & \hbox{(a)} \\
        \partial_{\mu} p^\mu & = &
        - \frac{1}{2}|p^{\gos\gos}|^2 & \hbox{(b)}
        \end{array}\right.
\end{equation}

\subsubsection{Cancellation of the sources}
We deduced from (\ref{ELvarpi},a) that
$\textbf{F}^{\mu\nu}
:= \textsf{b}^{\mu\mu'}\textsf{b}^{\nu\nu'}\textbf
F_{\mu'\nu'} = - p^{\mu\nu}$. However we also deduced from
(\ref{FmunuMaxwell}) that the coefficients
$\textbf{F}^{\mu\nu}$ are functions of
$x \in \R^4$ only.
Hence we deduce by averaging both sides of
(\ref{ELvarphi2}a)
over a fiber $\mathcal{F}_x$ that
\[
 \partial_{\nu} \textbf{F}^{\mu\nu} = \frac{\int_{\mathcal{F}_x}\partial_{\nu} \textbf{F}^{\mu\nu}\dR s}{\int_{\mathcal{F}_x}\dR s}
 = \frac{\int_{\mathcal{F}_x}-\partial_{\nu} p^{\mu\nu}\dR s}{\int_{\mathcal{F}_x}\dR s}
 =  \frac{\int_{\mathcal{F}_x}\partial_sp^\mu \dR s}{\int_{\mathcal{F}_x}\dR s}
 =  \frac{\int_{\mathcal{F}_x}\dR p^\mu}{\int_{\mathcal{F}_x}\dR s} = 0
\]
and we conclude that the Maxwell equation in vacuum holds
\begin{equation}\label{Maxwell}
 \frac{\partial\textbf{F}^{\mu\nu}}{\partial x^\nu} = 0.
\end{equation}

\subsubsection{Gauge symmetries}
We consider the transformation:
\begin{equation}\label{hypotheticsymmetry}
 (\theta,\pi)\longmapsto (\theta + \alpha,\pi+\psi)
\end{equation}
and we look at sufficient conditions for this transformation to provide us with
a gauge symmetry of the action
$\mathcal{A}[\theta,\pi] =
 \int_\mathcal{F} \pi\wedge \dR\theta +
 \frac{1}{4}|\pi^{\gos\gos}|^2_{\R^4}\dR x^{(4)}\wedge \theta$. We have the \emph{a priori}
decompositions
$\alpha = \alpha_\mu(x,y) \dR x^\mu + \alpha_4(x,y)\dR y$
and
$\psi = \frac{1}{2}\psi^{\mu\nu}(x,y)\dR x_{\mu\nu}^{(2)} \wedge \dR y
- \psi^\mu(x,y) \dR x_\mu^{(3)}$.
In order to keep the quantity
$|\pi^{\gos\gos}|^2_{\R^4}:= \frac{1}{2} \textsf{b}_{\mu\mu'}\textsf{b}_{\nu\nu'}
\pi^{\mu\nu}\pi^{\mu'\nu'}$ invariant, we assume that the coefficients $\psi^{\mu\nu}$ vanish, so that
$\psi = - \psi^\mu(x,y) \dR x_\mu^{(3)}$.

Then the computation of $\mathcal{A}[\theta+\alpha,\pi+\psi]$
gives us
\[
 \mathcal{A}[\theta+\alpha,\pi+\psi]
 = \mathcal{A}[\theta,\pi]
 + \int_\mathcal{F}(\pi+\psi)\wedge \dR\alpha + \psi\wedge \dR\theta
 + \frac{1}{2}|\pi^{\gos\gos}|^2_{\R^4}\dR x^{(4)}\wedge \alpha
\]
We note that
$\dR x^{(4)}\wedge \alpha = \alpha_4\dR x^{(4)}\wedge \dR y$, thus,
in order for the last term in the r.h.s. to cancel, we
need to assume $\alpha_4=0$.
Hence $\alpha = \alpha_\mu(x,y) \dR x^\mu$.
Then we observe that we need to require that $\dR\alpha = 0$
for $(\pi+\psi)\wedge \dR\alpha$ to vanish and, if so, we need
to assume that $\int_\mathcal{F}\psi\wedge \dR\theta = 0$
for having
$\mathcal{A}[\theta+\alpha,\pi+\psi]
 = \mathcal{A}[\theta,\pi]$.

For that purpose we assume that $\psi$ has compact support or decays at infinity so that
\[
 \int_\mathcal{F}\psi\wedge \dR\theta
 = \int_\mathcal{F}\dR(\theta\wedge \psi) + \theta\wedge \dR\psi
 = \int_\mathcal{F} \theta\wedge \dR\psi
\]
Then it suffices to choose $\psi$ so that $\dR\psi = 0$ for
(\ref{hypotheticsymmetry}) to be a symmetry of $\mathcal{A}$. Hence, to summarize, if
\begin{enumerate}
 \item $\alpha = \alpha_\mu(x,y) \dR x^\mu\in \Omega^1(\mathcal{F})$ is \emph{closed};
 \item $\psi = - \psi^\mu(x,y) \dR x_\mu^{(3)} \in \Omega^3(\mathcal{F})$ is \emph{closed} and \emph{decays at infinity},
\end{enumerate}
then $\mathcal{A}[\theta+\alpha,\pi+\psi]
 = \mathcal{A}[\theta,\pi]$.

However since $\dR\alpha = \frac{1}{2}
 (\frac{\partial \alpha_\nu}{\partial x^\mu}
 - \frac{\partial \alpha_\mu}{\partial x^\nu})\dR x^\mu\wedge \dR x^\nu
 - \frac{\partial \alpha_\mu}{\partial y}\dR x^\mu\wedge \dR y$
and
$\dR\psi = -\frac{\partial \psi^\mu}{\partial x^\mu}\dR x^{(4)}
+ \frac{\partial \psi^\mu}{\partial y}\dR x^{(3)}_\mu \wedge \dR y$,
the previous conditions imply that coefficients $\alpha_\mu$
and $\psi^\mu$ are independant of $y$. Hence
$\alpha = \alpha_\mu(x) \dR x^\mu$ and
$\psi = - \psi^\mu(x) \dR x_\mu^{(3)}$, with
\[
 \frac{\partial \alpha_\nu}{\partial x^\mu}
 - \frac{\partial \alpha_\mu}{\partial x^\nu} = 0
 \quad\hbox{ and }\quad
 \frac{\partial \psi^\mu}{\partial x^\mu} = 0.
\]
The first equation is equivalent to the existence of a function $V\in \mathscr{C}^\infty(\R^4)$ such that $\alpha = \dR V$.

\subsubsection{Invariance by fiber bundle diffeomorphisms}
Let us consider a diffeomorphism $T:\mathcal{F}\longrightarrow \mathcal{F}$ such that $P\circ T = P$, i.e. of the form
\[
 \begin{array}{cccc}
  T: & \mathcal{F} & \longrightarrow
  & \mathcal{F} \\
  & (x,y) & \longmapsto & (x,f(x,y))
 \end{array}
\]
and such that $\frac{\partial f}{\partial y}>0$.
It acts on the fields by pull-back
\[
 (\theta,\pi) \longmapsto (T^*\theta,T^*\pi)
\]
We note that $\textsf{Q}^{\gos\gos}[T^*\theta,T^*\pi] = \textsf{Q}^{\mu\nu}[T^*\theta,T^*\pi]E_\mu\otimes E_\nu$
is defined implicitely by using (\ref{implicitpimumnu0}), i.e.
\[
(T^*\pi)\wedge \dR x^\mu\wedge \dR x^\nu = \textsf{Q}^{\mu\nu}[T^*\theta,T^*\pi]\ \dR x^{(n)}\wedge T^*\theta.
\]
On the other hand the pull-back by $T$ of both sides of the relation
$\pi\wedge \dR x^\mu\wedge \dR x^\nu = \textsf{Q}^{\mu\nu}[\theta,\pi] \dR x^{(n)}\wedge \theta$ gives us
\[
 (T^*\pi)\wedge \dR x^\mu\wedge \dR x^\nu = \left(\textsf{Q}^{\mu\nu}[\theta,\pi]\circ T\right) \dR x^{(n)}\wedge T^*\theta
\]
By comparing both relations we deduce that $\textsf{Q}^{\gos\gos}[T^*\theta,T^*\pi]
= \textsf{Q}^{\gos\gos}[\theta,\pi]\circ T$.
This implies that
$|\pi^{\gos\gos}|^2_{\R^4}$ is transformed into $|\pi^{\gos\gos}|^2_{\R^4}\circ T$. Thus
the Lagrangian density transforms as
\[
 \pi\wedge \dR\theta
 + \frac{1}{2}\ |\pi^{\gos\gos}|^2_{\R^4}\dR x^{(4)}\wedge \theta
   \longmapsto
  T^*\left(\pi\wedge \dR\theta
 + \frac{1}{2}|\pi^{\gos\gos}|^2_{\R^4}\dR x^{(4)}\wedge \theta\right)
\]
Hence the action
$\mathcal{A}[\theta,\pi] = \int_\mathcal{F}
\frac{1}{2}|\pi^{\gos\gos}|^2_{\R^4}\dR x^{(4)}\wedge \theta
 + \pi\wedge \dR\theta$ is invariant by this transformation.

This invariance by fiber bundle diffeomorphisms may
be fixed as follows.
Consider some $(\theta,\pi)\in \mathscr{E}_{\textsf{Maxwell}}$
 and, for any $x\in \R^4$, let
$u(x):= \frac{1}{2\pi} \int_{\mathcal{F}_x}\theta
$ and
$f(x,y):= \frac{1}{u(x)}\int_0^y\theta_4(x,y')dy'$ and
define the map
\[
 \begin{array}{cccc}
  T: & \mathcal{F} & \longrightarrow & \R^4\times (\R/2\pi\Z) \\
  & (x,y) & \longmapsto & \left(x,f(x,y)\hbox{ mod }[2\pi]\right),
 \end{array}
\]
which is a diffeomorphism.
Then
$\dR f = \frac{\partial f}{\partial x^\mu}\dR x^\mu
+ \frac{\theta_4(x,y)}{u(x)}\dR y$
and thus
\[
 \theta = \theta_\mu \dR x^\mu + \theta_4\dR y
 = \left(\theta_\mu - u\frac{\partial f}{\partial x^\mu}\right)\dR x^\mu + u\,\dR f.
\]
Hence by defining
$\phi_\mu:=
  \left(\theta_\mu - u\frac{\partial f}{\partial x^\mu}\right)\circ T^{-1}$
and $\phi:= \phi_\mu \dR x^\mu + u\,\dR s$ and by observing that $u\circ T=u$, we have
\begin{equation}\label{gaugediffeocercle}
T^*\phi
= (\phi_\mu\circ T)\dR x^\mu + (u\circ T)\dR f
= \theta
\end{equation}
Thus the image of the transformation
$(\theta,\pi)\longmapsto ((T^{-1})^*\theta,(T^{-1})^*\pi)$
is $(\phi,p)$, so that
$\phi$ has the form $\phi = \phi_\mu \dR x^\mu + \phi_4 \dR s$,
where $\phi_4(x,s) = u(x)$ is independant on $s$.

This show that, by such a 'gauge transform', which does not change the action as seen in the previous paragraph, we can
assume that the coefficient $\frac{\partial}{\partial y}\iN \theta$ is independant of the coordinate on the fiber.

\subsection{Study of the Yang--Mills case}\label{backtoYM}

We now prove Theorem \ref{statementTheoYMgeneral} in Case (ii), i.e. for $\gog$ which is the Lie algebra of a compact, simply connected structure group $\widehat{\goG}$ and on a curved base pseudo Riemannian manifold $(\mathcal{X},\textbf{g})$. Recall that, since $\widehat{\goG}$ is compact, its Lie algebra $\gog$ is unimodular.
We endow $\gou:= \gos\oplus \gog$  with the metric $\textsf{h}$ such that its restriction to $\gos$ coincides with $\textsf{b}$, its restriction to $\gog$ coincides with $\textsf{k}$ and $\gos\perp \gog$. We let $\beta^\gos$ be a $\textbf{g}$-orthonormal coframe and we set $f^\gou:= \beta^\gos + \theta^\gog$. Note that, by hypothesis, $f^\gou$ is a coframe on $\mathcal{F}$.

Abusing notation we denote by $\hbox{Ad}:\widehat{\goG}\longrightarrow GL(\gou)$ and $\hbox{ad}:\gog \longrightarrow gl(\gou)$ the representations which extends trivially, respectively, the adjoint representations $\hbox{Ad}:\widehat{\goG}\longrightarrow GL(\gog)$ and $\hbox{ad}:\gog \longrightarrow gl(\gog)$, i.e. such that: $\forall g\in \widehat{\goG},\forall \xi\in \gog$,
\begin{equation}\label{definitionActionderho}
 \forall (X,\zeta)\in \gos\times \gog,\quad
 \hbox{Ad}_g(X+\zeta) = X+\hbox{Ad}_g\zeta,\quad
 \hbox{ad}_\xi(X+\zeta) = X+ [\xi,\zeta]
\end{equation}
In other words, $\gos$ and $\gog$ are stable by $\hbox{Ad}_\goG$ and $\hbox{ad}_\gog$ and their restrictions to
$\gos$ is trivial whereas their restrictions to $\gog$ coincides with, respectively, $\hbox{Ad}_\goG$
and $\hbox{ad}_\gog$.

Letting $\textbf{c}^i_{jk}$ be the structure coefficients of $\gog$ in the basis $(\textbf{t}_1,\cdots,\textbf{t}_r)$ and using the notation $\textbf{c}{^\gog}_{\gog\gog}:= \textbf{c}^i_{jk}\textbf{t}_i\otimes
\textbf{t}^j\otimes \textbf{t}^k \in \gog\otimes \gog^*\otimes \gog^*$ (see (\ref{cgggref})), we can write that, $\forall \xi^\gog,\eta^\gog\in \gog$, $[\xi^\gog,\eta^\gog] = \textbf{c}{^\gog}_{\celg\celg}\xi^{\celg}\eta^{\celg}$.

\subsubsection{First variation}

\noindent
\textbf{First variation with respect to $\pi_\gog$} ---
We write that the action is stationary with respect to variations $(\theta^\gog,\pi_\gog)\longmapsto (\theta^\gog,\pi_\gog+\varepsilon \delta\pi_\gog)$, for $\varepsilon$ small, where $\delta\pi_\gog = \chi_\gog = \frac{1}{2}\chi{_\gog}^{\celu\celu}f^{(N-2)}_{\celu\celu}$
(by using Convention (\ref{ConventionV})),  so that $\delta\pi_\gog$ is induced by $\delta\pi{_\gog}^{\gou\gou}$. Similarly
the curvature 2-form $\Theta^\gog:= \dR\theta^\gog+\frac{1}{2}[\theta^\gog\wedge \theta^\gog]$ decomposes as
\[
 \Theta^\gog = \frac{1}{2}\Theta{^\gog}_{\celu\celu}
 f^{\celu\celu}
 = \frac{1}{2}\Theta{^\gog}_{\cels\cels}
 f^{\cels\cels}
 + \Theta{^\gog}_{\cels\celg}
 f^{\cels\celg}
 + \frac{1}{2}\Theta{^\gog}_{\celg\celg}
 f^{\celg\celg}
\]
Hence $\pi_{\celg}\wedge
 \Theta^{\celg}
 = \frac{1}{2}\pi{_\celg}^{\celu\celu}\
 \Theta{^\celg}_{\celu\celu}\;f^{(N)}$. By using (\ref{decompositionofpiwithindices}) and (\ref{fifi}) we obtain the condition
\[
 \int_\mathcal{F}
 \left(\frac{1}{2}
 \chi{_\celg}^{\cels\cels}
 \left( \pi{^\celg}_{\cels\cels}
 + \Theta{^\celg}_{\cels\cels}\right)
 + \chi{_\celg}^{\cels\celg}
  \Theta{^\celg}_{\cels\celg}
 + \frac{1}{2}
 \chi{_\celg}^{\celg\celg}
  \Theta{^\celg}_{\celg\celg}
 \right)f^{(N)} = 0,
 \quad \forall \chi{_\gog}^{\gou\gou}
\]
which gives us the relations
\begin{equation}\label{elymp}
 \left\{\begin{array}{cclc}
       \pi{^\gog}_{\gos\gos}
 + \Theta{^\gog}_{\gos\gos} & = & 0 & \hbox{(a)} \\
       \Theta{^\gog}_{\gos\gog} & = & 0 & \hbox{(b)}  \\
       \Theta{^\gog}_{\gog\gog} & = & 0 & \hbox{(c)} .
        \end{array}
\right.
\end{equation}

\noindent
\textbf{First variation with respect to $\theta^\gog$} ---
We now look at the first variation of $\mathcal{A}$ through variations $(\theta^\gog,\pi_\gog)\longmapsto (\theta^\gog + \varepsilon \delta\theta^\gog,\pi_\gog)$, where $\delta\theta^\gog$ has a compact support.
It is useful
to decompose $\delta\theta^\gog$ as:
\[
\delta\theta^\gog = \tau{^\gog} = \tau{^\gog}_{\cels} \beta^{\cels} +
\tau{^\gog}_{\celg}
\theta^{\celg}
\]
This induces the variation $f^{(N)} \longmapsto f^{(N)} + \varepsilon\delta f^{(N)} + o(\varepsilon)$ with $\delta f^{(N)} = \tau{^{\celg}}_{\celg} f^{(N)}$. From (\ref{implicitpimumnu0}), which
implies $0 =\delta(\pi{_\gog}^{\gos\gos}
f^{(N)})$,
we deduce that the induced variation of
$\pi{_\gog}^{\gos\gos}$ is equal to $\delta\pi{_\gog}^{\gos\gos} = -\tau{^{\celg}}_{\celg}\pi{_\gog}^{\gos\gos}$ and thus $|\pi{_\gog}^{\gos\gos}|^2 \longmapsto |\pi{_\gog}^{\gos\gos}|^2 + \varepsilon \delta|\pi{_\gog}^{\gos\gos}|^2 + o(\varepsilon)$ with
$\delta |\pi{_\gog}^{\gos\gos}|^2 = -2\tau{^{\celg}}_{\celg}|\pi{_\gog}^{\gos\gos}|^2$. Hence
\[
 \delta\left(\frac{1}{2}|\pi{_\gog}^{\gos\gos}|^2f^{(N)}\right)
  =   \left(-\tau{^{\celg}}_{\celg} +\frac{\tau{^{\celg}}_{\celg}}{2}\right) |\pi{_\gog}^{\gos\gos}|^2f^{(N)}
 = - \frac{|\pi{_\gog}^{\gos\gos}|^2}{2}\tau{^{\celg}}_{\celg} f^{(N)}
= -\frac{|\pi{_\gog}^{\gos\gos}|^2}{2}
 \tau^{\celg}\wedge f^{(N-1)}_{\celg}
\]
Let us set
$\dR^\theta(\tau^\gog):=
\dR(\tau^\gog) +
[\theta^\gog\wedge\tau^\gog]$ and $\dR^\theta(\pi_\gog):=
\dR \pi_\gog +
\hbox{ad}^*_{\theta}\wedge \pi_\gog$.
We remark that
$\delta\Theta^\gog
= \delta\left(\dR\theta^\gog
 + \frac{1}{2}[\theta^\gog\wedge
 \theta^\gog]\right)
 = \dR^\theta \tau^\gog$ and
thus by (\ref{twistedLeibniz})
\[
 \delta\left(\pi_{\celg}\wedge \Theta^{\celg}\right) =  \delta\Theta^{\celg}\wedge
 \pi_{\celg}
 = \left(\dR^\theta\tau^{\celg}\right)\wedge
 \pi_{\celg}
  \\
= \dR^\theta
 \left(\tau^{\celg}\wedge \pi_{\celg}\right)
  + \tau^{\celg}\wedge
   \dR^\theta \pi_{\celg}
\]
Lastly we observe that
$\int_\mathcal{F}\dR^\theta\left(\tau^{\celg}\wedge \pi_{\celg}\right)
= \int_\mathcal{F}\dR\left(\tau^{\celg}\wedge \pi_{\celg}\right)$
since the coefficients of $\tau^\celg\wedge \pi_\celg$ are in $\R$, a trivial representation of $\gog$.
Thus the first variation of the action vanishes iff
\[
 \int_\mathcal{F} \tau^{\celg}\wedge
 \left(\dR^\theta\pi_{\celg}-
 \frac{|\pi{_\gog}^{\gos\gos}|^2}{2}
 f^{(N-1)}_{\celg}\right),
 \quad \forall \tau^\gog\hbox{ with compact support,}
\]
which give us the equation
\begin{equation}\label{elymphi}
\dR^\theta\pi_{\gog} =
 \frac{|\pi{_\gog}^{\gos\gos}|^2}{2}
 f^{(N-1)}_{\gog}
\end{equation}

\subsubsection{Principal bundle structure and equivariance of the connection}

We first exploit Equation (c) in system (\ref{elymp}), i.e. $\hbox{d}\theta^\gog +\frac{1}{2}[\theta^\gog\wedge\theta^\gog] = 0$.

Consider on the product manifold $\widehat{\goG} \times \mathcal{Y} = \{(h,\textsf{y})\in \widehat{\goG}\times \mathcal{Y}\}$ the $\gog$-valued 1-form
$\tau^\gog := \theta^\gog - h^{-1}\hbox{d}h$.
It satisfies the identity $\hbox{d}\tau^\gog = \hbox{d}\theta^\gog +\frac{1}{2}[\theta^\gog\wedge\theta^\gog]
- [\theta^\gog\wedge \tau^\gog] + \frac{1}{2}[\tau^\gog\wedge \tau^\gog]$ and its rank is clearly equal to $r$.
However Equation (c) in (\ref{elymp}) implies that,
for any fiber $\mathcal{F}_x$,
$\hbox{d}\theta^\gog + \frac{1}{2}[\theta^\gog\wedge \theta^\gog]|_{\mathcal{F}_x} = 0$
and thus
$\hbox{d}(\tau^\gog|_{\mathcal{F}_x\times \widehat{\goG}}) = 0\hbox{ mod}[\tau^\gog]$.
Hence, by Frobenius' theorem,
for any $(g_0,\textsf{y}_0)\in \widehat{\goG}\times \mathcal{F}_x$,
there exists a unique $r$-dimensional submanifold $\Gamma\subset \mathcal{F}_x\times \widehat{\goG}$ which is
a maximal solution of $\tau^\gog|_\Gamma=0$  and which contains $(g_0,\textsf{y}_0)$.

It is clear also that, $\forall (g,\textsf{y})\in \widehat{\goG} \times \mathcal{F}_x$, $\forall (\xi,v) \in T_{g}\widehat{\goG}\times T_{\textsf{y}}\mathcal{F}_x$, the equation $g^{-1}\hbox{d}g(\xi) = \theta^\gog(v)$ defines the graph of a vector space isomorphism between $T_{g}\widehat{\goG}$ and $T_{\textsf{y}}\mathcal{F}_x$. This implies that, around each point $(g,\textsf{y})\in \Gamma$, $\Gamma$ is locally the graph of a diffeomorphism between a neighbourhood of $g$ in $\widehat{\goG}$ and a neighbourhood of $\textsf{y}$ in $\mathcal{F}_x$.
But we have more: since each fiber $\mathcal{F}_\textsf{x}$ is actually a maximal solution of the system $\theta^\gos|_\textsf{f} = 0$, we can apply the following lemma to deduce that there exists a map from $\widehat{\goG}$ to $\mathcal{F}_\textsf{x}$, the graph of which is $\Gamma$, and thus this map is a universal cover of $\mathcal{F}_\textsf{x}$.

\begin{lemm}\label{lemmaGcoversleaves}
 Assume that $\widehat{\goG}$ is simply connected and that $(\mathcal{Y},\theta^\gos,\theta^\gog)$ is $\gog$-complete. Let $\textsf{\emph{f}}$ be a maximal integral solution of the system $\theta^\gos|_{\textsf{\emph{f}}} = 0$ of dimension $r$. Then $\widehat{\goG}$ is a \textbf{universal cover} of $\emph{\textsf{f}}$.

As a corollary, $\emph{\textsf{f}}$ is diffeomorphic to a quotient of $\widehat{\goG}$ by a finite subgroup and, if $\widehat{\goG}$ is furthermore \textbf{compact}, then $\emph{\textsf{f}}$ is compact.
\end{lemm}
\emph{Proof} --- Fix any base point $\textsf{y}_0\in \textsf{f}$ and consider:
\begin{itemize}
 \item the set $\mathscr{P}_{\widehat{\goG},1_{\widehat{\goG}}}$ of based paths $\gamma\in \mathscr{C}^1([0,1],\widehat{\goG})$ such  that $\gamma(0) = 1_{\widehat{\goG}}$ and
 \item the set $\mathscr{P}_{\textsf{f},\textsf{y}_0}$  of based paths $u\in \mathscr{C}^1([0,1],\textsf{f})$ such that $u(0) = \textsf{y}_0$.
\end{itemize}
We define an operator $\mathscr{T}$ from $\mathscr{P}_{\widehat{\goG},1_{\widehat{\goG}}}$ to $\mathscr{P}_{\textsf{f},\textsf{y}_0}$ as follows: to any $\gamma\in \mathscr{P}_{\widehat{\goG},1_{\widehat{\goG}}}$ we associate the unique path $u = \mathscr{T}(\gamma)\in \mathscr{P}_{\textsf{f},\textsf{y}_0}$ such that
\[
u(0) = \textsf{y}_0\quad \hbox{ and }\quad
\left[\ \gamma^{-1}\dR \gamma = u^*\theta^\gog
 \quad \Longleftrightarrow \quad (\gamma,u)^*\tau^\gog = 0 \ \right]
\]
We will show that, for any $\gamma\in \mathscr{P}_{\widehat{\goG},1_{\widehat{\goG}}}$, the end point $\mathscr{T}(\gamma)(1)$ of $u=\mathscr{T}(\gamma)$ depends uniquely on the end point $\gamma(1)$ of $\gamma$, i.e.,
\begin{equation}\label{implicationbyhomotopy}
\forall \gamma_0,\gamma_1\in \mathscr{P}_{\goG,1_{\widehat{\goG}}},\quad
\gamma_0(1) = \gamma_1(1) \quad
\Longrightarrow \quad \mathscr{T}(\gamma_0)(1) = \mathscr{T}(\gamma_1)(1)
\end{equation}
Since $\widehat{\goG}$ is connected, for any $g\in \widehat{\goG}$, there exists a path $\gamma\in \mathscr{P}_{\widehat{\goG},1_{\widehat{\goG}}}$ such that $\gamma(1) = g$, thus (\ref{implicationbyhomotopy}) shows the existence of a unique map $T:\widehat{\goG} \longrightarrow \textsf{f}$ such that, for any $\gamma\in \mathscr{P}_{\widehat{\goG},1_{\widehat{\goG}}}$, $T(\gamma(1)) = \mathscr{T}(\gamma)(1)$. The graph of $T$ clearly coincides with the integral leaf of $\tau^\gog$ in $\widehat{\goG} \times \textsf{f}$ passing through $(1_{\widehat{\goG}},\textsf{y}_0)$ and thus $T$ is a smooth cover of $\textsf{f}$, which is actually the universal cover since $\widehat{\goG}$ is simply connected.

Let us prove (\ref{implicationbyhomotopy}).
Let $\gamma_0$ and $\gamma_1$ be in $\mathscr{P}_{\widehat{\goG},1_{\widehat{\goG}}}$ and assume that $\gamma_0(1) = \gamma_1(1)$. Since $\widehat{\goG}$ is simply connected there exists a smooth homotopy $\Gamma\in\mathscr{C}^1([0,1]^2,\widehat{\goG})$ such that, $\forall t,s\in [0,1]$,
\[
 \xymatrix{ 
 \Gamma(0,1) = 1_{\widehat{\goG}} \ar@{-}[r] & \Gamma(t,1) = \gamma_1(t) \ar@{-}[r] & \Gamma(1,1) = \gamma_0(1)  \\
 \Gamma(0,s) = 1_{\widehat{\goG}} \ar@{-}[u] & & \Gamma(1,s) = \gamma_0(1) \ar@{-}[u] & \\
 \Gamma(0,0) = 1_{\widehat{\goG}} \ar@{-}[u] \ar@{-}[r] & \Gamma(t,0) = \gamma_0(t) \ar@{-}[r] & \Gamma(1,0) = \gamma_0(1) \ar@{-}[u] 
 }
\]
To this map we associate the unique map $U\in\mathscr{C}^1([0,1]^2,\textsf{f})$ defined by
\[
 \left\{
 \begin{array}{rcl}
  U(0,0) & = & \textsf{y}_0 \\
  \theta^\gog_{U(t,0)}\left(\frac{\partial U}{\partial t}(t,0)\right) & = & \left(\Gamma^{-1}\frac{\partial \Gamma}{\partial t}\right)(t,0),\quad \forall t\in [0,1] \\
  \theta^\gog_{U(t,s)}\left(\frac{\partial U}{\partial s}(t,s)\right) & = & \left(\Gamma^{-1}\frac{\partial \Gamma}{\partial s}\right)(t,s),\quad \forall (t,s)\in [0,1]^2
 \end{array}\right.
\]
Thus if we set $F:= (\Gamma,U)\in \mathscr{C}^1([0,1]^2,\widehat{\goG}\times \textsf{f})$, the previous relations read $F(0,0) = (1_{\widehat{\goG}},\textsf{f}_0)$ and
\begin{equation}\label{homotopyjumelee}
  (F^*\tau)_{(t,0)}\left(\frac{\partial}{\partial t}\right) = 0 \quad \hbox{ and }\quad (F^*\tau)_{(t,s)}\left(\frac{\partial}{\partial s}\right) = 0,\quad \forall t,s\in [0,1]
\end{equation}
Set $\sigma:= \frac{1}{2}\hbox{ad}_{\tau^\gog} - \hbox{ad}_{\theta^\gog}|_{\widehat{\goG}\times \textsf{f}} \in \hbox{End}(\gog)\otimes \Omega^1(\widehat{\goG}\times \textsf{f})$, so that $\hbox{d}\tau^\gog|_{\widehat{\goG}\times \textsf{f}} = \sigma \wedge \tau^\gog|_{\widehat{\goG}\times \textsf{f}}$, and set $\alpha:= F^*\tau^\gog$ and $\beta:= F^*\sigma$. Then $\hbox{d}\alpha = \beta\wedge \alpha$ and the second relation in (\ref{homotopyjumelee}) translates as $\alpha\left(\frac{\partial}{\partial s}\right) = 0$. We now use Cartan's formula
\[
 \hbox{d}\alpha\left(\frac{\partial}{\partial t},\frac{\partial}{\partial s}\right) + \alpha\left(\left[\frac{\partial}{\partial t},\frac{\partial}{\partial s}\right]\right)
=  \frac{\partial}{\partial t}\left(\alpha\left(\frac{\partial}{\partial s}\right)\right) -  \frac{\partial}{\partial s}\left(\alpha\left(\frac{\partial}{\partial t}\right)\right)
\]
which simplifies to
\[
 \beta\wedge\alpha\left(\frac{\partial}{\partial t},\frac{\partial}{\partial s}\right) + 0 = 0  - \frac{\partial}{\partial s}\left(\alpha\left(\frac{\partial}{\partial t}\right)\right)
\]
and thus
\[
 \frac{\partial}{\partial s}\left(\alpha\left(\frac{\partial}{\partial t}\right)\right)(t,s) = \beta\left(\frac{\partial}{\partial s}\right)\alpha\left(\frac{\partial}{\partial t}\right)(t,s),\quad
 \forall (t,s)\in [0,1]^2
\]
Since by (\ref{homotopyjumelee}) we also have the initial condition $\alpha\left(\frac{\partial}{\partial t}\right)(t,0) = 0$, $\forall t\in [0,1]$, we deduce that
\[
 \alpha\left(\frac{\partial}{\partial t}\right)(t,s) = 0,\quad \forall (t,s)\in [0,1]^2
\]
This means that, $\forall (t;s)\in [0,1]^2$, $(F^*\tau)_{(t,s)}\left(\frac{\partial}{\partial t}\right) = 0$, i.e., $\theta^\gog_{U(t,s)}\left(\frac{\partial U}{\partial t}(t,s)\right)= \left(\Gamma^{-1}\frac{\partial \Gamma}{\partial t}\right)(t,s)$. This can also be translated by defining the maps $\gamma_s\in \mathscr{P}_{\textsf{f},\textsf{y}_0}$ and $u_s\in \mathscr{P}_{\widehat{\goG},1_{\widehat{\goG}}}$ such that, respectively, $\forall (t,s)\in [0,1]^2$, $\Gamma(t,s) = \gamma_s(t)$ and  $U(t,s) = u_s(t)$, by writing $u_s^*\theta^\gog = \gamma_s^{-1}\hbox{d}\gamma_s$, $\forall s\in [0,1]$. Since, $\forall s\in [0,1]$, $u_s(0) = \textsf{y}_0$, we conclude that $u_s = \mathscr{T}(\gamma_s)$.

But we also have, by the definition of $\Gamma$, $\Gamma^{-1}\frac{\partial \Gamma}{\partial s}(1,s) = 0$, $\forall s\in [0,1]$, and hence $\frac{\partial U}{\partial s}(1,s) = 0$, $\forall s\in [0,1]$. This implies that $u_s(1) = u_0(1)$, i.e., $\mathscr{T}(\gamma_s)(1) = \mathscr{T}(\gamma_0)(1)$, $\forall s\in [0,1]$ and, in particular $\mathscr{T}(\gamma_1)(1) = \mathscr{T}(\gamma_0)(1)$.\hfill $\square$\\

A consequence of this Lemma \ref{lemmaGcoversleaves} is that, if $\widehat{\goG}$ is compact, all fibers are compact. Hence by a result of Ehresmann \cite{Ehresmann50} we deduce that $\mathcal{F}$ has a structure of fiber bundle over $\mathcal{X}$. In particular all fibers are diffeomorphic to a quotient $\goG$ of $\widehat{\goG}$. (Note that the latter conclusion can also be achieved by applying a straightforward variant of Lemma \ref{lemmaifGnoncompact} below.)

Thus, by choosing some (possibly local) section $\Sigma$ of $\mathcal{F}$, there exists a unique
map $g:\mathcal{F}\longrightarrow \goG$ such that, for any $x$,
\begin{equation}\label{YMnormalisationAtteinte}
\theta^\gog - g^{-1}\dR g|_{\mathcal{F}_\textsf{x}} = 0\quad \Longleftrightarrow \quad
\textbf{A}{^\gog}|_{\mathcal{F}_\textsf{x}} = 0,\hbox{ where }
\textbf{A}{^\gog}:= g\theta^\gog g^{-1} - \dR g\ g^{-1}.
\end{equation}
and such that
$g$ is equal to $1_\goG$ on $\Sigma$. Condition (\ref{YMnormalisationAtteinte}) implies that the 1-form $\textbf{A}{^\gog}\in \gog\otimes \Omega^1(\mathcal{F})$ decomposes as
$\textbf{A}{^\gog} = \textbf{A}{^\gog}_{\cels} \beta^{\cels}$. It also means that
$\theta^\gog = g^{-1}\textbf{A}{^\gog}g + g^{-1}\dR g $
is \textbf{normalized} and implies that
$\dR\theta^\gog + \frac{1}{2}[\theta^\gog\wedge \theta^\gog] = g^{-1}(\dR\textbf{A}^\gog+\frac{1}{2}[\textbf{A}^\gog\wedge \textbf{A}^\gog])g$, i.e. by defining $\textbf{F}^\gog:= \dR\textbf{A}{^\gog}+\frac{1}{2}[\textbf{A}^\gog\wedge \textbf{A}^\gog]$,
\begin{equation}\label{PhiegalAdgmoinsunF}
 \textbf{F}^\gog = \hbox{Ad}_{g}\Theta^\gog
\end{equation}
For any function $\alpha$ on $\mathcal{F}$, let us denote by
$\partial_\gos \alpha$ and $\partial_\gog \alpha$  the coefficients
in the decomposition $\dR\alpha = \partial_\cels\alpha\ \beta^{\cels}
+ \partial_\celg \alpha\ \theta^{\celg}$.
Then through the decomposition $\textbf{A}^\gog = \textbf{A}{^\gog}_{\cels}\beta^{\cels}$, $\textbf{F}^\gog$ decomposes as
\[
\textbf{F}^\gog
= \frac{1}{2}(\partial_{\cels_1}\textbf{A}{^\gog}_{\cels_2} - \partial_{\cels_2}\textbf{A}{^\gog}_{\cels_1} +[\textbf{A}{^\gog}_{\cels_1},\textbf{A}{^\gog}_{\cels_2}])\beta^{\cels_1\cels_2}
- \partial_\celg \textbf{A}{^\gog}_{\cels}\ \beta^{\cels}\wedge\theta^{\celg}.
\]
Equations (\ref{elymp}--b) and (\ref{PhiegalAdgmoinsunF}) now
imply $\partial_\gog \textbf{A}{^\gog}_\gos = 0$, which means that $\textbf{A}{^\gog}_\gos$ is constant on each fiber (i.e. the coefficients $\textbf{A}{^\gog}_\gos$ depends only on $\textsf{x}\in \mathcal{X}$). Equivalentely $\theta^\gog$ is \textbf{equivariant}.
Hence the coefficients $\textbf{F}{^\gog}_{\gos\gos}$ in the decomposition
$\textbf{F}^\gog = \frac{1}{2}
\textbf{F}{^\gog}_{{\cels_1}{\cels_2}}\beta^{\cels_1\cels_2}$ are also independent of $g$.

We next introduce the frame $e^\gou:= \hbox{Ad}_g
f^\gou$.
This implies in particular  by (\ref{YMnormalisationAtteinte}) that
\begin{equation}\label{eegalAdgphi}
 e^\gog:= \hbox{Ad}_g\theta^\gog = \textbf{A}^\gog
+ dg\,g^{-1}
\end{equation}
We also set
\begin{equation}\label{pegalAdetoilepi}
 p_\gog:= \hbox{Ad}^*_g\pi_\gog
 = \hbox{Ad}^*_g\left(\frac{1}{2} \pi{_\gog}^{\celu\celu}
 f^{(N-2)}_{\celu\celu}\right)
\end{equation}
and its decomposition by using the
$(N-2)$-form $e^{(N-2)}_{\gou\gou}:=
\hbox{Ad}_g^*\otimes \hbox{Ad}_g^*\left(f^{(N-2)}_{\gou\gou}\right)$:
\[
 p_\gog = \frac{1}{2} p{_\gog}^{\celu\celu}
 e^{(N-2)}_{\celu\celu}
\]
where, according to (\ref{eetf56}), $p{_\gog}^{\gou\gou}:=
 \hbox{Ad}_g^*\otimes \hbox{Ad}_g\otimes \hbox{Ad}_g\left(\pi{_\gog}^{\gou\gou}\right)$.

In particular (since the action of $\hbox{Ad}_g$ on $\gos$ is trivial)
 $p{_\gog}^{\gos\gos} =
 \hbox{Ad}_g^*\otimes 1_\gos\otimes 1_\gos
 \left(\pi{_\gog}^{\gos\gos}\right)$.
At this point we exploit Equation (\ref{elymp}--a) that we translate as $\pi{_\gog}^{\gos\gos}
 + \Theta{_\gog}^{\gos\gos}
 = 0$, where
$\Theta{_\gog}^{\gos\gos}
= \left(\textsf{k}_{\gog\celg}\otimes
\textsf{b}^{\gos\cels}\otimes
\textsf{b}^{\gos\cels}\right) \Theta{^\celg}_{\cels\cels}$.
Thus actually $p{_\gog}^{\gos\gos} =
- \hbox{Ad}_g^*\otimes 1_\gos\otimes 1_\gos
 \left(\Theta{_\gog}^{\gos\gos}\right)$. Hence by using (\ref{57})
\[
\begin{array}{ccl}
 -p{_\gog}^{\gos\gos}
 & = & \hbox{Ad}_g^*\otimes 1_\gos\otimes 1_\gos
 \left(\Theta{_\gog}^{\gos\gos}\right)
 = \hbox{Ad}_g^*\otimes 1_\gos\otimes 1_\gos
 \left(
 \left(
\textsf{k}_{\gog\celg}\otimes
\textsf{b}^{\gos\cels}\otimes
\textsf{b}^{\gos\cels}\right) \Theta{^\celg}_{\cels\cels}
 \right) \\
&  \overset{\star}{=} & \left(
\textsf{k}_{\gog\celg}\otimes
\textsf{b}^{\gos\cels}\otimes
\textsf{b}^{\gos\cels}\right)
\left(\hbox{Ad}_g\otimes 1_{\gos^*}\otimes 1_{\gos^*}
 \left(\Theta{^{\celg}}_{\cels\cels}\right)\right)
= \left(
\textsf{k}_{\gog\celg}\otimes
\textsf{b}^{\gos\cels}\otimes
\textsf{b}^{\gos\cels}\right)\textbf{F}{^\celg}_{\cels\cels}
\end{array}
\]
where in $\overset{\star}{=}$ we used the fact that $\textsf{k}$ is invariant by $\hbox{Ad}_g$, i.e.
$\hbox{Ad}_g^*\otimes \hbox{Ad}_g^*
\left(\textsf{k}_{\gog\gog}\right)
= \textsf{k}_{\gog\gog}$.
Hence by
setting $\textbf{F}{_\gog}^{\gos\gos}
:= \left(
\textsf{k}_{\gog\celg}\otimes
\textsf{b}^{\gos\cels}\otimes
\textsf{b}^{\gos\cels}\right)
\textbf{F}{^\celg}_{\cels\cels}$ (\ref{elymp}-a), translates as \begin{equation}\label{elymp1}
 p{_\gog}^{\gos\gos}
 = - \textbf{F}{_\gog}^{\gos\gos}
\end{equation}
Lastly we translate (\ref{elymphi}) as follows: by (\ref{eetf56}) we have $\hbox{Ad}_g^* f^{(N-1)}_{\gog} = e^{(N-1)}_{\gog}$. Moreover
by using (\ref{pegalAdetoilepi}) and (\ref{domegap=Adgdvarphivarpi})
we obtain that $\dR^\textbf{A}p_{\gog}:= \dR p_\gog + \hbox{ad}^*_{\textbf{A}}p_\gog = \hbox{Ad}_g^* \left(
 \dR^\theta\pi_{\gog} \right)$. Hence since $\textsf{k}$
is $\hbox{Ad}_g$-invariant (which implies
$|\pi{_\gog}^{\gos\gos}|^2 = |p{_\gog}^{\gos\gos}|^2$) and because of (\ref{elymphi}) and (\ref{elymp1}) we deduce
\begin{equation}\label{elymA}
 \dR^\textbf{A}p_{\gog} =  \hbox{Ad}_g^* \left(
 \frac{|\pi{_\gog}^{\gos\gos}|^2}{2}
 f^{(N-1)}_{\gog}\right)
 =\frac{|\textbf{F}{_\gog}^{\gos\gos}|^2}{2}
 e^{(N-1)}_{\gog}
\end{equation}


\noindent
\textbf{} ---

\subsubsection{Computation of the left hand side of (\ref{elymA})}It turns out that Equation (\ref{elymA}) implies that the connection $\textbf{A}^\gog$ is a solution of the Yang--Mills system of equations. However the proof of that fact requires a careful computation of the left hand side of (\ref{elymA}) using a decomposition of $p_\gog$ in the basis $e^{(N-2)}_{\gou\gou}$ obtained out of $e^\gou$.
(Note that an alternative method is possible, by using the coframe
$(\beta^\gos,g^{-1}\dR g)$ instead of $e^\gou$.)
This is the most delicate part.

Let $\underline{\gamma}^{so(\gos)} \in
so(\gos,\textsf{b})\otimes\Omega^1(\mathcal{X})$
be the connection 1-form of the Levi-Civita connection
$\nabla$
on $(\mathcal{X},\textbf{g})$ and
$\gamma:= \gamma^{so(\gos)}:= P^*\underline{\gamma}^{so(\gos)}
\in so(\gos,\textsf{b})\otimes\Omega^1(\mathcal{F})$.
The orthogonal splitting $\gou = \gos\oplus \gog$ induces an embedding of $so(\gos,\textsf{b})$ in $so(\gou,\textsf{h})$ so that actually $\gamma
\in so(\gou,\textsf{h})\otimes\Omega^1(\mathcal{F})$.
Similarly $\hbox{ad}_\textbf{A} \in so(\gog,\textsf{k})\otimes\Omega^1(\mathcal{F})
\subset so(\gou,\textsf{h})\otimes\Omega^1(\mathcal{F})$ and thus $\gamma + \hbox{ad}_\textbf{A}\in so(\gou,\textsf{h})\otimes \Omega^1(\mathcal{X})$. We then define the connection $\dR^{\gamma,\textbf{A}}$ acting on functions $\xi^\gou$ from $\mathcal{F}$ to $\gou$ by
\begin{equation}\label{nablaAxiU}
\begin{array}{c}
\dR^{\gamma,\textbf{A}}\xi^\gou:=
 \dR^\textbf{A} \xi^\gou + \gamma\, \xi^\gou
 = \dR^{\gamma}\xi^\gos + \dR^{\textbf{A}}\xi^\gog
 \\
 \hbox{with }\quad
 \dR^{\gamma}\xi^\gos := \dR \xi^\gos
 + \gamma\, \xi^{\gos}
 \quad\hbox{ and }\quad \dR^{\textbf{A}}\xi^\gog:= \dR \xi^\gog
 + \hbox{ad}_\textbf{A} \xi^\gog
\end{array}
\end{equation}
and extend it by using the graded Leibniz rule to any exterior differential form with coefficients in a tensor product of $\gou$ and $\gou^*$.

A key point is to observe that, since the action of $so(\gos,\textsf{b})$ on $\gog$ is trivial, $\dR^\textbf{A}p_{\gog} = \dR^{\gamma,\textbf{A}}p_{\gog}$. By using the decomposition $p_{\gog} = \frac{1}{2} p{_\gog}^{\celu\celu} e^{(N-2)}_{\celu\celu}$ and
the Leibniz rule (\ref{twistedLeibniz}) we deduce $\dR^\textbf{A}p_\gog =
 \frac{1}{2}\dR^{\gamma,\textbf{A}} p{_\gog}^{\celu\celu}\wedge e_{\celu\celu}^{(N-2)} +
 \frac{1}{2} p{_\gog}^{\celu\celu}
\dR^{\gamma,\textbf{A}}e^{(N-2)}_{\celu\celu}$. Moreover if we denote by $\partial^{\gamma,\textbf{A}}_{\gou} p{_\gog}^{\gou\gou}\in \gou^*\otimes \gog^*\otimes \gou\otimes\gou\otimes \mathscr{C}^\infty(\mathcal{F})$ the coefficients such that
$\dR^{\gamma,\textbf{A}} p{_\gog}^{\gou\gou}
= (\partial^{\gamma,\textbf{A}}_{\celu} p{_\gog}^{\gou\gou}) e^{\celu}$, then $\frac{1}{2}\dR^{\gamma,\textbf{A}} p{_\gog}^{\celu\celu}\wedge e_{\celu\celu}^{(N-2)}
= \frac{1}{2}(\partial^{\gamma,\textbf{A}}_{\celu} p{_\gog}^{\celu_1\celu_2}) e^{\celu} \wedge   e_{\celu_1\celu_2}^{(N-2)}
= \partial^{\gamma,\textbf{A}}_{\celu_2} p{_\gog}^{\celu\celu_2}e_{\celu}^{(N-1)}$. Hence
\begin{equation}\label{dApi}
 \dR^\textbf{A}p_\gog = \dR^{\gamma,\textbf{A}}p_{\gog} =
 \partial^{\gamma,\textbf{A}}_{\celu_2} p{_\gog}^{\celu\celu_2}e_{\celu}^{(N-1)} +
 \frac{1}{2} p{_\gog}^{\celu\celu}
\dR^{\gamma,\textbf{A}}e^{(N-2)}_{\celu\celu}
\end{equation}
By introducing the coefficients
$\partial_{\gou}p{_\gou}^{\gou\gou}$ such that
$\dR p{_\gou}^{\gou\gou}
= \left(\partial_{\celu}p{_\gou}^{\gou\gou}\right) e^{\celu}$ and the coefficients $\gamma_\gou$ such that
$\gamma = \gamma_\celu e^{\celu}$, the coefficients $\partial^{\gamma,\textbf{A}}_\gou p{_\gog}^{\gou\gou}$ read\footnote{Alternatively, for instance, the second relation in this system reads $\partial^{\gamma,\textbf{A}}_{\gou}p{_\gog}^{\gos\gog}
= \partial_{\gou} p{_\gog}^{\gos\gog}
- \mathbf{c}{^\celg}_{\celg_0\gog}\mathbf{A}{^{\celg_0}}_{\gou}p{_\celg}^{\gos\gog} + \gamma{^\gos}_{\cels\gou}p{_\gog}^{\cels\gog} + \mathbf{c}{^\gog}_{\celg_0\celg}\mathbf{A}{^{\celg_0}}_{\gou}p{_\gog}^{\gos\celg}$.}
\[
\left\{
 \begin{array}{ccccc}
\partial^{\gamma,\textbf{A}}_{\gou}p{_\gog}^{\gos\gos}
& = & \partial_{\gou} p{_\gog}^{\gos\gos}
& + & \left(\hbox{ad}^*_{\textbf{A}_{\gou}}\otimes 1\otimes 1
+ 1\otimes \gamma_\gou\otimes 1 + 1\otimes 1\otimes \gamma_\gou\right)p{_\gog}^{\gos\gos}
 \\
\partial^{\gamma,\textbf{A}}_{\gou}p{_\gog}^{\gos\gog}
& = & \partial_{\gou} p{_\gog}^{\gos\gog}
& + & \left(\hbox{ad}^*_{\textbf{A}_{\gou}}\otimes 1\otimes 1
+ 1\otimes \gamma_\gou\otimes 1 + 1\otimes 1\otimes \hbox{ad}_{\textbf{A}_{\gou}}\right)p{_\gog}^{\gos\gog}
 \\
\partial^{\gamma,\textbf{A}}_{\gou}p{_\gog}^{\gog\gos}
& = & \partial_{\gou} p{_\gog}^{\gog\gos}
& + & \left(\hbox{ad}^*_{\textbf{A}_{\gou}}\otimes 1\otimes 1
+ 1\otimes \hbox{ad}_{\textbf{A}_{\gou}}\otimes 1 + 1\otimes 1\otimes \gamma_\gou\right)p{_\gog}^{\gog\gos}
 \\
\partial^{\gamma,\textbf{A}}_{\gou}p{_\gog}^{\gog\gog}
& = & \partial_{\gou} p{_\gog}^{\gog\gog}
& + & \left(\hbox{ad}^*_{\textbf{A}_{\gou}}\otimes 1\otimes 1
+ 1\otimes \hbox{ad}_{\textbf{A}_{\gou}}\otimes 1 + 1\otimes 1\otimes \hbox{ad}_{\textbf{A}_{\gou}}\right)p{_\gog}^{\gog\gog}
 \end{array}\right.
\]
The sum $\partial^{\gamma,\textbf{A}}_{\celu} p{_\gog}^{\gou\celu}$ splits as $\partial^{\gamma,\textbf{A}}_{\celu} p{_\gog}^{\gou\celu} = \partial^{\gamma,\textbf{A}}_{\celu} p{_\gog}^{\gos\celu} + \partial^{\gamma,\textbf{A}}_{\celu} p{_\gog}^{\gog\celu}$, with
\[
 \partial^{\gamma,\textbf{A}}_{\celu} p{_\gog}^{\gos\celu}
  = \partial^{\gamma,\textbf{A}}_{\cels} p{_\gog}^{\gos\cels} + \partial^{\gamma,\textbf{A}}_{\celg} p{_\gog}^{\gos\celg}
  \quad \hbox{ and}\quad
  \partial^{\gamma,\textbf{A}}_{\celu}p{_\gog}^{\gog\celu}
  =  \partial^{\gamma,\textbf{A}}_{\cels}p{_\gog}^{\gog\cels} + \partial^{\gamma,\textbf{A}}_{\celg}p{_\gog}^{\gog\celg}
\]
where the first terms on the r.h.s are
\[
 \left\{\begin{array}{ccccc}
\partial^{\gamma,\textbf{A}}_{\cels} p{_\gog}^{\gos\cels}
& := &
\partial_{\cels} p{_\gog}^{\gos\cels}
& + & \left(\hbox{ad}^*_{\textbf{A}_{\cels}}\otimes 1\otimes 1
+ 1\otimes \gamma_{\cels}\otimes 1 + 1\otimes 1\otimes \gamma_{\cels}\right)p{_\gog}^{\gos\cels}
\\
\partial^{\gamma,\textbf{A}}_{\cels} p{_\gog}^{\gog\cels}
& := &
\partial_\cels p{_\gog}^{\gog\cels}
& + & \left(\hbox{ad}^*_{\textbf{A}_{\cels}}\otimes 1\otimes 1
+ 1\otimes \hbox{ad}_{\textbf{A}_{\cels}}\otimes 1 + 1\otimes 1\otimes \gamma_{\cels}\right)p{_\gog}^{\gog\cels}
\end{array}\right.
\]
The expressions of the second terms $\partial^{\gamma,\textbf{A}}_{\celg} p{_\gog}^{\gos\celg}$ and $\partial^{\gamma,\textbf{A}}_{\celg}p{_\gog}^{\gog\celg}$ simplify because of  the observations that
$\textbf{A}{^\gog}_{\gou} = \textbf{A}{^\gog}_{\gos}$ and $\gamma_\gou = \gamma_\gos$
(i.e. $\textbf{A}{^\gog}_{\gog} = \gamma_\gog = 0$):
\[
 \partial^{\gamma,\textbf{A}}_{\celg} p{_\gog}^{\gos\celg} = \partial_\celg p{_\gog}^{\gos\celg}
  \quad \hbox{ and}\quad
  \partial^{\gamma,\textbf{A}}_{\celg}p{_\gog}^{\gog\celg} = \partial_\celg p{_\gog}^{\gog \celg}
\]
Thus
\begin{equation}\label{dApimi1}
 \left(\partial^{\gamma,\textbf{A}}_{\celu_1} p{_\gog}^{\celu\celu_1}\right) e_{\celu}^{(N-1)}
 =
 \left(\partial^{\gamma,\textbf{A}}_{\cels_1} p{_\gog}^{\cels\cels_1}
   + \partial_\celg p{_\gog}^{\cels\celg}\right)e_{\cels}^{(N-1)}
    + \left( \partial^{\gamma,\textbf{A}}_{\cels} p{_\gog}^{\celg\cels}
   + \partial_{\celg_1} p{_\gog}^{\celg\celg_1}\right)e_{\celg}^{(N-1)}
\end{equation}

In order to compute the second term $\frac{1}{2} p{_\gog}^{\celu\celu}
\dR^{\gamma,\textbf{A}}e^{(N-2)}_{\celu\celu}$ we use (\ref{nablaLeibniz}), i.e.
$\dR^{\gamma,\textbf{A}}e^{(N-2)}_{\gou\gou}
= \dR^{\gamma,\textbf{A}} e^{\celu}\wedge
e^{(N-2)}_{\gou\gou\celu}$ and hence
we need to compute $\dR^{\gamma,\textbf{A}} e^{\gou}$.
We recall that
$e^\gou = \hbox{Ad}_gf^\gou
= \hbox{Ad}_g(\beta^\gos + \theta^\gog)
= \beta^\gos + \textbf{A}^\gog + \dR g\ g^{-1}$,
i.e. $e^\gou = e^\gos+e^\gog$ with $e^\gos = \beta^\gos$ and
$e^\gog = \textbf{A}^\gog + \dR g\ g^{-1}$.
Thus  since by
(\ref{dAeagalFee}) $\dR^\textbf{A}e^\gog
= \textbf{F}^\gog +
\frac{1}{2}[e^\gog\wedge e^\gog]$ we have
\[
 \dR^{\gamma,\textbf{A}}e^\gou =
 \dR^{\gamma} e^\gos + \dR^\textbf{A}e^\gog
= \dR^{\gamma} e^\gos +
 \textbf{F}^\gog +
 \frac{1}{2}[e^\gog\wedge e^\gog]
\]
where $\dR^{\gamma} e^\gos:= \dR e^\gos + \gamma\wedge e^\cels$.
However the latter quantity is the torsion, which vanishes since $\gamma = \gamma^{so(\gos)}$ corresponds to the Levi-Civita connection.
Thus the previous identity reduces to
\begin{equation}\label{82bis}
 \dR^{\gamma,\textbf{A}}e^\gou =
\textbf{F}^\gog +
 \frac{1}{2}[e^\gog\wedge e^\gog]
 = \frac{1}{2} \textbf{F}{^\gog}_{\cels\cels}
 e^{\cels\cels}
 + \frac{1}{2}\textbf{c}{^\gog}_{\celg\celg}e^{\celg\celg}
\end{equation}
where we used the notation
$\textbf{c}{^\gog}_{\gog\gog}$ introduced in (\ref{cgggref}).
Hence
\[
 \begin{array}{ccl}
  \dR^{\gamma,\textbf{A}}e^{(N-2)}_{\gou_1\gou_2}
  & = & \dR^{\gamma,\textbf{A}}e^{\celu}\wedge e^{(N-3)}_{\gou_1\gou_2\celu}
  = \left(\frac{1}{2} \textbf{F}{^\celg}_{\cels_1\cels_2}
 e^{\cels_1\cels_2}
 + \frac{1}{2}\textbf{c}{^\celg}_{\celg_1\celg_2}e^{\celg_1\celg_2} \right) \wedge  e^{(N-3)}_{\gou_1\gou_2\celg}
 \\
 & = &
 \textbf{F}{^\celg}_{\gou_1\gou_2} e^{(N-1)}_{\celg}
 + \textbf{c}{^\celg}_{\gou_1\gou_2}e^{(N-1)}_{\celg}
 + \textbf{c}{^\celg}_{\celg\gou_1}e^{(N-1)}_{\gou_2}
 + \textbf{c}{^\celg}_{\gou_2\celg}e^{(N-1)}_{\gou_1}\\
 & = &
 \textbf{F}{^\celg}_{\gou_1\gou_2} e^{(N-1)}_{\celg}
 + \textbf{c}{^\celg}_{\gou_1\gou_2}e^{(N-1)}_{\celg}
 \end{array}
\]
where we used the hypothesis that $\gog$ is unimodular, i.e. $\textbf{c}{^\celg}_{\gou\celg} = \textbf{c}{^\celg}_{\celg\gou} = 0$. Thus
\begin{equation}\label{dApimi2}
\frac{1}{2}
p{_\gog}^{\celu_1\celu_2}
\dR^{\gamma,\textbf{A}} e_{\celu_1\celu_2}^{(N-2)}
=
\frac{1}{2} \left(\textbf{F}{^{\celg}}_{\cels_1\cels_2}
 p{_\gog}^{\cels_1\cels_2}
 + \textbf{c}{^{\celg}}_{\celg_1\celg_2}p{_\gog}^{\celg_1\celg_2}\right)
 e_{\celg}^{(N-1)}
\end{equation}

By collecting (\ref{dApimi1}) and (\ref{dApimi2}) in (\ref{dApi}) we obtain
\begin{equation}\label{dApisynth}
\begin{array}{ccl}
 \displaystyle \dR^\textbf{A}p_\gog
 & = &  \displaystyle
 \left(\partial^{\gamma,\textbf{A}}_{\cels_1} p{_\gog}^{\cels\cels_1}
   + \partial_\celg p{_\gog}^{\cels\celg} \right)e_{\cels}^{(N-1)} \\
   & & +\ \displaystyle \left( \partial^{\gamma,\textbf{A}}_{\cels} p{_\gog}^{\celg\cels}
   + \partial_{\celg_1} p{_\gog}^{\celg\celg_1}
   + \frac{1}{2} \textbf{F}{^{\celg}}_{\cels_1\cels_2}
 p{_\gog}^{\cels_1\cels_2}
 + \frac{1}{2}\textbf{c}{^{\celg}}_{\celg_1\celg_2}p{_\gog}^{\celg_1\celg_2}
   \right)e_{\celg}^{(N-1)}
\end{array}
\end{equation}
Note however that it follows from (\ref{82bis}) that
$\dR e^\gog = \frac{1}{2}\textbf{c}{^\gog}_{\celg_1\celg_2}e^{\celg_1\celg_2} + \textbf{F}^\gog - [A^\gog\wedge e^\gog]$, which implies that $\dR e^{(N-2)}_{\gog_1\gog_2} = \dR e^{\celg}\wedge e^{(N-3)}_{\gog_1\gog_2\celg} = \textbf{c}{^{\celg}}_{\gog_1\gog_2}e^{(N-1)}_{\celg}$, thus
\[
 \dR\left(\frac{1}{2}p{_\gog}^{\celg_1\celg_2}
e_{\celg_1\celg_2}^{(N-2)}\right)
= \left(\partial_{\celg_1} p{_\gog}^{\celg\celg_1} + \frac{1}{2}\textbf{c}{^{\celg}}_{\celg_1\celg_2}p{_\gog}^{\celg_1\celg_2}
   \right)e_{\celg}^{(N-1)}
\]
Hence (\ref{dApisynth}) can be written as
\begin{equation}\label{dApYM}
 \dR^\textbf{A}p_\gog
= \left(\partial^{\gamma,\textbf{A}}_{\cels_1} p{_\gog}^{\cels\cels_1}
   + \partial_\celg p{_\gog}^{\cels\celg} \right)e_{\cels}^{(N-1)} + \left( \partial^{\gamma,\textbf{A}}_{\cels} p{_\gog}^{\celg\cels}
   + \frac{1}{2} \textbf{F}{^{\celg}}_{\cels_1\cels_2}
 p{_\gog}^{\cels_1\cels_2}   \right)e_{\celg}^{(N-1)} + \dR\left(\frac{1}{2}p{_\gog}^{\celg_1\celg_2}
e_{\celg_1\celg_2}^{(N-2)}\right)
\end{equation}

\subsubsection{Cancellation of the sources}
We come back to Equation (\ref{elymA}) ($\dR^\textbf{A}p_\gog = \frac{1}{2}|\textbf{F}{_\gog}^{\gos\gos}|^2 e_\gog^{(N-1)}$) which is equivalent to the fact that the r.h.s. of (\ref{dApYM}) is equal to $\frac{1}{2}|\textbf{F}{_\gog}^{\gos\gos}|^2 e_\gog^{(N-1)}$. By using (\ref{elymp1}) ($p{_\gog}^{\gos\gos} = - \textbf{F}{_\gog}^{\gos\gos}$) we deduce the two following equations
\begin{equation}\label{finalYMs}
\partial^{\gamma,\textbf{A}}_{\cels} \textbf{F}{_\gog}^{\gos\cels}
= \partial_\celg p{_\gog}^{\gos\celg}
\end{equation}
and
\begin{equation}\label{finalYMg}
(\partial^{\gamma,\textbf{A}}_{\cels}  p{_\gog}^{\celg\cels})
 e_{\celg}^{(N-1)}
   + \dR\left(\frac{1}{2}p{_\gog}^{\celg_1\celg_2}
e_{\celg_1\celg_2}^{(N-2)}\right)
= \frac{1}{2}|\textbf{F}{_\gog}^{\gos\gos}|^2 e_\gog^{(N-1)}
- \frac{1}{2} \textbf{F}{^{\celg}}_{\cels_1\cels_2}
 \textbf{F}{_\gog}^{\cels_1\cels_2}e_{\celg}^{(N-1)}
\end{equation}
Here comes the conclusion about (\ref{finalYMs}).
Let $(\textbf{t}^1,\cdots,\textbf{t}^r)$ be a basis of $\gog^*$ and set $\textbf{t}^{(r)}:= \textbf{t}^1\wedge \cdots\wedge \textbf{t}^r$, $(e^\gog)^{(r)}:= (e^\gog)^*\textbf{t}^{(r)}$ and
$(e^\gog)^{(r-1)}_\gog:= (e^\gog)^*\textbf{t}^{(r-1)}_\gog$.
A \textbf{first observation} here is that, for any $\textsf{x}\in \mathcal{X}$,
since the fiber $\mathcal{F}_\textsf{x}$ is \emph{compact},
the integration of both sides of
(\ref{finalYMs}) on $\mathcal{F}_\textsf{x}$ gives us (note that $\dR e_\gog^{(N-1)} = 0$ because $\gog$ is unimodular)
\[
\int_{\mathcal{F}_\textsf{x}}
\partial^{\gamma,\textbf{A}}_\cels \textbf{F}{_\gog}^{\gos\cels}
\  (e^\gog)^{(r)}
=  \int_{\mathcal{F}_\textsf{x}}\left(\partial_\celg p{_\gog}^{\gos \celg}\right)(e^\gog)^{(r)}
=  \int_{\mathcal{F}_\textsf{x}}\dR\left(p{_\gog}^{\gos \celg}\,(e^\gog)^{(r-1)}_\celg\right) = 0.
\]
A \textbf{second observation} is that the left hand side of (\ref{finalYMs})
is \emph{constant on any fiber}.
Thus, again since $\mathcal{F}_\textsf{x}$ is \emph{compact},
\begin{equation}\label{fundamentalELforYM1bis}
\partial^{\gamma,\textbf{A}}_\cels \textbf{F}{_\gog}^{\gos\cels}
=
\frac{\int_{\mathcal{F}_\textsf{x}} \partial^{\gamma,\textbf{A}}_\cels \textbf{F}{_\gog}^{\gos\cels}
 \ (e^\gog)^{(r)}}{\int_{\mathcal{F}_\textsf{x}} (e^\gog)^{(r)}} = 0.
\end{equation}
And this relation exactly means that $\textbf{A}^\gog$ is a solution of the (pure) Yang--Mills equations.

\subsubsection{A conservation law for the current}
Let us introduce the notation
\begin{equation}
 J{_\gog}^\gos:= \partial_{\celg}p{_\gog}^{\gos\celg}
\end{equation}
for the right hand side of (\ref{finalYMs}).
As seen previously (\ref{finalYMs}) implies that $J{_\gog}^\gos$ is constant on each fiber and hence is a function of $\textsf{x}\in \mathcal{X}$. However it may not vanish in the case where, in the previous model, $\goG$ is not compact, because (\ref{fundamentalELforYM1bis}) would not hold in general. However Equation (\ref{finalYMg}) still implies a conservation law on $J{_\gog}^\gos$, as shown by the following.
\begin{prop}
 Let $(p_\gog,\mathbf{F}_\gog)$ be a solution (\ref{finalYMg}). Then
 \begin{equation}\label{conservationlawforYM}
  \partial_\cels^\mathbf{A}J{_\gog}^\cels = 0
 \end{equation}
\end{prop}
\emph{Proof} --- By computing the exterior differential of both sides of (\ref{finalYMg}) and by using the facts that $\partial_\gog\textbf{F}{^\gog}_{\gos\gos} = 0$ and $\dR e^{(N-1)}_\gog = 0$
one obtains that
$\partial_\celg\left(\partial^{\gamma,\textbf{A}}_\cels p{_\gog}^{\celg\cels}\right) = 0$.
Recall that 
\[
 \partial_{\gos}^{\gamma,\mathbf{A}} p{_{\gog}}^{\gog\gos}
 = \partial_{\gos} p{_{\gog}}^{\gog\gos} - \mathbf{c}{^{\celg_1}}_{\celg_0\gog}\mathbf{A}{^{\celg_0}}_{\gos} p{_{\celg_1}}^{\gog\gos} 
 + \mathbf{c}{^{\gog}}_{\celg_0\celg_2}\mathbf{A}{^{\celg_0}}_{\gos}p{_{\gog}}^{\celg_2\gos} + \gamma{^\gos}_{\gos\cels}p{_{\gog}}^{\gog\cels}
\]
and hence, since $\gamma{^\cels}_{\cels\gos} = 0$ because the coefficients of $\gamma$ are in $so(\gos,\textsf{b})$,
\[
 \partial_{\cels}^{\gamma,\mathbf{A}} p{_{\gog}}^{\gog\cels}
 = \partial_{\cels} p{_{\gog}}^{\gog\cels} - \mathbf{c}{^{\celg_1}}_{\celg_0\gog}\mathbf{A}{^{\celg_0}}_{\cels} p{_{\celg_1}}^{\gog\cels} 
 + \mathbf{c}{^{\gog}}_{\celg_0\celg_2}\mathbf{A}{^{\celg_0}}_{\cels}p{_{\gog}}^{\celg_2\cels}
\]
Thus
\begin{equation}\label{etapepourcourantYM}
 \partial_\celg\left(\partial^{\gamma,\textbf{A}}_\cels p{_\gog}^{\celg\cels}\right) 
 = \partial_{\celg}\left(\partial_{\cels} p{_{\gog}}^{\celg\cels}\right)
 - \mathbf{c}{^{\celg_1}}_{\celg_0\gog}\mathbf{A}{^{\celg_0}}_{\cels} \left(\partial_{\celg}p{_{\celg_1}}^{\celg\cels} \right)
 + \mathbf{c}{^{\celg}}_{\celg_0\celg_2}\mathbf{A}{^{\celg_0}}_{\cels}\left(\partial_{\celg}p{_{\gog_1}}^{\celg_2\cels}\right)
\end{equation}
However by using Cartan's formula (\ref{82bis}) implies that
\[
 e^\gou([\partial_\gog,\partial_\gos])
= - \dR e^\gou(\partial_\gog,\partial_\gos)
= \gamma \wedge e^\gos(\partial_\gog,\partial_\gos) + [\textbf{A}^\gog\wedge e^\gog](\partial_\gog,\partial_\gos)
= -\textbf{c}{^\gog}_{\celg_0\celg_2}
\textbf{A}{^{\celg_0}}_\gos e^{\celg_2}(\partial_\gog)
\]
and hence $[\partial_\gog,\partial_\gos] = - \textbf{c}{^\celg}_{\celg_0\gog}
\textbf{A}{^{\celg_0}}_\gos\partial_{\celg}$, which implies 
\[
\partial_{\celg}\left(\partial_{\cels} p{_{\gog}}^{\celg\cels}\right) = \partial_{\cels}\left(\partial_{\celg} p{_{\gog}}^{\celg\cels}\right)- \textbf{c}{^\celg}_{\celg_0\celg_2}
\textbf{A}{^{\celg_0}}_\gos \left( \partial_{\celg}p{_{\gog}}^{\celg_2\cels}\right) 
\]
This leads to the following simplification in (\ref{etapepourcourantYM})
\[
 \partial_\celg\left(\partial^{\gamma,\textbf{A}}_\cels p{_\gog}^{\celg\cels}\right) 
 = \partial_{\cels}\left(\partial_{\celg} p{_{\gog}}^{\celg\cels}\right)
 - \mathbf{c}{^{\celg_1}}_{\celg_0\gog}\mathbf{A}{^{\celg_0}}_{\cels} \left(\partial_{\celg}p{_{\celg_1}}^{\celg\cels} \right)
 = \partial_{\cels}J{_\gog}^\cels
 - \mathbf{c}{^{\celg_1}}_{\celg_0\gog}\mathbf{A}{^{\celg_0}}_{\cels} J{_{\celg_1}}^\cels
\]
The right hand side of the latter equation is equal to $\partial_{\cels}^\mathbf{A} J{_\gog}^\cels$. Since we know from the beginning that the left hand side is zero, we deduce (\ref{conservationlawforYM}).\hfill $\square$

\subsubsection{Standard gauge symmetries}
The Yang--Mills action (\ref{actionYMnaivegeneral}) is invariant by several types of gauge symmetries, which generalizes the gauge symmetries of the Maxwell model seen previously.

For any $g\in \mathcal{C}^\infty(\mathcal{X},\goG)$ the action
$\mathcal{A}[\theta^\gog,\pi_\gog] = \int_\mathcal{F} \frac{1}{2}
|\pi{_\gog}^{\gos\gos}|^2\beta^{(n)}\wedge \theta^{(r)}
+ \pi_{\celg}\wedge
(\dR\theta^{\celg} + \frac{1}{2}[\theta^\gog\wedge\theta^\gog]^{\celg})$
is invariant by the gauge transformation
\[
 \left\{
 \begin{array}{ccl}
  \theta^\gog & \longmapsto & \hbox{Ad}_g\theta^\gog -\dR g\,g^{-1}
  = g\theta^\gog g^{-1} -\dR g\,g^{-1} \\
  \pi_\gog & \longmapsto & \hbox{Ad}_g^*\pi_\gog
 \end{array}\right.
\]
meaning that
\begin{equation}\label{standardgaugesymmetry1}
\mathcal{A}[\hbox{Ad}_g\theta^\gog -\dR g\,g^{-1},\hbox{Ad}_g^*\pi_\gog]
= \mathcal{A}[\theta^\gog,\pi_\gog]
\end{equation}
Indeed on the one hand since the scalar product $\textsf{k}$ on $\gog$
is invariant by the adjoint action of $\goG$,
we have
$|\hbox{Ad}^*_g\otimes 1_\gos\otimes 1_\gos (\pi{_\gog}^{\gos\gos})|^2
= |\pi{_\gog}^{\gos\gos}|^2$.
On the other hand the relations
\[
 \dR\left(\hbox{Ad}_g\theta^\gog -\dR g\,g^{-1}\right)
+ \frac{1}{2}\left[\left(\hbox{Ad}_g\theta^\gog -\dR g\,g^{-1}\right)\wedge\left(\hbox{Ad}_g\theta^\gog -\dR g\,g^{-1}\right)\right]
=\hbox{Ad}_g\left(\dR\theta^\gog + \frac{1}{2}[\theta^\gog\wedge \theta^\gog]\right)
\]
and
$\hbox{Ad}_g^*\pi_{\celg}\wedge \hbox{Ad}_g(\dR\theta^{\celg} + \frac{1}{2}[\theta^\gog\wedge \theta^\gog]^{\celg})
= \pi_{\celg}\wedge (\dR\theta^\celg + \frac{1}{2}[\theta^\gog\wedge \theta^\gog]^\celg)$
imply that the integral
$\int_\mathcal{F} \pi_{\celg}\wedge
(\dR\theta^{\celg}
+ \frac{1}{2}[\theta^\gog\wedge\theta^\gog]^{\celg})$ is invariant by this transformation.
Hence (\ref{standardgaugesymmetry1}) follows.

\subsubsection{Gauge symmetries of the dual fields}\label{paragraphgaugesymmetryYMpi}
Let $\chi_\gog \in \gog^*\otimes \Omega^{N-2}(\mathcal{F})$ and assume that we replace $\pi_\gog$ by
$\pi_\gog + \chi_\gog$. Then by observing that
$\Theta^\gog:= \dR\theta^\gog + \frac{1}{2}[\theta^\gog\wedge \theta^\gog] = \dR^{\theta/2}\theta^\gog$
and by using (\ref{twistedLeibniz})
\[
 (\pi_{\celg} + \chi_{\celg})\wedge \Theta^{\celg}
= \pi_{\celg}\wedge \Theta^{\celg} + (\dR^{\theta/2}\theta^\gog)\wedge \chi_{\celg}
=  \pi_{\celg}\wedge \Theta^{\celg} + \dR^{\theta/2}\left(\theta^{\celg} \wedge
\chi_{\celg}\right)
+ \theta^{\celg} \wedge \dR^{\theta/2}\chi_{\celg}
\]
But since $\theta^{\celg} \wedge \chi_{\celg}$ has real coefficients (hence in a trivial representation of $\gog$), we have actually
$\dR^{\theta/2}\left(\theta^{\celg} \wedge \chi_{\celg}\right)
 = \dR\left(\theta^{\celg} \wedge \chi_{\celg}\right)$, so that
\begin{equation}\label{74}
 (\pi_{\celg} + \chi_{\celg})\wedge \Theta^{\celg}
=  \pi_{\celg}\wedge \Theta^{\celg} + \dR\left(\theta^{\celg} \wedge
\chi_{\celg}\right)
+ \theta^{\celg} \wedge \dR^{\theta/2}\chi_{\celg}
\end{equation}
Assume further that
\begin{equation}\label{hypogaugesymmetryYM}
 \chi_\gog \wedge \beta^{\mathfrak{ss}} = 0,
\end{equation}
i.e. $\chi_\gog$ decomposes as
$\chi_\gog = \chi{_\gog}^{\cels\celg} f^{(N-2)}_{\cels\celg}
+ \frac{1}{2}\chi{_\gog}^{\celg\celg} f^{(N-2)}_{\celg\celg}$ or
$\chi{_\gog}^{\mathfrak{ss}} = 0$. Then
$|\pi{_\gog}^{\mathfrak{ss}} + \chi{_\gog}^{\mathfrak{ss}}|^2 =
|\pi{_\gog}^{\mathfrak{ss}}|^2$.
Hence if we assume that
$\chi_\gog \in \gog^*\otimes \Omega^{N-2}(\mathcal{F})$ satisfies (\ref{hypogaugesymmetryYM}) and
\emph{decreases at infinity} (or is compactly supported),
so that
$\int_\mathcal{F}\dR\left(\theta^{\celg} \wedge
\chi_{\celg}\right) = 0$, it follows then from (\ref{74}) that
\[
 \mathcal{A}[\theta^\gog,\pi_\gog+ \chi_\gog]
 = \mathcal{A}[\theta^\gog,\pi_\gog]
 + \int_\mathcal{F} \theta^{\celg} \wedge \dR^{\theta/2}\chi_{\celg}
\]
Thus the action satisfies
$\mathcal{A}[\theta^\gog,\pi_\gog+ \chi_\gog]
 = \mathcal{A}[\theta^\gog,\pi_\gog]$
if $\chi_\gog \in \gog^*\otimes \Omega^{N-2}(\mathcal{F})$
satisfies (\ref{hypogaugesymmetryYM}),
 \emph{decreases at infinity}
and satisfies
$\theta^{\celg} \wedge \dR^{\theta/2}\chi_{\celg} = 0$.
As a conclusion:
\begin{lemm}
 Let $\chi_\gog \in \gog^*\otimes \Omega^{N-2}(\mathcal{F})$. Assume that
\begin{enumerate}
 \item $\chi_\gog$ decays at infinity or has compact support;
 \item $\chi{_\gog}^{\mathfrak{ss}} = 0$, i.e. $\chi_\gog$ decomposes as
 \begin{equation}\label{contraintesurchi}
  \chi_\gog = \chi{_\gog}^{\cels\celg} f^{(N-2)}_{\cels\celg}
+ \frac{1}{2}\chi{_\gog}^{\celg\celg} f^{(N-2)}_{\celg\celg}
\end{equation}
\item
\begin{equation}\label{conditiononchiYM}
 \theta^{\celg} \wedge \emph{\dR}^{\theta/2}\chi_{\celg}  = 0
\end{equation}
\end{enumerate}
then we have $\mathcal{A}[\theta^\gog,\pi_\gog+ \chi_\gog]
 = \mathcal{A}[\theta^\gog,\pi_\gog]$.
\end{lemm}
Note that Condition (\ref{contraintesurchi}) is actually sufficient for $\chi_\gog$ to be an \emph{on shell} gauge symmetry. Indeed if the Euler--Lagrange equations (\ref{elymp}) are satisfied then
$\Theta^\gog = \frac{1}{2}\Theta{^\gog}_{\cels\cels}\theta^{\cels\cels}$ and thus the action $\int_\mathcal{F}\frac{1}{2}|\pi{_\gog}^{\gos\gos}|^2 + \pi_{\celg}\wedge \Theta^{\celg}$ is obviously invariant by the transformation $(\theta^\gog,\pi_\gog)\longmapsto (\theta^\gog,\pi_\gog+\chi_\gog)$ if $\chi_\gog $ satisfies (\ref{contraintesurchi}).

\subsubsection{Invariance by fiber bundle diffeomorphisms}

Let $T:\mathcal{F}\longrightarrow\mathcal{F}$
be a diffeomorphism such that $P\circ T=P$ (i.e. which preserves each fiber of the fibration $P:\mathcal{F}\longrightarrow\mathcal{X}$). Then
our action enjoys the symmetry
\begin{equation}\label{verticaldiffeo}
 \mathcal{A}[T^*\theta^\gog,T^*\pi_\gog]
 = \mathcal{A}[\theta^\gog,\pi_\gog]
\end{equation}
Indeed remind that
$\pi{_\gog}^{\gos\gos} = \textsf{Q}{_\gog}^{\gos\gos}(\theta^\gog,\pi_\gog)$ (\ref{Qmap0})
is characterized by
$\pi{_\gog}^{\gos\gos}\beta^{(n)}\wedge \theta^{(r)} = \pi_\gog\wedge \beta^{\gos\gos}$   (\ref{implicitpimumnu0}). Hence since $\beta^{(n)}$ and $\beta^{\gos\gos}$ are invariant by $T^*$, this implies
$(\pi{_\gog}^{\gos\gos}\circ T)\beta^{(n)}\wedge T^*\theta^{(r)} = (T^*\pi_\gog)\wedge \beta^{\gos\gos}$, so that $\pi{_\gog}^{\gos\gos}\circ T$ satisfies the same relation as
$\textsf{Q}{_\gog}^{\gos\gos}(T^*\theta^\gog,T^*\pi_\gog)$. Hence
$\textsf{Q}{_\gog}^{\gos\gos}(T^*\theta^\gog,T^*\pi_\gog) = \textsf{Q}{_\gog}^{\gos\gos}(\theta^\gog,\pi_\gog)\circ T$. It follows that
\[
\begin{array}{l}
 \frac{1}{2}\left|\textsf{Q}{_\gog}^{\gos\gos}(T^*\theta^\gog,T^*\pi_\gog)\right|^2 \beta^{(n)}\wedge
 T^*\theta^{(r)}
 + T^*\pi_\celg\wedge
 \left(\dR T^*\theta^\celg + \frac{1}{2}
 [T^*\theta^\gog\wedge T^*\theta^\gog]^\celg\right) \\
 = T^*\left[\frac{1}{2}\left|\textsf{Q}{_\gog}^{\gos\gos}(\theta^\gog,\pi_\gog)\right|^2 \beta^{(n)}\wedge
 \theta^{(r)}
 + \pi_\celg\wedge
 \left(\dR\theta^\celg + \frac{1}{2}
 [\theta^\gog\wedge \theta^\gog]^\celg\right)\right]
\end{array}
\]
and by integration over $\mathcal{F}$ we deduce
(\ref{verticaldiffeo}).

\section{Kaluza--Klein theories}\label{SectionKK}
A Kaluza--Klein action functional can be obtained by adding a quantity of the kind $\int \pi_\celg\wedge \Theta^\celg$, where $\Theta^\gog:= \dR \theta^\gog + \frac{1}{2}[\theta^\gog\wedge\theta^\gog]$, from the Yang--Mills action to a higher dimensional version of the Palatini action functional as defined in \S \ref{paragraphPalatinidebase}.

Starting from the Palatini action described in \S \ref{paragraphPalatinidebase} we replace $\gos$ by a larger space $\gou:= \gos\oplus \gog$, where $(\gog,[\cdot,\cdot])$ is a Lie algebra of dimension $r$ and, in the role of $\gol$, we replace $so(\gos,\textsf{b})$ by $so(\gou,\textsf{h})$. Hence
\[
\gou:= \gos\oplus \gog\quad
\hbox{ and }\quad
\gol = so(\gou,\textsf{h})
\]
so that $\hbox{dim}\gou = N:= n+r$.
We extend the Lie bracket of $\gog$ on $\gou$ in such a way that $\gos$ is in the center of $(\gou,[\cdot,\cdot])$. In a way similar to the Yang--Mills theory (see \S \ref{backtoYM}) we assume that $\gog$ is the Lie algebra of a simply connected Lie group $\widehat{\goG}$ (but not necessarily compact in the following). We also assume that $\gou$ is endowed with a symmetric nondegenerate bilinear form $\textsf{h}$ which is invariant by the adjoint action of $\widehat{\goG}$ (see (\ref{definitionActionderho})) and such that $\gos\perp \gog$. We denote by $\textsf{b}$ and $\textsf{k}$ the restriction of $\textsf{h}$ on, respectively, $\gos$ and $\gog$.

Let $\mathcal{Y}$ be a smooth oriented manifold of dimension $N$. The dynamical fields on $\mathcal{Y}$ will be a pair $(\theta^\gou,\varphi^\gol)$, where $\theta^\gou\in \gou\otimes \Omega^1(\mathcal{Y})$ and $\varphi^\gol\in \gol\otimes \Omega^1(\mathcal{Y})$, for the 'Palatini' part of the action plus an extra field $\pi_\gou\in \gou^*\otimes\Omega^{N-2}(\mathcal{Y})$ which satisfies the constraint $\theta^{\gos\gos}\wedge \pi_\gou = 0$.
Hence the space of fields is:
\[
\begin{array}{r}
\mathscr{E}:=
 \{\, (\theta^\gou,\varphi^\gol,\pi_\gou)\in \left(\gou \otimes\Omega^1(\mathcal{Y})\right)\times \left(\gol \otimes\Omega^1(\mathcal{Y})\right) \times \left(\gou^*\otimes\Omega^{N-2}(\mathcal{Y})\right)\,; \\
(\theta^\gou,\varphi^\gol,\pi_\gou)\hbox{ are of class }\mathscr{C}^2\hbox{ and }
\theta^\mathfrak{s}\wedge \theta^\mathfrak{s} \wedge \pi_\gou =0 \}
\end{array}
\]
We let $\kappa{_\gol}^{\gou\gou}\in \gol^*\otimes\gou\otimes\gou$
be defined as in (\ref{kappasimple}) (this tensor is invariant by $\hbox{Ad}_{\widehat{\goG}}$ as expounded in \S \ref{paragraph1.2.5}) and we set
$\Phi^\gol:= \dR\varphi^\gol+
 \frac{1}{2}[\varphi^\gol\wedge \varphi^\gol]$ and, for shortness,
 $\Phi^{\gou\gou}:= \kappa{_\cell}^{\gou\gou}\Phi^\cell$.
Then with the same conventions as before, we define on $\mathscr{E}$ the action functional $\mathscr{A}$ by:
\begin{equation}\label{actionKK}
 \mathscr{A}[\theta^\gou,\varphi^\gol,\pi_\gou]:= \int_{\mathcal{Y}}  \pi_\celu \wedge \Theta^\celu + \frac{1}{2}\theta_{\celu\celu}^{(N-2)} \wedge \Phi^{\celu\celu} - \Lambda_0 \theta^{(N)}
\end{equation}
\begin{theo}
\label{theorem0}
Let $\widehat{\goG}$ be a simply connected Lie group with Lie algebra $\gog$, of dimension $r$, and $\gos$ be a vector space of dimension $n$.
Let $\gou = \gos\oplus \gog$, endowed with the Lie bracket $[\cdot,\cdot]$ which extends the Lie algebra structure on $\gog$ in such a way such that $\gos$ belongs to the center of $(\gou,[\cdot,\cdot])$. Assume that $\gou$ is endowed with a symmetric bilinear form $\emph{\textsf{h}} = \emph{\textsf{b}}\oplus \emph{\textsf{k}}$ which is invariant by the adjoint action of $\gog$ on $\gou$ and such that $\gos\perp \gog$.

Let $\mathcal{Y}$ be a connected, oriented manifold of dimension $N:= n+r$.
Let $(\theta^\gou,\varphi^\gol,\pi_\gou)\in \mathscr{E}$ be a critical point of $\mathscr{A}$ and let $\emph{\textbf{h}}:= (\theta^\gou)^*\emph{\textsf{h}}$. Assume that the rank of $\theta^\gou$ is equal to $N$ everywhere and that $(\mathcal{Y},\theta^\gos,\theta^\gog)$ is $\gog$-\emph{complete} (see Definition \ref{defigcomplete}).  
Then
\begin{enumerate}
 \item the exterior differential system $\theta^\gos|_\textsf{\emph{f}}=0$, for $r$-dimensional submanifolds $\textsf{\emph{f}}\subset \mathcal{Y}$, is completely integrable and $\mathcal{Y}$ is foliated by the integral leaves $\textsf{\emph{f}}$.
 \item there exists a Lie group $\goG$, which is a quotient of $\widehat{\goG}$ by a finite subgroup such that all integral leaves $\textsf{\emph{f}}$ are diffeomorphic to $\goG$.
\end{enumerate}
Assume the \textbf{additional hypothesis that $\goG$ is compact}. Then the foliation forms actually a fibration and the following holds.
\begin{enumerate}
\item[(iii)] the manifold $\mathcal{Y}$ acquires the structure of a principal bundle over
an $n$-dimensional manifold $\mathcal{X}$ with structure group $\mathfrak{G}$:
\[
   \goG\longrightarrow \mathcal{Y}
   \xrightarrow{\hbox{ }P\hbox{ }}
   \mathcal{X}
\]
  \item[(iv)] $\emph{\textbf{g}}:= (\theta^\gos)^*\emph{\textsf{b}} =  \emph{\textsf{b}}_{\cels\cels}\theta^\cels\otimes \theta^\cels$ is constant on each fiber of $P$ and induces a pseudo metric (also denoted by) $\emph{\textbf{g}}$ on $\mathcal{X}$;
  \item[(v)] in any local trivialization
  $\mathcal{Y}_U \simeq U\times \goG$ (where $U\subset \mathcal{X}$ is an open subset and $\mathcal{Y}_U:= {P}^{-1}(U)$) we can write $\theta^\gog = g^{-1}\emph{\textbf{A}}^\gog g + g^{-1}\dR g$, where $g\in \goG$ and $\emph{\textbf{A}}^\gog$ is depends only on $\emph{\textsf{x}}\in \mathcal{X}$;
\item[(vi)] $\emph{\textbf{g}}$ and $\emph{\textbf{A}}^\gog$ are solutions of the Einstein--Yang--Mills system
  \[
   \left\{
 \begin{array}{ccl}
  \mathbf{R}(\emph{\textbf{g}}){^\gos}_\gos
  - \frac{1}{2}\mathbf{R}\delta{^\gos}_\gos
  + \Lambda\delta{^\gos}_\gos & = &
\frac{1}{2}\emph{\textbf{F}}{_\celg}^{\gos\cels}\emph{\textbf{F}}{^\celg}_{\gos\cels}
-\frac{1}{4} |\emph{\textbf{F}}|^2\delta{^\gos}_\gos\\
\nabla^{T\mathcal{X},\emph{\textbf{A}}}_\cels \emph{\textbf{F}}_{\gog}{^{\gos\cels}} & = & 0
 \end{array}\right.
\]
with a
cosmological constant equal to $\Lambda = \Lambda_0 + \frac{1}{4}\langle \textsf{\emph{B}},\textsf{\emph{k}}\rangle$, where $\textsf{\emph{B}}_{\gog\gog}:= \textbf{\emph{c}}{^{\celg_1}}_{\celg_2\gog} \textbf{\emph{c}}{^{\celg_2}}_{\celg_1\gog}$ is the Killing form on $\gog$ and
$\langle \textsf{\emph{B}},\textsf{\emph{k}}\rangle:=
 \frac{1}{2}\textsf{\emph{B}}_{\celg\celg}\textsf{\emph{k}}^{\celg\celg}$.
 \end{enumerate}
\end{theo}
A straightforward corollary of Theorem \ref{theorem0} is the following.
\begin{coro}
\label{coroobvious}
Assume exactly the same Hypotheses as in Theorem \ref{theorem0} and, in addition, that $\widehat{\goG}$ is compact. Then Conclusions (iii) to (vi) in Theorem \ref{theorem0} hold.
\end{coro}
\begin{rema}
 One may replace $\kappa{_\gol}^{\gou\gou}$
given by (\ref{kappasimple}) by any tensor which is invariant by $\hbox{\emph{Ad}}_{\widehat{\goG}}$, as expounded in \S \ref{paragraph1.2.5}, and such that the map $\gou^*\otimes \gou^*\ni \xi_{\gou\gou}\longmapsto
\kappa{_\gol}^{\celu\celu}\xi_{\celu\celu}$ has a \emph{non trivial} kernel. Most computations still holds, however the interpretation of the resulting system of equations would be different.
\end{rema}

\begin{rema}
 The action $\mathscr{A}$ in (\ref{actionKK}) and the constraint $\pi_\gou\wedge \theta^{\gos\gos} = 0$ are obviously invariant under the action
 $(\theta^\gou,\varphi^\gol,\pi_\gou) \longmapsto (T^*\theta^\gou,T^*\varphi^\gol,T^*\pi_\gou)$
 of orientation preserving diffeomorphisms $T:\mathcal{Y}\longrightarrow\mathcal{Y}$.
 It is also invariant through the transformation $(\theta^\gou,\varphi^\gol,\pi_\gou) \longmapsto
 (\emph{Ad}_g\theta^\gou,(\emph{Ad}_g)\varphi^\gol(\emph{Ad}_g)^{-1},\emph{Ad}^*_g\pi_\gou)$, where $g\in \goG$ is \textbf{constant}. However there is apparentely no way to extend this finite symmetry to a gauge group  action, because the curvature form $\Phi^\gol = \hbox{\emph{d}}\varphi^\gol + \frac{1}{2}[\varphi^\gol\wedge\varphi^\gol]$ does not transform in a simple way.
\end{rema}

The next sections are devoted to the proof of Theorem \ref{theorem0}.

\subsection{The Euler--Lagrange equations}
In the following we assume that $(\theta^\gou,\varphi^\gol,\pi_\gou) \in \mathscr{E}$
is a critical point of $\mathscr{A}$ such that $\hbox{rank}\theta^\gou = N$.
We denote by $\textbf{h} = \textsf{b}_{\cels\cels}\theta^\cels\theta^\cels
+ \textsf{k}_{\celg\celg}\theta^\celg\theta^\celg$ the induced metric on $\mathcal{Y}$
and we assume that $(\mathcal{Y},\theta^\gos,\theta^\gog)$ is  $\gog$-complete. Recall that $\Theta^\gou:= \hbox{d}\theta^\gou + \frac{1}{2}[\theta^\gou\wedge \theta^\gou]$
and $\Phi^\gol:= d\varphi^\gol + \frac{1}{2}[\varphi^\gol\wedge \varphi^\gol]$ and the a priori decompositions
$\Theta^\gou = \frac{1}{2}\Theta{^\gou}_{\celu\celu}\,\theta^{\celu\celu}
 = \frac{1}{2}\Theta{^\gou}_{\cels\cels}\,\theta^{\cels\cels} + \Theta{^\gou}_{\cels\celg}\,\theta^{\cels\celg}
 + \frac{1}{2}\Theta{^\gou}_{\celg\celg}\,\theta^{\celg\celg}$
and $\pi_\gou = \frac{1}{2}\pi{_\gou}^{\celu\celu}\,\theta^{(N-2)}_{\celu\celu}
 = \frac{1}{2}\pi{_\gou}^{\cels\cels}\,\theta^{(N-2)}_{\cels\cels} +
\pi{_\gou}^{\cels\celg}\,\theta^{(N-2)}_{\cels\celg}
 + \frac{1}{2}\pi{_\gou}^{\celg\celg}\,\theta^{(N-2)}_{\celg\celg}$.
The constraint $\pi_\gou\wedge \theta^{\gos\gos} = 0$ in the definition of $\mathscr{E}$ then reads $\pi{_\gou}^{\gos\gos} = 0$ or
\begin{equation}\label{constraint1}
 \pi_\gou  =
\pi{_\gou}^{\cels\celg}\,\theta^{(N-2)}_{\cels\celg}
 + \frac{1}{2}\pi{_\gou}^{\celg\celg}\,\theta^{(N-2)}_{\celg\celg}
\end{equation}

\subsubsection{Study of the first variation}

\noindent
\textbf{First variation with respect to coefficients of $\pi_\gou$} ---
We write that the action functional is stationary with respect to first order variations $(\theta^\gou,\varphi^\gol,\pi_\gou)\longmapsto (\theta^\gou,\varphi^\gol,\pi_\gou + \varepsilon \delta\pi_\gou)$, where $\delta \pi_\gou  = \chi_\gou =
\chi{_\gou}^{\cels\celg}\,\theta^{(N-2)}_{\cels\celg}
 + \frac{1}{2}\chi{_\gou}^{\celg\celg}\,\theta^{(N-2)}_{\celg\celg}$, so that it respects (\ref{constraint1}). It gives us:
\[
 \forall \chi{_\gou}^{\gos\gol}, \chi{_\gou}^{\gol\gol},\quad
0 = \int_\mathcal{Y} \chi_\celu \wedge \Theta{^\celu}
= \int_\mathcal{Y} \left(\chi{_\celu}^{\cels\celg}\,\Theta{^\celu}_{\cels\celg}
 + \frac{1}{2}\chi{_\celu}^{\celg\celg}\,\Theta{^\celu}_{\celg\celg}\right)\theta^{(N)}
\]
This is equivalent to $\Theta{^\gou}_{\gos\gog} = \Theta{^\gou}_{\gog\gog} = 0$. Hence $\Theta^\gou
 = \frac{1}{2}\Theta{^\gou}_{\cels\cels}\,\theta^{\cels\cels}$, which reads
\begin{equation}\label{ELFrobeniusdetaillee}
 \left\{
 \begin{array}{ccl}
\hbox{d}\theta^\gos & = &  \frac{1}{2}\Theta{^\gos}_{\cels\cels}\,\theta^\cels\wedge \theta^\cels \\
\hbox{d}\theta^\gog +\frac{1}{2}[\theta^\gog\wedge \theta^\gog] & = &  \frac{1}{2}\Theta{^\gog}_{\cels\cels}\,\theta^\cels\wedge \theta^\cels
 \end{array}\right.
\end{equation}

\noindent
\textbf{First variation with respect to $\varphi^\gol$} ---
We look at first order variations
$(\theta^\gou,\varphi^\gol,\pi_\gou)\longmapsto (\theta^\gou,\varphi^\gol + \varepsilon \delta\varphi^\gol,\pi_\gou)$, for any
$\delta\varphi^\gol = \lambda^\gol$ with compact support. This induces the transformation $\Phi^{\gou\gou} \longmapsto \Phi^{\gou\gou} + \dR^\varphi \lambda^{\gou\gou}$ (where $\lambda^{\gou\gou} := \kappa{_\cell}^{\gou\gou} \lambda^\cell$) since $\kappa{_\gol}^{\gou\gou}$ is constant and $\hbox{ad}_\gol$-invariant. This yields the condition that, $\forall \lambda^\gol$ with compact support,
\[
 0 = \int_\mathcal{Y}
 \frac{1}{2}\dR^\varphi\lambda^{\celu\celu} \wedge\theta^{(N-2)}_{\celu\celu}
 = \int_\mathcal{Y}
 \frac{1}{2} \dR^\varphi \left(
 \lambda^{\celu\celu}  \wedge
 \theta^{(N-2)}_{\celu\celu}
 \right)
 + \frac{1}{2} \lambda^{\celu\celu} \wedge \dR^\varphi \theta^{(N-2)}_{\celu\celu}
\]
However $\dR^\varphi \left(
 \lambda^{\celu\celu} \wedge
 \theta^{(N-2)}_{\celu\celu}
 \right)
 = \dR \left(
 \lambda^{\celu\celu} \wedge
\theta^{(N-2)}_{\celu\celu}
 \right)$, since the adjoint action of $\gol$ on this quantity is trivial. We thus obtain the condition
\[
 0 = \int_\mathcal{Y} \dR \left(
 \frac{1}{2}\theta^{(N-2)}_{\celu\celu}\wedge \lambda^{\celu\celu}
 \right)
 + \frac{1}{2} \lambda^{\celu\celu}\dR^\varphi \theta^{(N-2)}_{\celu\celu}
\]
from which deduce (since $\gou^*\otimes \gou^*\ni \xi_{\gou\gou}\longmapsto
\kappa{_\gol}^{\celu\celu}\xi_{\celu\celu}$ has a trivial kernel) that $\dR^\varphi \theta^{(N-2)}_{\gou\gou} = 0$. Lastly since $\dR^\varphi \theta^{(N-2)}_{\gou\gou} = \dR^\varphi \theta^\celu \wedge \theta^{(N-3)}_{\gou\gou\celu}$ and $N>2$ we deduce $\dR^\varphi \theta^\gou = 0$ (a similar result is derived in (\ref{inverseThetaChapeau})).
This means that the connection on $T\mathcal{Y}$ defined by $\varphi^\gol$ is torsion free, i.e. coincides with
the Levi-Civita connection of $(\mathcal{Y},\textbf{h})$.\\

\noindent
\textbf{First variation with respect to $\theta^\gou$} ---
Lastly we look at variations $(\theta^\gou,\varphi^\gol,\pi_\gou)\longmapsto (\theta^\gou + \varepsilon\delta\theta^\gou, \varphi^\gol, \pi_\gou + \varepsilon\delta\pi_\gou)$, for any $\delta\theta^\gou = \tau^\gou$ with compact support, where $\delta\pi_\gou$ is chosen in such a way that the coefficients $\pi{_\gou}^{\gou\gou}$ are fixed (in particular we preserve the constraint $\pi_\gou\wedge \theta^{\gos\gos} = 0$).
Through these variations,
\[
 \left\{
 \begin{array}{cclll}
 \pi_\gou & \longmapsto & \pi_\gou + \varepsilon\delta\pi_\gou + o(\varepsilon) & \hbox{with} & \delta\pi_\gou = \pi{_\gou}^{\cels\celg}\
 \tau^\celu \wedge \theta^{(N-3)}_{\cels\celg\celu}
 + \frac{1}{2}
 \pi{_\gou}^{\celg\celg}\
 \tau^\celu \wedge \theta^{(N-3)}_{\celg\celg\celu} \\
  \theta_{\gou\gou}^{(N-2)} & \longmapsto &\theta_{\gou\gou}^{(N-2)} + \varepsilon\delta\theta_{\gou\gou}^{(N-2)} + o(\varepsilon) & \hbox{with} & \delta\theta_{\gou\gou}^{(N-2)} = \tau^\celu\wedge \theta^{(N-3)}_{\gou\gou\celu} \\
  \Theta^\gou & \longmapsto & \Theta^\gou + \varepsilon\delta\Theta^\gou+ o(\varepsilon) & \hbox{with} & \delta\Theta^\gou = \dR^\theta \tau^\gou \\
  \theta^{(N)} & \longmapsto & \theta^{(N)} + \varepsilon\delta\theta^{(N)}+ o(\varepsilon) & \hbox{with} & \delta\theta^{(N)} = \tau^\celu\wedge \theta^{(N-1)}_\celu
 \end{array}
\right.
\]
Hence in particular, on the one hand, by using (\ref{ELFrobeniusdetaillee}) and (\ref{fifi})
\[
\begin{array}{ccl}
 \displaystyle \delta\pi_\celu \wedge \Theta{^\celu}
& = & \displaystyle \left(\pi{_\celu}^{\cels_1\celg_2}\
 \tau^{\celu_3} \wedge \theta^{(N-3)}_{\cels_1\celg_2\celu_3}
 + \frac{1}{2}
 \pi{_\celu}^{\celg_1\celg_2}\
 \tau^{\celu_3} \wedge \theta^{(N-3)}_{\celg_1\celg_2\celu_3}\right)
 \wedge \frac{1}{2}\Theta{^\celu}_{\cels\cels}\theta^{\cels\cels} \\
 & = & \tau^\celu \wedge \left(\pi{_{\celu_1}}^{\cels\celg} \Theta{^{\celu_1}}_{\celu\cels}\ \theta^{(N-1)}_{\celg}\right)
\end{array}
\]
On the other hand
\[
\pi_\celu\wedge \delta \Theta^\celu
= \dR^\theta \tau^\celu \wedge\pi_\celu
= \dR^\theta \left( \tau^\celu \wedge\pi_\celu\right) + \tau^\celu \wedge\dR^\theta\pi_\celu
= \dR \left( \tau^\celu \wedge\pi_\celu\right) + \tau^\celu \wedge\dR^\theta\pi_\celu
\]
Thus by using the fact that $\tau^\gou$ has a compact support, we deduce the condition
\[
 \forall \tau^\gou,\quad
 0 = \int_\mathcal{Y}
 \tau^\celu \wedge \left(
 \dR^\theta\pi_\celu -
 \Theta{^{\celu_1}}_{\celu\cels}\pi{_{\celu_1}}^{\celg\cels}\ \theta^{(N-1)}_{\celg}
 + \frac{1}{2}\theta^{(N-3)}_{\celu_1\celu_2\celu}\wedge \Phi^{\celu_1\celu_2}
- \Lambda_0 \theta^{(N-1)}_\celu\right)
\]
which gives us the equation
\begin{equation}\label{ELEinsteinbrut}
 \dR^\theta\pi_\gou +
 \frac{1}{2}\theta^{(N-3)}_{\celu\celu\gou}\wedge \Phi^{\celu\celu}
- \Lambda_0 \theta^{(N-1)}_\gou
= \Theta{^\celu}_{\gou\cels}\pi{_\celu}^{\celg\cels}\ \theta^{(N-1)}_{\celg}
\end{equation}

\subsection{Geometric consequences of the Euler--Lagrange equations}

\subsubsection{Existence of a foliation}

From the first equation in (\ref{ELFrobeniusdetaillee})
we deduce that $\hbox{d}\theta^\gos \!=\! 0\hbox{ mod}[\theta^\gos]$.
Since the rank of $\theta^\gos$ is equal to $n$ everywhere, we deduce from
Frobenius' theorem that $\mathcal{Y}$ is foliated by integral leaves $\textsf{f}$ which are solutions of the system $\theta^\gos|_\textsf{f} = 0$ of dimension $r$. We denote by $\mathcal{X}$ the set of integral leaves.

\subsubsection{The structure of the leaves}\label{traitementcompact}

Consider on the product manifold $\widehat{\goG} \times \mathcal{Y} = \{(h,\textsf{y})\in \widehat{\goG}\times \mathcal{Y}\}$ the $\gog$-valued 1-form
\[
 \tau^\gog := \theta^\gog - h^{-1}\hbox{d}h
\]
It satisfies the identity $\hbox{d}\tau^\gog = \hbox{d}\theta^\gog +\frac{1}{2}[\theta^\gog\wedge\theta^\gog]
- [\theta^\gog\wedge \tau^\gog] + \frac{1}{2}[\tau^\gog\wedge \tau^\gog]$ and its rank is clearly equal to $r$.
However the second equation in (\ref{ELFrobeniusdetaillee}) implies that,
for any integral leaf $\textsf{f}$,
$\hbox{d}\theta^\gog + \frac{1}{2}[\theta^\gog\wedge \theta^\gog]|_\textsf{f} = 0$
and thus
$\hbox{d}(\tau^\gog|_{\textsf{f}\times \goG}) = 0\hbox{ mod}[\tau^\gog]$.
Hence, again by Frobenius' theorem,
for any $(g_0,\textsf{y}_0)\in \widehat{\goG}\times \textsf{f}$,
there exists a unique $r$-dimensional submanifold $\Gamma\subset \widehat{\goG}\times\textsf{f}$ which is
a solution of $\tau^\gog|_\Gamma=0$  and which contains $(g_0,\textsf{y}_0)$.

As for the Yang--Mills theory, $\forall (g,\textsf{y})\in \widehat{\goG} \times \textsf{f}$, $\forall (\xi,v) \in T_{g}\widehat{\goG}\times T_{\textsf{y}}\textsf{f}$, the equation $g^{-1}\hbox{d}g(\xi) = \theta^\gog(v)$ defines the graph of a vector space isomorphism between $T_{g}\widehat{\goG}$ and $T_{\textsf{y}}\textsf{f}$. This implies that, around each point $(g,\textsf{y})\in \Gamma$, $\Gamma$ is locally the graph of a diffeomorphism between a neighbourhood of $g$ in $\widehat{\goG}$ and a neighbourhood of $\textsf{y}$ in $\textsf{f}$.

But we have more: since $(\mathcal{Y},\theta^\gos,\theta^\gog)$ is $\gog$-complete by hypothesis, by applying Lemma \ref{lemmaGcoversleaves} we deduce that $\widehat{\goG}$ is a universal cover of each leaf $\textsf{f}$.

\subsubsection{All integral leaves are diffeomorphic}\label{traitementnoncompact}
In the following result we still assume the hypotheses of Theorem \ref{theorem0}.
\begin{lemm}\label{lemmaifGnoncompact}
Assume that $(\mathcal{Y},\theta^\gos,\theta^\gog)$ is $\gog$-complete and that $\mathcal{Y}$ is connected. Then for any pair $\textsf{\emph{f}}_0$, $\textsf{\emph{f}}_1$ of integral leaves, $\textsf{\emph{f}}_0$ is diffeomorphic to $\textsf{\emph{f}}_1$.
\end{lemm}
\emph{Proof} Let $\check{\textsf{f}}\subset \mathcal{Y}$ be any fixed integral leaf and consider
\[
 \check{\mathcal{Y}}:= \{ \textsf{y}\in \mathcal{Y}\ ;\ \hbox{ the leaf which contains }\textsf{y}\hbox{ is diffeomorphic to }\check{\textsf{f}}\}
\]
We will show that $\check{\mathcal{Y}}$
is open and closed. It is clear that $\check{\mathcal{Y}}\neq \emptyset$ since $\check{\textsf{f}} \subset\check{\mathcal{Y}}$. Thus it will prove that $\check{\mathcal{Y}} = \mathcal{Y}$ since $\mathcal{Y}$ is connected. \\

\noindent
\textbf{(i)} We first prove that $\check{\mathcal{Y}}$ is open.
Let $\textsf{y}_0\in \check{\mathcal{Y}}$ and let us denote by $\textsf{f}_0$ the leaf which contains $\textsf{y}_0$ (which is hence diffeomorphic to $\check{\textsf{f}}$).

To any \emph{fixed} $\xi^\gou\in \gou$ we associated the vector field $\textsf{X}(\xi^\gou)$ on $\mathcal{Y}$
defined by $\textsf{X}(\xi^\gou) = \xi^\celu\frac{\partial}{\partial\theta^\celu}$ (in an equivalent way, $\theta^\gou(X(\xi^\gou)) = \xi^\gou$). For any $(\textsf{y},\xi^\gou)\in \mathcal{Y}\times \gou$, denote by, if it exists, $e^{X(\xi^\gou)}(\textsf{y})$ the value at time $t=1$ of the solution $\gamma \in \mathscr{C}^1([0,1],\mathcal{Y})$ of the equation $\frac{d\gamma}{dt} = X(\xi^\gou)(\gamma)$, with the initial condition $\gamma(0) = \textsf{y}$.
We consider the open subset $\Delta_\Phi\subset \mathcal{Y}\times \gou$ and the map $\Phi:\Delta_\Phi \longrightarrow \mathcal{Y}$ such that $\Phi(\textsf{y},\xi^\gou) = e^{X(\xi^\gou)}(\textsf{y})$ and $\Delta_\Phi$ ('life set') is the maximal open subset of $\mathcal{Y}\times \gou$ on which $\Phi$ can be defined.

For any value $r>0$ we let $B_\gou(r)$ be the ball of radius $r$ centered at 0 in $\gou$ (for any norm on $\gou$). For $r>0$ sufficiently small, we define the map $\Psi:B_\gou(r) \longrightarrow \mathcal{Y}$ as follows. For any $\xi^\gou \in B_\gou(r)$, we use the unique splitting $\xi^\gou = \xi^\gos + \xi^\gog$ according to the decomposition $\gou = \gos\oplus \gog$ and we set
\begin{equation}\label{definitionPsipourfibration}
  \Psi(\xi^\gou) = \Phi(\Phi(\textsf{y}_0,\xi^\gos), \xi^\gog)
\end{equation}
The differential of $\Psi$ at $0$ is the inverse map of $\theta^\gou_{\textsf{y}_0}$ and hence is invertible. Thus, thanks to the inverse mapping theorem, by choosing $r$ sufficiently small we can assume that $\Psi$ is a diffeomorphism between $B_\gou(r)$ and its image $\mathcal{O}$ in $\mathcal{Y}$, which is a neighbourhood of $\textsf{y}_0$. Let $\textsf{z}\in \mathcal{O}$ be an arbitrary point and let $\textsf{f}$ be the integral leaf which contains $\textsf{z}$. We will show that $\textsf{f}$ is  diffeomorphic to $\textsf{f}_0$ and hence to $\check{\textsf{f}}$.
For that purpose we will show that the flow map $e^{X(\xi^\gos)}$ is defined on $\textsf{f}_0$ and maps $\textsf{f}_0$ to $\textsf{f}$ in a diffeomorphic way.

We remark that, since any integral leaf is invariant by $\textsf{y}\longmapsto\Phi(\textsf{y},\xi^\gog)$, $\textsf{f}$ contains also $\textsf{y}_1:= \Phi(\textsf{y}_0,\xi^\gos)$ and hence is characterized by this property.

Let $\textsf{y}_0'\in \textsf{f}_0$. There exists a path $\gamma\in \mathscr{C}^1([0,1],\textsf{f}_0)$ based on $\gamma(0) = \textsf{y}_0$ and with end point $\gamma(1) = \textsf{y}_0'$. We build the map $U: [0,1]^2 \longrightarrow \mathcal{Y}$ by:
\[
\left\{  \begin{array}{cccc}
U(t,0) & = & \Phi(\textsf{y}_0,t\xi^\gos) & \forall t\in [0,1] \\
\left(U^*\theta^\gou\left(\frac{\partial}{\partial s}\right)\right)(t,s) & = & \left(\gamma^*\theta^\gou\left(\frac{d}{ds}\right)\right)(s) & \forall (t,s)\in [0,1]\times [0,1]
 \end{array}\right.
\]
A key point is that, since $\gamma$ takes value in the leaf $\textsf{f}_0$, $\gamma^*\theta^\gos = 0$, which implies $\left( U^*\theta^\gos\left(\frac{\partial}{\partial s}\right)\right)(t,s) = 0$. This has the first consequence that the existence of $U$ is guaranteed by the hypothesis (iii), i.e. that the manifold is $\gog$-complete.
From Equations (\ref{ELFrobeniusdetaillee}) we deduce
\[
 \left\{
 \begin{array}{lcc}
  \hbox{d}\left(U^*\theta^\gos\right) & = & \frac{1}{2} U^*\left(\Theta{^\gos}_{\cels\cels}\theta^{\cels}\wedge \theta^{\cels}\right) \\
  \hbox{d}\left(U^*\theta^\gog\right) + \frac{1}{2} U^*\left([\theta^\gog\wedge \theta^\gog]\right) & = & \frac{1}{2} U^*\left(\Theta{^\gog}_{\cels\cels}\theta^{\cels}\wedge \theta^{\cels}\right)
 \end{array}\right.
\]
This implies, since $\left(U^*\theta^\gos\left(\frac{\partial}{\partial s}\right)\right)(t,s) = 0$, that
\[
 \left\{
 \begin{array}{lcc}
 \displaystyle \hbox{d}\left(U^*\theta^\gos\right)\left(\frac{\partial }{\partial t}, \frac{\partial }{\partial s}\right) & = & 0 \\
 \displaystyle \hbox{d}\left(U^*\theta^\gog\right)\left(\frac{\partial }{\partial t}, \frac{\partial }{\partial s}\right)
 + \left[U^*\theta^\gog\left(\frac{\partial }{\partial t}\right), U^*\theta^\gog\left(\frac{\partial }{\partial s}\right)\right]
& = & 0
 \end{array}\right.
\]
On the other hand Cartan's formula
\[
\hbox{d}\left(U^*\theta^\gou\right)\left(\frac{\partial }{\partial t}, \frac{\partial }{\partial s}\right) + U^*\theta^\gou \left(\left[\frac{\partial }{\partial t}, \frac{\partial }{\partial s}\right]\right) =
\frac{\partial }{\partial t}\left(U^*\theta^\gou\left( \frac{\partial }{\partial s}\right) \right) - \frac{\partial }{\partial s}\left(U^*\theta^\gou\left( \frac{\partial }{\partial t}\right) \right)
\]
simplifies to $\hbox{d}\left(U^*\theta^\gou\right)\left(\frac{\partial }{\partial t}, \frac{\partial }{\partial s}\right)  + 0=  0- \frac{\partial }{\partial s}\left(U^*\theta^\gou\left( \frac{\partial }{\partial t}\right) \right)$. We hence deduce that, for all $t\in [0,1]$, $s\longmapsto U^*\theta^\gou\left( \frac{\partial }{\partial t}\right)(t,s)$ is solution of the system of differential equations
\[
 \left\{
 \begin{array}{lcc}
 \displaystyle \frac{\partial }{\partial s}\left(U^*\theta^\gos\left( \frac{\partial }{\partial t}\right) \right) & = & 0 \\
 \displaystyle \frac{\partial }{\partial s}\left(U^*\theta^\gog\left( \frac{\partial }{\partial t}\right) \right)  & = & \displaystyle \left[U^*\theta^\gog\left(\frac{\partial }{\partial t}\right), U^*\theta^\gog\left(\frac{\partial }{\partial s}\right)\right]
\end{array}\right.
\]
However we also have the following initial conditions at $s=0$:
\[
 U^*\theta^\gos\left( \frac{\partial }{\partial t}\right)(t,0) = \xi^\gos\quad \hbox{ and }\quad
 U^*\theta^\gog\left( \frac{\partial }{\partial t}\right)(t,0) = 0
\]
We thus conclude that $U^*\theta^\gos\left( \frac{\partial }{\partial t}\right)(t,s) = \xi^\gos$ and $U^*\theta^\gog\left( \frac{\partial }{\partial t}\right)(t,s) = 0$, $\forall (t,s)\in [0,1]^2$. This is equivalent to the relation $\frac{\partial U}{\partial s} = X(\xi^\gos)(U)$. This shows that the flow map of $X(\xi^\gos)$ is well defined at least for all time in $[0,1]$ on $\textsf{f}_0$ and maps $\textsf{f}_0$ to $\textsf{f}$. Since the reasoning can be reversed (by exchanging $\textsf{f}_0$ and $\textsf{f}$) this map is actually a diffeomorphism and, in particular, $\textsf{f}$ is compact. Thus $\textsf{z}\in \check{\mathcal{Y}}$.\\

\noindent
\textbf{(ii)} We show that $\check{\mathcal{Y}}$ is closed. Let $\textsf{y}$ be in the closure of $\check{\mathcal{Y}}$. In a way similar to the previous step, for $r>0$ sufficiently small, we define the map $\Psi:B_\gou(r) \longrightarrow \mathcal{Y}$ by $\Psi(\xi^\gou) = \Phi(\Phi(\textsf{y},\xi^\gos), \xi^\gog)$, where, $\forall \xi^\gou \in B_\gou(r)$, $\xi^\gou = \xi^\gos + \xi^\gog$.
For $r>0$ sufficiently small, we can assume that $\Psi$ is a diffeomorphism between $B_\gou(r)$ and its image $\mathcal{O}$ in $\mathcal{Y}$ and $\mathcal{O}$ is a neighbourhood of $\textsf{y}$.

Since $\textsf{y}$ belongs to the closure of $\check{\mathcal{Y}}$, there exists a sequence $\left(\textsf{y}_n\right)_{n\in \N}$ of points in $\check{\mathcal{Y}}$ which converges to $\textsf{y}$. We can fix a value of $n$ sufficiently large so that $\textsf{y}_n\in \mathcal{O}$. Since $\textsf{y}_n\in \check{\mathcal{Y}}$, the leaf $\textsf{f}_n$ which contains $\textsf{y}_n$ is diffeomorphic to $\check{\textsf{f}}$. We can then repeat the arguments of the previous step by replacing $\textsf{y}_0$ by $\textsf{y}_0':= \Phi(\textsf{y}_n,-\xi^\gog_n)$, where $\xi^\gog_n$ is such that $\Psi(\xi^\gos_n + \xi^\gog_n) = \textsf{y}_n$. (Note that $\Phi(\textsf{y}_0',-\xi^\gos_n) = \textsf{y}$.) We thus obtain that $\textsf{f}$ is diffeomorphic to $\check{\textsf{f}}$. \hfill $\square$

\subsubsection{Intermediate conclusion}
By using Lemmas \ref{lemmaGcoversleaves} and \ref{lemmaifGnoncompact} we  immediately obtain Conclusions (i) and (ii) in Theorem \ref{theorem0} holds, i.e. that all integral fibers are diffeomorphic to a Lie group $\goG:= \widehat{\goG}/\pi_1(\textsf{f})$, where $\pi_1(\textsf{f})$ is the fundamental group of any integral leaf $\textsf{f}$.

\subsubsection{Construction of a principal fiber bundle structure}
In the following we exploit Lemmas \ref{lemmaGcoversleaves} and \ref{lemmaifGnoncompact} by assuming furthermore that \textbf{$\goG:= \widehat{\goG}/\pi_1(\textsf{f})$ is compact}. Then all integral leaves $\textsf{f}$ are compact and we will prove that these leaves are actually the fibers of a principal bundle with structure group $\goG$.

As in the proof of Lemma \ref{lemmaifGnoncompact}, to any $\xi^\gou\in \gou$ we associate the vector field $\textsf{X}(\xi^\gou)$ on $\mathcal{Y}$ such that $\theta^\gou(X(\xi^\gou)) = \xi^\gou$. A useful property is
\begin{equation}\label{gosetgogcommutent}
 \forall (\xi^\gos,\xi^\gog)\in \gos\times \gog, \quad [X(\xi^\gos),X(\xi^\gog)] = 0
\end{equation}
The proof of (\ref{gosetgogcommutent}) follows again from Cartan's formula $\hbox{d}\theta^\gou(X,Y) + \theta^\gou([X,Y]) = X\cdot \theta^\gou(Y) - Y\cdot \theta^\gou(X)$, with $X=X(\xi^\gos)$ and $Y=X(\xi^\gog)$, which gives $\theta^\gou([X,Y]) = - \hbox{d}\theta^\gou(X,Y)$. This implies by using (\ref{ELFrobeniusdetaillee}) that $\theta^\gou([X,Y]) = 0$ and hence $[X,Y]=0$.

For any integral leaf $\textsf{f}$ and any point $\textsf{y}_0\in \textsf{f}$ we define the map
\[
 \begin{array}{ccc}
  \goG & \longrightarrow & \textsf{f} \\
  g & \longmapsto & g\cdot \textsf{y}_0
 \end{array}
\]
as follows. Let $\hat{g}\in \widehat{\goG}$ be any point which is mapped to $g$ through the projection mapping $\widehat{\goG} \longrightarrow \goG = \widehat{\goG}/\pi_1(\textsf{f})$. We then set $g\cdot \textsf{y}_0 = T(\hat{g})$, where $T$ is the map constructed in the proof of Lemma \ref{lemmaGcoversleaves}. It follows from the definition of the action of $\pi_1(\textsf{f})$ on $\widehat{\goG}$ that this value does not depend on the choice of $\hat{g}$.

For any $r\in (0,+\infty)$ let $B_\gos(r)$ be the open ball of radius $r$ and of center 0 in $\gos$. We fix an arbitrary point $\textsf{y}_0\in \mathcal{Y}$ and we define the map
\[
 \begin{array}{cccc}
  A_r: & B_\gos(r)\times \goG & \longrightarrow & \mathcal{Y} \\
  & (\xi^\gos,g) & \longmapsto &
  A_r(\xi^\gos,g) = g\cdot \left(e^{X(\xi^\gog)}(\textsf{y}_0)\right)
 \end{array}
\]
Note that, for $g=\hbox{exp}\,\xi^\gog$, we have $A_r(\xi^\gos,\hbox{exp}\,\xi^\gog) = \Psi(\xi^\gos+\xi^\gog)$ (where $\Psi$ is defined by (\ref{definitionPsipourfibration})). For $r$ sufficiently small, it is clear that $A_r$ is well-defined and is a \emph{local} diffeomorphism. However it is not clear a priori whether $A_r$ is a \emph{global} diffeomorphism between $B_\gos(r)\times \goG$ and its image since $A_r$ may not be one-to-one in general. Indeed although, for any $\xi^\gos\in B_\gos(r)$, the restriction of $A_r$ to $\{\xi^\gos\}\times \goG$ is a diffeomorphism whose image is an integral leaf, it may happen that there exists two different values $\xi^\gos,\zeta^\gos\in B_\gos(r)$ such that $A_r(\{\xi^\gos\}\times \goG) = A_r(\{\zeta^\gos\}\times \goG)$.

For $h\in (0,+\infty)$ let $B_\gog(h)$ be the open ball of center 0 and of radius $h$ in $\gog$ and let $\Psi_{r,h}$ be the restriction of $\Psi$ (defined by (\ref{definitionPsipourfibration})) to $B_\gos(r)\times B_\gog(h)$. Since $\hbox{d}\Psi_{r,h}$ is invertible, we may choose $(r,h)$ in such a way that $\Psi_{r,h}$ is a diffeomorphism onto its image $\mathcal{O}_{r,h}:= \Psi_{r,h}(B_\gos(r)\times B_\gog(h))$.

Let $\check{\textsf{f}}$ be the integral leaf which contains $\textsf{y}_0$. Since $\check{\textsf{f}}$ is compact the intersection $\check{\textsf{f}}\ \cap \mathcal{O}_{r,h}$ is composed of a \emph{finite} number $N+1$ of connected components. We denote by $\check{\textsf{f}}_0,\check{\textsf{f}}_1,\cdots, \check{\textsf{f}}_N$ these connected components, where $\check{\textsf{f}}_0$ is the image of $\{0\}\times \goG$ by $\Psi_{r,g}$.

For any pair $\textsf{f}',\textsf{f}''$ of submanifolds of $\mathcal{O}_{r,h}$ which are open subsets of integrals leaves, define
\[
 d(\textsf{f}',\textsf{f}''):=
 \hbox{inf}\{ \|\zeta^\gos\|\ ;\ \zeta^\gos\in \gos, e^{X(\zeta^\gos)}(\textsf{f}') \cap \textsf{f}''\neq \emptyset \}
\]
It is clear that $\exists \delta\in (0,+\infty)$ such that $d(\check{\textsf{f}}_0,\check{\textsf{f}}_j) >2\delta$, $\forall j = 1,\cdots,N$. (This means in particular that the inverse image of $\check{\textsf{f}}$ by $\Psi_{2\delta,h}$ is reduced to $\{0\}\times B_\gog(h)$.)

Now we observe that, by the proof of Lemma \ref{lemmaifGnoncompact}, for all $\xi^\gos$ in a neighbourhood of 0 in $\gos$ and for any $j\in \{0,\cdots,N\}$, $e^{X(\xi^\gos)}(\check{\textsf{f}}_j)$ is well defined and depends in a \emph{continuous} way on $\xi^\gos$.
Thus in particular, $\exists \rho\in (0,\delta)$ such that $\forall \xi^\gos\in B_\gos(\rho)$, $\forall j \in \{1,\cdots,N\}$, $d\left(\check{\textsf{f}}_0,e^{X(\xi^\gos)}(\check{\textsf{f}}_j)\right) > \delta$. Hence, if $\xi^\gos\in B_\gos(\rho)$, on the one hand, $e^{X(\xi^\gos)}(\check{\textsf{f}}_0) = \Psi_{\rho,h}(\{\xi^\gos\}\times B_\gog(h)) \subset \Psi_{\rho,h}(B_\gos(\rho)\times B_\gog(h)) =:\mathcal{O}_{\rho,h}$ and, on the other hand, all the other connected components $e^{X(\xi^\gos)}(\check{\textsf{f}}_j)$ (for $1\leq j\leq N$) are outside $\mathcal{O}_{\delta,h}$. Since $\rho<\delta$, this ensures that the inverse image by $\Psi_{\rho,h}$ of the intersection of any integral leaf with $\mathcal{O}_{\rho,h}$ is reduced to $\{\xi^\gos\}\times B_\gog(h)$.

As a consequence the map $A_\rho$ is a diffeomorphism between $B_\gos(\rho)\times \goG$ and its image. This shows that $\mathcal{Y}$ has a principal bundle structure, with structure group $\goG$, the map $A_\rho$ providing us with a local trivialization. Hence
the set $\mathcal{X}$ of integral leaves has the structure of an $n$-dimensional manifold. We denote by $P:\mathcal{Y}\longrightarrow \mathcal{X}$
the quotient map.

Set $e^\gos:= \theta^\gos$.
From $\frac{\partial}{\partial \theta^\gog}\iN e^\gos = \frac{\partial}{\partial \theta^\gog}\iN \hbox{d}e^\gos = 0$
we deduce that there exists a coframe $\underline{e}^\gos$ on $\mathcal{X}$
such that $e^\gos = P^*\underline{e}^\gos$. Thus we can equipp $\mathcal{X}$ with the
pseudo Riemannian metric $\underline{\textbf{g}}:= \textsf{b}_{ab}\underline{e}^a\otimes\underline{e}^b$.

\subsubsection{Working in a local trivialization of the bundle}

In the following we choose an $n$-dimensional submanifold $\Sigma \subset \mathcal{Y}$
transverse to the fibration. Without loss of generality
(replacing $\mathcal{Y}$ by an open subset of $\mathcal{Y}$ if necessary)
we can assume that $\Sigma$ intersects all fibers of $P$ (i.e. defines a section of $P:\mathcal{Y}\longrightarrow \mathcal{X}$) and
we define the map $g:\mathcal{Y}\longrightarrow \goG$
which is constant equal to $1_\goG$ on $\Sigma$ and
such that
\[
\hbox{d}g-g\theta^\gog|_\textsf{f}=0 
\]
for any
integral leaf $\textsf{f}$.
We then define
\[
 \textbf{A}^\gou:= \hbox{Ad}_g\theta^\gou - \hbox{d}g\cdot g^{-1}
\]
which means that $\textbf{A}^\gos = \theta^\gos$ and $\textbf{A}^\gog:= \hbox{Ad}_g\theta^\gog - \hbox{d}g\cdot g^{-1}$. Obviously $\textbf{A}^\gos|_\textsf{f} = 0$ and moreover
the relation $\hbox{d}g-g\theta^\gog|_\textsf{f}=0$ translates as
$\textbf{A}^\gog|_\textsf{f} = 0$. Thus $\textbf{A}^\gou|_\textsf{f} = 0$ so that we have the decomposition
$\textbf{A}^\gou = \textbf{A}{^\gou}_\cels\theta^\cels$ (with $\textbf{A}{^\gos}_\gos = \delta{^\gos}_\gos$).
Moreover since
\begin{equation}\label{thetaenfonctiondeA}
\theta^\gou = g^{-1}\textbf{A}^\gou g + g^{-1}\hbox{d}g,
\end{equation}
we have
$\hbox{d}\theta^\gou +\frac{1}{2}[\theta^\gou\wedge\theta^\gou] = g^{-1}(\hbox{d}\textbf{A}^\gou +\frac{1}{2}[\textbf{A}^\gou\wedge \textbf{A}^\gou])g
= g^{-1}\textbf{F}^\gou g$, where $\textbf{F}^\gou:= \hbox{d}\textbf{A}^\gou + \frac{1}{2}[\textbf{A}^\gou\wedge \textbf{A}^\gou]$.
From (\ref{ELFrobeniusdetaillee}) we deduce $\frac{\partial}{\partial \theta^\gog}\iN
\Theta^\gou = 0$ which is equivalent to $\frac{\partial}{\partial \theta^\gou}\iN
\textbf{F}^\gog = 0$. But since $\frac{\partial}{\partial \theta^\gog}\iN \textbf{A}^\gou = 0$, this implies furthermore that $\frac{\partial}{\partial \theta^\gog}\iN \dR\textbf{A}^\gou = 0$ and thus the coefficients $\textbf{A}{^\gou}_\gos$ are constants on the fibers
$\textsf{f}$. 
Hence
\begin{equation}\label{decompositiondeF}
 \textbf{F}^\gou = \frac{1}{2}\textbf{F}{^\gou}_{\cels\cels}\theta^{\cels\cels}
\end{equation}
where the coefficients $\textbf{F}{^\gou}_{\gos\gos} = \hbox{Ad}_g\otimes 1_{\gos^*}\otimes 1_{\gos^*} \Theta{^\gou}_{\gos\gos}$ are constant on the fibers.

\subsection{The Euler--Lagrange system in a local trivialization}\label{transformationdejauge}
We proceed similarly as for Yang--Mills in \S \ref{backtoYM}.
Consider the map $\hbox{Ad}_g:\mathcal{Y}\longrightarrow \hbox{End}(\gou)$, where $g:\mathcal{Y}\longrightarrow \goG$ is the map defined previously. Actually $\hbox{Ad}_g$ takes values in $SO(\gou,\textsf{h})$ since $\textsf{h}$ is invariant by $\hbox{Ad}_\goG$.
We define the coframe $e^\gou:= \hbox{Ad}_g \theta^\gou = \textbf{A}^\gou + \dR g\ g^{-1}$. Note that
\begin{equation}\label{formulepoure}
 \begin{array}{ccccl}
  e^\gos & = & \theta^\gos & = & \textbf{A}^\gos \\
e^\gog  & = & \hbox{Ad}_g \theta^\gog
& = & \textbf{A}^\gog + \dR g\ g^{-1}
 \end{array}
\end{equation}
and (\ref{decompositiondeF}) becomes
\[
\textbf{F}^\gou = \frac{1}{2}\textbf{F}{^\gou}_{\cels\cels}e^{\cels\cels},
\]
By using (\ref{dAeagalFee}) we get
\begin{equation}\label{deextra}
 \dR^{\textbf{A}} e^\gog = \textbf{F}^\gog
  + \frac{1}{2} [e^\gog\wedge e^\gog]
\end{equation}
We also define $p_\gou:= \hbox{Ad}^*_g \pi_\gou$ and
\[
\omega^\gol:= (\hbox{Ad}_g)\varphi^\gol(\hbox{Ad}_g)^{-1} - \hbox{d}(\hbox{Ad}_g)(\hbox{Ad}_g)^{-1}
\in \gol\otimes \Omega^1(\mathcal{Y})
\]
We note then that
\[
 \Omega^\gol:= \dR \omega^\gol + \frac{1}{2}[\omega^\gol\wedge \omega^\gol]
 = (\hbox{Ad}_g)\Phi^\gol(\hbox{Ad}_g)^{-1}
 \hbox{ and set }\Omega^{\gou\gou}:=
 \kappa{_\cell}^{\gou\gou}\Omega^\cell
\]
We translate Equation (\ref{ELEinsteinbrut}) by computing the images of its both sides by $\hbox{Ad}^*_g$ in terms of these new variables. From Lemma \ref{mainlemma} we deduce
$\hbox{Ad}_g^*(\dR^{\theta^\gou} \pi_\gou) =
\dR^{\hbox{\scriptsize{Ad}}_g\theta^\gou - dg\, g^{-1}} \hbox{Ad}^*_\gog \pi_\gou
= \dR^{\theta^\gos + \textbf{A}^\gog} p_\gou$.
However since $\gos$ belongs to the center of $(\gou,[\cdot,\cdot])$, this relation reduces to $\hbox{Ad}_g^*(\dR^{\theta^\gou} \pi_\gou)
= \dR^{\textbf{A}^\gog} p_\gou$. Hence by using the fact that $\kappa{_\gol}^{\gou\gou}$ is invariant by $\hbox{Ad}_\goG$, we get
\begin{equation}\label{ELEinsteinbrut2}
 \dR^\textbf{A}p_\gou +
 \frac{1}{2}e^{(N-3)}_{\celu\celu\gou}\wedge \Omega^{\celu\celu}
- \Lambda_0 e^{(N-1)}_\gou
= \mathbf{F}{^\celu}_{\gou\cels}p{_\celu}^{\celg\cels}\ e^{(N-1)}_{\celg}
\end{equation}
We note that the second term on the l.h.s. is nothing but (minus) the Einstein tensor
$\mathbf{E}(\textbf{h}){_\gou}^\gou$ (see (\ref{notationEinstein}))
on $(\mathcal{Y},\textbf{h})$:
\begin{equation}\label{kappaeOmegaEin}
 \frac{1}{2} e_{\celu\celu\gou}^{(N-3)}\wedge \Omega^{\celu\celu}
 = - \mathbf{E}(\textbf{h}){_\gou}^\celu e^{(N-1)}_\celu
\end{equation}
Thus we obtain
\begin{equation}\label{ELEinsteinbrut3}
\mathbf{E}(\textbf{h}){_\gou}^\celu e^{(N-1)}_\celu
+ \Lambda_0 e^{(N-1)}_\gou
=  \dR^\textbf{A}p_\gou  -
 \mathbf{F}{^\celu}_{\gou\cels}p{_\celu}^{\celg\cels}\ e^{(N-1)}_{\celg}
\end{equation}
The computation of $\dR^\textbf{A} p_\gou$ follows the same steps as for the Yang--Mills case (see (\ref{dApisynth})), by using $\dR^\textbf{A}$ given by (\ref{deextra}) instead of $\dR^{\gamma,\textbf{A}}$ given by (\ref{nablaAxiU}) and with the simplification that $p{_\gou}^{\gos\gos} = 0$:
\begin{equation}\label{dApKK}
 \dR^\textbf{A}p_\gou
= \partial_\celg p{_\gou}^{\cels\celg} e_{\cels}^{(N-1)} + \left( \partial^\textbf{A}_{\cels} p{_\gou}^{\celg\cels} +
  \partial_{\celg_1}p{_\gou}^{\celg\celg_1}
  + \frac{1}{2}\textbf{c}{^\celg}_{\celg_1\celg_2}p{_\gou}^{\celg_1\celg_2}
  \right)
  e_{\celg}^{(N-1)}
\end{equation}
Hence, by using $p{_\gou}^{\gos\gos} = 0$, we can write (\ref{ELEinsteinbrut3}) as the system
\begin{equation}\label{ELEinsteinbrut4}
\left\{
 \begin{array}{ccl}
  \mathbf{E}(\textbf{h}){_\gou}^\gos
+ \Lambda_0 \delta{_\gou}^\gos & = & \partial_\celg p{_\gou}^{\gos\celg}  \\
\mathbf{E}(\textbf{h}){_\gou}^\gog
+ \Lambda_0 \delta{_\gou}^\gog & =  &
 \partial^\textbf{A}_{\cels} p{_\gou}^{\gog\cels} +
  \partial_{\celg}p{_\gou}^{\gog\celg}
  + \frac{1}{2} p{_\gou}^{\celg\celg} \textbf{c}{^\gog}_{\celg\celg}
  - \mathbf{F}{^\celu}_{\gou\cels}p{_\celu}^{\gog\cels}
  \end{array}\right.
\end{equation}
Equivalentely by using the splitting $\gou^* = \gos^* + \gog^*$ and with the simplification $\mathbf{F}{^\gou}_{\gog\gos} = 0$,
\begin{equation}\label{ELEBrutSysteme}
\begin{array}{l}
 \left(
 \begin{array}{cc}
  \mathbf{E}(\textbf{h}){_\gos}^\gos
+ \Lambda_0 \delta{_\gos}^\gos
& \mathbf{E}(\textbf{h}){_\gog}^\gos \\
\mathbf{E}(\textbf{h}){_\gos}^\gog
& \mathbf{E}(\textbf{h}){_\gog}^\gog
+ \Lambda_0 \delta{_\gog}^\gog
 \end{array}\right)
\\  =
 \left(
 \begin{array}{ccc}
  \partial_\celg p{_\gos}^{\gos\celg}
  & & \partial_\celg p{_\gog}^{\gos\celg} \\
  \partial^\textbf{A}_{\cels} p{_\gos}^{\gog\cels} +
  \partial_{\celg}p{_\gos}^{\gog\celg}
  + \frac{1}{2} p{_\gos}^{\celg\celg} \textbf{c}{^\gog}_{\celg\celg}
  - \mathbf{F}{^\celu}_{\gos\cels}p{_\celu}^{\gog\cels}
& & \partial^\textbf{A}_{\cels} p{_\gog}^{\gog\cels} +
  \partial_{\celg}p{_\gog}^{\gog\celg}
  + \frac{1}{2} p{_\gog}^{\celg\celg} \textbf{c}{^\gog}_{\celg\celg}
 \end{array}\right)
\end{array}
 \end{equation}
Observe here that, because of the symmetry of the Einstein tensor and since $\gos \perp \gog$, we have
$\textsf{h}^{\gos\cels}\mathbf{E}(\textbf{h}){_\cels}^\gog = \textsf{h}^{\gog\celg}\mathbf{E}(\textbf{h}){_\celg}^\gos$

Again a crucial point is to observe that
the l.h.s. $\mathbf{E}(\mathbf{h}){_\gou}^\gos
+ \Lambda_0 \delta{_\gou}^\gos$ of the first equation in (\ref{ELEinsteinbrut4}) is constant on any fiber of the fibration $P:\mathcal{Y}\longrightarrow \mathcal{X}$. By setting
$(e^\gog)^{(r)}:= e^{n+1}\wedge \cdots \wedge e^N$ and $(e^\gog)^{(r-1)}_\gog:= \frac{\partial}{\partial e^i}\iN e^{(r)}\textbf{t}^i$ and by using the fact that the fibers are compact we deduce from (\ref{ELEinsteinbrut4}) that the cancellation phenomenon holds:
\begin{equation}\label{cancellationKK}
  \mathbf{E}(\textbf{h}){_\gou}^\gos
+ \Lambda_0 \delta{_\gou}^\gos = \frac{\int_{\mathcal{Y}_{\textsf{x}}}\left(\mathbf{E}(\textbf{h}){_\gou}^\gos
+ \Lambda_0 \delta{_\gou}^\gos\right)(e^\gog)^{(r)}}{\int_{\mathcal{Y}_{\textsf{x}}}(e^\gog)^{(r)}}= \frac{\int_{\mathcal{Y}_{\textsf{x}}}\dR\left(p{_\gou}^{\gos\celg}(e^\gog)^{(r-1)}_\celg\right)}{\int_{\mathcal{Y}_{\textsf{x}}}(e^\gog)^{(r)}} = 0
\end{equation}
Hence by taking into account the symmetry of the Einstein tensor we deduce that (\ref{ELEBrutSysteme}) reduces to
\begin{equation}\label{ELESysteme}
 \left(
 \begin{array}{cc}
  \mathbf{E}(\textbf{h}){_\gos}^\gos
+ \Lambda_0 \delta{_\gos}^\gos
& \mathbf{E}(\textbf{h}){_\gog}^\gos \\
\mathbf{E}(\textbf{h}){_\gos}^\gog
& \mathbf{E}(\textbf{h}){_\gog}^\gog
+ \Lambda_0 \delta{_\gog}^\gog
 \end{array}\right)
 =
 \left(
 \begin{array}{cc}
  0 & 0 \\
  0 & \partial^\textbf{A}_{\cels} p{_\gog}^{\gog\cels} +
  \partial_{\celg}p{_\gog}^{\gog\celg}
  + \frac{1}{2} p{_\gog}^{\celg\celg} \textbf{c}{^\gog}_{\celg\celg}
 \end{array}\right)
 \end{equation}
Beware that it does mean that $(\mathcal{Y},\textbf{h})$ is a solution of the Einstein equation with a cosmological constant since $\mathbf{E}(\textbf{h}){^\gog}_\gog
+ \Lambda_0 \delta{^\gog}_\gog$ does not vanish in general.

\subsection{The Einstein--Yang--Mills system}\label{calculfinal}
Lastly we translate equations $\mathbf{E}(\textbf{h}){_\gos}^\gos
+ \Lambda_0 \delta{_\gos}^\gos = 0$ and $\mathbf{E}(\textbf{h}){_\gos}^\gog = 0$ as equations on fields defined on $\mathcal{X}$. We introduce a basis $(\textbf{u}_1,\cdots,\textbf{u}_N)$ of $\gou$ such that $(\textbf{u}_1,\cdots,\textbf{u}_b)$ is a basis of $\gos$ and $(\textbf{u}_{n+1},\cdots,\textbf{u}_N)$ is a basis of $\gog$.

From $e^\gou$ we build the metric $\textbf{g}:= (\theta^\gos)^*\textsf{b} = (e^\gos)^*\textsf{b}$ on $\mathcal{X}$ and the associated Levi-Civita connection $\nabla^{T\mathcal{X}}$. The connection form $\gamma^{so(\gos)}\in so(\gos,\textsf{b})\otimes \Omega^1(\mathcal{F})$ of $\nabla^{T\mathcal{X}}$
can be computed by comparing (\ref{ELFrobeniusdetaillee}), which gives us $\dR e^\gos = \frac{1}{2}\Theta{^\gos}_{\cels\cels}e^{\cels\cels}$, with the zero torsion condition
$\dR e^\gos + \gamma^{so(\gos)}\wedge e^\gos = 0$: by using the notations $\gamma{^a}_b$ for the matrix coefficients of $\gamma^{so(\gos)}$ in the basis $(\textbf{u}_1,\cdots,\textbf{u}_n)$ and $\gamma{^a}_{bc}$ for its coefficients (see (\ref{exportKK}), we have $\gamma{^a}_{bc}=
 \frac{1}{2}\left(\Theta{^a}_{bc}
 - \textsf{b}^{aa'}\textsf{b}_{bb'}\Theta{^{b'}}_{a'c}
 - \textsf{b}^{aa'}\textsf{b}_{cc'}\Theta{^{c'}}_{a'b}\right)$.

Let also $(\omega{^A}_B)_{1\leq A,B\leq N}$ be the matrix coefficients of the Levi-Civita connection 1-form $\omega^\gol$.
In \cite{helein2020} $\omega^\gol$ is computed in function of $\gamma^{so(\gos)}$ and of $\textbf{A}^\gog$ and $\textbf{F}^\gog$. The result is the following: let
$\omega{^\gou}_\gou:= \omega{^A}_B\textbf{u}_A\otimes\textbf{u}^B$ and
$\gamma{^\gos}_\gos = \gamma{^a}_b\textbf{u}_a\otimes\textbf{u}^b = \gamma{^a}_{bc} \textbf{u}_a\otimes \textbf{u}^b\otimes e^c$.
By setting $\textbf{F}{_\gog}{^\gos}_\gos:=  (\textsf{k}_{\gog\celg}\otimes \textsf{b}^{\gos\cels}\otimes 1) \textbf{F}{^{\celg}}_{\cels\gos}$
and $\textbf{F}{_{\gog\gos}}^\gos:= (\textsf{k}_{\gog\celg}\otimes 1 \otimes \textsf{b}^{\gos\cels}) \textbf{F}{^{\celg}}_{\gos\cels}$
\[
 \left(\begin{array}{cc}
\omega{^\gos}_\gos & \omega{^\gos}_\gog \\
\omega{^\gog}_\gos & \omega{^\gog}_\gog
       \end{array}\right)
 = \left(\begin{array}{cc}
        \gamma{^\gos}_\gos
- \frac{1}{2}\textbf{F}{_{\celg}}{^\gos}_{\gos}e^\celg
& \frac{1}{2}\textbf{F}{_{\gog\cels}}^{\gos}e^\cels \\
 \frac{1}{2}\textbf{F}{^\gog}_{\gos\cels}e^\cels & \frac{1}{2} \textbf{c}{^\gog}_{\gog\celg}(e^\celg - 2\textbf{A}^\celg)
       \end{array}\right)
\]
We deduce the curvature 2-form $\Omega{^\gou}_\gou = \hbox{d}\omega{^\gou}_\gou +  \omega{^\gou}_\celu\wedge \omega{^\celu}_\gou$ and the components of the Ricci tensor $\mathbf{R}(\textbf{h}){_\gos}^\gos$ and of the Einstein tensor $\mathbf{E}(\textbf{h}){_\gos}^\gos$.
By setting $|\textbf{F}|^2:= \frac{1}{2}
\textbf{F}{_\celg}^{\cels_1\cels_2}\textbf{F}{^\celg}_{\cels_1\cels_2}$ and $\langle \textsf{B},\textsf{k}\rangle:=
 \frac{1}{2}\textbf{c}{^{\celg_1}}_{\celg_2\celg_3} \textbf{c}{^{\celg_2}}_{\celg_1\celg_4}\textsf{k}^{\celg_3\celg_4}$ (here $\textsf{B}_{\gog\gog}:= \textbf{c}{^{\celg_1}}_{\celg_2\gog} \textbf{c}{^{\celg_2}}_{\celg_1\gog}$ is the Killing form on $\gog$), the scalar curvature reads $\mathbf{R}(\textbf{h}) = \mathbf{R}(\gamma)
 - \frac{1}{2} |\textbf{F}|^2
 - \frac{1}{2} \langle \textsf{B},\textsf{k}\rangle$ and
\begin{equation}
 \mathbf{E}(\textbf{h}){_\gos}^\gos = \mathbf{E}(\textbf{g}){_\gos}^\gos - \frac{1}{2}\left(\textbf{F}{^\celg}_{\gos\cels}\textbf{F}{_\celg}^{\gos\cels}
 - \frac{1}{2} |\textbf{F}|^2\delta{_\gos}^\gos\right)
+ \frac{1}{4}\langle \textsf{B},\textsf{k}\rangle \delta{_\gos}^\gos
\end{equation}
\begin{equation}\label{Ricciadelta}
\mathbf{E}(\textbf{h}){_\gog}^\gos = \frac{1}{2}\left(\partial_\cels \textbf{F}_{\gog}{^{\gos\cels}}
 +  \textbf{F}{_\gog}^{\cels_1\cels} \gamma{^\gos}_{\cels_1\cels} + \gamma{^\cels}_{\cels_1\cels} \textbf{F}{_\gog}^{\gos\cels_1}
 - \textbf{c}{^{\celg_1}}_{\celg_2\gog}\textbf{A}{^{\celg_2}}_\cels \textbf{F}{_{\celg_1}}^{\gos\cels}\right)
\end{equation}
\begin{equation}\label{Riccialphad}
 \mathbf{E}(\textbf{h}){_\gog}^\gog = \frac{1}{4} \textbf{F}{_\gog}{^{\cels\cels}}\textbf{F}{^\gog}{_{\cels\cels}}
 - \frac{1}{4} \textbf{c}{^{\celg_1}}_{\gog\celg_3} \textbf{c}{^{\gog}}_{\celg_1\celg_2} \textsf{k}^{\celg_2\celg_3}
 - \frac{1}{2}\mathbf{R}(\textbf{h})\delta{_\gog}^\gog
\end{equation}
We note that (\ref{Ricciadelta}) can be written $\mathbf{E}(\textbf{h}){_\gog}^\gos = \frac{1}{2}\nabla^{T\mathcal{X},\textbf{A}}_\cels \textbf{F}_{\gog}{^{\gos\cels}}$, where $\nabla^{T\mathcal{X},\textbf{A}} = \nabla^{T\mathcal{X}} + \hbox{ad}^*_\textbf{A}\wedge$.
In conclusion, by setting $\Lambda := \frac{1}{4}\langle \textsf{B},\textsf{k}\rangle + \Lambda_0$, we get a solution of the system
\begin{equation}\label{EYMfinal}
 \left\{
 \begin{array}{ccl}
  \mathbf{E}(\textbf{g}){_\gos}^\gos
  + \Lambda\delta{_\gos}^\gos & = &
\frac{1}{2}\textbf{F}{^\celg}_{\gos\cels} \textbf{F}{_\celg}^{\gos\cels}
- \frac{1}{4} |\textbf{F}|^2\delta{_\gos}^ \gos\\
\nabla^{T\mathcal{X},\textbf{A}}_\cels \textbf{F}_{\gog}{^{\gos\cels}} & = & 0
 \end{array}\right.
\end{equation}
i.e. the Einstein--Yang--Mills system on $(\mathcal{X},\textbf{g})$ with the connection $\textbf{A}^\gog$ on $\mathcal{Y}\longrightarrow \mathcal{X}$ and the cosmological constant $ \frac{1}{4}\langle \textsf{B},\textsf{k}\rangle + \Lambda_0$.

\section{Gravity theory}\label{sectionGravity}

We now turn to generalized gravity theories the formulation of which takes place on manifolds which look \emph{locally} as principal bundles. For general solutions the corresponding space-time will be built as a set of leaves of a foliation (hence non separable in general). In special circumstances this quotient space is a true manifold and we recover usual gravity theories on this manifold.

We let $\widehat{\goL}$ and $\widehat{\goP}$ be two simply connected unimodular Lie groups and we assume that $\widehat{\goL}$ is a subgroup of $\widehat{\goP}$. As a motivation we may think that $\widehat{\goL}$ is the connected component of the identity of the Spin group $Spin_0(1,3)$ and that $\widehat{\goP}$ is the corresponding Spin Poincaré group $Spin_0(1,3)\ltimes \R^4$.
We let $\gol$ and $\gop$ be, respectively, the Lie algebras of $\widehat{\goL}$ and $\widehat{\goP}$.

The unknown fields will be a $\gop$-valued 1-form $\varphi^\gop$ which is a coframe on an oriented manifold $\mathcal{F}$
(where $\hbox{dim}\mathcal{F} = \hbox{dim}\goP =:N$) and a dual field $\pi_\gop$,
which is an $(N-2)$-form with
coefficients in $\gop^*$.
Then by looking at the Euler--Lagrange equations of the action functional $\int_\mathcal{F}\pi_{\celp}\wedge
(\dR\varphi^{\celp} + \frac{1}{2}[\varphi^\gop\wedge \varphi^\gop]^{\celp})$
on a class of fields satisfying a particular constraint we find dynamical
equations which implies the existence of a foliation of $\mathcal{F}$ which, under some \emph{extra topological
hypotheses}, gives rise to a principal bundle structure
on $\mathcal{F}$ with a structure group $\goL$, which is a quotient of $\widehat{\goL}$ by a finite subgroup. The space of leaves $\mathcal{X}$ has the same dimension as $\widehat{\goP}/\widehat{\goL}$ and can be interpreted as
the space-time $\mathcal{X}$. The dynamical equations then imply that on can extract some
fields defined on $\mathcal{X}$ out of $\varphi^\gop$, which satisfy an
Einstein--Cartan system of equations.

\subsection{General setting for gravity}

\subsubsection{Hypotheses on
the structure groups}
We denote by $\gol$ and $\gop$ the Lie algebras of, respectively, $\widehat{\goL}$ and $\widehat{\goP}$.
Our hypotheses are:
\begin{enumerate}
 \item
$\gop$ is \emph{reductive}, i.e.
there exists some vector subspace $\gos\subset \gop$
such that
\begin{equation}\label{hyposplustegalp}
 \gol\oplus \gos = \gop.
\end{equation}
and $\gos$ is stable under the ajdoint action of $\widehat{\goL}$, i.e.
\begin{equation}\label{hypotstablesouss}
 \hbox{Ad}_{\widehat{\goL}}\gos\subset \gos,
 \quad \hbox{ i.e.: }\quad \forall g\in \widehat{\goL},
 \forall \xi\in \gos, \hbox{Ad}_g\xi
 = g\xi g^{-1}\in \gos.
\end{equation}
\item $\widehat{\goP}/\widehat{\goL}$ is a \emph{symmetric space}, which amounts to assume that
\begin{equation}\label{hypottincluss}
 [\gos,\gos]\subset \gol,
 \quad \hbox{ meaning that }\quad \forall \xi,\zeta\in \gos,
 [\xi,\zeta]\in \gol
\end{equation}
\item the Lie algebras $\gop$ and $\gol$ are unimodular.
\end{enumerate}
Note that the fact that $\widehat{\goL}$ is a subgroup of $\widehat{\goP}$,
(\ref{hyposplustegalp}) and (\ref{hypottincluss}) imply respectively that:
\begin{equation}\label{stabilitiesofLiesubalgebras}
 [\gol,\gol]\subset \gol,\quad
 [\gol,\gos]\subset \gos \hbox{ and }
 [\gos,\gos]\subset \gol.
\end{equation}
The latter property is equivalent to the fact that the linear map $\tau: \gop \longrightarrow \gop$ such that $\gol$ and $\gos$ are the eigenspaces of $\tau$ for the eigenvalues 1 and $-1$, respectively, is a Lie algebra automorphism.

We define
$\gos^\perp
:= \{\alpha\in \gop^*;\langle\alpha,\xi\rangle = 0,
\forall \xi\in \gos\}$ and
similarly
$\gol^\perp
:= \{\alpha\in \gop^*;\langle\alpha,\xi\rangle = 0,
\forall \xi\in \gol\}$ and
we will systematically use the identifications
\[
 \gol^*:= \gos^\perp
 \hbox{ and }\gos^*:= \gol^\perp.
\]
We hence have $\gop^* = \gol^* \oplus \gos^*$.

Note that, if $\alpha\in \gol^*=\gos^\perp$,
then $\forall (\xi,\zeta)\in
(\gol\times \gos)\cup (\gos\times \gol)$,
$[\xi,\zeta]\in \gos$ because of
(\ref{stabilitiesofLiesubalgebras}),
and hence $\langle \hbox{ad}_\xi^*\alpha,\zeta\rangle =
\langle \alpha,[\xi,\zeta]\rangle = 0$. Hence
$(\alpha,\xi)\in \gol^*\times \gol$
implies $\hbox{ad}_\xi^*\alpha\in \gos^\perp = \gol^*$ and
$(\alpha,\xi)\in \gol^*\times \gos$
implies $\hbox{ad}_\xi^*\alpha\in \gol^\perp = \gos^*$.
A similar reasonning shows that
$(\alpha,\xi)\in \gos^*\times \gol$
implies $\hbox{ad}_\xi^*\alpha\in \gos^*$ and
$(\alpha,\xi)\in \gos^*\times \gos$
implies $\hbox{ad}_\xi^*\alpha\in \gol^*$.
To summarize:
\begin{equation}\label{adjstability}
\begin{array}{cc}
 \hbox{ad}_\gol^*\gol^* \subset \gol^*, &
 \hbox{ad}_\gol^*\gos^* \subset \gos^* \\
 \hbox{ad}_\gos^*\gol^* \subset \gos^*, &
 \hbox{ad}_\gos^*\gos^* \subset \gol^*
\end{array}
\end{equation}

\subsubsection{The space of fields and action functional}
We assume that $\widehat{\goP}$ and $\widehat{\goL}$
satisfy Hypotheses (\ref{hyposplustegalp},\ref{hypotstablesouss},\ref{hypottincluss}). We suppose that there exists
some tensor
$\kappa{_\gop}^{\gos\gos} \ \in\
 \gop^*\otimes \Lambda^2\gos
 \ \subset \
 \gop^*\otimes \gos\otimes \gos$
which is invariant by the adjoint action of $\widehat{\goL}$:
\begin{equation}\label{invarianceOfKappa}
 \hbox{Ad}^*_g\otimes\hbox{Ad}_g\otimes\hbox{Ad}_g
 (\kappa{_\gop}^{\gos\gos})
 = \kappa{_\gop}^{\gos\gos},\quad \forall g\in \widehat{\goL}.
\end{equation}
A fundamental example of a tensor $\kappa{_\gop}^{\gos\gos}$ is presented in \S \ref{paragraph1.2.5}.
We fix a non vanishing volume form $\hbox{vol}_\gop\in \Lambda^{N}\gop^*$ and we consider a $N$-dimensional oriented manifold
$\mathcal{F}$.
We then consider the class of fields
\begin{equation}\label{constraintaxiomatic}
\begin{array}{ccr}
 \mathscr{E}_{\textsf{E}} & := & \left\{(\pi_\gop,\varphi^\gop)\in (\gop^*\otimes \Omega^{N-2}(\mathcal{F}))\times
(\gop\otimes \Omega^1(\mathcal{F})) \hbox{ of class }\mathscr{C}^2 \right. \\
& & \left. \hbox{rank}\varphi^\gop_\textsf{z} = N, \forall \textsf{z}\in \mathcal{F} \hbox{ and }
\pi_\gop\wedge \varphi^{\gos\gos} = \kappa{_\gop}^{\gos\gos}\;(\varphi^\gop)^*\hbox{vol}_\gop \right\}
\end{array}
\end{equation}
we set $\Phi^\gop:= \dR \varphi^\gop + \frac{1}{2}[\varphi^\gop\wedge \varphi^\gop]$
and we define on $\mathscr{E}_{\textsf{E}}$ the functional
\begin{equation}\label{functionalAEinstein}
 \mathscr{A}[\pi_\gop,\varphi^\gop] = \int_{\mathcal{F}}
 \pi_{\celp}\wedge \Phi^{\celp}
\end{equation}
\begin{theo}\label{theoBigOne}
Let $\widehat{\goP}$ be a simply connected Lie group of  finite dimension $N$ and $\widehat{\goL}\subset \widehat{\goP}$ a simply connected Lie subgroup. Let $\gop$ and $\gol$ be their respective Lie algebras. Assume that $\gop$ and $\gol$ are \emph{unimodular}, that there exists a vector subspace $\gos\subset \gop$ which is stable by $\hbox{\emph{Ad}}_{\widehat{\goL}}$ and such that $\gop = \gos\oplus \gol$ and that $\widehat{\goP}/\widehat{\goL}$ is a \emph{symmetric space} (\ref{hypottincluss}).

Let $\mathcal{F}$ be smooth oriented manifold of dimension $N$ and consider the functional $\mathscr{A}$ defined by (\ref{functionalAEinstein}) on $\mathscr{E}_{\textsf{\emph{E}}}$.
Assume that $\kappa{_\gop}^{\gos\gos}$ (in the definition of $\mathscr{E}_{\textsf{\emph{E}}}$) satisfies the \emph{additional hypothesis}:
\begin{equation}\label{additionalhypothesis}
 \kappa{_\gop}^{\gos\gos} =
 \kappa{_\gol}^{\gos\gos} \in \gol^*\otimes \Lambda^2\gos
 \quad \hbox{ i.e. }\kappa{_\gos}^{\gos\gos} = 0.
\end{equation}
Let $(\pi_\gop,\varphi^\gop)\in \mathscr{E}_{\textsf{\emph{E}}}$ be a smooth critical point of $\mathscr{A}$. Then
\begin{enumerate}
 \item $\mathcal{F}$ is foliated by smooth leaves $\textsf{\emph{f}}$ of dimension $r:= \hbox{\emph{dim}}\gol$, which are solutions of the exterior differential system $\theta^\gos|_{\textsf{\emph{f}}} = 0$.
 \item for any point in $\mathcal{F}$ there exists an open neighbourhood $\mathcal{O}\subset \mathcal{F}$ of this point such we can endow the set of 
intersections $\mathcal{X}_\mathcal{O}:= \{\textsf{\emph{f}}\cap \mathcal{O}\,;\, \textsf{\emph{f}} \hbox{ is an integral leaf}\}$ with a structure of manifold $\mathcal{X}_\mathcal{O}$ of dimension $n:= \hbox{\emph{dim}}\gos$.
 \item there exist local charts $\mathcal{O}\ni \textsf{\emph{z}} \longmapsto (\textsf{\emph{x}},g)\in \mathcal{X}_\mathcal{O}\times \widehat{\goL}$, such that the projection map $\mathcal{O}\xrightarrow{\hbox{ }P_\mathcal{O}\hbox{ }} \mathcal{X}_\mathcal{O}$ is a submersion and we have the decompositions $\varphi^\gos = g^{-1}\theta^\gos g$ and $\varphi^\gol = g^{-1}\omega^\gol g + g^{-1}\emph{\dR}g$, where $\theta^\gos$ and $\omega^\gol$ are pull-backs by $\mathcal{O}\xrightarrow{\hbox{ }P_\mathcal{O}\hbox{ }} \mathcal{X}_\mathcal{O}$ of 1-forms on $\mathcal{X}_\mathcal{O}$. Moreover $\textbf{\emph{g}}:=(\varphi^\gos)^*\textsf{\emph{b}}:=
 \textsf{\emph{b}}_{\cels\cels}\varphi^\cels\otimes \varphi^\cels$ is the pull-back  by $\mathcal{O}\xrightarrow{\hbox{ }P_\mathcal{O}\hbox{ }} \mathcal{X}_\mathcal{O}$ of a pseudo metric (also denoted by) $\textbf{\emph{g}}$ on $T\mathcal{X}_\mathcal{O}$ and $\theta^\gos$ provides us with an orthonormal coframe for $\textbf{\emph{g}}$ and $\omega^\gog$ defines a connection on $T\mathcal{X}_\mathcal{O}$ which respects  $\textbf{\emph{g}}$.
 \item $\theta^\gos$, $\omega^\gol$ and $p_\gop:= \hbox{\emph{Ad}}^*_g\pi_\gop$ are solutions of the following equations
\begin{equation}\label{newugly}
 \begin{array}{c}
  \frac{1}{2}\kappa{_\gog}^{\cels_1\cels_2}
  \Omega{^\gol}_{\cels_1\cels_2}
+ (\partial_{\cels}^\omega
+ \Theta{^{\ast}}_{\cels\ast})
  p{_\gop}^{\gol\cels}
  + \textbf{c}{^{\cels_0}}_{\gol\cels}\ p{_{\cels_0}}^{\gol\cels}
  + \partial_{\cell_1}p{_\gop}^{\gol\cell_1}
  + \frac{1}{2}\textbf{c}{^\gol}_{\cell_1\cell_2}p{_\gop}^{\cell_1\cell_2}
    \\
 = \Theta{^{\cels_0}}_{\gos\cels}\,
p{_{\cels_0}}^{\gol\cels}
+ \Omega{^{\cell}}_{\gos\cels} p{_{\cell}}^{\gol\cels}
- \frac{1}{2} Q\delta{_\gol}^\gol
 \end{array}
\end{equation}
and
\begin{equation}\label{newnice}
\frac{1}{2}\kappa{_\gog}^{\cels_1\cels_2}
 \mathring{\Theta}{^\gos}_{\cels_1\cels_2}
 + \partial_{\cell}p{_\gop}^{\gos\cell}
= \Omega{^{\cell}}_{\gos\cels_1}
\kappa{_{\cell}}^{\gos\cels_1}
-\frac{1}{2} (\Omega{^{\cell}}_{\cels_1\cels_2}
+ \mathbf{c}{^{\cell}}_{\cels_1\cels_2})
\kappa{_{\cell}}^{\cels_1\cels_2}\delta{_\gos}^\gos 
\end{equation}
where we set $\Theta{^\ast}_{\gos\ast}:= \Theta{^\cels}_{\gos\cels}$ and $\mathring{\Theta}{^\gos}_{\gos_1\gos_2}
=\Theta{^\gos}_{\gos_1\gos_2} - \delta^\gos_{\gos_2}\Theta{^\ast}_{\gos_1\ast} + \delta^\gos_{\gos_1}\Theta{^\ast}_{\gos_2\ast}$ and
\begin{equation}\label{premierQ}
Q = \Omega{^{\cell}}_{\cels\cels}\kappa{_{\cell}}^{\cels\cels}
+ \mathbf{c}{^{\cell}}_{\cels\cels}
\kappa{_{\cell}}^{\cels\cels}
\end{equation}
\item if we assume furthermore that the integral leaves $\textsf{\emph{f}}$ are the fibers of a global fibration $\mathcal{F}\xrightarrow{\hbox{ }P\hbox{ }} \mathcal{X}$, then the previous equations make sense on this fiber bundle
\end{enumerate}
\end{theo}
Comments on Equations (\ref{newugly}) and (\ref{newnice}) may be welcome. By defining the \emph{generalized Cartan tensor} $\tilde{\mathbf{C}}{_\gog}^\gos:= - \frac{1}{2} \kappa{_\gog}^{\cels_1\cels_2} \mathring{\Theta}{^\gos}_{\cels_1\cels_2}$  (equivalent to the torsion tensor $\Theta{^\gos}_{\gos\gos}$ in most situations), the \emph{generalized Einstein tensor} $\tilde{\mathbf{E}}{_\gos}^\gos = \Omega{^\celg}_{\cels_1\gos}\kappa{_\celg}^{\cels_1\gos}- \frac{1}{2}(\Omega{^\celg}_{\cels_1\cels_2}\kappa{_\celg}^{\cels_1\cels_2})\delta{_\gos}^\gos$ (see (\ref{deltass}) for the definition of $\delta{_\gos}^\gos$) and by setting $T{_\gop}^\gos:= \partial_\cell p{_\gop}^{\gos\cell}$, Equation (\ref{newnice}) has the form of a generalized Einstein--Cartan system
\begin{equation}
 \left\{
 \begin{array}{ccl}
 \tilde{\mathbf{C}}{_\gol}^\gos
 & = & T{_\gol}^\gos
 \\
\tilde{\mathbf{E}}{_\gos}^\gos + \Lambda\delta{_\gos}^\gos
& = &
T{_\gos}^\gos
 \end{array}\right.
\end{equation}
where  $\Lambda := - \frac{1}{2}\mathbf{c}{^\celg}_{\cels_1\cels_2}\kappa{_\celg}^{\cels_1\cels_2}$. Hence $T{_\gol}^\gos$ can be interpreted as an angular momentum tensor and $T{_\gos}^\gos$ as a stress-energy tensor.

Equation (\ref{newugly}) does not look that friendly but leads however to interesting open questions. We prove in Lemma \ref{lemmedonnepartialT} that, independently of (\ref{newnice}), Equation (\ref{newugly}) implies that $T{_\gop}^\gos$ is a solution of
\begin{equation}\label{conservationdesdeuxtenseurs}
 \partial_\cels^\omega T{_\gop}^\cels + \Theta{^{\ast}}_{\cels\ast}
T{_\gop}^\cels
+ \textbf{c}{^{\cels_0}}_{\gol\cels}\ T{_{\cels_0}}^\cels
  =
  \Theta{^{\cels_0}}_{\gos\cels_1}\,
T{_{\cels_0}}^{\cels_1}
+ \Omega{^{\cell_0}}_{\gos\cels_1}
T{_{\cell_0}}^{\cels_1}
\end{equation}
which expresses the conservation of the angular and the stres-energy momentum tensors. We will derive in Proposition \ref{propoBianchi} the constraint equations on the Cartan and the Einstein tensors which derive from the Bianchi identities and check that they are compatible with (\ref{conservationdesdeuxtenseurs}).

\begin{coro}\label{corooftheBigOne}
Assume all the hypotheses of Theorem \ref{theoBigOne} and that furthermore the integral leaves $\textsf{\emph{f}}$ are the fibers of a principal bundle structure $\mathcal{F}\xrightarrow{\hbox{ }P\hbox{ }} \mathcal{X}$ with structure group $\goL$, where $\goL$ is a quotient of $\widehat{\goL}$ by a finite subgroup.
 
Assume in addition that $\goL$ is compact, or that the first derivatives of $p{_\gop}^{\gog\gos}$ decay to zero at infinity in each fiber. Then the fields $\theta^\gos$ and $\omega^\gog$ are solutions of a generalized Einstein--Cartan system of equations in vacuum, i.e.
\[
 \tilde{\mathbf{C}}{_\gol}^\gos
 =
\tilde{\mathbf{E}}{_\gos}^\gos + \Lambda\delta{_\gos}^\gos = 0
\]

\end{coro}
The next paragraphs until \S \ref{theend} are devoted to the proof of Theorem \ref{theoBigOne}.
Most computations will be performed without assuming Hypothesis (\ref{additionalhypothesis}). The latter hypothesis will be used only in the conclusion. The proof of Corollary \ref{corooftheBigOne} will be given in Section \ref{sectionexploitation}.

\subsection{Study of the critical points}
We let $(\textbf{t}_a)_{1\leq a\leq n}$ be a basis of $\gos$ and
let $(\textbf{t}_i)_{n+1\leq i\leq N}$ be a basis of $\gol$.
Then $(\textbf{t}_I)_{1\leq I\leq N}:= (\textbf{t}_a)_{1\leq a\leq n}\cup (\textbf{t}_i)_{n+1\leq i\leq N}$
is a basis of $\gop$. Here we make the following
implicit assumptions on the indices:
$1\leq I,J,K,\ldots \leq N$,
$1\leq a,b,c,\ldots \leq n$
and $n+1\leq i,j,k,\ldots \leq N$.
We denote by $(\textbf{t}^I)_{1\leq I\leq N}$ the basis of $\gop^*$
which is dual of $(\textbf{t}_I)_{1\leq I\leq N}$.
Note that $(\textbf{t}^a)_{1\leq a\leq n}$ is a basis of
$\gos^*:=\gol^{\perp}$
and $(\textbf{t}^i)_{n< i\leq N}$ is a basis of
$\gol^*:= \gos^{\perp}\subset \gop$.
We denote by $\textbf{c}^K_{IJ}$ the structure coefficients of $\gop$ such
that $[\textbf{t}_I,\textbf{t}_J] = \textbf{t}_K\textbf{c}^K_{IJ}$ and $\hbox{ad}^*_{\textbf{t}_I}\textbf{t}^J
 = - \textbf{c}^J_{IK}\textbf{t}^K$.
We can thus decompose
$\kappa{_\gop}^{\gos\gos}
 = \frac{1}{2}\kappa{_I}^{bc}\;
 \textbf{t}^I\otimes(\textbf{t}_b\wedge\textbf{t}_c)\
 \simeq \frac{1}{2}\kappa{_I}^{bc}
 \textbf{t}^I(\textbf{t}_b\wedge\textbf{t}_c)$.

Without loss of generality we assume that $\hbox{vol}_\gop = \textbf{t}^1\wedge \cdots \wedge \textbf{t}^{N}$.
Hence the constraint
$\pi_\gop\wedge \varphi^\gos\wedge \varphi^\gos  = \kappa{_\gop}^{\gos\gos}\varphi^*\hbox{vol}_\gop$ reads
\begin{equation}\label{constraintaxiomaticindex}
 \pi_{\gop} \wedge \varphi^{\gos}\wedge \varphi^{\gos} = \kappa{_{\gop}}^{\gos\gos}\; \varphi^{(N)},\quad \hbox{ where }\varphi^{(N)}:= \varphi^1\wedge \cdots\wedge \varphi^{N}.
\end{equation}
Since $\varphi^\gop\in \gop\otimes \Omega^1(\mathcal{F})$ is
a coframe on $\mathcal{F}$ we can decompose
\begin{equation}\label{18gravity}
 \Phi^\gop:= \dR\varphi^\gop + \frac{1}{2}[\varphi^\gop\wedge \varphi^\gop]
  = \frac{1}{2}\Phi{^{\gop}}_{\celp\celp}\,\varphi^{\celp\celp}
\end{equation}
(see (\ref{ConventionV}))
and
$\pi_\gop\in \gop^* \otimes \Omega^{N-2}(\mathcal{F})$ as
(see (\ref{ConventionV*}))
$\pi_{\gop} = \frac{1}{2}
 \pi{_\gop}^{\celp\celp}
 \ \varphi^{(N-2)}_{\celp\celp}
 = \frac{1}{2}
 \pi{_\gop}^{\cels\cels}
 \ \varphi^{(N-2)}_{\cels\cels}
 + \pi{_\gop}^{\cels\cell}
 \ \varphi^{(N-2)}_{\cels\cell}
 + \frac{1}{2}  \pi{_\gop}^{\cell\cell}
 \ \varphi^{(N-2)}_{\cell\cell}$. Condition (\ref{constraintaxiomaticindex}) reads
\begin{equation}\label{constraintcoordinates}
  \pi{_\gop}^{\gos\gos} = \kappa{_\gop}^{\gos\gos}.
\end{equation}

\noindent
\textbf{First variation with respect to the coefficients of $\pi^\gop$} ---
We look at infinitesimal variations of the form
\[
(\varphi^\gop,\pi_\gop)\longmapsto
(\varphi^\gop ,\pi_\gop + \varepsilon \delta\pi_\gop)
\]
where $\delta\pi_\gop = \chi_\gop$ has the form $\chi_\gop = \chi{_\gop}^{\cels\cell}\varphi^{(N-2)}_{\cels\cell} + \frac{1}{2}\chi{_\gop}^{\cell\cell}\varphi^{(N-2)}_{\cell\cell}$, so that the constraint  (\ref{constraintaxiomaticindex}) is preserved. The first variation of the action vanishes under such variations iff
\[
 \forall \chi{_\gop}^{\gos\gol}, \chi{_\gop}^{\gol\gol}\varphi^{(N-2)}_{\cell\cell},\quad
 \int_\mathcal{F}\left(\chi{_\gop}^{\cels\cell}\Phi{^\celp}_{\cels\cell} + \frac{1}{2}\chi{_\gop}^{\cell\cell}\Phi{^\celp}_{\cell\cell}
 \right)\varphi^{(N)} = 0
\]
which leads to the equations
\begin{equation}\label{gravEL1}
\Phi{^{\gop}}_{\gol\gol} =
\Phi{^{\gop}}_{\gos\gol} = 0
\end{equation}

\noindent
\textbf{First variation with respect to $\varphi^\gop$} ---
We compute the first variation of the action under an infinitesimal variation of $(\varphi^\gop,\pi_\gop)$ of the form
\[
(\varphi^\gop,\pi_\gop)\longmapsto
(\varphi^\gop + \varepsilon \delta\varphi^\gop,\pi_\gop + \varepsilon \delta\pi_\gop + o(\varepsilon)),
\]
where $\delta \varphi^\gop = \lambda^\gop = \lambda{^\gop}_{\celp}\varphi^{\celp}$ has a compact support and by keeping the coefficients $\pi{_\gop}^{\gop\gop}$ constant, so that
$\delta\pi_\gop = \frac{1}{2}\pi{_\gop}^{\celp_1\celp_2}\lambda{^{\celp_3}}\wedge \varphi^{(N-3)}_{\celp_1\celp_2\celp_3}$ (we hence preserve the constraint (\ref{constraintaxiomaticindex})). The first order variation of $\pi_\celp\wedge \Phi^\celp$ splits as
$\delta (\pi_\celp\wedge \Phi^\celp ) =
\delta\pi_\celp\wedge \Phi^\celp + \pi_\celp\wedge \delta\Phi^\celp$. On the one hand:
\[
\begin{array}{ccl}
 \delta\pi_\celp\wedge \Phi^\celp
 & = & \left(\frac{1}{2}\pi{_\celp}^{\celp_1\celp_2}\,\lambda{^{\celp_3}}\wedge \varphi^{(N-3)}_{\celp_1\celp_2\celp_3}\right) \wedge
 \left(\frac{1}{2}\Phi{^\celp}_{\celp_4\celp_5}\varphi^{\celp_4\celp_5}\right) \\
 & = & \frac{1}{2}\Phi{^\celp}_{\celp_4\celp_5}
 \left(\pi{_\celp}^{\celp_1\celp_4}\lambda{^{\celp_5}}\wedge \varphi^{(N-1)}_{\celp_1}
 + \pi{_\celp}^{\celp_5\celp_2}\lambda{^{\celp_4}}\wedge \varphi^{(N-1)}_{\celp_2}
 + \pi{_\celp}^{\celp_4\celp_5}\lambda{^{\celp_3}}\wedge \varphi^{(N-1)}_{\celp_3}
 \right) \\
 & = & \frac{1}{2}\Phi{^\celp}_{\celp_4\celp_5}
 \left(\pi{_\celp}^{\celp_1\celp_4}\lambda{^{\celp_5}}_{\celp_1}
 + \pi{_\celp}^{\celp_5\celp_2}\lambda{^{\celp_4}}_{\celp_2}
 + \pi{_\celp}^{\celp_4\celp_5}\lambda{^{\celp_3}}_{\celp_3}
 \right)\varphi^{(N)}
\end{array}
\]
It is then convenient to introduce the following notations
\begin{equation}\label{definitionfamillePsi}
  \left\{\begin{array}{ccccc}
         \Psi{_{\gop\gop}}^{\gop\gop} & := & & &
         \Phi{^{\celp}}_{\gop\gop}\pi{_{\celp}}^{\gop\gop} \\
         \Psi{_\gop}^\gop & := &
         \Psi{_{\gop\celp}}^{\gop\celp} & = &
         \Phi{^{\celp}}_{\gop\celp_2}\pi{_{\celp}}^{\gop\celp_2} \\
         \Psi & := & \Psi{_{\celp}}^{\celp}\quad & = &
         \Phi{^{\celp}}_{\celp_1\celp_2}\pi{_{\celp}}^{\celp_1\celp_2}
        \end{array}\right.
\end{equation}
so that we obtain
$\delta\pi_\celp\wedge \Phi^\celp
 = \frac{1}{2}
 \left(
 - \Psi{_{\celp_5}}^{\celp_1}
 \lambda{^{\celp_5}}_{\celp_1}
 - \Psi{_{\celp_4}}^{\celp_2}
 \lambda{^{\celp_4}}_{\celp_2}
 + \Psi \lambda{^{\celp_3}}_{\celp_3}
 \right)\varphi^{(N)}$, i.e.
\[
 \delta\pi_\celp\wedge \Phi^\celp
 = - \left(
  \Psi{_{\celp_1}}^{\celp_2}
 \lambda{^{\celp_1}}_{\celp_2}
 - \frac{1}{2}\Psi \lambda{^{\celp}}_{\celp}
 \right)\varphi^{(N)}
 = -\left( \Psi{_{\celp_1}}^{\celp}
 \lambda{^{\celp_1}}- \frac{1}{2}\Psi \lambda{^{\celp}}
 \right)\wedge \varphi^{(N-1)}_{\celp}
\]
On the other hand $\delta\Phi^\gop = \dR\lambda^{\gop} + [\varphi^\gop\wedge \lambda^\gop] = \dR^\varphi \lambda^\gop$ and thus, by (\ref{twistedLeibniz}),
\[
 \pi_{\celp}\wedge
\delta\Phi^{\celp} =
(\dR^\varphi\lambda^{\celp})
\wedge\pi_{\celp}
 = \dR^\varphi(\lambda^{\celp}
\wedge\pi_{\celp})
+ \lambda^{\celp}
\wedge \dR^\varphi\pi_{\celp}
\]
where actually, since the coefficients of $\lambda^{\celp}
\wedge\pi_{\celp}$ are in a trivial represent of $\gop$, $\dR^\varphi(\lambda^{\celp}
\wedge\pi_{\celp}) = \dR(\lambda^{\celp}
\wedge\pi_{\celp})$.

In conclusion
$\delta (\pi_\celp\wedge \Phi^\celp ) =
\dR(\lambda^{\celp}
\wedge\pi_{\celp})
+ \lambda^{\celp}
\wedge \left(\dR^\varphi\pi_{\celp} - \Psi{_\celp}^{\celp_1}\varphi^{(N-1)}_{\celp_1} + \frac{1}{2}\Psi\varphi^{(N-1)}_\celp\right)$. Thus since $\lambda^\gop$ has a compact support the action is stationary with respect to these variations iff
\[
 \forall \lambda^\gop,\quad
 \int_\mathcal{F}\lambda^{\celp}
\wedge \left(\dR^\varphi\pi_{\celp} - \Psi{_\celp}^{\celp_1}\varphi^{(N-1)}_{\celp_1} + \frac{1}{2}\Psi\varphi^{(N-1)}_\celp\right)
= 0
\]
which leads to the equation (see (\ref{deltass}) for the definition of $\delta{_\gop}^\gop$)
\begin{equation}\label{gravEL2}
 \dR^\varphi\pi_{\gop} = \Psi{_\gop}^{\celp}\varphi^{(N-1)}_{\celp} - \frac{1}{2}\Psi\varphi^{(N-1)}_\gop
 = \left(\Psi{_\gop}^{\celp}- \frac{1}{2}\Psi\delta{_\gop}^\celp\right)\varphi^{(N-1)}_{\celp}
\end{equation}
We observe that direct consequences of  (\ref{gravEL1}) and (\ref{definitionfamillePsi})
are $\Psi{_{\gol\gol}}^{\gop\gop} = \Psi{_{\gos\gol}}^{\gop\gop} = 0$
and hence
\[
 \Psi{_\gol}^\gop = 0,\quad
 \Psi{_\gos}^\gop = \Psi{_{\gos\cels}}^{\gop\cels} = \Phi{^{\celp}}_{\gos\cels}\pi{_{\celp}}^{\gop\cels},
 \quad \Psi = \Psi{_{\cels}}^{\cels} = \Phi{^{\celp}}_{\cels\cels}\pi{_{\celp}}^{\cels\cels}
\]
This implies that (\ref{gravEL2}) can be written
\[
 \left\{\begin{array}{ccc}
\dR^\varphi\pi{_\gos} &  = & \Psi{_\gos}^\celp\varphi^{(N-1)}_\celp -
\frac{1}{2}\Psi\,\varphi^{(N-1)}_\gos\\
\dR^\varphi\pi{_\gol} & = &
- \frac{1}{2}\Psi\,\varphi^{(N-1)}_\gol
 \end{array}\right.
\]
or $\dR^\varphi\pi{_{\gop}} = \Psi{_{\gop}} - \frac{1}{2}\Psi\,
\varphi^{(N-1)}_{\gop}$,
where
$\Psi{_{\gop}} :=
\Psi{_{\gop}}^{\celp}
\, \varphi^{(N-1)}_{\celp}$.
In conclusion the Euler--Lagrange system is
\begin{equation}\label{gravELcompact}
\left\{\begin{array}{ccl}
    \Phi{^\gop}_{\gol\gol} & = & \Phi{^\gop}_{\gos\gol} = 0 \\
    \dR^\varphi\pi{_{\gop}} & = &
    \displaystyle
    \Psi{_{\gop}} - \frac{1}{2}\Psi\,
\varphi^{(N-1)}_{\gop}
\end{array}\right.
\end{equation}
or, by splitting $\gop = \gol \oplus \gos$ and by using the relation $\Psi{_\gol}^\gop = 0$,
\begin{equation}\label{gravELa}
\left\{\begin{array}{ccl}
    \dR\varphi^\gop +
    \frac{1}{2}[\varphi^\gop\wedge \varphi^\gop] & = &
 \displaystyle \frac{1}{2}\Phi{^\gop}_{\cels\cels} \ \varphi^{\cels\cels} \\
 \dR^\varphi\pi_\gop & = & 
\Psi{_\gos}^\cell\  \varphi^{(N-1)}_\cell
+ \Psi{_\gos}^\cels\  \varphi^{(N-1)}_\cels
-\frac{1}{2}\Psi\ \varphi^{(N-1)}_\gop
\end{array}
\right.
\end{equation}

\subsubsection{Spontaneous foliation}
We first exploit the Euler--Lagrange equation
$\dR\varphi^\gop + \frac{1}{2}[\varphi^\gop\wedge \varphi^\gop] = \frac{1}{2}\Phi{^\gop}_{\cels\cels} \varphi^{\cels\cels}$.
For that purpose we split $\varphi^\gop = \varphi^\gol + \varphi^\gos$ and similarly
$[\varphi^\gop\wedge \varphi^\gop] = [\varphi^\gop\wedge \varphi^\gop]^\gol +
[\varphi^\gop\wedge \varphi^\gop]^\gos$, according to the decomposition.
${\gop} = \gol \oplus \gos$
We have
\[
 [\varphi^\gop\wedge \varphi^\gop] = [(\varphi^\gol + \varphi^\gos)\wedge
 (\varphi^\gol + \varphi^\gos)]
 = [\varphi^\gol\wedge \varphi^\gol]
 + 2[\varphi^\gol\wedge \varphi^\gos]
 + [\varphi^\gos\wedge \varphi^\gos].
\]
Thus by using the hypotheses $[\gol,\gol]\subset \gol$,
$[\gos,\gos]\subset \gol$
and $[\gol,\gos]\subset \gos$, we deduce
\[
\begin{array}{cclll}
 [\varphi^\gop\wedge \varphi^\gop]^\gol & = & [\varphi^\gol\wedge \varphi^\gol]
& & +\; [\varphi^\gos\wedge \varphi^\gos] \\
 {[\varphi^\gop\wedge \varphi^\gop]}^\gos & = & & 2[\varphi^\gol\wedge \varphi^\gos]
&
\end{array}
\]
Hence the relation $\dR\varphi^\gop + \frac{1}{2}[\varphi^\gop\wedge \varphi^\gop] = \frac{1}{2}\Phi{^\gop}_{\cels\cels} \varphi^\cels\wedge \varphi^\cels$
is equivalent to
\begin{equation}\label{splittingEL}
 \left\{\begin{array}{ccl}
         \dR\varphi^\gol + \frac{1}{2}[\varphi^\gol\wedge \varphi^\gol]
+ \frac{1}{2}[\varphi^\gos\wedge \varphi^\gos]
         & = & \frac{1}{2}\Phi{^\gol}_{\cels\cels} \varphi^\cels\wedge \varphi^\cels \\
\dR\varphi^\gos + [\varphi^\gol\wedge \varphi^\gos]
& = & \frac{1}{2}\Phi{^\gos}_{\cels\cels} \varphi^\cels\wedge \varphi^\cels
        \end{array}
\right.
\end{equation}
In order to apply these relations
let us look for $r$-dimensional submanifolds $\textsf{f}$ of $\mathcal{F}$ which are
solutions of the Pfaffian system
\begin{equation}\label{AbstractPfaffianSystem0}
 \varphi^\gos|_\textsf{f} = 0
\end{equation}
By using the second equation in (\ref{splittingEL})
($\dR\varphi^\gos + [\varphi^\gol\wedge\varphi^\gos]
= \frac{1}{2}\Phi{^\gos}_{\cels\cels}\varphi^{\cels\cels}$, which implies $\hbox{d}\varphi^\gos = 0 \hbox{ mod }[\varphi^\gos]$)
we deduce from Frobenius' theorem that, for any point $\textsf{z}\in \mathcal{F}$, there exists
a neighbourhood of $\textsf{z}$ in $\mathcal{F}$ such that there exists an unique solution $\textsf{f}$
to (\ref{AbstractPfaffianSystem0}) that passes through $\textsf{z}$.
We hence deduce the existence of a foliation of $\mathcal{F}$ by leaves $\mathsf{f}$ of
dimension $r$ and codimension $n$. For any $\textsf{z}\in \mathcal{F}$, we denote by $\textsf{f}_\textsf{z}$ the unique integral leaf which contains $\textsf{z}$.

We denote by $\mathcal{X}:= \{ \mathsf{f}_\textsf{z};\textsf{z}\in \mathcal{F}\}$ the set of leaves and
\begin{equation}
 \begin{array}{cccl}
  \textsf{x}: & \mathcal{F} & \longrightarrow & \mathcal{X} \\
  & \textsf{z} & \longmapsto & \textsf{x}(\textsf{z})\hbox{ such that }\textsf{z}\in \textsf{f}_{\textsf{x}(\textsf{z})}
 \end{array}
\end{equation}
the quotient map. Note that in general
$\mathcal{X}$ is just a topological space and may not be a manifold, unless
it is a separated (Hausdorff) space.

In the following we restrict ourself to some open subset $\mathcal{O}$ of $\mathcal{F}$ such that
there exists an $n$-dimensional submanifold $\Sigma$ which crosses transversally each leaf in $\mathcal{O}$ at one and only one point. Then the image of the restriction $\textsf{x}|_\mathcal{O}$
has the structure of an $n$-dimensional manifold, which may be identified with an
open subset of $\Sigma$.

\subsubsection{Local  principal bundle structure and trivialization}
Consider the product manifold $\mathcal{O} \times \widehat{\goL}:= \{(\textsf{z},h);\; \textsf{z}\in \mathcal{O},h\in \widehat{\goL}\}$
and the $\gol$-valued 1-form $\psi^\gol\in \gol\otimes \Omega^1(\mathcal{O}\times \goL)$
defined by $\psi^\gol:= \dR h -h\varphi^\gol$.
Observe that
\[
 \dR\psi^\gol =  - h\left(\dR\varphi^\gol + \frac{1}{2}[\varphi^\gol\wedge \varphi^\gol]\right)
 - \psi^\gol\wedge \varphi^\gol.
\]
However the first equation in (\ref{splittingEL}) 
implies that the restriction of $\dR\varphi^\gol +\frac{1}{2}[\varphi^\gol\wedge \varphi^\gol]$
on any leaf $\textsf{f}$ vanishes:
$\dR\varphi^\gol +\frac{1}{2}[\varphi^\gol\wedge \varphi^\gol]|_\textsf{f} = 0$.
Thus, for any leaf $\textsf{f}$, $\dR\psi^\gol|_{\textsf{f}\times \widehat{\goL}}= 0\hbox{ mod }[\psi^\gol]$,
which implies by Frobenius' theorem that
the Pfaffian system $\psi^\gol|_{\textsf{f}\times \widehat{\goL}} = 0$ is integrable on each fiber
and thus there exists a map $g:\mathcal{O}\longrightarrow \widehat{\goL}$ such that
$(\hbox{Id}\times g)^*\psi^\gol|_\textsf{f}=0$, i.e.
\begin{equation}\label{GstructureAtteinte}
 \dR g|_{\textsf{f}} = g\varphi^\gol |_{\textsf{f}}
\quad \Longleftrightarrow \quad
\varphi^\gol|_{\textsf{f}} = g^{-1}\dR g|_{\textsf{f}}.
\end{equation}
Moreover by requiring that $g$ is equal to $1_{\widehat{\goL}}$ on $\Sigma$, $g$ is unique.

Note that since the rank of $\varphi^\gol$ is equal to the dimension of the leaves,
each restriction map $g|_\textsf{f}$ is a local diffeomorphism between $\textsf{f}$
and an open neighbourhood of $1_{\widehat{\goL}}$ in $\widehat{\goL}$.
Thus, by replacing $\mathcal{O}$ by another open subset if necessary,
we deduce that there exists a neighbourhood
$\mathcal{V}_{\widehat{\goL}}$ of $1_{\widehat{\goL}}$ in $\widehat{\goL}$
such that the map
\[
 \begin{array}{ccl}
  \mathcal{O} & \longrightarrow & \Sigma \times \mathcal{V}_{\widehat{\goL}} \subset \Sigma \times \widehat{\goL} \\
  \textsf{z} & \longmapsto & (\textsf{x}(\textsf{z}),g)
 \end{array}
\]
is a diffeomorphism.

Let us define the ${\gop}$-valued 1-forms
\begin{equation}\label{introductione&A}
 e^\gop:= \hbox{Ad}_g\varphi^\gop
 \quad \hbox{ and }\quad
 \textbf{A}^\gop:= \hbox{Ad}_g\varphi^\gop - \dR g\, g^{-1} = e^\gop - \dR g\, g^{-1}
\end{equation}
or equivalentely
\begin{equation}\label{introductione&Abis}
 \left\{ \begin{array}{ccl}
  e^\gos & := & \hbox{Ad}_g\varphi^\gos  \\
  e^\gol & := & \hbox{Ad}_g\varphi^\gol
 \end{array}\right.
 \quad \hbox{ and }\quad
 \left\{ \begin{array}{cclcl}
  \textbf{A}^\gos & := & \hbox{Ad}_g\varphi^\gos & = & e^\gos \\
  \textbf{A}^\gol & := & \hbox{Ad}_g\varphi^\gol - \dR g\, g^{-1}
  & = & e^\gol - \dR g\, g^{-1}
 \end{array}\right.
\end{equation}
and the ${\gop}$-valued 2-form
\begin{equation}\label{definitionF}
 \textbf{F}^\gop:= \dR\textbf{A}^\gop + \frac{1}{2}[\textbf{A}^\gop\wedge \textbf{A}^\gop]
\end{equation}
A direct computation of $\textbf{F}^\gop$ gives the following. We denote by $\textbf{A}{^\gop}_\gop$ the coefficients in the decomposition $\textbf{A}{^\gop} = \textbf{A}{^\gop}_\celp e^\celp$ and by $\partial_\gop \textbf{A}{^\gop}_\gop$ the coefficients such that $\dR\textbf{A}{^\gop}_\gop = \partial_\celp\textbf{A}{^\gop}_\gop e^\celp$. We obtain
\begin{equation}\label{calculdirectF}
 \textbf{F}^\gop
 = \frac{1}{2} \left(\partial_{\celp_1}\textbf{A}{^\gop}_{\celp_2}
 - \partial_{\celp_2}\textbf{A}{^\gop}_{\celp_1} + [\textbf{A}{^\gop}_{\celp_1}, \textbf{A}{^\gop}_{\celp_1}]\right)e^{\celp_1\celp_2}
\end{equation}
By (\ref{introductione&Abis}) Equation  (\ref{AbstractPfaffianSystem0}) translates as $e^\gos|_\textsf{f} = 0$.
Still by (\ref{introductione&Abis})
we get $\varphi^\gol = \hbox{Ad}_{g^{-1}}\textbf{A}^\gol + g^{-1}\dR g$, so that Relation (\ref{GstructureAtteinte}) reads $\textbf{A}^\gol|_\textsf{f} = 0$.
The latter relation is thus equivalent to $\textbf{A}^\gol = \textbf{A}{^\gol}_\cels e^\cels$ (i.e. $\textbf{A}{^\gol}_\gol = 0$). But we also have $\textbf{A}^\gos = e^\gos$ and thus we conclude that $\textbf{A}^\gop = \textbf{A}{^\gop}_\cels e^\cels$ (i.e. $\textbf{A}{^\gop}_\gol = 0$). Hence (\ref{calculdirectF}) reduces to
\begin{equation}\label{calculF2}
  \textbf{F}^\gop
 = \frac{1}{2} \left(\partial_{\cels_1}\textbf{A}{^\gop}_{\cels_2}
 - \partial_{\cels_2}\textbf{A}{^\gop}_{\cels_1} + [\textbf{A}{^\gop}_{\cels_1}, \textbf{A}{^\gop}_{\cels_1}]\right)e^{\cels_1\cels_2}
 + \partial_{\cell}\textbf{A}{^\gop}_{\cels}e^{\cell\cels}
\end{equation}
A consequence of (\ref{definitionF}) and (\ref{introductione&A}) is
\begin{equation}\label{FegalAdgPhi}
 \textbf{F}^\gop= \hbox{Ad}_g
 \left(\dR\varphi^\gop + \frac{1}{2}[\varphi^\gop\wedge \varphi^\gop]\right)  = \hbox{Ad}_g\Phi^\gop
\end{equation}
Thus by letting $\textbf{F}{^\gop}_{\gop\gop}:= \hbox{Ad}_g\otimes \hbox{Ad}^*_g\otimes \hbox{Ad}^*_g\,\Phi{^\gop}_{\gop\gop}$ (see (\ref{18gravity})) we have $\textbf{F}{^\gop} = \frac{1}{2}\textbf{F}{^\gop}_{\celp\celp}e^{\celp\celp}$ by Lemma \ref{lemme2point3}. However (\ref{gravEL1}) translates as
\begin{equation}\label{translation114}
\textbf{F}{^\gop}_{\gos\gol} = \textbf{F}{^\gop}_{\gol\gos} = \textbf{F}{^\gop}_{\gol\gol} = 0
\end{equation}
and thus
$\textbf{F}{^\gop} = \frac{1}{2}\textbf{F}{^\gop}_{\cels\cels}e^{\cels\cels}$.
By comparing with (\ref{calculF2}) we deduce that $\partial_\gol\textbf{A}{^\gop}_\gos = 0$. As a consequence:
\begin{equation}\label{conclusionsurAetF}
 \hbox{The coefficients }\textbf{A}{^\gop}_\gos\hbox{ and }\textbf{F}{^\gop}_{\gos\gos}\hbox{ are constant on each fiber }\textsf{f}
\end{equation}
The next step is to look at the image
by $\hbox{Ad}_g$ of both sides of the relation
$\dR^\varphi\pi{_{\gop}} = \Psi{_{\gop}} - \frac{1}{2}\Psi\,
\varphi^{(N-1)}_{\gop}$
in
(\ref{gravELcompact})
(recall that $\Psi{_{\gop}} :=
\Psi{_{\gop}}^{\celp}
\, \varphi^{(N-1)}_{\celp}$), i.e.
to compute both sides of
\begin{equation}\label{equationatraduire}
\hbox{Ad}_g^*\left( \dR^\varphi\pi_{\gop}\right) =
\hbox{Ad}_g^*\left( \Psi{_{\gop}} - \frac{1}{2}\Psi\,
\varphi^{(N-1)}_{\gop}\right)
\end{equation}

\subsubsection{Translation of Equation (\ref{equationatraduire})}

We recall that $e^{\gop}:=
 \hbox{Ad}_g\varphi^{\gop}$.
We also introduce
\begin{equation}\label{introductionp=Adgvarpi}
 p_{\gop}:= \hbox{Ad}_g^*\pi_{\gop} .
\end{equation}
and we set $p{_\gop}^{\gop\gop}
 := \hbox{Ad}_g^*\otimes\hbox{Ad}_g\otimes\hbox{Ad}_g\ \pi{_\gop}^{\gop\gop}$. Since $\gop$ is
unimodular we have the decomposition
$p_\gop = \frac{1}{2}
 p{_\gop}^{\celp\celp}\ e^{(N-2)}_{\celp\celp}$, by (\ref{eetf56}) and (\ref{58}).

Let us define
\begin{equation}\label{definitionQ}
\left\{\begin{array}{ccccc}
         Q{_{\gop\gop}}^{\gop\gop} & := & & &
         \textbf{F}{^\celp}_{\gop\gop}\
         p{_\celp}^{\gop\gop} \\
         Q{_\gop}^\gop & := &
         Q{_{\gop\celp}}^{\gop\celp} &
         =  &
         \textbf{F}{^\celp}_{\gop\celp}\
         p{_\celp}^{\gop\celp} \\
         Q & := & Q{_\celp}^\celp & = &
         \textbf{F}{^\celp}_{\celp\celp} \
         p{_\celp}^{\celp\celp}
       \end{array}\right.
\end{equation}
It follows from these definitions that
\[
\begin{array}{ccl}
 Q{_{\gop\gop}}^{\gop\gop} & = &
 \left(\hbox{Ad}_g\otimes\hbox{Ad}_g^*\otimes\hbox{Ad}_g^*\ \Phi{^\celp}_{\gop\gop}\right)
 \left(\hbox{Ad}_g^*\otimes\hbox{Ad}_g\otimes\hbox{Ad}_g\ \pi{_\celp}^{\gop\gop}\right) \\
 & = & \hbox{Ad}_g^*\otimes\hbox{Ad}_g^*\otimes\hbox{Ad}_g\otimes\hbox{Ad}_g\
 (\Phi{^\celp}_{\gop\gop}\pi{_\celp}^{\gop\gop})\\
 & = & \hbox{Ad}_g^*\otimes\hbox{Ad}_g^*\otimes\hbox{Ad}_g\otimes\hbox{Ad}_g\
 \Psi{_{\gop\gop}}^{\gop\gop}.
\end{array}
\]
(see (\ref{definitionfamillePsi}))
and hence
\begin{equation}\label{Qicci}
 Q{_\gop}^\gop =
 \hbox{Ad}_g^*\otimes\hbox{Ad}_g\ \Psi{_\gop}^{\gop}
 \quad \hbox{ and }\quad Q = \Psi
\end{equation}
Thus
\[
 Q_\gop := Q{_\gop}^\celp\ e^{(N-1)}_\celp
=  \left(\hbox{Ad}_g^*\otimes\hbox{Ad}_g\ \Psi{_\gop}^{\celp}\right)
\left(\hbox{Ad}_g^*\varphi^{(N-1)}_\celp\right)
 = \hbox{Ad}_g^*\left(\Psi{_\gop}^{\celp}
 \varphi^{(N-1)}_\celp\right)
 = \hbox{Ad}_g^*\Psi_{\gop}
\]
Hence
\begin{equation}\label{Qeinstein}
 \hbox{Ad}_g^*\left( \Psi{_\gop} - \frac{1}{2}\Psi\,
\varphi^{(N-1)}_{\gop}\right)
= Q{_\gop}- \frac{1}{2}Qe^{(N-1)}_{\gop}
\end{equation}
Thus Equation (\ref{equationatraduire}) is equivalent to
$\hbox{Ad}_g^*\left( \dR^\varphi\pi_{\gop}\right) =
Q{_\gop}- \frac{1}{2}Qe^{(N-1)}_{\gop}$.
We can conclude by using (\ref{domegap=Adgdvarphivarpi})
which says that $\hbox{Ad}_g^*\left( \dR^\varphi\pi_{\gop}\right) = \dR^\textbf{A}p_\gop$ that Equation (\ref{equationatraduire})
is equivalent to the \emph{fundamental equation}
\begin{equation}\label{equationfondamentaleGrave}
\begin{array}{|c|}
 \hline
 \dR^\textbf{A}p_\gop =  Q_\gop -
 \frac{1}{2}Q\,e_\gop^{(N-1)} \\
\hline
\end{array}
\end{equation}

\subsection{The dynamical equation (\ref{equationfondamentaleGrave}) in a local trivialization}

\subsubsection{Remarks on the dual fields and computation of the right hand side}

First the facts that $p{_\gop}^{\gop\gop}
 = (\hbox{Ad}_g^*\otimes\hbox{Ad}_g\otimes\hbox{Ad}_g)\pi{_\gop}^{\gop\gop}$ and that the adjoint (respectively
 coadjoint) action of $\widehat{\goL}$ on $\gop$ (respectively
 $\gop^*$) leaves the decomposition
 $\gop = \gol \oplus \gos$ (respectively
 $\gop^* = \gol^* \oplus \gos^*$)
invariant imply in particular
 that $p{_\gop}^{\gos\gos}
 = (\hbox{Ad}_g^*\otimes\hbox{Ad}_g\otimes\hbox{Ad}_g)\pi{_\gop}^{\gos\gos}$. Hence since
 $\pi{_\gop}^{\gos\gos} = \kappa{_\gop}^{\gos\gos}$
is $\widehat{\goL}$-invariant,
 we deduce that
\begin{equation}\label{constraintcoordinates1}
p{_\gop}^{\gos\gos} = \kappa{_\gop}^{\gos\gos}.
\end{equation}
Second (\ref{gravELcompact}) and (\ref{translation114}) imply
$\textbf{F}{^\gop}_{\gos\gol}
= \textbf{F}{^\gop}_{\gol\gos} = 0$.
Hence $Q{_{\gos\gol}}^{\gop\gop}
= Q{_{\gol\gos}}^{\gop\gop}
= Q{_{\gol\gol}}^{\gop\gop} =0$ and thus
\begin{equation}\label{reductionQourbure}
 Q{_{\gop\gop}}^{\gop\gop}
 = Q{_{\gos\gos}}^{\gop\gop}
 + Q{_{\gos\gol}}^{\gop\gop}
 + Q{_{\gol\gos}}^{\gop\gop}
 + Q{_{\gol\gol}}^{\gop\gop}
 = Q{_{\gos\gos}}^{\gop\gop}.
\end{equation}
This implies also that
$Q{_\gol}^\gop = Q{_{\gol\celp}}^{\gop\celp}
= Q{_{\gol\cels}}^{\gop\cels}
+ Q{_{\gol\cell}}^{\gop\cell} = 0$
and
$Q{_\gos}^\gop = Q{_{\gos\celp}}^{\gop\celp}
= Q{_{\gos\cels}}^{\gop\cels}
+ Q{_{\gos\cell}}^{\gop\cell} =
Q{_{\gos\cels}}^{\gop\cels}$.
To summarize
\[
 \left(\begin{array}{cc}
    Q{_\gol}^\gol & Q{_\gos}^\gol \\
    Q{_\gol}^\gos & Q{_\gos}^\gos
\end{array}\right)
=
 \left(\begin{array}{cc}
    0 & \textbf{F}{^\celp}_{\gos\cels}\,
    p{_\celp}^{\gol\cels} \\
    0 & \textbf{F}{^\celp}_{\gos\cels}\,
    \kappa{_\celp}^{\gos\cels}
\end{array}\right)
\quad \hbox{ and } \quad
Q = \textbf{F}{^\celp}_{\cels\cels}\,
    \kappa{_\celp}^{\cels\cels}
\]
Thus by setting
$\delta{_\gol}^\gol:= \delta^i_j\textbf{t}^j\otimes \textbf{t}_i$ and
$\delta{_\gos}^\gos:= \delta^a_b\textbf{t}^b\otimes \textbf{t}_a$,
\begin{equation}\label{detailQ}
\begin{array}{|c|}
 \hline
\displaystyle Q_\gop - \frac{1}{2}Q\,e_\gop^{(N-1)} =
\left(\textbf{F}{^\celp}_{\gos\cels_1}\,
    p{_\celp}^{\cell\cels_1}
- \frac{1}{2} Q\delta{_\gol}^\cell \right) e^{(N-1)}_\cell
 + \left( \textbf{F}{^\celp}_{\gos\cels_1}\,
    \kappa{_\celp}^{\cels\cels_1}
-\frac{1}{2} Q\delta{_\gos}^\cels\right) e^{(N-1)}_\cels \\
\hline
\end{array}
\end{equation}

\subsubsection{Introducing a solder form and a connection form}
Recall that by (\ref{introductione&Abis})
\begin{equation}\label{relationeA}
 e^{\gop} = \dR g\,g^{-1} + \textbf{A}^\gop
\end{equation}
and, by decomposing $\textbf{A}^\gop = \textbf{A}^\gos + \textbf{A}^\gol$,
that $e^\gos = \textbf{A}^\gos$ and
$e^\gol = \dR g\,g^{-1} + \textbf{A}^\gol$. 
For later interpretation, we give
special names to these two forms:
\begin{equation}\label{141bis}
\begin{array}{|c|}
 \hline
 \theta^\gos:= \textbf{A}^\gos = e^\gos
 \\
 \hline
\end{array}
\quad\hbox{ and } \quad
\begin{array}{|c|}
 \hline
 \omega^\gol:= \textbf{A}^\gol
 = e^\gol - \hbox{d}g\,g^{-1},
 \\
 \hline
\end{array}
\end{equation}
so that
\begin{equation}\label{decompositionAomegatheta}
\textbf{A}^\gop = \theta^\gos + \omega^\gol  .
\end{equation}
We will see later that $\omega^\gol$ plays the role of a connection 1-form
and $\theta^\gos$ the role of a soldering 1-form (meaning that the components
$\theta^a = e^a$ forms a coframe over the space-time).
We also define
\begin{equation}\label{definitionOmega}
\Omega^\gol:= \dR \omega^\gol + \frac{1}{2}[\omega^\gol\wedge\omega^\gol] \in \gol\otimes \Omega^2(\mathcal{F})
\end{equation}
which can be interpreted as a curvature form,
and
\begin{equation}\label{definitionTheta}
 \Theta^\gos:= \dR^\omega e^\gos = \dR^\omega\theta^\gos = \dR\theta^\gos
 + [\omega^\gol\wedge \theta^\gos]
 \in \gos\otimes \Omega^2(\mathcal{F})
\end{equation}
which can be interpreted as a
torsion form. It follows from (\ref{definitionF}) that
\begin{equation}\label{FenfonctionThetaetOmega}
 \mathbf{F}^\gop:= \hbox{d}\mathbf{A}^\gop + \frac{1}{2}[\mathbf{A}^\gop\wedge \mathbf{A}^\gop]
 = \Theta^\gos + \Omega^\gol + \frac{1}{2}[\theta^\gos\wedge \theta^\gos]
\end{equation}

\subsubsection{Computation of the left hand side of (\ref{equationfondamentaleGrave})}

Since $\gos$, $\gol$ and $\kappa{_\gop}^{\gos\gos}$ are not stable by $\hbox{Ad}_{\widehat{\goP}}$ but are stable by $\hbox{Ad}_\goG$ it will be convenient to define
\begin{equation}\label{domegap}
 \dR^\omega p_\gop:= \dR p_\gop + \hbox{ad}_\omega^* \wedge p_\gop,
\end{equation}
to split
\begin{equation}\label{dApEgaldomegapplusadthetap}
 \dR^\textbf{A} p_\gop = \dR^\omega p_\gop + \hbox{ad}_\theta^* \wedge p_\gop.
\end{equation}
and to compute separately $\dR^\omega p_\gop$ and $\hbox{ad}_\theta^* \wedge p_\gop$.
Using (\ref{twistedLeibniz}) with $\dR^\omega$
and the decomposition
$p_\gop:= \frac{1}{2}
p{_\gop}^{\celp\celp}\
e_{\celp\celp}^{(N-2)}$,
we get $\dR^\omega p_\gop = \frac{1}{2}
 \dR^\omega p{_\gop}^{\celp\celp}
 \wedge e_{\celp\celp}^{(N-2)}
 + \frac{1}{2}p{_\gop}^{\celp\celp}\
\dR^\omega e_{\celp\celp}^{(N-2)}$.
Summarizing with (\ref{dApEgaldomegapplusadthetap}) we see
that we need to compute each term in the r.h.s. of
\begin{equation}\label{dApEgaltroistermes}
\begin{array}{|c|}
 \hline
\dR^\textbf{A} p_\gop = \frac{1}{2}
 \dR^\omega p{_\gop}^{\celp\celp}
 \wedge e_{\celp\celp}^{(N-2)}
 + \frac{1}{2}p{_\gop}^{\celp\celp}\
\dR^\omega e_{\celp\celp}^{(N-2)}
 + \hbox{ad}_\theta^* \wedge p_\gop  \\
 \hline
\end{array}
\end{equation}

\noindent
\textbf{Computation of $\frac{1}{2}
 \dR^\omega p{_\gop}^{\celp\celp}\wedge e_{\celp\celp}^{(N-2)}$} --- 
Let us introduce the coefficients $\partial_\gop p{_\gop}^{\gop\gop}$ and $\partial_\gop^\omega p{_\gop}^{\gop\gop}$ such that $\dR p{_\gop}^{\gop\gop} = \left(\partial_\celp p{_\gop}^{\gop\gop} \right) e^\celp$ and $\dR^\omega p{_\gop}^{\gop\gop} = \left(\partial_\celp^\omega p{_\gop}^{\gop\gop} \right) e^\celp$.
 They are related by
 \begin{equation}\label{relation96}
  \partial_\gop^\omega p{_\gop}^{\gop\gop}
  = \partial_\gop p{_\gop}^{\gop\gop}
  + \left(\hbox{ad}_{\omega_\gos}^*\otimes 1 \otimes 1
  + 1 \otimes \hbox{ad}_{\omega_\gos} \otimes 1 + 1 \otimes 1 \otimes\hbox{ad}_{\omega_\gos}\right)
  p{_\gop}^{\gop\gop}
 \end{equation}
which, through the decomposition (see (\ref{cgggref}) for the notation)
\begin{equation}\label{settingforcstructure}
\textbf{c}{^\gop}_{\gop\gop}:=
\textbf{c}^I_{JK} \textbf{t}_I \otimes\textbf{t}^J\otimes\textbf{t}^K = \textbf{c}{^\gol}_{\gol\gol} + \textbf{c}{^\gos}_{\gol\gos} + \textbf{c}{^\gos}_{\gos\gol} + \textbf{c}{^\gol}_{\gos\gos}
\end{equation}
means:
\[
 \partial^\omega_\gop p{_{\gop_0}}^{\gop_1\gop_2}
 = \partial_\gop p{_{\gop_0}}^{\gop_1\gop_2}
- \textbf{c}{^\celp}_{\cell\gop_0}\,\omega{^\cell}_\gos\, p{_\celp}^{\gop_1\gop_2}
+ \textbf{c}{^{\gop_1}}_{\cell\celp}\,\omega{^\cell}_\gos\, p{_{\gop_0}}^{\celp\gop_2}
+ \textbf{c}{^{\gop_2}}_{\cell\celp}\,\omega{^\cell}_\gos\, p{_{\gop_0}}^{\gop_1\celp}
\]
Then
\begin{equation}\label{relation97}
  \frac{1}{2}
 \dR^\omega p{_\gop}^{\celp\celp}\wedge e_{\celp\celp}^{(N-2)}
= \frac{1}{2}\left(\partial_\celp^\omega p{_\gop}^{\celp_1\celp_2}  \right)e^\celp \wedge e_{\celp_1\celp_2}^{(N-2)}
 = \partial_{\celp}^\omega p{_\gop}^{\celp_1\celp}\ e_{\celp_1}^{(N-1)}
 \end{equation}
We have $\partial_{\celp}^\omega p{_\gop}^{\gop\celp} = \partial_{\cels}^\omega p{_\gop}^{\gop\cels} + \partial_{\cell}^\omega p{_\gop}^{\gop\cell}$ but, since $\omega^\gol = \omega{^\gol}_\cels e^\cels$
(i.e. $\omega{^\gol}_\gol = 0$), actually $\partial_{\cell}^\omega p{_\gop}^{\gop\cell} = \partial_{\cell} p{_\gop}^{\gop\cell}$. Hence
$\partial_{\celp}^\omega p{_\gop}^{\gop\celp} = \partial_{\cels}^\omega p{_\gop}^{\gop\cels} + \partial_{\cell} p{_\gop}^{\gop\cell}$,
i.e.
\[
 \left\{ \begin{array}{ccl}
\partial_{\celp}^\omega p{_\gop}^{\gos\celp} & = & \partial_{\cels}^\omega p{_\gop}^{\gos\cels} + \partial_{\cell} p{_\gop}^{\gos\cell} \\
\partial_{\celp}^\omega p{_\gop}^{\gol\celp} & = & \partial_{\cels}^\omega p{_\gop}^{\gol\cels} + \partial_{\cell} p{_\gop}^{\gol\cell}
 \end{array}\right.
\]
Moreover $p{_\gop}^{\gos\gos} = \kappa{_\gop}^{\gos\gos}$ as observed in (\ref{constraintcoordinates1}). Thus since $\kappa{_\gop}^{\gos\gos}$
 is constant and $\hbox{ad}_\gol$-invariant,
 $\partial_{\gos}^\omega p{_\gop}^{\gos\gos} =\partial_{\gos}^\omega \kappa{_\gop}^{\gos\gos} = 0$. Hence $\partial_{\celp}^\omega p{_\gop}^{\gos\celp} = 0 + \partial_{\cell} p{_\gop}^{\gos\cell}$. In conclusion (\ref{relation97}) gives us
\begin{equation}\label{domegapcd}\begin{array}{|c|}
 \hline
  \frac{1}{2}
 \dR^\omega p{_\gop}^{\celp\celp}\wedge e_{\celp\celp}^{(N-2)}
= \partial_{\cell} p{_\gop}^{\cels\cell}\,e^{(N-1)}_\cels
+ \left(\partial_{\cels}^\omega p{_\gop}^{\cell\cels} + \partial_{\cell_1} p{_\gop}^{\cell\cell_1}\right)e^{(N-1)}_\cell  \\
 \hline
\end{array}
\end{equation}

\noindent
\textbf{Computation of $\frac{1}{2}p{_\gop}^{\celp\celp}\ \dR^\omega e_{\celp\celp}^{(N-2)}$} --- 
By applying (\ref{nablaLeibniz}) we get $\dR^\omega e^{(N-2)}_{\gop\gop} = \dR^\omega e^\celp \wedge \omega^{(N-2)}_{\gop\gop\celp}$.
We thus need to compute  $\dR^\omega e^\gop$. For that purpose we split
$e^\gop = e^\gos + e^\gol$ and we use (\ref{dAeagalFee}),
i.e. $\dR^\omega e^\gol = \Omega^\gol
 + \frac{1}{2}[e^\gol\wedge e^\gol]$. Hence
\begin{equation}\label{117bis}
 \dR^\omega e^\gop = \dR^\omega e^\gos
 + \dR^\omega e^\gol
 = \Theta^\gos + \Omega^\gol
 + \frac{1}{2}[e^\gol\wedge e^\gol]
\end{equation}
or by using the notation (\ref{settingforcstructure}),
$\dR^\omega e^\gop = \frac{1}{2}\Theta{^\gos}_{\cels\cels} e^{\cels\cels} + \frac{1}{2}\Omega{^\gol}_{\cels\cels}e^{\cels\cels}
+ \frac{1}{2}\textbf{c}{^\gol}_{\cell\cell}e^{\cell\cell}$.
Thus
\[
 \dR^\omega e_{\gop\gop}^{(N-2)}
 = \frac{1}{2}\Theta{^\cels}_{\cels_1\cels_2} e^{\cels_1\cels_2}\wedge e^{(N-2)}_{\gop\gop\cels} + \frac{1}{2}\Omega{^\cell}_{\cels_1\cels_2}e^{\cels_1\cels_2}\wedge e^{(N-2)}_{\gop\gop\cell}
+ \frac{1}{2}\textbf{c}{^\cell}_{\cell_1\cell_2}e^{\cell_1\cell_2}
\wedge e^{(N-2)}_{\gop\gop\cell}
\]
We compute the r.h.s. by using (\ref{fifi}) (e). The first term is (recall (\ref{conclusionsurAetF}) and (\ref{141bis}))
\begin{equation}\label{reasoningp}
 \frac{1}{2}\Theta{^\cels}_{\cels_1\cels_2} e^{\cels_1\cels_2}\wedge e^{(N-2)}_{\gop_1\gop_2\cels}
 = \Theta{^\cels}_{\gop_1\gop_2} e^{(N-1)}_{\cels}
 + \Theta{^\cels}_{\gop_2\cels} e^{(N-1)}_{\gop_1}
 +\Theta{^\cels}_{\cels\gop_1} e^{(N-1)}_{\gop_2}
\end{equation}
(we use here the fact that $\Theta{^\gos}_{\gop\gop} = \Theta{^\gos}_{\gos\gos}$, because of (\ref{conclusionsurAetF}) and (\ref{141bis}))
the second term is
$\frac{1}{2}\Omega{^\cell}_{\cels_1\cels_2}e^{\cels_1\cels_2}\wedge e^{(N-2)}_{\gop_1\gop_2\cell}
= \Omega{^\cell}_{\gop_1\gop_2} e^{(N-1)}_{\cell} +  \Omega{^\cell}_{\gop_2\cell} e^{(N-1)}_{\gop_1}
 +\Omega{^\cell}_{\cell\gop_1} e^{(N-1)}_{\gop_2}$, which is equal to $\Omega{^\cell}_{\gop_1\gop_2} e^{(N-1)}_{\cell}$ because $\Omega{^\gol}_{\gop_2\gol} = \Omega{^\gol}_{\gol\gop_1}  = 0$.
The last term is $\frac{1}{2}\textbf{c}{^\cell}_{\cell_1\cell_2}e^{\cell_1\cell_2}\wedge e^{(N-2)}_{\gop_1\gop_2\cell}
= \textbf{c}{^\cell}_{\gop_1\gop_2} e^{(N-1)}_{\cell} + \textbf{c}{^\cell}_{\gop_2\cell} e^{(N-1)}_{\gop_1}
 + \textbf{c}{^\cell}_{\cell\gop_1} e^{(N-1)}_{\gop_2}$ which simplifies to $\textbf{c}{^\cell}_{\gop_1\gop_2} e^{(N-1)}_{\cell}$, because $\textbf{c}{^\cell}_{\gop\cell} = \textbf{c}{^\cell}_{\cell\gop} = 0$ since $\gol$ is unimodular. In conclusion by setting
\begin{equation}\label{Thetaast}
 \Theta{^\ast}_{\gos\ast}:= \Theta{^\cels}_{\gos\cels}
\end{equation}
we get
\[
 \dR^\omega e_{\gop_1\gop_2}^{(N-2)}
= \Theta{^\cels}_{\gop_1\gop_2} e^{(N-1)}_{\cels}
 + \Theta{^\ast}_{\gop_2\ast} e^{(N-1)}_{\gop_1}
 - \Theta{^\ast}_{\gop_1\ast} e^{(N-1)}_{\gop_2}
 + \Omega{^\cell}_{\gop_1\gop_2} e^{(N-1)}_{\cell}
 + \textbf{c}{^\cell}_{\gop_1\gop_2} e^{(N-1)}_{\cell}
\]
This implies by using the fact that $\Theta{^\gos}_{\gop\gop}
= \Theta{^\gos}_{\gos\gos}$,
$\Omega{^\gol}_{\gop\gop} = \Omega{^\gol}_{\gos\gos}$ and
$\textbf{c}{^\cell}_{\gop\gop} = \textbf{c}{^\cell}_{\gol\gol}$, that
\[
\begin{array}{ll}
\frac{1}{2}p{_\gop}^{\celp\celp}\ \dR^\omega e_{\celp\celp}^{(N-2)} 
= & \frac{1}{2}p{_\gop}^{\cels_1\cels_2}\Theta{^\cels}_{\cels_1\cels_2} e^{(N-1)}_{\cels}
 - p{_\gop}^{\cels_1\celp_2}\Theta{^\ast}_{\cels_1\ast} e^{(N-1)}_{\celp_2} \\
& \ +\ \frac{1}{2}\left( p{_\gop}^{\cels_1\cels_2}\Omega{^\cell}_{\cels_1\cels_2}
 + p{_\gop}^{\cell_1\cell_2}\textbf{c}{^\cell}_{\cell_1\cell_2}\right) e^{(N-1)}_{\cell}
 \end{array}
\]
By splitting $p{_\gop}^{\cels_1\celp_2}
\Theta{^\ast}_{\cels_1\ast}\,e^{(N-1)}_{\celp_2} = p{_\gop}^{\cels_1\cels_2}
\Theta{^\ast}_{\cels_1\ast}\,e^{(N-1)}_{\cels_2} + p{_{\gop_2}}^{\cels_1\cell_2}
\Theta{^\ast}_{\cels_1\ast}\,e^{(N-1)}_{\cell_2}$ and by grouping together we obtain
\begin{equation}\label{prepdomegaeeN-3}
\begin{array}{ccc}
 \frac{1}{2}p{_\gop}^{\celp\celp}\ \dR^\omega e_{\celp\celp}^{(N-2)}
 & = & \left(\frac{1}{2}p{_\gop}^{\cels_1\cels_2}
\Theta{^\cels}_{\cels_1\cels_2}
- p{_\gop}^{\cels_1\cels}
\Theta{^\ast}_{\cels_1\ast}
\right)e^{(N-1)}_\cels
\\
&  & +\
\left(\frac{1}{2}p{_\gop}^{\cels_1\cels_2}
\Omega{^\cell}_{\cels_1\cels_2}
- p{_\gop}^{\cels_1\cell}
\Theta{^\ast}_{\cels_1\ast}
+ \frac{1}{2}p{_\gop}^{\cell_1\cell_2}
\textbf{c}{^\cell}_{\cell_1\cell_2}
\right)e^{(N-1)}_\cell
\end{array}
\end{equation}

\subsubsection{Introducing the generalized Cartan and Einstein tensors}
For further use we introduce the notations
\[
 \mathring{\Theta}{^\gos}_{\gos_1\gos_2}:=
 \Theta{^\gos}_{\gos_1\gos_2} + \delta^\gos_{\gos_1}\Theta{^\ast}_{\gos_2\ast} - \delta^\gos_{\gos_2}\Theta{^\ast}_{\gos_1\ast}
 = \Theta{^\gos}_{\gos_1\gos_2} - \delta^\gos_{\gos_1}\Theta{^\cels}_{\cels\gos_2} - \delta^\gos_{\gos_2}\Theta{^\cels}_{\gos_1\cels}
\]
and
\[
 \mathring{\Omega}{^\gog}{_{\gos_1\gos_2\gos_3}}^\gos := 
 \Omega{^\gog}_{\gos_1\gos_2}\delta^\gos_{\gos_3} + \Omega{^\gog}_{\gos_2\gos_3}\delta^\gos_{\gos_1} + \Omega{^\gog}_{\gos_3\gos_1}\delta^\gos_{\gos_2}
\]
We observe (through a computation similar to (\ref{reasoningp})) that the first coefficient in the right hand side of (\ref{prepdomegaeeN-3}) can be written
\[
 \frac{1}{2}p{_\gop}^{\cels_1\cels_2}
\Theta{^\gos}_{\cels_1\cels_2}
- p{_\gop}^{\cels_1\gos}
\Theta{^\ast}_{\cels_1\ast}
= \frac{1}{2}p{_\gop}^{\cels_1\cels_2}\mathring{\Theta}{^\gos}_{\cels_1\cels_2}
\]
so that (\ref{prepdomegaeeN-3}) reads
\begin{equation}\label{pdomegaeeN-3}\begin{array}{|c|}
 \hline
  \frac{1}{2}p{_\gop}^{\celp\celp}\ \dR^\omega e_{\celp\celp}^{(N-2)}
 = \frac{1}{2}p{_\gop}^{\cels_1\cels_2}\mathring{\Theta}{^\cels}_{\cels_1\cels_2}e^{(N-1)}_\cels
+ \left(\frac{1}{2}p{_\gop}^{\cels\cels}\Omega{^\cell}_{\cels\cels}
- p{_\gop}^{\cels\cell}\Theta{^\ast}_{\cels\ast}
+ \frac{1}{2}p{_\gop}^{\cell_1\cell_2}\textbf{c}{^\cell}_{\cell_1\cell_2}\right) e_{\cell}^{(N-1)}  \\
 \hline
\end{array}
\end{equation}
Actually $\mathring{\Theta}{^\gos}_{\gos\gos}$ and $\mathring{\Omega}{^\gog}{_{\gos\gos\gos}}^\gos$ are defined implicitely by
\begin{equation}\label{defimathringTheta}
 \mathring{\Theta}{^\cels}_{\gos\gos}\ e^{(N-1)}_\cels:= \Theta{^\cels}\wedge e^{(N-3)}_{\gos\gos\cels}
\end{equation}
and
\begin{equation}\label{defimathringOmega}
 \mathring{\Omega}{^\gog}{_{\gos_1\gos_2\gos}}^\cels\  e^{(N-1)}_\cels:= \Omega{^\gog}\wedge e^{(N-3)}_{\gos_1\gos_2\gos}
\end{equation}
We note that these relations imply that $\mathring{\Theta}{^\gos}_{\gos\gos} = \mathring{\Omega}{^\gog}{_{\gos\gos\gos}}^\gos = 0$ whenever $N\leq 2$.

We further define 
\begin{equation}\label{defiCtenseurCartan}
 \tilde{\mathbf{C}}{_\gog}^\gos:= - \frac{1}{2} \kappa{_\gog}^{\cels_1\cels_2} \mathring{\Theta}{^\gos}_{\cels_1\cels_2}
\end{equation}
and 
\begin{equation}\label{defiEtenseurEinstein}
 \tilde{\mathbf{E}}{_\gos}^\gos:=  -\frac{1}{2} \kappa{_\celg}^{\cels_1\cels_2} \mathring{\Omega}{^\celg}{_{\cels_1\cels_2\gos}}^\gos 
\end{equation}
which can also be defined implicitely by
\begin{equation}\label{defiCtenseurCartanimplicit}
 \tilde{\mathbf{C}}{_\gog}^\cels\, \ e^{(N-1)}_\cels:= - \frac{1}{2} \kappa{_\gog}^{\cels_1\cels_2}\Theta{^\cels}\wedge e^{(N-3)}_{\cels_1\cels_2\cels}
\end{equation}
and
\begin{equation}\label{defiEtenseurEinsteinimplicit}
 \tilde{\mathbf{E}}{_\gos}^\cels\, \ e^{(N-1)}_\cels:= - \frac{1}{2} \kappa{_\celg}^{\cels_1\cels_2}\Omega{^\celg}\,\wedge e^{(N-3)}_{\cels_1\cels_2\gos}
\end{equation}
We shall see that, in standard situations, $\tilde{\mathbf{C}}{_\gog}^\gos$ (or equivalentely $\mathring{\Theta}{^\gos}_{\gos\gos}$) corresponds to the \emph{Cartan tensor} and $\tilde{\mathbf{E}}{_\gos}^\gos$ to the \emph{Einstein tensor}. Indeed
\begin{equation}\label{EssestEinstein}
  \begin{array}{ccl}
  \tilde{\mathbf{E}}{_\gos}^\gos & = & - \frac{1}{2} \left(\Omega{^\celg}_{\cels_1\cels_2}\delta^\gos_{\gos} + \Omega{^\celg}_{\cels_2\gos}\delta^\gos_{\cels_1} + \Omega{^\celg}_{\gos\cels_1}\delta^\gos_{\cels_2}\right) \kappa{_\celg}^{\cels_1\cels_2} \\
  & = & - \frac{1}{2}(\Omega{^\celg}_{\cels_1\cels_2}\kappa{_\celg}^{\cels_1\cels_2})\delta^\gos_{\gos}
  - \frac{1}{2}\Omega{^\celg}_{\cels_2\gos}\kappa{_\celg}^{\gos\cels_2}
  - \frac{1}{2}\Omega{^\celg}_{\gos\cels_1}\kappa{_\celg}^{\cels_1\gos} \\
  & = & \Omega{^\celg}_{\cels_1\gos}\kappa{_\celg}^{\cels_1\gos}- \frac{1}{2}(\Omega{^\celg}_{\cels_1\cels_2}\kappa{_\celg}^{\cels_1\cels_2})\delta^\gos_{\gos} 
 \end{array}
\end{equation}
i.e., by denoting $\tilde{\mathbf{R}}{_\gos}^\gos:= \Omega{^\celg}_{\cels_1\gos}\kappa{_\celg}^{\cels_1\gos}$ and $\tilde{\mathbf{R}}: = \tilde{\mathbf{R}}{_\cels}^\cels$ (generalized versions of, respectively, the Ricci tensor and the scalar curvature), we have
$\tilde{\mathbf{E}}{_\gos}^\gos= \tilde{\mathbf{R}}{_\gos}^\gos - \frac{1}{2}\tilde{\mathbf{R}}\delta^\gos_{\gos}$.\\

\noindent
\textbf{Computation of $\hbox{ad}_\theta^* \wedge p_\gol$} ---
Our last task consists in computing
$\hbox{ad}_\theta^*\wedge p_\gop$, i.e., since
$\theta = e^\gos$,
$\hbox{ad}_\theta^*\wedge p_\gop
= - \textbf{c}{^{\celp}}_{\cels\gop}e^{\cels}\wedge p_{\celp}$.
Since $p_\gop =
\frac{1}{2}p{_\gop}^{\cels\cels}
 \ e_{\cels\cels}^{(N-2)}
+ p{_\gop}^{\cels\cell}e_{\cels\cell}^{(N-2)}
+ \frac{1}{2}p{_\gop}^{\cell\cell}e_{\cell\cell}^{(N-2)}$, this quantity is the sum of three terms.
The first one is
$- \frac{1}{2}\textbf{c}{^\celp}_{\cels\gop}
 \ p{_\celp}^{\cels_1\cels_2}\ e^{\cels} \wedge
 \ e_{\cels_1\cels_2}^{(N-2)}$.
Since by (\ref{constraintcoordinates1}), $p{_\gop}^{\gos\gos} = \kappa{_\gop}^{\gos\gos}$, it is equal to
\[
 - \frac{1}{2}\textbf{c}{^\celp}_{\cels\gop}
 \ \kappa{_\celp}^{\cels_1\cels_2}\ e^{\cels} \wedge
 \ e_{\cels_1\cels_2}^{(N-2)}
  = \textbf{c}{^\celp}_{\cels_1\gop}\ \kappa{_{\celp}}^{\cels_1\cels_2}
  \ e_{\cels_2}^{(N-1)}
\]
The second term is
\[
 - \textbf{c}{^\celp}_{\cels\gop}
 \ p{_\celp}^{\cels_1\cell_2}\ e^{\cels} \wedge
 \ e_{\cels_1\cell_2}^{(N-2)}
 = \textbf{c}{^\celp}_{\cels\gop}\ p{_\celp}^{\cels\cell}
\ e_{\cell}^{(N-1)}
\]
Lastly since $\theta^\gos\wedge e_{\gol\gol}^{(N-2)} = 0$, the last term
in the r.h.s. vanishes. Hence
we get
\begin{equation}\label{adthetap}\begin{array}{|c|}
 \hline
 \hbox{ad}_\theta^*\wedge p_\gop
 = \textbf{c}{^\celp}_{\cels_1\gop}\kappa{_\celp}^{\cels_1\cels}   e_\cels^{(N-1)} +
 \textbf{c}{^\celp}_{\cels\gop}\ p{_\celp}^{\cels\cell}
\ e_{\cell}^{(N-1)} \\
 \hline
\end{array}
\end{equation}

\subsubsection{Conclusion}\label{theend}
We go back to (\ref{dApEgaltroistermes}), by collecting
(\ref{domegapcd}), (\ref{pdomegaeeN-3}) and
(\ref{adthetap}):
\begin{equation}\label{dApfinal}
  \begin{array}{|ccl|}
  \hline
  \dR^\textbf{A} p{_\gop} & = &
   \left(\frac{1}{2}p{_\gop}^{\cels_1\cels_2}
  \Omega{^\cell}_{\cels_1\cels_2}
+ (\partial_{\cels}^\omega
+ \Theta{^{\ast}}_{\cels\ast})
  p{_\gop}^{\cell\cels}
  - \textbf{c}{^\celp}_{\cels\gop}\ p{_\celp}^{\cell\cels}
  + \partial_{\cell_1}p{_\gop}^{\cell\cell_1}
  + \frac{1}{2}\textbf{c}{^\cell}_{\cell_1\cell_2}p{_\gop}^{\cell_1\cell_2}
  \right)
  e^{(N-1)}_\cell
   \\
  & & +\ \left(\frac{1}{2}p{_\gop}^{\cels_1\cels_2}
 \mathring{\Theta}{^\cels}_{\cels_1\cels_2}
 + \partial_{\cell}p{_\gop}^{\cels\cell}
 + \textbf{c}{^\celp}_{\cels_1\gop}\kappa{_\celp}^{\cels_1\cels} \right) e_\cels^{(N-1)}  \\
 \hline
 \end{array}
\end{equation}
Summarizing with (\ref{detailQ}) and taking into account that $p{_\gop}^{\gos\gos} = \kappa{_\gop}^{\gos\gos}$, the fundamental equation $\dR^\textbf{A} p{_\gop} = Q_\gop - \frac{1}{2}Qe_\gop^{(N-1)}$ (\ref{equationfondamentaleGrave}) is equivalent to the system
\begin{equation}\label{bigECsystem}
\left\{
\begin{array}{r}
 \frac{1}{2}\kappa{_\gop}^{\cels_1\cels_2}
  \Omega{^\gol}_{\cels_1\cels_2}
+ (\partial_{\cels}^\omega
+ \Theta{^{\ast}}_{\cels\ast})
  p{_\gop}^{\gol\cels}
  - \textbf{c}{^\celp}_{\cels\gop}\ p{_\celp}^{\gol\cels}
  + \partial_{\cell_1}p{_\gop}^{\gol\cell_1}
  + \frac{1}{2}\textbf{c}{^\gol}_{\cell_1\cell_2}p{_\gop}^{\cell_1\cell_2}
\quad\quad\quad \\
 = \textbf{F}{^\celp}_{\gos\cels_1}\,
    p{_\celp}^{\gol\cels_1}
- \frac{1}{2} Q\delta{_\gol}^\gol \\
\frac{1}{2}\kappa{_\gop}^{\cels_1\cels_2}
 \mathring{\Theta}{^\gos}_{\cels_1\cels_2}
 + \partial_{\cell}p{_\gop}^{\gos\cell}
 + \textbf{c}{^\celp}_{\cels\gop}\kappa{_\celp}^{\cels\gos}
= \textbf{F}{^\celp}_{\gos\cels_1}\,
    \kappa{_\celp}^{\gos\cels_1}
-\frac{1}{2} Q\delta{_\gos}^\gos
\end{array}
\right.
\end{equation}
where $Q = \textbf{F}{^\celp}_{\cels\cels}\,
\kappa{_\celp}^{\cels\cels}$. Note that,
by (\ref{FenfonctionThetaetOmega}),
$\textbf{F}{^\gop}_{\gos\gos}
= \Theta{^\gos}_{\gos\gos} + \Omega{^\gol}_{\gos\gos} + \textbf{c}{^\gol}_{\gos\gos}$. Hence the first equation of (\ref{bigECsystem}) reads
\[
 \begin{array}{r}
 \frac{1}{2}\kappa{_\gop}^{\cels_1\cels_2}
  \Omega{^\gol}_{\cels_1\cels_2}
+ (\partial_{\cels}^\omega
+ \Theta{^{\ast}}_{\cels\ast})
  p{_\gop}^{\gol\cels}
  - \textbf{c}{^\celp}_{\cels\gop}\ p{_\celp}^{\gol\cels}
  + \partial_{\cell_1}p{_\gop}^{\gol\cell_1}
  + \frac{1}{2}\textbf{c}{^\gol}_{\cell_1\cell_2}p{_\gop}^{\cell_1\cell_2}
   \quad\quad\quad \\
 = \Theta{^{\cels_0}}_{\gos\cels}\,
p{_{\cels_0}}^{\gol\cels}
+ (\Omega{^{\cell}}_{\gos\cels}
+ \mathbf{c}{^{\cell}}_{\gos\cels})
    p{_{\cell}}^{\gol\cels}
- \frac{1}{2} Q\delta{_\gol}^\gol
\end{array}
\]
However we observe that the term $- \textbf{c}{^\celp}_{\cels\gop}\ p{_\celp}^{\gol\cels}$ on the left hand side is equal to $\textbf{c}{^\celp}_{\gop\cels}\ p{_\celp}^{\gol\cels} = \textbf{c}{^{\cels_0}}_{\gop\cels}\ p{_{\cels_0}}^{\gol\cels} + \textbf{c}{^\cell}_{\gop\cels}\ p{_\cell}^{\gol\cels} = \textbf{c}{^{\cels_0}}_{\gol\cels}\ p{_{\cels_0}}^{\gol\cels} + \textbf{c}{^\cell}_{\gos\cels}\ p{_\cell}^{\gol\cels}$, whereas the term $\textbf{c}{^\cell}_{\gos\cels}\ p{_\cell}^{\gol\cels}$ appears also on the right hand side. Hence the first equation in (\ref{bigECsystem}) simplifies to
\begin{equation}\label{bigECsystem1}
 \begin{array}{|r|}
 \hline 
 \frac{1}{2}\kappa{_\gop}^{\cels_1\cels_2}
  \Omega{^\gol}_{\cels_1\cels_2}
+ (\partial_{\cels}^\omega
+ \Theta{^{\ast}}_{\cels\ast})
  p{_\gop}^{\gol\cels}
  + \textbf{c}{^{\cels_0}}_{\gol\cels}\ p{_{\cels_0}}^{\gol\cels}
  + \partial_{\cell_1}p{_\gop}^{\gol\cell_1}
  + \frac{1}{2}\textbf{c}{^\gol}_{\cell_1\cell_2}p{_\gop}^{\cell_1\cell_2}
   \quad\quad\quad \\
 = \Theta{^{\cels_0}}_{\gos\cels}\,
p{_{\cels_0}}^{\gol\cels}
+ \Omega{^{\cell}}_{\gos\cels} p{_{\cell}}^{\gol\cels}
- \frac{1}{2} Q\delta{_\gol}^\gol \\
 \hline 
\end{array}
\end{equation}
with
\begin{equation}\label{Qmenfin}
Q = \Theta{^{\cels}}_{\cels_1\cels_2}\,
\kappa{_{\cels}}^{\cels_1\cels_2}
+ (\Omega{^{\cell}}_{\cels_1\cels_2}
+ \mathbf{c}{^{\cell}}_{\cels_1\cels_2})
\kappa{_{\cell}}^{\cels_1\cels_2}
\end{equation}
We call Equation (\ref{bigECsystem1}) the \textbf{(dynamical) equation on hidden fields}.\\

\noindent
The second equation reads
\[
 \begin{array}{r}
\frac{1}{2}\kappa{_\gop}^{\cels_1\cels_2}
 \mathring{\Theta}{^\gos}_{\cels_1\cels_2}
 + \partial_{\cell}p{_\gop}^{\gos\cell}
 + \textbf{c}{^\celp}_{\cels\gop}\kappa{_\celp}^{\cels\gos}
= \Theta{^{\cels_0}}_{\gos\cels_1}\,
\kappa{_{\cels_0}}^{\gos\cels_1}
+ (\Omega{^{\cell}}_{\gos\cels_1}
+ \mathbf{c}{^{\cell}}_{\gos\cels_1})
\kappa{_{\cell}}^{\gos\cels_1}
-\frac{1}{2} Q\delta{_\gos}^\gos
\end{array}
\]
Similarly the left hand side contains the term $\textbf{c}{^\celp}_{\cels\gop}\kappa{_\celp}^{\cels\gos}
= \textbf{c}{^{\cels_0}}_{\cels\gop}\kappa{_{\cels_0}}^{\cels\gos} + \textbf{c}{^\cell}_{\cels\gop}\kappa{_\cell}^{\cels\gos}
= \textbf{c}{^{\cels_0}}_{\cels\gol}\kappa{_{\cels_0}}^{\cels\gos} + \textbf{c}{^\cell}_{\cels\gos}\kappa{_\cell}^{\cels\gos}$, whereas the right hand side contains $\mathbf{c}{^{\cell}}_{\gos\cels_1}
\kappa{_{\cell}}^{\gos\cels_1} = \textbf{c}{^{\cell}}_{\cels\gos}\kappa{_{\cell}}^{\cels\gos}$. This leads to the simplification of the second equation in (\ref{bigECsystem}):
\begin{equation}\label{bigECsystem2}
\begin{array}{|r|}
 \hline 
\frac{1}{2}\kappa{_\gop}^{\cels_1\cels_2}
 \mathring{\Theta}{^\gos}_{\cels_1\cels_2}
 + \partial_{\cell}p{_\gop}^{\gos\cell}
 + \textbf{c}{^{\cels_0}}_{\cels\gol}\kappa{_{\cels_0}}^{\cels\gos}
= \Theta{^{\cels_0}}_{\gos\cels_1}\,
\kappa{_{\cels_0}}^{\gos\cels_1}
+ \Omega{^{\cell}}_{\gos\cels_1}
\kappa{_{\cell}}^{\gos\cels_1}
-\frac{1}{2} Q\delta{_\gos}^\gos \\
 \hline 
\end{array}
\end{equation}
We call Equation (\ref{bigECsystem2}) the \textbf{(dynamical) equation on physical fields}.\hfill $\square$\\

\noindent
This completes the proof of Theorem \ref{theoBigOne}.
It is important to keep in mind that, in all Equations (\ref{bigECsystem1}) and (\ref{bigECsystem2}),
\begin{itemize}
 \item quantities $\kappa{_\gop}^{\gos\gos}$ and $\mathbf{c}{^\gop}_{\gop\gop}$ are constant,
 \item  quantities $\omega^\gol$, $\mathbf{F}{^\gop}_{\gos\gos}$  and hence $\Theta{^\gos}_{\gos\gos}$, $\mathring{\Theta}{^\gos}_{\gos\gos}$, $\Omega{^\gol}_{\gos\gos}$ and $Q$ are constant along the integral leaves.
\end{itemize}

\subsection{Analysis of the dynamical equations}

\subsubsection{The equation on hidden fields (\ref{bigECsystem1})}\label{themysterious}
This equation can be written as:
\begin{equation}\label{bigECsystem1form}
 \begin{array}{r}
 \left(\frac{1}{2}\kappa{_\gop}^{\cels_1\cels_2}
  \Omega{^\cell}_{\cels_1\cels_2}
+ (\partial_{\cels}^\omega
+ \Theta{^{\ast}}_{\cels\ast})
  p{_\gop}^{\cell\cels}
+ \textbf{c}{^{\cels_0}}_{\gol\cels}\ p{_{\cels_0}}^{\cell\cels}
  + \partial_{\cell_1}p{_\gop}^{\cell\cell_1}
  + \frac{1}{2}\textbf{c}{^\cell}_{\cell_1\cell_2}p{_\gop}^{\cell_1\cell_2}
  \right)
  e^{(N-1)}_\cell
   \quad\quad\quad \\
 = \left( \Theta{^{\cels_0}}_{\cels\cels_1}\,
p{_{\cels_0}}^{\gol\cels_1}
+ \Omega{^{\cell_0}}_{\gos\cels_1}
    p{_{\cell_0}}^{\cell\cels_1} \right)e^{(N-1)}_\cell
- \frac{1}{2}Q e^{(N-1)}_\gol
\end{array}
\end{equation}
We observe that the quantity $\partial_{\cell_1}p{_\gop}^{\gol\cell_1}
+ \frac{1}{2}\textbf{c}{^\gol}_{\cell_1\cell_2}p{_\gop}^{\cell_1\cell_2}$ in the left hand side represents an exact term. Indeed, since $\omega^\gol = \omega{^\gol}_\cels e^\cels$, one has $\hbox{d} e_{\gol\gol}^{(N-2)} = \dR^\omega e_{\gol\gol}^{(N-2)}$, which implies because of (\ref{117bis}) that
$\hbox{d} e_{\gol\gol}^{(N-2)} = \dR^\omega e_{\gol\gol}^{(N-2)} = \textbf{c}{^\cell}_{\gol\gol} e_\cell^{(N-1)}$. Hence
\[
 \hbox{d}\left(\frac{1}{2}p{_\gop}^{\cell\cell}e_{\cell\cell}^{(N-2)}\right)
 = \frac{1}{2}\hbox{d} p{_\gop}^{\cell\cell}\wedge e_{\cell\cell}^{(N-2)}
 + \frac{1}{2}p{_\gop}^{\cell\cell}\hbox{d} e_{\cell\cell}^{(N-2)}
 = \left(\partial_{\cell_1} p{_\gop}^{\cell\cell_1}
 + \frac{1}{2}p{_\gop}^{\cell_1\cell_2}\textbf{c}{^\cell}_{\cell_1\cell_2}\right)
 e_\cell^{(N-1)}
\]
Thus (\ref{bigECsystem1}) or (\ref{bigECsystem1form}) is equivalent to
\begin{equation}\label{bigECsystem1form1}
 \begin{array}{r}
 \left(\frac{1}{2}\kappa{_\gop}^{\cels_1\cels_2}
  \Omega{^\cell}_{\cels_1\cels_2}
+ (\partial_{\cels}^\omega
+ \Theta{^{\ast}}_{\cels\ast})
  p{_\gop}^{\cell\cels}
+ \textbf{c}{^{\cels_0}}_{\gol\cels}\ p{_{\cels_0}}^{\cell\cels} \right)
  e^{(N-1)}_\cell
+ \hbox{d}\left(\frac{1}{2}p{_\gop}^{\cell\cell}e_{\cell\cell}^{(N-1)}\right)
   \quad\quad\quad \\
 = \left( \Theta{^{\cels_0}}_{\gos\cels_1}\,
p{_{\cels_0}}^{\cell\cels_1}
+ \Omega{^{\cell_0}}_{\gos\cels_1}
p{_{\cell_0}}^{\cell\cels_1} \right)e^{(N-1)}_\cell
- \frac{1}{2}Q e^{(N-1)}_\gol
\end{array}
\end{equation}
Hence $p{_\gop}^{\gol\gol}$ enters into play in the system only through an exact term.

\subsubsection{A conservation law}
By applying the exterior derivative to both sides of (\ref{bigECsystem1form1}) and by using the facts that $\partial_\gol\Omega{^\gol}_{\gos\gos} = \partial_\gol\Theta{^\gos}_{\gos\gos} = \partial_\gol Q = 0$ and that $\hbox{d}e_\gol^{(N-1)} = 0$ because $\gol$ is unimodular, we obtain
\begin{equation}\label{bigECsystem1derivee}
 \partial_\cell (\partial_{\cels}^\omega p{_\gop}^{\cell\cels}) + \Theta{^{\ast}}_{\cels\ast}
  \partial_\cell p{_\gop}^{\cell\cels}
+ \textbf{c}{^{\cels_0}}_{\gol\cels}\ \partial_\cell p{_{\cels_0}}^{\cell\cels}  =
  \Theta{^{\cels_0}}_{\gos\cels_1}\,
\partial_\cell p{_{\cels_0}}^{\cell\cels_1}
+ \Omega{^{\cell_0}}_{\gos\cels_1}
\partial_\cell p{_{\cell_0}}^{\cell\cels_1}
\end{equation}
\begin{lemm}\label{lemmedonnepartialT}
 Let us define
 \begin{equation}\label{notationTps}
  T{_\gop}^\gos:= \partial_\cell p{_\gop}^{\gos\cell}
 \end{equation}
and set $\partial_\gos^\omega T{_\gop}^\gos = \partial_\gos T{_\gop}^\gos - \mathbf{c}{^\celp}_{\cell\gop} \omega{^\cell}_\gos \, T{_\celp}^\gos + \mathbf{c}{^\gos}_{\cell\cels} \omega{^\cell}_\gos \, T{_\gop}^\cels$. Then
\[
 \partial_\cell (\partial_{\cels}^\omega p{_\gop}^{\cels\cell}) = \partial_\cels^\omega T{_\gop}^\cels
\]
\end{lemm}
A consequence of Lemma \ref{lemmedonnepartialT} is that Equation (\ref{bigECsystem1derivee}) is equivalent to
\begin{equation}\label{bigECsystem1derivee1}
\begin{array}{|c|}
 \hline
 \displaystyle
  \partial_\cels^\omega T{_\gop}^\cels + \Theta{^{\ast}}_{\cels\ast}
T{_\gop}^\cels
+ \textbf{c}{^{\cels_0}}_{\gol\cels}\ T{_{\cels_0}}^\cels
  =
  \Theta{^{\cels_0}}_{\gos\cels_1}\,
T{_{\cels_0}}^{\cels_1}
+ \Omega{^{\cell_0}}_{\gos\cels_1}
T{_{\cell_0}}^{\cels_1}     \\
\hline
\end{array}
\end{equation}
Equation (\ref{bigECsystem1derivee1}) can be splitted into the system
\begin{equation}\label{deuxombres}
  \left\{ 
 \begin{array}{ccc}
\partial_\cels^\omega T{_\gol}^\cels + \Theta{^{\ast}}_{\cels\ast}
T{_\gol}^\cels
+ \textbf{c}{^{\cels_0}}_{\gol\cels}\ T{_{\cels_0}}^\cels
& = & 0 \\
\partial_\cels^\omega T{_\gos}^\cels + \Theta{^{\ast}}_{\cels\ast}
T{_\gos}^\cels
& = &
\Theta{^{\cels_0}}_{\gos\cels_1}\,
T{_{\cels_0}}^{\cels_1}
+ \Omega{^{\cell_0}}_{\gos\cels_1}
T{_{\cell_0}}^{\cels_1}
 \end{array}\right.
\end{equation}
We will see later on that $T{_\gol}^\gos$ and $T{_\gos}^\gos$ can be interpreted as, respectively, an angular momentum tensor and a stress-energy tensor. Hence Relations (\ref{deuxombres}) express the conservation of these tensors.
The proof of Lemma \ref{lemmedonnepartialT} rests on the following result.
\begin{lemm}\label{lemmeCommutateursGraves}
 The vector fields $\partial_\gop$ satisfy the following commutation relation:
\begin{equation}\label{collectiondecommutateurs}
 \left\{
 \begin{array}{ccc}
[\partial_{\gos_1},\partial_{\gos_2}] & = & -\ \Theta{^\cels}_{\gos_1\gos_2}\,\partial_\cels - \Omega{^\cell}_{\gos_1\gos_2}\,\partial_\cell
+ (\mathbf{c}{^\cels}_{\cell\gos_2}\,\omega{^{\cell}}_{\gos_1} - \mathbf{c}{^\cels}_{\cell\gos_1}\,\omega{^{\cell}}_{\gos_2})\partial_\cels
\\
{[} \partial_\gos ,\partial_\gol  {]} & = & \mathbf{c}{^\cell}_{\cell_1\gol}\,\omega{^{\cell_1}}_\gos\,\partial_\cell \\
{[} \partial_{\gol_1},\partial_{\gol_2} {]} & = & -\ \mathbf{c}{^\cell}_{\gol_1\gol_2}\,\partial_\cell
 \end{array}\right.
\end{equation}
\end{lemm}
\emph{Proof of Lemma \ref{lemmeCommutateursGraves}} --- We first deduce from (\ref{117bis}), which reads $\hbox{d}e^\gop + [\omega^\gol \wedge e^\gop] = \Theta^\gos + \Omega^\gol + \frac{1}{2}[e^\gol\wedge e^\gol]$, that
\begin{equation}\label{lebeaucalculdedep}
 \hbox{d}e^\gop = \Theta^\gos + \Omega^\gol + \frac{1}{2}[e^\gol\wedge e^\gol]- [\omega^\gol \wedge e^\gop]
\end{equation}
By using Cartan's formula $\hbox{d}e^\gop(\partial_{\gop_1},\partial_{\gop_1}) + e^\gop([\partial_{\gop_1},\partial_{\gop_2}]) = \partial_{\gop_1}(e^\gop(\partial_{\gop_2})) - \partial_{\gop_2}(e^\gop(\partial_{\gop_1})) = 0$ we get $e^\gop([\partial_{\gop_1},\partial_{\gop_2}]) = - \hbox{d}e^\gop(\partial_{\gop_1},\partial_{\gop_1})$ and hence, by (\ref{lebeaucalculdedep}) we get that
\[
 e^\gop([\partial_{\gop_1},\partial_{\gop_2}]) =  -\Theta^\gos(\partial_{\gop_1},\partial_{\gop_2}) - \Omega^\gol(\partial_{\gop_1},\partial_{\gop_2}) - \frac{1}{2}[e^\gol\wedge e^\gol](\partial_{\gop_1},\partial_{\gop_2}) + [\omega^\gol \wedge e^\gop](\partial_{\gop_1},\partial_{\gop_2})
\]
Hence we deduce (\ref{collectiondecommutateurs}). \hfill $\square$\\

\noindent
\emph{Proof of Lemma \ref{lemmedonnepartialT}} --- From
\[
 \partial_\cels^\omega p{_\gop}^{\cels\gol} = \partial_\cels p{_\gop}^{\cels\gol} - \mathbf{c}{^\celp}_{\cell\gop} \omega{^\cell}_\cels \, p{_\celp}^{\cels\gol} + \mathbf{c}{^\gol}_{\cell\cell_2} \omega{^\cell}_\cels \, p{_\gop}^{\cels\cell_2}  + \mathbf{c}{^\cels}_{\cell\cels_2} \omega{^\cell}_\cels \, p{_\gop}^{\cels_2\gol}
\]
we get first
\[
\begin{array}{ccl}
\partial_\cell (\partial_{\cels}^\omega p{_\gop}^{\cels\cell})
& = &
\partial_{\cell}\left(\partial_\cels p{_\gop}^{\cels\cell}
- \mathbf{c}{^\celp}_{\cell_1\gop} \omega{^{\cell_1}}_\cels \, p{_\celp}^{\cels\cell}
+ \mathbf{c}{^\cell}_{\cell_1\cell_2} \omega{^{\cell_1}}_\cels \, p{_\gop}^{\cels\cell_2}
+ \mathbf{c}{^\cels}_{\cell_1\cels_2} \omega{^{\cell_1}}_\cels \, p{_\gop}^{\cels_2\cell}\right) \\
& = & \partial_{\cell}\left(\partial_\cels p{_\gop}^{\cels\cell}\right)
- \mathbf{c}{^\celp}_{\cell_1\gop} \omega{^{\cell_1}}_\cels \, T{_\celp}^{\cels}
+ \mathbf{c}{^\cell}_{\cell_1\cell_2} \omega{^{\cell_1}}_\cels \, \partial_{\cell}p{_\gop}^{\cels\cell_2}
+ \mathbf{c}{^\cels}_{\cell_1\cels_2} \omega{^{\cell_1}}_\cels \, T{_\gop}^{\cels_2}
\end{array}
\]
On the other hand we deduce from (\ref{collectiondecommutateurs}) that
\[
 \partial_\gol\left(\partial_\gos p{_\gop}^{\gos\gol}\right) =
 \partial_\gos\left(\partial_\gol p{_\gop}^{\gos\gol}\right)
 + [\partial_\gol,\partial_\gos]p{_\gop}^{\gos\gol}
  =
 \partial_\gos\left(\partial_\gol p{_\gop}^{\gos\gol}\right)
 - \mathbf{c}{^{\cell}}_{\cell_1\gol}
 \omega{^{\cell_1}}_\gos
 \partial_{\cell}p{_\gop}^{\gos\gol}
\]
hence 
\[
 \partial_\cell\left(\partial_\cels p{_\gop}^{\cels\cell}\right)
 =
 \partial_\cels\left(\partial_\cell p{_\gop}^{\cels\cell}\right)
 - \mathbf{c}{^{\cell}}_{\cell_1\cell_2}
 \omega{^{\cell_1}}_\cels
 \partial_{\cell}p{_\gop}^{\cels\cell_2}
 = \partial_\cels T{_\gop}^{\cels}
 - \mathbf{c}{^{\cell}}_{\cell_1\cell_2}
 \omega{^{\cell_1}}_\cels
 \partial_{\cell}p{_\gop}^{\cels\cell_2}
\]
We thus deduce by replacing this expression for $\partial_\cell\left(\partial_\cels p{_\gop}^{\cels\cell}\right)$ in the formula giving $\partial_\cell (\partial_{\cels}^\omega p{_\gop}^{\cels\cell})$ that
\[
 \partial_\cell (\partial_{\cels}^\omega p{_\gop}^{\cels\cell})
 = \partial_\cels T{_\gop}^{\cels}
- \mathbf{c}{^\celp}_{\cell_1\gop} \omega{^{\cell_1}}_\cels \, T{_\celp}^{\cels}
+ \mathbf{c}{^\cels}_{\cell_1\cels_2} \omega{^{\cell_1}}_\cels \, T{_\gop}^{\cels_2}
\]
The equivalence between (\ref{bigECsystem1derivee}) and (\ref{bigECsystem1derivee1}) is then straightforward.
\hfill $\square$\\

\subsubsection{The equation on physical fields (\ref{bigECsystem2})}\label{sectionphysicaleqgrav}
We can split  (\ref{bigECsystem2}) into the system
\begin{equation}\label{bigECsystem2split}
 \left\{
 \begin{array}{ccc}
 - \frac{1}{2}\kappa{_\gol}^{\cels_1\cels_2}
 \mathring{\Theta}{^\gos}_{\cels_1\cels_2}
 & = & \partial_{\cell}p{_\gol}^{\gos\cell}
 + \textbf{c}{^{\cels_0}}_{\cels\gol}\ \kappa{_{\cels_0}}^{\cels\gos}
 \\
\Omega{^{\cell}}_{\gos\cels_1}
\kappa{_{\cell}}^{\gos\cels_1}
-\frac{1}{2} Q\delta{_\gos}^\gos
& = &
\partial_{\cell}p{_\gos}^{\gos\cell}
-(\Theta{^{\cels_0}}_{\gos\cels_1}\,
\kappa{_{\cels_0}}^{\gos\cels_1}
- \frac{1}{2}
 \mathring{\Theta}{^\gos}_{\cels_1\cels_2}\kappa{_\gos}^{\cels_1\cels_2})
 \end{array}\right.
\end{equation}
It can be rephrased by using the notations $\tilde{\mathbf{C}}{_\gol}^\gos$ and $\tilde{\mathbf{E}}{_\gos}^\gos$ given in (\ref{defiCtenseurCartan}) (\ref{defiEtenseurEinstein}) and (\ref{EssestEinstein}) and the relation (\ref{Qmenfin}) for $Q$ for the left hand sides and by using the notation $T{_\gop}^\gos$ (\ref{notationTps}) for the right hand sides:
\[
 \left\{
 \begin{array}{ccl}
 \tilde{\mathbf{C}}{_\gol}^\gos
 & = & T{_\gol}^\gos
 + \textbf{c}{^{\cels_0}}_{\cels\gol}\ \kappa{_{\cels_0}}^{\cels\gos}
 \\
\tilde{\mathbf{E}}{_\gos}^\gos
+ \Lambda\delta{_\gos}^\gos
& = &
T{_\gos}^\gos
+ \frac{1}{2} (\Theta{^\cels}_{\cels_1\cels_2}\kappa{_\cels}^{\cels_1\cels_2})\delta{_\gos}^\gos
\Theta{^{\cels_0}}_{\gos\cels_1}\,
\kappa{_{\cels_0}}^{\gos\cels_1}
+ \frac{1}{2}
 \mathring{\Theta}{^\gos}_{\cels_1\cels_2}\kappa{_\gos}^{\cels_1\cels_2}
 \end{array}\right.
\]
where $\Lambda:= - \frac{1}{2}\mathbf{c}{^{\cell}}_{\cels_1\cels_2}
\kappa{_{\cell}}^{\cels_1\cels_2}$.

In the following we will assume the \textbf{additional Hypothesis} $\kappa{_\gos}^{\gos\gos} =0$ (\ref{additionalhypothesis}), which is satisfied in all usual gravity theories. The previous system simplifies to
\begin{equation}\label{bigECsystem2splitsimple}
 \left\{
 \begin{array}{ccl}
 \tilde{\mathbf{C}}{_\gol}^\gos
 & = & T{_\gol}^\gos
 \\
\tilde{\mathbf{E}}{_\gos}^\gos + \Lambda\delta{_\gos}^\gos
& = &
T{_\gos}^\gos
 \end{array}\right.
\end{equation}
The first equation can be interpreted as a generalization of the \emph{Cartan equation} involving the torsion of the connection $\omega^\gol$, whereas in the second equation the left hand side can be interpreted as the \emph{Einstein tensor} with a cosmological constant $\Lambda$.

Indeed recall that $\tilde{\mathbf{C}}{_\gog}^\gos:= - \frac{1}{2} \kappa{_\gog}^{\cels_1\cels_2} \mathring{\Theta}{^\gos}_{\cels_1\cels_2}$ (\ref{defiCtenseurCartan}) and assume for simplicity that $T{_\gol}^\gos = 0$. Then, by assuming that the tensor $\kappa{_\gol}^{\gos\gos}$ is non degenerate (i.e. the map $\Lambda^2\gos^*\ni \xi_{\gos\gos}\longmapsto \xi_{\cels\cels}\kappa{_\gol}^{\cels\cels} \in \gol^*$ is invertible), which is true in all standard situations, then the equation $\mathbf{C}{_\gog}^\gos =0$ is equivalent to $\mathring{\Theta}{^\gos}_{\gos\gos} =0$. Since 
$\mathring{\Theta}{^\gos}_{\gos\gos}
= \mathring{\Theta}{^c}_{ab}\textbf{t}_c\otimes \textbf{t}^a\otimes \textbf{t}^b$, with
$\mathring{\Theta}{^c}_{ab}
=\Theta{^c}_{ab} - \delta^c_b\Theta{^d}_{ad} - \delta^c_a\Theta{^d}_{db}$, we have
$\mathring{\Theta}{^c}_{ac} = (2-n)\Theta{^d}_{ad}$, from which
we deduce that, if $n\neq 2$, $\Theta{^\gos}_{\gos\gos} = \frac{1}{2}\Theta{^c}_{ab}\textbf{t}_c\otimes \textbf{t}^a\otimes \textbf{t}^b$, with
\begin{equation}\label{inverseThetaChapeau}
 \Theta{^c}_{ab} = \mathring{\Theta}{^c}_{ab} - \frac{1}{n-2}
 \left(\delta^c_b\mathring{\Theta}{^d}_{ad} + \delta^c_a\mathring{\Theta}{^d}_{db}\right)
\end{equation}
and thus, if $n>2$,
$\Theta{^\gos}_{\gos\gos} = 0$ if and only if $\mathring{\Theta}{^\gos}_{\gos\gos} = 0$. Alternatively this conclusion can also be deduced from (\ref{defimathringTheta}).
Similarly the second equation in  (\ref{bigECsystem2splitsimple}) relates the  sum of the generalized \emph{Einstein} tensor and a cosmological constant on the left hand side to $T{_\gos}^\gos = - \partial_{\cell}p{_\gos}^{\gos\cell}$ which plays here the role of a stress-energy tensor.

We see that the only way coupling between the fields $(\theta^\gos,\omega^\gol)$ and the fields $p{_\gop}^{\gop\gop}$ in the generalized Einstein--Cartan system is operated by $T{_\gop}^\gos = \partial_{\cell}p{_\gop}^{\gos\cell} $. Moreover System (\ref{bigECsystem2splitsimple}) tells us that $T{_\gop}^\gos$ is constant on each leaf of the integral foliation, since the left hand side of this system is so.

\subsubsection{Constraints on the generalized Einstein and Cartan tensors}
It is well-known in General Relativity that, for any connection without torsion, the Einstein tensor $\mathbf{E}{_\gos}^\gos$ satisfies the constraint $\partial^\omega_\cels\mathbf{E}{_\gos}^\cels = 0$. Thus a necessary condition for the Einstein equation $\mathbf{E}{_\gos}^\gos = T{_\gos}^\gos$ to admit solutions is that the stress-energy tensor $T{_\gos}^\gos$ satisfies the  similar relation $\partial^\omega_\cels T{_\gos}^\cels = 0$. The latter constraint expresses the conservation of the energy-momentum tensor.

Similarly the left hand side of the Einstein--Cartan system (\ref{bigECsystem2splitsimple}) satisfies constraints which match with the conservation laws (\ref{deuxombres}).
Indeed if the generalized Einstein--Cartan system   (\ref{bigECsystem2splitsimple}) and the conservation law (\ref{bigECsystem1derivee1}) are satisfied, then, by replacing $T{_\gou}^\gos$ in (\ref{bigECsystem1derivee1}) by its expression in function of $\tilde{\mathbf{C}}{_\gol}^\gos$ and $\tilde{\mathbf{E}}{_\gos}^\gos$ given by (\ref{bigECsystem2splitsimple}), one obtains
\[
  \left\{
 \begin{array}{ccc}
 \partial_\cels^\omega \tilde{\mathbf{C}}{_\gol}^\cels + \Theta{^{\ast}}_{\cels\ast}
\tilde{\mathbf{C}}{_\gol}^\cels
+ \textbf{c}{^{\cels_0}}_{\gol\cels}\left(\tilde{\mathbf{E}}{_{\cels_0}}^\cels
+ \Lambda\delta{_{\cels_0}}^\cels\right) & = & 0 \\
\partial_\cels^\omega \left(\tilde{\mathbf{E}}{_\gos}^\cels + \Lambda \delta{_\gos}^\cels\right) + \Theta{^{\ast}}_{\cels\ast}
\left(\tilde{\mathbf{E}}{_\gos}^\cels + \Lambda \delta{_\gos}^\cels\right)
& = &
\Theta{^{\cels_0}}_{\gos\cels_1}\left(\tilde{\mathbf{E}}{_{\cels_0}}^{\cels_1} + \Lambda \delta{_{\cels_0}}^{\cels_1}\right)
+ \Omega{^{\cell_0}}_{\gos\cels_1}
\tilde{\mathbf{C}}{_{\cell_0}}^{\cels_1}
 \end{array}\right.
\]
This system can easily be simplified by observing that $\mathbf{c}{^{\celp}}_{\celp\gol} - \mathbf{c}{^{\cell}}_{\cell\gol} = 0 - 0 = 0$ (because $\gop$ are $\gol$ unimodular),  $\textbf{c}{^{\cels_0}}_{\gol\cels}\,\Lambda\delta{_{\cels_0}}^\cels = \Lambda\, \textbf{c}{^{\cels}}_{\gol\cels} = 0$, $\partial_\cels^\omega \left(\Lambda \delta{_\gos}^\cels\right) = 0$ and $\Theta{^{\ast}}_{\cels\ast}
\left( \Lambda \delta{_\gos}^\cels\right) = \Theta{^{\cels_0}}_{\gos\cels_1}\left( \Lambda \delta{_{\cels_0}}^{\cels_1}\right) = \Lambda \Theta{^\ast}_{\gos\ast}$, leading to the following
 \begin{equation}\label{contraintedeBianchi}
\left\{
 \begin{array}{ccc}
 \partial_\cels^\omega \tilde{\mathbf{C}}{_\gol}^\cels + \Theta{^{\ast}}_{\cels\ast}
\tilde{\mathbf{C}}{_\gol}^\cels
+ \textbf{c}{^{\cels_0}}_{\gol\cels}\,\tilde{\mathbf{E}}{_{\cels_0}}^\cels & = & 0 \\
\partial_\cels^\omega \tilde{\mathbf{E}}{_\gos}^\cels  + \Theta{^{\ast}}_{\cels\ast}
\tilde{\mathbf{E}}{_\gos}^\cels
& = &
\Theta{^{\cels_0}}_{\gos\cels_1}\tilde{\mathbf{E}}{_{\cels_0}}^{\cels_1}
+ \Omega{^{\cell_0}}_{\gos\cels_1}
\tilde{\mathbf{C}}{_{\cell_0}}^{\cels_1}
 \end{array}\right.
\end{equation}
In Proposition \ref{propoBianchi} (see the Appendix) we prove that (\ref{contraintedeBianchi}) is actually a consequence of the very definition of $\tilde{\mathbf{C}}{_\gol}^\gos$ and $\tilde{\mathbf{E}}{_\gos}^\gos$, confirming that (\ref{bigECsystem1}) is a \emph{necessary} condition for  (\ref{bigECsystem2splitsimple}) to have solutions.

\subsection{Exploitation of the equations}\label{sectionexploitation}
In the following we prove Corollary \ref{corooftheBigOne}. We start by assuming the following\\
\begin{equation}\label{fibrationhypothesis}
\begin{array}{c}
 \hbox{\textbf{Fibration hypothesis :}}\\
 \hbox{\emph{The integral leaves of the exterior differential system (\ref{AbstractPfaffianSystem0}) }} \\
 \hbox{\emph{form a fibration of $\mathcal{F}$.}}
\end{array}
\end{equation}

Hence the manifold $\mathcal{F}$ is the total space of a principal bundle over some manifold $\mathcal{X}$, the structure group of which is $\widehat{\mathfrak{G}}$ or a quotient of $\widehat{\mathfrak{G}}$ by a finite subgroup. Note that, even if the group is compact and $\mathcal{F}$ is $\gol$-complete (see Definition \ref{defigcomplete}), so that one can prove that each leaf is compact, there may be some obstructions for the \emph{Fibration hypothesis} to be true since the topology of the leaves may vary (see \cite{pierard2}).

A consequence of (\ref{fibrationhypothesis}) is that we can interpret $\theta^\gos$ and $\omega^\gol$ as respectively a soldering and a connection form on this bundle defining a geometric structure on $\mathcal{X}$.

The key point in the Einstein--Cartan equations (\ref{bigECsystem2split}) or (\ref{bigECsystem2splitsimple}) is that the left hand sides are constant
along the fibers. Thus these Einstein--Cartan equations imply that the restrictions on any fiber of
$T{_\gol}^\gos = \partial_\cell  p{_\gol}^{\cell\gos}$ and
$T{_\gos}^\gos = \partial_\cell p{_\gos}^{\cell\gos}$ are also constant. Hence
if we assume that one of the two following hypotheses holds, we will deduce that
these right hand sides vanish.
\begin{enumerate}
 \item (non compactness) the fibers are non compact (which occurs if, for instance,
 $\widehat{\goL}$ is isomorphic to $Spin(1,n)$) and the first derivatives of
 the field $p{_\gop}^{\gos\gol}$ decay at infinity in each fiber.
 \item (compactness) each fiber $\textsf{f}$ is compact. This case occurs if, for instance,
 $\widehat{\goL}$ is isomorphic to $Spin(n)$ (or its spin group).
\end{enumerate}
Indeed if we assume (i), then we deduce that each
$\partial_\cell p{_\gop}^{\cell\gos}$
decays at infinity on each fiber, but on the other hand such a
quantity is constant along the fibers, hence it vanishes.

If we assume (ii) then the value of $\partial_\cell p{_\gop}^{\cell\gos}$ at any point is equal to its average
on the fiber $\textsf{f}$ which contains this point, hence, by
setting $(e^\gol)^{(r)}:= e^{n+1}\wedge \cdots \wedge e^N$ and
$(e^\gol)^{(r-1)}_\gol:= \frac{\partial}{\partial e^\gol}\iN (e^\gol)^{(r)}$,
\begin{equation}
T{_\gop}^\gos =
 \partial_\cell p{_\gop}^{\cell\gos}
 = \frac{\displaystyle \int_\textsf{f}
 \partial_\cell p{_\gop}^{\cell\gos}(e^\gol)^{(r)}}{\displaystyle \int_\textsf{f}(e^\gol)^{(r)}}
 = \frac{\displaystyle \int_\textsf{f}
 \dR\left( p{_\gop}^{\cell\gos}(e^\gol)^{(r-1)}_\cell\right)}{\displaystyle \int_\textsf{f}(e^\gol)^{(r)}} = 0
\end{equation}
and we achieve the same conclusion.
Hence assuming either (i) or (ii) and assuming also $\kappa{_\gos}^{\gos\gos} = 0$ for simplicity we deduce from (\ref{bigECsystem2splitsimple}) that
our fields are solutions of the system
\begin{equation}\label{brutbeautiful}
 \left\{
 \begin{array}{ccc}
  \frac{1}{2}\kappa{_\gol}^{\cels_1\cels_2}
 \mathring{\Theta}{^\gos}_{\cels_1\cels_2}
 & = & 0 \\
\Omega{^{\cell}}_{\cels_1\gos}\,
\kappa{_{\cell}}^{\cels_1\gos}
-\frac{1}{2} (\Omega{^{\cell}}_{\cels_1\cels_2}
\kappa{_{\cell}}^{\cels_1\cels_2})\delta{_\gos}^\gos + \Lambda\delta{_\gos}^\gos
& = & 0
 \end{array}\right.
\end{equation}
The first equation will imply that the generalized torsion
$\Theta$ vanishes, provided that the kernel of the linear map
$\mathring{\Theta}{^\gos}_{\gos\gos}\longmapsto \frac{1}{2}\kappa{_\gol}^{\cels\cels}\mathring{\Theta}{^\gos}_{\cels\cels}$ is $\{0\}$,
which will be the case in the following examples. The second equation is a generalization of the Einstein
equation in vacuum with a cosmological constant (and it will be so in basic examples).

\subsection{Gauge symmetries}

\subsubsection{Invariance by diffeomorphisms}
$\mathcal{A}$ is invariant under the transformation
$(\pi_\gop,\varphi^\gop)\longmapsto (T^*\pi_\gop,T^*\varphi^\gop)$,
where $T:\mathcal{F}\longrightarrow\mathcal{F}$
is a diffeomorphism which
preserves the orientation $
 \mathcal{A}[T^*\pi_\gop,T^*\varphi^\gop] = \mathcal{A}[\pi_\gop,\varphi^\gop]$.
Moreover the constraint $\pi_\gop\wedge
\varphi^\mathfrak{s}\wedge \varphi^\mathfrak{s}= \kappa{_\gop}^{\gos\gos} \varphi^{(N)}$
is invariant by such transformations.

\subsubsection{Adjoint action of the gauge group}
For any $g\in \mathcal{C}^\infty(\mathcal{F},\widehat{\goL})$ the action is clearly invariant by the gauge transformation
\begin{equation}\label{standardgaugesymmetry}
 \left\{
 \begin{array}{ccl}
  \varphi^\gop & \longmapsto & \hbox{Ad}_g\varphi^\gop -\dR g\,g^{-1} = g\varphi^\gop g^{-1} -\dR g\,g^{-1} \\
  \pi_\gop & \longmapsto & \hbox{Ad}_g^*\pi_\gop
 \end{array}\right.
\end{equation}
Moreover, since $\hbox{Ad}_g\kappa{_\gop}^{\gos\gos} = \kappa{_\gop}^{\gos\gos}$, the constraint $\pi_\gop\wedge\varphi^\mathfrak{s}\wedge \varphi^\mathfrak{s} = \kappa{_\gop}^{\gos\gos} \varphi^{(N)}$
is also invariant by this action.

\subsubsection{Gauge symmetries of the dual gauge fields} 
Using exactly the same arguments as in \S
\ref{paragraphgaugesymmetryYMpi} for Yang--Mills theory
we may write, for any $\chi_\gop\in \gop^*\otimes \Omega^{N-2}(\mathcal{F})$,
\[
 \mathcal{A}[\pi_\gop+\chi_\gop,\varphi^\gop] =
 \mathcal{A}[\pi_\gop,\varphi^\gop] +
 \int_{\mathcal{F}} \dR\left(\varphi^\celp\wedge \chi_\celp\right)
  + \varphi^\celp\wedge \dR^{\varphi/2}\chi_\celp
\]
Thus if $\chi_\gop\in \gop^*\otimes \Omega^{N-2}(\mathcal{F})$ decays rapidly at infinity and is a solution of
\begin{equation}\label{conditiondejaugelineaireE}
 \varphi^\celp\wedge \dR^{\varphi/2}\chi_\celp = 0
\end{equation}
then
\[
 \mathcal{A}[\pi_\gop+\chi_\gop,\varphi^\gop] =
 \mathcal{A}[\pi_\gop,\varphi^\gop].
 \]
If furthermore $\chi_\gop$ satisfies
\begin{equation}\label{conditiondejaugecontrainte}
\chi_\gop \wedge
 \varphi^\mathfrak{s}\wedge \varphi^\mathfrak{s} = 0,
\end{equation}
then $(\pi_\gop+\chi_\gop,\varphi^\gop)$ satisfies the constraint (\ref{constraintaxiomaticindex}).
Hence these three conditions are sufficient for having a gauge symmetry of the variational problem.

Moreover as for the Yang--Mills model (see \S \ref{paragraphgaugesymmetryYMpi})) any field $\chi_\gop$ which satisfies (\ref{conditiondejaugecontrainte}) provides us with an \emph{on shell} symmetry of the action.

\subsection{Appendix: proof of the generalized Bianchi identities}
We prove here that identities (\ref{contraintedeBianchi}), i.e.
\[
\left\{
 \begin{array}{ccc}
 \partial_\cels^\omega \tilde{\mathbf{C}}{_\gol}^\cels + \Theta{^{\ast}}_{\cels\ast}
\tilde{\mathbf{C}}{_\gol}^\cels
+ \textbf{c}{^{\cels_0}}_{\gol\cels}\,\tilde{\mathbf{E}}{_{\cels_0}}^\cels & = & 0 \\
\partial_\cels^\omega \tilde{\mathbf{E}}{_\gos}^\cels  + \Theta{^{\ast}}_{\cels\ast}
\tilde{\mathbf{E}}{_\gos}^\cels
& = &
\Theta{^{\cels_0}}_{\gos\cels_1}\tilde{\mathbf{E}}{_{\cels_0}}^{\cels_1}
+ \Omega{^{\cell_0}}_{\gos\cels_1}
\tilde{\mathbf{C}}{_{\cell_0}}^{\cels_1}
 \end{array}\right.
\]
are structure equations, hence automatically satisfied.
\begin{prop}\label{propoBianchi}
Let $\tilde{\mathbf{C}}{_\gol}^\gos$ and $\tilde{\mathbf{E}}{_\gos}^\gos$ be defined by, respectively, (\ref{defiCtenseurCartan}) and (\ref{defiEtenseurEinstein}) (or, equivalentely, (\ref{defiCtenseurCartanimplicit}) and (\ref{defiEtenseurEinsteinimplicit})). Assume that the tensor $\kappa{_\gop}^{\gos\gos}$ is invariant by the adjoint action of $\gol$. Then the relations (\ref{contraintedeBianchi}) hold.
\end{prop}

\noindent
\emph{Proof of Proposition \ref{propoBianchi}} --- \emph{Step 1: preliminary results} --- We first prove the Bianchi relations
\begin{equation}\label{Bianchisimple}
 \left\{
 \begin{array}{ccl}
  \hbox{d}^\omega\Theta^\gos + \left[e^\gos \wedge \Omega^\gol\right]  & = & 0 \\
  \hbox{d}^\omega\Omega^\gol & = & 0
 \end{array}\right.
\end{equation}
These relations follows from the relation
\begin{equation}\label{Bianchipure}
 \hbox{d}\mathbf{F}^\gop + [\mathbf{A}^\gop\wedge \mathbf{F}^\gop] = 0
\end{equation}
where we recall that $\mathbf{A}^\gop = \theta^\gos + \omega^\gog = e^\gos + \omega^\gog$ and $\mathbf{F}^\gop
:= \hbox{d}\mathbf{A}^\gop
+\frac{1}{2}[\mathbf{A}^\gop\wedge \mathbf{A}^\gop]
=  \Theta^\gos + \Omega^\gol + \mathbf{c}^\gol$, where we set $\mathbf{c}^\gol:= \frac{1}{2}\mathbf{c}{^\gol}_{\cels\cels}e^\cels\wedge e^\cels$.
Identity (\ref{Bianchipure}) thus reads
\[
 \hbox{d}\left(\Theta^\gos + \Omega^\gol + \mathbf{c}^\gol\right) + \left[\left(\theta^\gos + \omega^\gog\right)\wedge \left(\Theta^\gos + \Omega^\gol + \mathbf{c}^\gol\right)\right] = 0
\]
which, through the decomposition $\gop = \gos \oplus \gol$, splits into~:
\[
 \left\{
 \begin{array}{ccl}
  \hbox{d}\Theta^\gos + \left[\theta^\gos \wedge (\Omega^\gol+\mathbf{c}^\gol)\right] + \left[ \omega^\gog\wedge \Theta^\gos \right] & = & 0 \\
  \hbox{d}\left(\Omega^\gol + \mathbf{c}^\gol\right) + \left[\theta^\gos\wedge \Theta^\gos\right]  + \left[ \omega^\gog\wedge \left(\Omega^\gol + \mathbf{c}^\gol\right) \right] & = & 0
 \end{array}\right.
\]
or
\begin{equation}\label{Bianchipasencoresimple}
 \left\{
 \begin{array}{ccl}
  \hbox{d}^\omega\Theta^\gos + \left[\theta^\gos \wedge (\Omega^\gol+\mathbf{c}^\gol)\right]  & = & 0 \\
  \hbox{d}^\omega\left(\Omega^\gol + \mathbf{c}^\gol\right) + \left[\theta^\gos\wedge \Theta^\gos\right]   & = & 0
 \end{array}\right.
\end{equation}
However, a consequence of the Jacobi identity is that
\[
 [\theta^\gos\wedge \mathbf{c}^\gol] = [e^\gos\wedge \mathbf{c}^\gol] = \mathbf{c}{^\gos}_{\cels\celg} e^\cels\wedge \left(\frac{1}{2}\mathbf{c}{^\celg}_{\cels_1\cels_2}e^{\cels_1}\wedge e^{\cels_2}\right)
 = \frac{1}{2}\mathbf{c}{^\gos}_{\cels\celp}\mathbf{c}{^\celp}_{\cels_1\cels_2}\
 e^{\cels\cels_1\cels_2} = 0
\]
and on the other hand
\[
 \begin{array}{ccl}
\hbox{d}^\omega \mathbf{c}^\gol
  & = &
\hbox{d}^\omega \left(\frac{1}{2}\mathbf{c}{^\gol}_{\cels_1\cels_2}e^{\cels_1}\wedge e^{\cels_2} \right)
= \frac{1}{2}\mathbf{c}{^\gol}_{\cels_1\cels_2} \left(\hbox{d}^\omega e^{\cels_1} \wedge e^{\cels_2} - e^{\cels_1}\wedge \hbox{d}^\omega e^{\cels_2}\right) \\
& = &
\frac{1}{2}\mathbf{c}{^\gol}_{\cels_1\cels_2} \left( \Theta^{\cels_1}\wedge e^{\cels_2} - e^{\cels_1}\wedge \Theta^{\cels_2}\right) = [\Theta^\gos\wedge e^\gos]
 \end{array}
\]
Hence (\ref{Bianchipasencoresimple}) implies (\ref{Bianchisimple}).

We also need the two following lemmas.
\begin{lemm}\label{lemmedeNmoinsun}
By using Notation (\ref{Thetaast}) for $\Theta{^\ast}_{\gos\ast}$,
we have
\begin{equation}\label{domegaesparlemme}
   \hbox{\emph{d}}^\omega e^{(N-1)}_\gos
  = \Theta{^\ast}_{\gos\ast}\ e^{(N)}
\end{equation}
\end{lemm}
\emph{Proof} ---
$\hbox{d}^\omega e^{(N-1)}_\gos
 = \hbox{d}^\omega e^{\cels} \wedge e^{(N-2)}_{\gos\cels}
 = \Theta^\cels \wedge e^{(N-2)}_{\gos\cels}
 = \Theta{^\cels}_{\gos\cels} \ e^{(N)} = \Theta{^\ast}_{\gos\ast}\ e^{(N)}$.
\hfill $\square$\\
\begin{lemm}\label{aucoroetonnant}
 \begin{equation}\label{desss}
\begin{array}{ccc}
  \hbox{\emph{d}}^\omega e^{(N-3)}_{\gos_1\gos_2\gos_3}
& = & \Theta{^\ast}
_{\gos_3\ast}e^{(N-2)}_{\gos_1\gos_2}
+ \Theta{^\ast}_{\gos_1\ast}e^{(N-2)}_{\gos_2\gos_3}
+ \Theta{^\ast}_{\gos_2\ast}e^{(N-2)}_{\gos_3\gos_1} \\
& &
+\ \Theta{^\cels}_{\gos_3\gos_2}e^{(N-2)}_{\cels\gos_1}
+ \Theta{^\cels}_{\gos_1\gos_3}e^{(N-2)}_{\cels\gos_2}
+ \Theta{^\cels}_{\gos_2\gos_1}e^{(N-2)}_{\cels\gos_3}
\end{array}
\end{equation}
\end{lemm}
\emph{Proof} --- A computation gives
\[
\begin{array}{ccc}
 e^{\gos_4\gos_5}\wedge e^{(N-4)}_{\gos_1\gos_2\gos_3\gos}
 & = &
  \delta^{\gos_4\gos_5}_{\gos_3\gos}e^{(N-2)}_{\gos_1\gos_2} + \delta^{\gos_4\gos_5}_{\gos_1\gos}e^{(N-2)}_{\gos_2\gos_3} + \delta^{\gos_4\gos_5}_{\gos_2\gos}e^{(N-2)}_{\gos_3\gos_1} \\
& &  +\ \delta^{\gos_4\gos_5}_{\gos_3\gos_2}e^{(N-2)}_{\gos\gos_1}
  + \delta^{\gos_4\gos_5}_{\gos_1\gos_3}e^{(N-2)}_{\gos\gos_2}
  + \delta^{\gos_4\gos_5}_{\gos_2\gos_1}e^{(N-2)}_{\gos\gos_3}
\end{array}
\]
and we deduce (\ref{desss}) by developping  $\hbox{d}^\omega e^{(N-3)}_{\gos_1\gos_2\gos_3}
 = \hbox{d}^\omega e^\cels \wedge e^{(N-4)}_{\gos_1\gos_2\gos_3\cels} = \frac{1}{2}\Theta{^\cels}_{\cels_4\cels_5}e^{\cels_4\cels_5} \wedge e^{(N-4)}_{\gos_1\gos_2\gos_3\cels}$. \hfill $\square$\\

\noindent
\emph{Step 2: the proof of the first relation in (\ref{contraintedeBianchi})} --- We first compute the term $\partial_\cels^\omega \tilde{\mathbf{C}}{_\gol}^\cels + \Theta{^{\ast}}_{\cels\ast}
\tilde{\mathbf{C}}{_\gol}^\cels$. We start by observing that, by Lemma \ref{lemmedeNmoinsun},
\[
\left(\partial^\omega_\cels \tilde{\mathbf{C}}{_\gol}^\cels + \Theta{^\ast}_{\cels\ast}\tilde{\mathbf{C}}{_\gol}^\cels\right)e^{(N)}
= \hbox{d}^\omega \tilde{\mathbf{C}}{_\gol}^\cels \wedge e^{(N-1)}_\cels + \tilde{\mathbf{C}}{_\gol}^\cels \ \hbox{d}^\omega e^{(N-1)}_\cels
= \hbox{d}^\omega\left(\tilde{\mathbf{C}}{_\gol}^\cels\ e^{(N-1)}_\cels\right)
\]
This implies by using first (\ref{defiCtenseurCartanimplicit}), then (\ref{Bianchisimple}), that
\[
\begin{array}{ccl}
 \left(\partial^\omega_\cels \tilde{\mathbf{C}}{_\gol}^\cels + \Theta{^\ast}_{\cels\ast}\tilde{\mathbf{C}}{_\gol}^\cels\right)e^{(N)} & = &  - \frac{1}{2}\kappa{_\gol}^{\cels_1\cels_2}\hbox{d}^\omega\left(\Theta{^\cels} \wedge e^{(N-3)}_{\cels_1\cels_2\cels}\right) \\
& = &  - \frac{1}{2}\kappa{_\gol}^{\cels_1\cels_2} \left(\hbox{d}^\omega\Theta{^\cels} \wedge e^{(N-3)}_{\cels_1\cels_2\cels}
 +\Theta{^\cels} \wedge\hbox{d}^\omega e^{(N-3)}_{\cels_1\cels_2\cels}\right)\\
& = &  - \frac{1}{2}\kappa{_\gol}^{\cels_1\cels_2} \left(- \left[\theta^\gos \wedge \Omega^\gol\right]^\cels \wedge e^{(N-3)}_{\cels_1\cels_2\cels}
 +\Theta{^\cels} \wedge\hbox{d}^\omega e^{(N-3)}_{\cels_1\cels_2\cels}\right)
\end{array}
\]
However we get from (\ref{desss}) that
\[
 \begin{array}{ccl}
  \Theta{^{\cels_3}} \wedge\hbox{d}^\omega e^{(N-3)}_{\gos_1\gos_2\cels_3}
  & = & \Theta{^\ast}
_{\cels_3\ast}\Theta{^{\cels_3}}_{\gos_1\gos_2}e^{(N)}
+ \Theta{^\ast}_{\gos_1\ast}\Theta{^{\cels_3}}_{\gos_2\cels_3}e^{(N)}
+ \Theta{^\ast}_{\gos_2\ast}\Theta{^{\cels_3}}_{\cels_3\gos_1}e^{(N)} \\
& &
+\ \Theta{^\cels}_{\cels_3\gos_2}\Theta{^{\cels_3}}_{\cels\gos_1}e^{(N)}
+ \Theta{^\cels}_{\gos_1\cels_3}\Theta{^{\cels_3}}_{\cels\gos_2}e^{(N)}
+ \Theta{^\cels}_{\gos_2\gos_1}\Theta{^{\cels_3}}_{\cels\cels_3}e^{(N)} \\
  & = & \Theta{^\ast}
_{\cels\ast}\Theta{^{\cels}}_{\gos_1\gos_2}e^{(N)}
+ \Theta{^\ast}_{\gos_1\ast}\Theta{^{\ast}}_{\gos_2\ast}e^{(N)}
- \Theta{^\ast}_{\gos_2\ast}\Theta{^{\ast}}_{\gos_1\ast}e^{(N)} \\
& &
+\ \Theta{^\cels}_{\cels_3\gos_2}\Theta{^{\cels_3}}_{\cels\gos_1}e^{(N)}
+ \Theta{^\cels}_{\gos_1\cels_3}\Theta{^{\cels_3}}_{\cels\gos_2}e^{(N)}
+ \Theta{^\cels}_{\gos_2\gos_1}\Theta{^{\ast}}_{\cels\ast}e^{(N)} \\
& = & 0
 \end{array}
\]
(indeed the first term and the last term cancel together, the second and the third ones also, the fourth and the fifth ones also). It follows that
\[
 \begin{array}{ccl}
\left(\partial^\omega_\cels \tilde{\mathbf{C}}{_\gol}^\cels + \Theta{^\ast}_{\cels\ast}\tilde{\mathbf{C}}{_\gol}^\cels\right)e^{(N)}
& = &  \frac{1}{2}\kappa{_\gol}^{\cels_1\cels_2} \left[\theta^\gos \wedge \Omega^\gol\right]^\cels\wedge e^{(N-3)}_{\cels_1\cels_2\cels} \\
& = &
\frac{1}{2}\kappa{_\gol}^{\cels_1\cels_2}\mathbf{c}{^{\cels_0}}_{\cels\cell} e^\cels \wedge \Omega^\cell \wedge e^{(N-3)}_{\cels_1\cels_2\cels_0} \\
  & = &
 \frac{1}{2}\kappa{_\gol}^{\cels_1\cels_2} \Omega^\cell \wedge
 \left(\mathbf{c}{^{\cels}}_{\cels\cell} e^{(N-2)}_{\cels_1\cels_2}+ \mathbf{c}{^{\cels_0}}_{\cels_1\cell}e^{(N-2)}_{\cels_2\cels_0} + \mathbf{c}{^{\cels_0}}_{\cels_2\cell}e^{(N-2)}_{\cels_0\cels_1}\right)
\end{array}
\]
But since $\mathbf{c}{^{\cels}}_{\cels\gol} = 0$ this gives us
\[
 \begin{array}{ccl}
 \left(\partial^\omega_\cels \tilde{\mathbf{C}}{_\gol}^\cels + \Theta{^\ast}_{\cels\ast}\tilde{\mathbf{C}}{_\gol}^\cels\right)e^{(N)}
 & = & \frac{1}{2}\kappa{_\gol}^{\cels_1\cels_2}
 \mathbf{c}{^{\cels_0}}_{\cels_1\cell}\Omega^\cell \wedge e^{(N-2)}_{\cels_2\cels_0} + \frac{1}{2}\kappa{_\gol}^{\cels_1\cels_2}\mathbf{c}{^{\cels_0}}_{\cels_2\cell}\Omega^\cell \wedge e^{(N-2)}_{\cels_0\cels_1} \\
 & = &
 \frac{1}{2}\left(\mathbf{c}{^{\cels_0}}_{\cels_1\cell}\kappa{_\gol}^{\cels_1\cels_2}
 \Omega{^\cell}_{\cels_2\cels_0} + \mathbf{c}{^{\cels_0}}_{\cels_2\cell}\kappa{_\gol}^{\cels_1\cels_2}\Omega{^\cell}_{\cels_0\cels_1}\right) e^{(N)}
 \end{array}
\]
By exchanging indices $\cels_0\leftrightarrow\cels_1$ in the first term and indices $\cels_0\leftrightarrow\cels_2$ in the second term, we obtain that the two first terms in the left hand side of the first equation in (\ref{contraintedeBianchi}) are equal to
\begin{equation}\label{etapedCgsplusTheta}
 \partial^\omega_\cels \tilde{\mathbf{C}}{_\gol}^\cels + \Theta{^\ast}_{\cels\ast}\tilde{\mathbf{C}}{_\gol}^\cels
 =
 \frac{1}{2}\left(\mathbf{c}{^{\cels_1}}_{\cels_0\cell}\kappa{_\gol}^{\cels_0\cels_2}
 \Omega{^\cell}_{\cels_2\cels_1} + \mathbf{c}{^{\cels_2}}_{\cels_0\cell}\kappa{_\gol}^{\cels_1\cels_0}\Omega{^\cell}_{\cels_2\cels_1}\right)
\end{equation}
We now compute the term $\textbf{c}{^{\cels_0}}_{\gol\cels}\tilde{\mathbf{E}}{_{\cels_0}}^\cels
$. For that purpose we use (\ref{defiEtenseurEinsteinimplicit}):
\[
 \begin{array}{ccl}
  \textbf{c}{^{\cels_0}}_{\gol\cels}\tilde{\mathbf{E}}{_{\cels_0}}^\cels e^{(N)}
  & = & \textbf{c}{^{\cels_0}}_{\gol\cels_1} e^{\cels_1} \wedge \tilde{\mathbf{E}}{_{\cels_0}}^\cels e^{(N)}_\cels \\
  & = & \textbf{c}{^{\cels_0}}_{\gol\cels} e^{\cels} \wedge
  \left(-\frac{1}{2}\kappa{_\cell}^{\cels_1\cels_2}\Omega^\cell \wedge e^{(N-3)}_{\cels_1\cels_2\cels_0} \right) \\
  & = &
  -\frac{1}{2}\kappa{_\cell}^{\cels_1\cels_2}\Omega^\cell \wedge
  \left( \textbf{c}{^{\cels_0}}_{\gol\cels_0} e^{(N-2)}_{\cels_1\cels_2} + \textbf{c}{^{\cels_0}}_{\gol\cels_1} e^{(N-2)}_{\cels_2\cels_0} + \textbf{c}{^{\cels_0}}_{\gol\cels_2} e^{(N-2)}_{\cels_0\cels_1}\right)
 \end{array}
\]
and thus since $\textbf{c}{^{\cels_0}}_{\gol\cels_0} = 0$
\[
\begin{array}{ccl}
 \textbf{c}{^{\cels_0}}_{\gol\cels}\tilde{\mathbf{E}}{_{\cels_0}}^\cels e^{(N)}
 & = &
 -\frac{1}{2}\kappa{_\cell}^{\cels_1\cels_2}\Omega^\cell \wedge
\left( \textbf{c}{^{\cels_0}}_{\gol\cels_1} e^{(N-2)}_{\cels_2\cels_0} + \textbf{c}{^{\cels_0}}_{\gol\cels_2} e^{(N-2)}_{\cels_0\cels_1}\right) \\
& = & -\frac{1}{2}
\left(\textbf{c}{^{\cels_0}}_{\gol\cels_1}\kappa{_\cell}^{\cels_1\cels_2}\Omega{^\cell}_{\cels_2\cels_0} + \textbf{c}{^{\cels_0}}_{\gol\cels_2}\kappa{_\cell}^{\cels_1\cels_2}\Omega{^\cell}_{\cels_0\cels_1} \right)e^{(N)}
\end{array}
\]
By exchanging indices $\cels_0\leftrightarrow\cels_1$ in the first term and $\cels_0\leftrightarrow\cels_2$ in the second term, this gives us
\begin{equation}\label{etapecE}
 \textbf{c}{^{\cels_0}}_{\gol\cels}\tilde{\mathbf{E}}{_{\cels_0}}^\cels
 = -\frac{1}{2}
\left(\textbf{c}{^{\cels_1}}_{\gol\cels_0}\kappa{_\cell}^{\cels_0\cels_2}\Omega{^\cell}_{\cels_2\cels_1} + \textbf{c}{^{\cels_2}}_{\gol\cels_0}\kappa{_\cell}^{\cels_1\cels_0}\Omega{^\cell}_{\cels_2\cels_1} \right)
\end{equation}
Now by gathering (\ref{etapedCgsplusTheta}) and (\ref{etapecE}) we obtain
\begin{equation}\label{lesdeuxetapesCgsreunies}
\begin{array}{l}
 \partial^\omega_\cels \tilde{\mathbf{C}}{_\gol}^\cels + \Theta{^\ast}_{\cels\ast}\tilde{\mathbf{C}}{_\gol}^\cels
+ \textbf{c}{^{\cels_0}}_{\gol\cels}\tilde{\mathbf{E}}{_{\cels_0}}^\cels \\
= \frac{1}{2}\left(\mathbf{c}{^{\cels_1}}_{\cell\cels_0}\kappa{_\gol}^{\cels_0\cels_2}
+ \mathbf{c}{^{\cels_2}}_{\cell\cels_0}\kappa{_\gol}^{\cels_1\cels_0}\right)\Omega{^\cell}_{\cels_1\cels_2}
 +\frac{1}{2}
\left(\textbf{c}{^{\cels_1}}_{\gol\cels_0}\kappa{_\cell}^{\cels_0\cels_2} + \textbf{c}{^{\cels_2}}_{\gol\cels_0}\kappa{_\cell}^{\cels_1\cels_0} \right)\Omega{^\cell}_{\cels_1\cels_2}
\end{array}
\end{equation}
To conclude we use the fact that $\kappa{_\gol}^{\gos\gos}$ is invariant by the adjoint action of $\gol$, i.e. $\hbox{ad}_{\gol_1}\kappa_{\gol_2}^{\gos_1\gos_2} =0$. Since $\hbox{ad}_{\gol_1}\kappa_{\gol_2}^{\gos_1\gos_2} =
- \mathbf{c}{^{\celg}}_{\gol_1\gol_2}\kappa{_{\cell}}^{\gos_1\gos_2} +
\mathbf{c}{^{\gos_1}}_{\gol_1\cels_0}\kappa{_{\gol_2}}^{\cels_0\gos_2} + \mathbf{c}{^{\gos_2}}_{\gol_1\cels_0}\kappa{_{\gol_2}}^{\gos_1\cels_0}$,
this implies
\begin{equation}
\mathbf{c}{^{\gos_1}}_{\gol_1\cels_0}\kappa{_{\gol_2}}^{\cels_0\gos_2} + \mathbf{c}{^{\gos_2}}_{\gol_1\cels_0}\kappa{_{\gol_2}}^{\gos_1\cels_0} = \mathbf{c}{^{\cell_0}}_{\gol_1\gol_2}\kappa{_{\cell_0}}^{\gos_1\gos_2}
\end{equation}
By applying this relation for $(\gol_1,\gol_2) = (\cell,\gol)$ in the first term of the right hand side of (\ref{lesdeuxetapesCgsreunies}) and for $(\gol_1,\gol_2) = (\gol,\cell)$ in the second term of the right hand side, this gives us
\[
 \partial^\omega_\cels \tilde{\mathbf{C}}{_\gol}^\cels + \Theta{^\ast}_{\cels\ast}\tilde{\mathbf{C}}{_\gol}^\cels
 + \textbf{c}{^{\cels_0}}_{\gol\cels}\tilde{\mathbf{E}}{_{\cels_0}}^\cels = \frac{1}{2}\left(\mathbf{c}{^{\cell_0}}_{\cell\gol}\kappa{_{\cell_0}}^{\cels_1\cels_2} + \mathbf{c}{^{\cell_0}}_{\gol\cell}\kappa{_{\cell_0}}^{\cels_1\cels_2}\right)\Omega{^\cell}_{\cels_1\cels_2} = 0
\]

\noindent
\emph{Step 3: the proof of the second relation in (\ref{contraintedeBianchi})} --- It amounts to show that
\begin{equation}\label{contraintesurtenseurE}
  \partial_\cels^\omega \tilde{\mathbf{E}}{_\gos}^\cels  + \Theta{^{\ast}}_{\cels\ast}
\tilde{\mathbf{E}}{_\gos}^\cels
=
\Theta{^{\cels_0}}_{\gos\cels}\tilde{\mathbf{E}}{_{\cels_0}}^{\cels}
+ \Omega{^{\cell}}_{\gos\cels}
\tilde{\mathbf{C}}{_{\cell}}^{\cels}
\end{equation}
On the one hand by using first (\ref{domegaesparlemme}), then (\ref{defiEtenseurEinsteinimplicit}) we get
\[
 \begin{array}{ccl}
  \left(\partial_\cels^\omega \tilde{\mathbf{E}}{_\gos}^\cels + \Theta{^\ast}_{\cels\ast}\tilde{\mathbf{E}}{_\gos}^\cels\right) e^{(N)}
  & = & \hbox{d}^\omega \left(\mathbf{E}{_\gos}^\cels e^{(N-1)}_\cels\right)
  = \hbox{d}^\omega \left(- \frac{1}{2}\kappa{_\cell}^{\cels_1\cels_2} \Omega^\cell \wedge e^{(N-3)}_{\cels_1\cels_1\gos}\right) \\
  & = &
  - \frac{1}{2}\kappa{_\cell}^{\cels_1\cels_2} \left(
  \hbox{d}^\omega\Omega^\cell \wedge e^{(N-3)} + \Omega^\cell \wedge\hbox{d}^\omega e^{(N-3)}_{\cels_1\cels_1\gos}\right)
 \end{array}
\]
and since $\hbox{d}^\omega\Omega^\gol = 0$ by (\ref{Bianchisimple}) we deduce by using (\ref{desss}) that
\begin{equation}\label{contraintesurtenseurEgauche}
  \begin{array}{ccc}
\partial_\cels^\omega \tilde{\mathbf{E}}{_\gos}^\cels + \Theta{^\ast}_{\cels\ast}\tilde{\mathbf{E}}{_\gos}^\cels
& = &
- \frac{1}{2}\kappa{_\cell}^{\cels_1\cels_2} \left(
\Theta{^\ast}
_{\gos\ast}\Omega{^\cell}_{\cels_1\cels_2}
+ \Theta{^\ast}_{\cels_1\ast}\Omega{^\cell}_{\cels_2\gos}
+ \Theta{^\ast}_{\cels_2\ast}\Omega{^\cell}_{\gos\cels_1} \right. \\
& &
+ \left. \Theta{^\cels}_{\gos\cels_2}\Omega{^\cell}_{\cels\cels_1}
+ \Theta{^\cels}_{\cels_1\gos}\Omega{^\cell}_{\cels\cels_2}
+ \Theta{^\cels}_{\cels_2\cels_1}\Omega{^\cell}_{\cels\gos}
\right)
 \end{array}
\end{equation}
On the other hand
\[
  \left(\Theta{^{\cels_0}}_{\gos\cels}\tilde{\mathbf{E}}{_{\cels_0}}^{\cels}
+ \Omega{^{\cell}}_{\gos\cels}
\tilde{\mathbf{C}}{_{\cell}}^{\cels}\right) e^{(N)}
= \Theta{^{\cels_0}}_{\gos\cels_3}e^{\cels_3}\wedge \tilde{\mathbf{E}}{_{\cels_0}}^{\cels}e^{(N-1)}_\cels
+ \Omega{^{\cell}}_{\gos\cels_3}e^{\cels_3}\wedge
\tilde{\mathbf{C}}{_{\cell}}^{\cels}e^{(N-1)}_\cels
\]
and thus by using (\ref{defiCtenseurCartanimplicit}) and (\ref{defiEtenseurEinsteinimplicit})
\[
 \begin{array}{l}
  \left(\Theta{^{\cels_0}}_{\gos\cels}\tilde{\mathbf{E}}{_{\cels_0}}^{\cels}
+ \Omega{^{\cell}}_{\gos\cels}
\tilde{\mathbf{C}}{_{\cell}}^{\cels}\right) e^{(N)} \\
\quad = \Theta{^{\cels_0}}_{\gos\cels_3}e^{\cels_3}\wedge \left(-\frac{1}{2}\kappa{_\cell}^{\cels_1\cels_2}\Omega^\cell \wedge e^{(N-3)}_{\cels_1\cels_2\cels_0} \right) 
+ \Omega{^{\cell}}_{\gos\cels_3}e^{\cels_3}\wedge
\left(-\frac{1}{2}\kappa{_\cell}^{\cels_1\cels_2}\Theta^{\cels_0} \wedge e^{(N-3)}_{\cels_1\cels_2\cels_0} \right) \\
\quad =
-\frac{1}{2}\kappa{_\cell}^{\cels_1\cels_2}
\Omega^\cell \wedge \left(\Theta{^{\cels_0}}_{\gos\cels_0}e^{(N-2)}_{\cels_1\cels_2} + \Theta{^{\cels_0}}_{\gos\cels_2}e^{(N-2)}_{\cels_0\cels_1} +\Theta{^{\cels_0}}_{\gos\cels_1}e^{(N-2)}_{\cels_2\cels_0} \right) \\
\quad  \quad \quad \quad  \quad
-\frac{1}{2}\kappa{_\cell}^{\cels_1\cels_2}\Theta^{\cels_0} \wedge
\left(\Omega{^\cell}_{\gos\cels_0}e^{(N-2)}_{\cels_1\cels_2} + \Omega{^\cell}_{\gos\cels_2}e^{(N-2)}_{\cels_0\cels_1} + \Omega{^\cell}_{\gos\cels_1}e^{(N-2)}_{\cels_2\cels_0} \right) \\
\quad =
-\frac{1}{2}\kappa{_\cell}^{\cels_1\cels_2} \left(\Theta{^{\cels_0}}_{\gos\cels_0}\Omega{^\cell}_{\cels_1\cels_2} + \Theta{^{\cels_0}}_{\gos\cels_2}\Omega{^\cell}_{\cels_0\cels_1} +\Theta{^{\cels_0}}_{\gos\cels_1}\Omega{^\cell}_{\cels_2\cels_0} \right. \\
\quad   \quad \quad \quad  \quad
+ \left. \Theta{^{\cels_0}}_{\cels_1\cels_2}\Omega{^\cell}_{\gos\cels_0} + \Theta{^{\cels_0}}_{\cels_0\cels_1}\Omega{^\cell}_{\gos\cels_2} + \Theta{^{\cels_0}}_{\cels_2\cels_0}\Omega{^\cell}_{\gos\cels_1} \right)e^{(N)}
 \end{array}
\]
Hence
\begin{equation}\label{contraintesurtenseurEdroite}
\begin{array}{ccl}
 \Theta{^{\cels_0}}_{\gos\cels}\tilde{\mathbf{E}}{_{\cels_0}}^{\cels}
+ \Omega{^{\cell}}_{\gos\cels}
\tilde{\mathbf{C}}{_{\cell}}^{\cels}
& = & -\frac{1}{2}\kappa{_\cell}^{\cels_1\cels_2} \left(\Theta{^{\ast}}_{\gos\ast}\Omega{^\cell}_{\cels_1\cels_2} + \Theta{^{\cels}}_{\gos\cels_2}\Omega{^\cell}_{\cels\cels_1} +\Theta{^{\cels}}_{\gos\cels_1}\Omega{^\cell}_{\cels_2\cels} \right. \\
& &  \quad 
+ \left. \Theta{^{\cels}}_{\cels_1\cels_2}\Omega{^\cell}_{\gos\cels} + \Theta{^{\ast}}_{\ast\cels_1}\Omega{^\cell}_{\cels_2\gos} + \Theta{^{\ast}}_{\cels_2\ast}\Omega{^\cell}_{\gos\cels_1} \right)
\end{array}
\end{equation}
By comparing (\ref{contraintesurtenseurEgauche}) and (\ref{contraintesurtenseurEdroite}) we conclude that (\ref{contraintesurtenseurE}) is satisfied (in the right hand side of  (\ref{contraintesurtenseurEdroite}) the term ranked 1, 2, 3, 4, 5, 6 coincides with, respectively, the term ranked 1, 4, 5, 6, 2, 3 in the right hand side of  (\ref{contraintesurtenseurEgauche})). \hfill $\square$

\section{Applications}
\subsection{Gravity with a cosmological constant}\label{soussectioncosmo}
The most natural theory is obtained by choosing $\widehat{\goP}$ to be the universal cover of the Poincar{\'e} group
$Spin_0(1,n-1)\ltimes \R^n$ for $n\geq 2$. However it is also interesting to consider their deformations
$Spin_0(2,n-1)$ and $Spin_0(1,n)$. Since our description is local it can be given by using their quotients $SO_0(1,n-1)\ltimes \R^n$,
$SO_0(2,n-1)$ and $SO_0(1,n)$, respectively.
These Lie groups can be represented as subgroups of the matrix group $GL(n+1,\R)$ as follows. We define
$\underline{\textsf{h}}:= \left(\begin{array}{ccc}
   \textsf{h}_{11} & \cdots & \textsf{h}_{1n} \\
   \vdots &  & \vdots \\
   \textsf{h}_{n1} & \cdots & \textsf{h}_{nn}
          \end{array}\right)$ and $\overline{\textsf{h}} = (\underline{\textsf{h}})^{-1} =
\left(\begin{array}{ccc}
   \textsf{h}^{11} & \cdots & \textsf{h}^{1n} \\
\vdots &  & \vdots \\
\textsf{h}^{n1} & \cdots & \textsf{h}^{nn}
\end{array}\right)$
(a Minkowski metric on $\R^n$) and
$\overline{\textsf{H}}:=
\left(\begin{array}{cc}
\overline{\textsf{h}}  & 0 \\ 0 & k
  \end{array}\right)$, where $k\in \R$ (a metric on $\R^{n+1}$).
We let
\begin{equation}
 \goP_k(n):= \left\{G\in GL(n+1,\R);\;
 G\overline{\textsf{H}}G^t = \overline{\textsf{H}},
 \hbox{det}G=1\right\}.
\end{equation}
Assuming that the signature of $\overline{\textsf{h}}$ is
$(-,+,\cdots,+)$, we have the following identifications
\begin{itemize}
 \item if $k<0$, $\goP_k(n)$ is isomorphic to $SO(1,n)$;
 \item if $k=0$, $\goP_0(n)$ is isomorphic to the Poincar{\'e} group $\goP(n) = SO(1,n-1)\ltimes \R^n$;
 \item if $k>0$, $\goP_k(n)$ is isomorphic to $SO(2,n-1)$.
\end{itemize}
In each case we get a theory of gravitation with a cosmological
constant $\Lambda = \frac{n(n-1)k}{2}$.
The representation of the Lorentz subgroup $\goL_k(n)$ ($\simeq SO(1,n-1)$) is
\[
 \goL_k(n):= \left\{G = \left(\begin{array}{cc}
      g & 0 \\ 0 & 1
  \end{array}\right);\;
  g\in GL(n,\R),
 g\overline{\textsf{h}}g^t = \overline{\textsf{h}},
 \hbox{det}g=1\right\}.
\]
For $\left(\begin{array}{cc} g & 0 \\ 0 & 1\end{array}\right)\in \goL_k(n)$, we deduce from $g\overline{\textsf{h}}g^t = \overline{\textsf{h}}$ the following useful relations:
\begin{equation}\label{usefulA0}
 g\overline{\textsf{h}} = \overline{\textsf{h}}(g^{-1})^t
 \quad \hbox{and} \quad
 \underline{\textsf{h}}g^{-1} = g^t\underline{\textsf{h}}
\end{equation}

\subsubsection{Lie algebras}
The Lie algebra of $\goP_k(n)$ is
$\gop_k(n):= \left\{\xi\in M(n+1,\R);\;
 \xi\overline{\textsf{H}} + \overline{\textsf{H}}\xi^t =0
 \right\}$. Any element $\xi\in \gop_k(n)$ can be written
 \[
  \xi = \left(\begin{array}{cccc}
    \xi{^1}_1 & \cdots & \xi{^1}_n & \xi{^1}_{n+1} \\
   \vdots &  & \vdots & \vdots \\
   \xi{^n}_1 & \cdots & \xi{^n}_n & \xi{^n}_{n+1} \\
   \xi{^{n+1}}_1 & \cdots & \xi{^{n+1}}_n & 0
 \end{array}\right)
 = \left(\begin{array}{cccc}
    \xi^{1b}\textsf{h}_{b1} & \cdots & \xi^{1b}\textsf{h}_{bn} & \xi^1 \\
   \vdots &  & \vdots & \vdots \\
   \xi^{nb}\textsf{h}_{b1} & \cdots & \xi^{nb}\textsf{h}_{bn} & \xi^{n} \\
   -k\xi^b\textsf{h}_{b1} & \cdots & -k\xi^b\textsf{h}_{bn} & 0
 \end{array}\right)
 \]
 where $(\xi^{ab})_{1\leq a,b\leq n}$ and $(\xi^a)_{1\leq a\leq n}$ are real coefficients such that $\xi^{ab} + \xi^{ba} = 0$. Clearly there exists a unique family of matrices $\left(\textbf{t}_A\right)_{1\leq A\leq n(n+1)/2} =
\left((\textbf{t}_a)_{1\leq a\leq n},(\textbf{t}_{ab})_{1\leq a<b\leq n}\right)$ in
$\gop_k(n)$ such that, $\forall \xi\in \gop_k(n)$, $\xi=
 \sum_{1\leq a<b\leq n}\textbf{t}_{ab}\xi^{ab}
 + \sum_{1\leq a\leq n}\textbf{t}_a\xi^a$. Obviously this family forms a base of $\gop_k(n)$. It is convenient to define
$\textbf{t}_{ba}:= -\textbf{t}_{ab}$, for
$1\leq a\leq b\leq n$, and to write
\[
 \xi = \frac{1}{2}
 \sum_{1\leq a,b\leq n}\textbf{t}_{ab}\xi^{ab}
 + \sum_{1\leq a\leq n}\textbf{t}_a\xi^a
 = \frac{1}{2}\textbf{t}_{ab}\xi^{ab}
 + \textbf{t}_a\xi^a
\]
The Lie algebra of $\goL_k(n)$ is simply
$\gol_k(n):= \left\{\xi =\frac{1}{2}\textbf{t}_{ab}\xi^{ab};\
\xi^{ab}\in\R,
\xi^{ab} +\xi^{ba} = 0\right\}$.
and we have $\gop_k(n) = \gol_k(n) \oplus \mathfrak{s}_k(n)$, with
$\mathfrak{s}_k(n):=
 \left\{\xi=\textbf{t}_a\xi^a;\
 \xi^a\in \R\right\}$.
The Lie bracket in this basis reads
\[
\left[\begin{array}{cc}
       {[\textbf{t}_{ab},\textbf{t}_{cd}]} &
       {[\textbf{t}_{ab},\textbf{t}_c]}
       \\
       {[\textbf{t}_{a},\textbf{t}_{cd}]}
       &
        {[\textbf{t}_a,\textbf{t}_c]}
      \end{array}\right]
= \left[\begin{array}{cc}
       \textsf{h}_{bc}\textbf{t}_{ad}
    - \textsf{h}_{bd}\textbf{t}_{ac}
    - \textsf{h}_{ac}\textbf{t}_{bd}
    + \textsf{h}_{ad}\textbf{t}_{bc}
    & \quad \textsf{h}_{bc}\textbf{t}_a
 - \textsf{h}_{ac}\textbf{t}_b \quad
 \\
 \textsf{h}_{ac}\textbf{t}_d
 - \textsf{h}_{ad}\textbf{t}_c
 & -k \textbf{t}_{ac}
      \end{array}\right]
\]
Equivalentely the structure coefficients $\textbf{c}^I_{JK} =
\langle\textbf{t}^I,[\textbf{t}_B,\textbf{t}_K]\rangle$ of the Lie algebra $\gop_k(n)$ in the chosen basis
are given by
\[
 \begin{array}{|ccc|ccc|}
 \hline
\left(\begin{array}{c}
\textbf{c}^{[ef]}_{[ab][cd]} \\
\textbf{c}^{e}_{[ab][cd]}
\end{array}\right)
& = &
\left(\begin{array}{c}
\delta^{ef}_{ad}\textsf{h}_{bc}
    - \delta^{ef}_{ac}\textsf{h}_{bd}
    - \delta^{ef}_{bd}\textsf{h}_{ac}
    + \delta^{ef}_{bc}\textsf{h}_{ad} \\
0
\end{array}\right) 
& \left(\begin{array}{c}
\textbf{c}^{[ef]}_{[ab]c} \\
\textbf{c}^{e}_{[ab]c}
\end{array}\right)
& = & 
\left(\begin{array}{c}
0 \\
 \delta^{e}_a\textsf{h}_{bc} - \delta^{e}_b\textsf{h}_{ac}
\end{array}\right) 
\\
\hline 
\left(\begin{array}{c}
\textbf{c}^{[ef]}_{a\;\;[cd]} \\ \textbf{c}^{e}_{a[cd]}
\end{array}\right)
& = &
\left(\begin{array}{c}
0 \\ \delta^{e}_d\textsf{h}_{ac} - \delta^{e}_c\textsf{h}_{ad}
\end{array}\right)
&
\left(\begin{array}{c}
\textbf{c}^{[ef]}_{\;a\;\;\;c} \\ \textbf{c}^{e}_{ac}
\end{array}\right)
& = &
\left(\begin{array}{c}
-k \delta^{ef}_{ac} \\ 0
\end{array}\right) \\ 
\hline
 \end{array}
\]
where $\delta^{ef}_{ab}:=
\delta^{e}_a\delta^{f}_b - \delta^{e}_b\delta^{f}_a$.

The adjoint action of an element $ g\in\goL_k(n)$ on $\xi\in \gop_k(n)$ reads $\hbox{Ad}_g\left(\frac{1}{2}\textbf{t}_{ab}\xi^{ab}
 + \textbf{t}_a\xi^a\right)
 = \frac{1}{2}\textbf{t}_{ab}g^a_{a'}g^b_{b'}\xi^{a'b'}
 + \textbf{t}_ag^a_{a'}\xi^{a'}$
and the coadjoint action of $g\in \goL_k(n)$ on $\alpha  = \frac{1}{2}\alpha_{ab}\textbf{t}^{ab}
 + \alpha_a\textbf{t}^a\in \gop_k^*(n)$
 expresses as
$\hbox{Ad}_g^*\left(\frac{1}{2}\alpha_{ab}\textbf{t}^{ab}
 + \alpha_a\textbf{t}^a\right)
 = \frac{1}{2}\alpha_{a'b'}(g^{-1})^{a'}_a(g^{-1})^{b'}_b\textbf{t}^{ab}
 + \alpha_{a'}(g^{-1})^{a'}_a\textbf{t}^a$.

\subsubsection{Checking the hypotheses}
Hypothesis (i), that $\goP_k(n)$ is unimodular, can be checked by a direct computation:
on the one hand, for any $a,b$,
\[
\begin{array}{ccl}
 \frac{1}{2}c^{[cd]}_{[ab][cd]} + c^{c}_{[ab]c} & = & \frac{1}{2}\left(
 \delta^{cd}_{ad}\textsf{h}_{bc} - \delta^{cd}_{ac}\textsf{h}_{bd}
 - \delta^{cd}_{bd}\textsf{h}_{ac} + \delta^{cd}_{bc}\textsf{h}_{ad}\right)
 + \delta^c_a\textsf{h}_{bc} - \delta^c_b\textsf{h}_{ac} \\
 & = & \frac{1}{2}\left((n-1)\textsf{h}_{ba} - (1-n)\textsf{h}_{ba}
 - (n-1)\textsf{h}_{ab} + (1-n)\textsf{h}_{ab}\right)
 + \textsf{h}_{ba} - \textsf{h}_{ab} = 0
\end{array}
\]
and, on the other hand, for any $a$, we have obviously
$\frac{1}{2}c^{[cd]}_{a[cd]} + c^{c}_{ac} = 0 + 0 = 0$.

Hypothesis (ii), that $\hbox{Ad}_{\goL}\mathfrak{s} \subset \mathfrak{s}$ and $[\mathfrak{s},\mathfrak{s}]\subset \gol$, is straightforward.
We choose 
\[
\left(\begin{array}{cc}\kappa{_{[cd]}}^{ab} &
\kappa{_c}^{ab} \end{array}\right)
:= \left(\begin{array}{cc} \delta^{ab}_{cd} & 0 \end{array}\right)
\quad \Longleftrightarrow \quad
\kappa{_\gop}^{\gos\gos} := \frac{1}{2}\textbf{t}^{ab} \otimes
 \textbf{t}_a\wedge \textbf{t}_b
\]
We can check easily that $\hbox{Ad}_g\kappa{_\gop}^{\gos\gos} = \kappa{_\gop}^{\gos\gos}$, $\forall g\in \goL_k(n)$, i.e. that
Hypothesis (iii) is satisfied. Lastly $\kappa{_\gop}^{\gos\gos}$ satisfies obviously Hypothesis
(\ref{additionalhypothesis}), i.e. that $\kappa{_\gop}^{\gos\gos} \in \gol^*\otimes\Lambda^2\mathfrak{s}$.
Hence Theorem \ref{theoBigOne} can be applied: any smooth critical point $(\varphi^\gop,\pi_\gop)\in \mathscr{E}_\textsf{E}$ of $\mathscr{A}$ (given by (\ref{functionalAEinstein}) gives rise locally to a solution of the system (\ref{bigECsystem2splitsimple}).

\subsubsection{The equations of dynamics}
Assume $n>2$.
Let $(\pi_\gop,\varphi^\gop)\in \mathscr{E}_\textsf{E}$ be a critical point of the action (\ref{functionalAEinstein}) $\int_\mathcal{F}\pi_\celp\wedge \varphi^\celp$ and assume that it satisfies the \emph{Fibration hypothesis} (\ref{fibrationhypothesis}).
Then the manifold $\mathcal{F}$ is
fibered over an $n$-dimensional manifold $\mathcal{X}$ and $\mathcal{X}$
is equipped with a metric $\textbf{g}$ the pull-back by $\mathcal{F}\longrightarrow\mathcal{X}$ of which is $\textsf{h}_{ab}e^a\otimes e^b$ and $T\mathcal{X}$
is endowed with a metric connection $\nabla$ defined by $\omega$.

Let us assume furthermore that either (i) or (ii) in \S \ref{sectionexploitation} holds. Then the fields $(\pi_\gop,\varphi^\gop)$ give rise to a solution of the generalized Einstein--Cartan system in vacuum with a cosmological constant $\Lambda:= -\frac{1}{2}\mathbf{c}{^\celg}_{\cels\cels}\kappa{_\celg}^{\cels\cels}$
i.e. System (\ref{brutbeautiful}).

Since $\kappa{_\gog}^{\gos\gos}$ is given by $\kappa{_{[ab]}}^{cd} = \delta_{ab}^{cd}$  the first
equation in System (\ref{brutbeautiful}) is obviously equivalent to $\mathring{\Theta}{^\gos}_{\gos\gos} = 0$, which, since $n>2$, is itself equivalent to $\Theta{^\gos}_{\gos\gos} = 0$ as seen in \S \ref{sectionphysicaleqgrav}. This means that the connection $\nabla$ is torsion free, i.e. that it
is the Levi-Civita connection for the metric $\textbf{g}$.

The second equation in System (\ref{brutbeautiful}) reads
$\Omega{^\celg}_{\cels\gos}\kappa{_\celg}^{\cels\gos} - \frac{1}{2}\left(\Omega{^\celg}_{\cels\gos}\kappa{_\celg}^{\cels\gos}\right)\delta{_\gos}^\gos + \Lambda \delta{_\gos}^\gos = 0$. The computation in terms of the standard Riemann and Ricci tensors $\mathbf{R}{^{\gos\gos}}_{\gos\gos}$ and $\mathbf{R}{^\gos}_\gos$ is straightforward:
\[
 \Omega{^\celg}_{\cels a}\kappa{_\celg}^{\cels b} = \frac{1}{2}\Omega{^{[cd]}}_{ae}\kappa{_{[cd]}}^{be}
 = \frac{1}{2}\Omega{^{[cd]}}_{ae}\delta_{cd}^{be} = \Omega{^{[be]}}_{ae}
 = \mathbf{R}{^{be}}_{ae} = \mathbf{R}{^b}_a = \mathbf{R}{_a}^b
\]
(we use the symmetry of the Ricci tensor). Hence $\Omega{^\celg}_{\cels\gos}\kappa{_\celg}^{\cels\gos} = \mathbf{R}{^a}_a = \mathbf{R}$ is the scalar curvature. We also have
\[
 \Lambda = -\frac{1}{2}\mathbf{c}{^\celg}_{\cels\cels}\kappa{_\celg}^{\cels\cels}
 = - \frac{1}{4}\left(-k\delta_{ab}^{cd}\right)\delta^{ab}_{cd} = \frac{n(n-1)}{2}k
\]
Thus we obtain that $\tilde{\mathbf{E}}{_\gos}^\gos = \mathbf{E}{_\gos}^\gos$, so that the second equation in (\ref{brutbeautiful}) is exactly the Einstein equation
\begin{equation}\label{Einsteinequationcosmo}
 \mathbf{E}{_a}^b + \Lambda \delta{_a}^b = 0,
\end{equation}
with $\mathbf{E}{_a}^b:= \mathbf{R}{_a}^b - \frac{1}{2}\mathbf{R}\delta{_a}^b$ and  the \emph{cosmological constant} $\Lambda = \frac{n(n-1)}{2}k$.

\subsection{Gravity with a Barbero--Immirzi parameter}
This example is a variant of the previous case for $n=4$.
We use the groups $\goP_k(4)$ and $\goL_k(4)$.
Hence Hypotheses (i), (ii) and (iii) have been already checked. However the tensor
$\kappa$ is now
\begin{equation}\label{4KappaStructureCste}
\left(\begin{array}{cc}\kappa{_{[cd]}}^{ab} &
         \kappa{_c}^{ab} \end{array}\right)
     = \left(\begin{array}{cc} \delta^{ab}_{cd} - \frac{1}{\gamma}\epsilon{^{ab}}_{cd} & 0 \end{array}\right)
\end{equation}
where $\epsilon{^{ab}}_{cd}:= \epsilon_{a'b'cd}\textsf{h}^{a'a}\textsf{h}^{b'b}$
and $\epsilon_{abcd}$ is the completely antisymmetric tensor such that
$\epsilon_{1234}=1$ and where $\gamma\in \C^*$ is a constant (the \emph{Barbero--Immirzi parameter}).
Alternatively
\[
 \kappa{_\gop}^{\gos\gos}
 = \frac{1}{2}\textbf{t}^{ab}\otimes \textbf{t}_a\wedge \textbf{t}_b
 - \frac{1}{4\gamma}\textsf{h}^{aa'}\textsf{h}^{bb'}\epsilon_{a'b'cd}\textbf{t}^{[cd]}\otimes
 \textbf{t}_a\wedge \textbf{t}_b
 = \frac{1}{2}\textbf{t}^{ab} \otimes
 \textbf{t}_a\wedge \textbf{t}_b
 -\frac{1}{\gamma}\eta{_\gop}^{\gos\gos},
\]
where $\eta{_\gop}^{\gos\gos}:= \frac{1}{4}\textsf{h}^{aa'}\textsf{h}^{bb'}\epsilon_{a'b'cd}\textbf{t}^{[cd]} \otimes
\textbf{t}_a\wedge \textbf{t}_b$.

Hypothesis (\ref{additionalhypothesis}) is obviously satisfied.
In order to check that $\kappa{_\gop}^{\gos\gos}$ defined by (\ref{4KappaStructureCste}) is invariant by the adjoint action of $\goL_k(4)$, it suffices to check that $\eta{_\gop}^{\gos\gos}:= \frac{1}{4}\textsf{h}^{aa'}\textsf{h}^{bb'}\epsilon_{a'b'cd}\textbf{t}^{[cd]} \otimes
\textbf{t}_a\wedge \textbf{t}_b$ is so.
Using (\ref{usefulA0}) we get
\[
 \begin{array}{ccl}
  \hbox{Ad}_g^*\otimes\hbox{Ad}_g\otimes\hbox{Ad}_g\eta{_\gop}^{\gos\gos}
  & = & \hbox{Ad}_g^*\otimes\hbox{Ad}_g\otimes\hbox{Ad}_g\left(\frac{1}{4}\textsf{h}^{aa''}\textsf{h}^{bb''}\epsilon_{a''b''cd}\textbf{t}^{[cd]}\otimes
  \textbf{t}_a\wedge \textbf{t}_b
 \right) \\
  & = & \frac{1}{4}
 g^a_{a'}\textsf{h}^{a'a''}
 g^b_{b'}\textsf{h}^{b'b''}\epsilon_{a''b''c'd'}
 (g^{-1})^{c'}_c(g^{-1})^{d'}_d\textbf{t}^{[cd]}
 \otimes\textbf{t}_a\wedge \textbf{t}_b \\
 & = & \frac{1}{4}
 \textsf{h}^{aa'}\textsf{h}^{bb'}(g^{-1})^{a''}_{a'}(g^{-1})^{b''}_{b'}\epsilon_{a''b''c'd'}
 (g^{-1})^{c'}_c(g^{-1})^{d'}_d\textbf{t}^{[cd]}\otimes
 \textbf{t}_a\wedge \textbf{t}_b \\
 & = & \frac{1}{4}
 \textsf{h}^{aa'}\textsf{h}^{bb'} \hbox{det}(g^{-1})\epsilon_{a'b'cd}
 \textbf{t}^{[cd]}\otimes
 \textbf{t}_a\wedge \textbf{t}_b\\
 & = & \frac{1}{4}
 \epsilon{^{ab}}_{cd}\textbf{t}^{[cd]}\otimes\textbf{t}_a\wedge \textbf{t}_b
 = \eta{_\gop}^{\gos\gos}
 \end{array}
\]
where we have used $\hbox{det}(g^{-1}) = 1$.

Although the action takes complex values, this does not change the derivation of the Euler--Lagrange equations and, in particular, our conclusion about the local fibration of $\mathcal{F}$ over some 4-dimensional manifold $\mathcal{X}$. Thus if assume the \emph{Fibration hypothesis} (\ref{fibrationhypothesis}) and one of the two hypotheses (i) and (ii) in \S \ref{sectionexploitation}, we  get Equations (\ref{brutbeautiful}).

Let us prove that the first equation (c), i.e. $\frac{1}{2}\kappa{_\gol}^{\cels\cels}\mathring{\Theta}{^\gos}_{\cels\cels} = 0$, implies that $ \mathring{\Theta}{^\gos}_{\gos\gos} = 0$ (which is equivalent to the fact that the connection $\nabla$
is torsion free). The proof relies on two different arguments, according to the value of $\gamma$:
\begin{itemize}
\item if $\gamma = \pm i$ the condition $\frac{1}{2}\kappa{_\gol}^{\cels\cels}\mathring{\Theta}{^\gos}_{\cels\cels} = 0$ reads also
$\mathring{\Theta}{^c}_{ab} =
\frac{1}{2\gamma}\mathring{\Theta}{^c}_{a'b'} \epsilon{^{a'b'}}_{ab}$, which implies straightforwardly that $\mathring{\Theta}{^\gos}_{\gos\gos} = 0$ since this quantity is real (this case corresponds to the
Ashtekar action).
 \item in general, if
$\overline{\textsf{h}}$ is a Minkowski metric which is suitably normalized, the condition $\frac{1}{2}\kappa{_\gol}^{\cels\cels}\mathring{\Theta}{^\gos}_{\cels\cels} = 0$ is equivalent to
$\left(1+\frac{1}{\gamma^2}\right)
\mathring{\Theta}{^\gos}_{\gos\gos} = 0$ (see Lemma \ref{lemmeBarbero} below). This implies  $\mathring{\Theta}{^\gos}_{\gos\gos} = 0$ if $\gamma\neq \pm i$.
\end{itemize}
\begin{lemm}\label{lemmeBarbero} Assume that the metric $\textsf{\emph{h}}$ is Minkowski\footnote{In an Euclidean theory where we would assume
that the metric $\overline{\textsf{h}}$ has the signature $(+,\cdots,+)$, the
natural normalization would be $\hbox{det}\overline{\textsf{h}} = 1$, leading to
the relation $\left(1-\frac{1}{\gamma^2}\right)\mathring{\Theta}{^c}_{ab} = 0$.} and that $\hbox{\emph{det}}\,\overline{\textsf{\emph{h}}} = - 1$ Then the condition $\frac{1}{2}\kappa{_\gol}^{\cels\cels}\mathring{\Theta}{^\gos}_{\cels\cels} = 0$ implies $\left(1+\frac{1}{\gamma^2}\right)
\mathring{\Theta}{^c}_{ab} = 0$.
\end{lemm}
\emph{Proof} --- Condition $\frac{1}{2}\kappa{_\gol}^{\cels\cels}\mathring{\Theta}{^\gos}_{\cels\cels} = 0$ is equivalent to
$\mathring{\Theta}{^c}_{ab} = \frac{1}{2\gamma}\mathring{\Theta}{^c}_{a'b'}\epsilon{^{a'b'}}_{ab}$. By iterating this relation we
obtain $\mathring{\Theta}{^c}_{ab} = \frac{1}{4\gamma^2}\mathring{\Theta}{^c}_{a''b''}\epsilon{^{a''b''}}_{a'b'}\epsilon{^{a'b'}}_{ab}$.
But since
\[
\begin{array}{ccl}
\epsilon{^{ab}}_{c'd'} \epsilon{^{c'd'}}_{cd}
& = & \textsf{h}^{aa'}\textsf{h}^{bb'}\epsilon_{a'b'c'd'}
\textsf{h}^{c'c''}\textsf{h}^{d'd''}\epsilon_{c''d''cd}
= \left( \textsf{h}^{aa'} \textsf{h}^{bb'}\textsf{h}^{c''c'}\textsf{h}^{d''d'}\epsilon_{a'b'c'd'} \right)
\epsilon_{c''d''cd} \\
& = & \displaystyle
\sum_{1\leq\, c'',\ d''\leq 4}(\hbox{det}\,\overline{\textsf{h}})\epsilon_{abc''d''}\  \epsilon_{c''d''cd}
= 2(\hbox{det}\,\overline{\textsf{h}}) \delta^{ab}_{cd}
\end{array}
\]
we deduce
\[
 \mathring{\Theta}{^c}_{ab}
 = \frac{1}{2\gamma^2}(\hbox{det}\,\overline{\textsf{h}})\delta^{a'b'}_{ab}
 \mathring{\Theta}{^c}_{a'b'}
 = \frac{\hbox{det}\,\overline{\textsf{h}}}{\gamma^2}\mathring{\Theta}{^c}_{ab}
\]
Thus if we normalize
$\overline{\textsf{h}}$ such that $\hbox{det}\overline{\textsf{h}} = -1$ (which is always possible if
$\overline{\textsf{h}}$ is a Minkowski metric), we deduce the result. \hfill $\square$\\

\noindent
Recall that the fact that $\Theta^\gos=0$ implies through the Bianchi
identity that $\Omega^\gol\wedge \theta^\gos = \dR\Theta^\gos + [\omega^\gol\wedge \Theta^\gos] =0$, which
reads $\mathbf{R}_{abcd} + \mathbf{R}_{acdb} + \mathbf{R}_{adbc} = 0$. This implies in particular
$\mathbf{R}_{abcd} = \mathbf{R}_{cdab}$. 

Let us look at the second equation in System (\ref{brutbeautiful}). By using the computation of $\Omega{^\celg}_{\cels\gos}\kappa{_\celg}^{\cels\gos}$ in the previous paragraph we obtain
\[
 \Omega{^\celg}_{\cels a}\kappa{_\celg}^{\cels b} = \frac{1}{2}\Omega{^{[cd]}}_{ae}\left(\delta_{cd}^{be} - \frac{1}{\gamma}\epsilon{^{be}}_{cd}\right)
 = \mathbf{R}{_a}^b -\frac{1}{2\gamma}\epsilon{^{be}}_{cd}\mathbf{R}{^{cd}}_{ae}
\]
However
\[
 \frac{1}{2\gamma}\epsilon{^{be}}_{cd}\mathbf{R}{^{cd}}_{ae}
 = \frac{1}{2\gamma}\epsilon^{acef}\mathbf{R}_{efbc}
= \frac{1}{2\gamma}\epsilon^{acef}\mathbf{R}_{bcef}
= \frac{1}{6\gamma}\epsilon^{acef}\left(\mathbf{R}_{bcef} + \mathbf{R}_{befc} +\mathbf{R}_{bfce}\right) = 0
\]
Hence $\Omega{^\celg}_{\cels a}\kappa{_\celg}^{\cels b} = \mathbf{R}{_a}^b$, which implies $\Omega{^\celg}_{\cels\gos}\kappa{_\celg}^{\cels\gos} = \mathbf{R}$.
Similarly, 
\[
 \Lambda = -\frac{1}{2}\mathbf{c}{^\celg}_{\cels\cels}\kappa{_\celg}^{\cels\cels}
 = - \frac{1}{4}\left(-k\delta_{ab}^{cd}\right)\left(\delta_{cd}^{ab} - \frac{1}{\gamma}\epsilon{^{ab}}_{cd}\right) = \frac{n(n-1)}{2}k
 - \frac{k}{4\gamma}\epsilon{^{ab}}_{ab}
 = 6k
\]
Hence the equation $\Omega{^\celg}_{\cels\gos}\kappa{_\celg}^{\cels\gos} - \frac{1}{2}\left(\Omega{^\celg}_{\cels\gos}\kappa{_\celg}^{\cels\gos}\right)\delta{_\gos}^\gos + \Lambda \delta{_\gos}^\gos = 0$ gives us again the Einstein equation with a cosmological constant $\mathbf{E}{_\gos}^\gos + 6k \delta{_\gos}^\gos = 0$.

\bibliographystyle{plain}
\bibliography{bibneo}

\begin{thebibliography}{10}

\bibitem{appelquist}
T.~Appelquist, A.~Chodos, and P.G.O. Freund, editors.
\newblock {\em Modern Kaluza-Klein theories}.
\newblock Addison-Wesley Pub. Co., 1987.

\bibitem{catren}
G.~Catren.
\newblock Geometric foundations of cartan gauge gravity.
\newblock {\em International Journal of Geometric Methods in Modern Physics},
  12(4), 2015.

\bibitem{ddfr1980}
A.~D'Adda, R.~D'Auria, P.~Fr{\'e}, and T.~Regge.
\newblock Geometrical formulation of supergravity theories on orthosymplectic
  supergroup manifolds.
\newblock {\em La Rivista Del Nuovo Cimento}, 3(6):1--80, 1980.

\bibitem{df1982}
R.~D'Auria and P.~Fr{\'e}.
\newblock Geometric supergravity in d=11 and its hidden supergroup.
\newblock {\em Nuclear Physics B}, 201(1):101--140, 1982.

\bibitem{df1980}
R.~D'Auria, P.~Fr{\'e}, and T.~Regge.
\newblock Graded-lie-algebra cohomology and supergravity.
\newblock {\em La Rivista Del Nuovo Cimento}, 3(12):1--37, 1980.

\bibitem{pierard1}
J.~Pierard de~Maujouy.
\newblock Dirac spinors on generalised frame bundles: a frame bundle
  formulation for einstein-cartan-dirac theory.
\newblock \href{https://arxiv.org/abs/2201.01108}{arxiv:2201.01108}.

\bibitem{pierard2}
J.~Pierard de~Maujouy.
\newblock In preparation.

\bibitem{Ehresmann50}
C.~Ehresmann.
\newblock Les connexions infinit{\'e}simales dans un espace fibré
  diff{\'e}rentiable.
\newblock pages 153--168, 1952.
\newblock S{\'e}minaire Nicolas Bourbaki,
  \href{http://www.numdam.org/item/SB_1948-1951__1__153_0}{www.numdam.org}.

\bibitem{gielen-wise}
Steffen Gielen and Derek~K. Wise.
\newblock Lifting general relativity to observer space.
\newblock 54(5):052501, May 2013.
\newblock \href{https://arxiv.org/abs/1210.0019}{arxiv:1210.0019}.

\bibitem{helein14}
Fr{\'{e}}d{\'{e}}ric H{\'{e}}lein.
\newblock Multisymplectic formulation of {Y}ang-{M}ills equations and
  {E}hresmann connections.
\newblock {\em Advances in Theoretical and Mathematical Physics},
  19(4):805--835, 2015.
\newblock \href{https://arxiv.org/abs/1406.3641}{arXiv:1406.3641}.

\bibitem{helein2020}
Fr{\'{e}}d{\'{e}}ric H{\'{e}}lein.
\newblock A variational principle for kaluza{\textendash}klein types theories.
\newblock {\em Advances in Theoretical and Mathematical Physics},
  24(2):305--326, 2020.
\newblock \href{https://arxiv.org/abs/1809.03375}{arXiv:1809.03375}.

\bibitem{helein22}
Fr{\'{e}}d{\'{e}}ric H{\'{e}}lein.
\newblock Dynamical mechanisms for kaluza{\textendash}klein theories.
\newblock {\em Letters in Mathematical Physics}, 112(5), sep 2022.
\newblock \href{https://arxiv.org/abs/2201.01981}{arXiv:2201.01981}.

\bibitem{heleinvey15}
Fr{\'{e}}d{\'{e}}ric H{\'{e}}lein and Dimitri Vey.
\newblock Curved space-times by crystallization of liquid fiber bundles.
\newblock {\em Foundations of Physics}, 47(1):1--41, sep 2016.
\newblock \href{https://arxiv.org/abs/1508.07765}{arXiv:1508.07765}.

\bibitem{jordan}
P.~Jordan.
\newblock Erweiterung der projektiven {R}elativit\"{a}tstheorie.
\newblock {\em Annalen der Physik}, 436(4-5):219--228, 1947.

\bibitem{kaluza}
Th. Kaluza.
\newblock On the unification problem in physics.
\newblock {\em Sitzungsberichte Pruss. Acad. Sci.}, page 966, 1921.
\newblock Reprinted in English in \cite{appelquist} and in
  \href{https://arxiv.org/abs/1803.08616}{arXiv:1803.08616}.

\bibitem{kerner}
R.~Kerner.
\newblock Generalization of the {K}aluza-{K}lein theory for an arbitrary
  non-{A}belian gauge group.
\newblock {\em Ann. Inst. H. Poincar{\'e}}, 9(2):143, 1968.
\newblock in
  \href{http://www.numdam.org/item?id=AIHPA_1968__9_2_143_0}{www.numdam.org}.

\bibitem{klein}
Oskar Klein.
\newblock Quantentheorie und f{\"u}nfdimensionale {R}elativit{\"a}tstheorie.
\newblock {\em Zeitschrift f{\"u}r Physik}, 37(12):895--906, dec 1926.

\bibitem{lurcat}
Fran\ifmmode \mbox{\c{c}}\else~\c{c}\fi{}ois Lur\ifmmode~\mbox{\c{c}}\else
  \c{c}\fi{}at.
\newblock Quantum field theory and the dynamical role of spin.
\newblock {\em Physics Physique Fizika}, 1:95--106, Sep 1964.

\bibitem{macdowellmansouri}
S.~W. MacDowell and F.~Mansouri.
\newblock Unified geometric theory of gravity and supergravity.
\newblock {\em Physical Review Letters}, 38(23):1376--1376, June 1977.

\bibitem{neemanregge}
Y.~Ne'eman and T.~Regge.
\newblock Gauge theory of gravity and supergravity on a group manifold.
\newblock {\em La Rivista Del Nuovo Cimento}, 1(5):1--43, 1978.

\bibitem{sharpe}
R.W. Sharpe.
\newblock {\em Differential geometry: Cartan's generalization of Klein's
  Erlangen program}, volume 166 of {\em Graduate Texts in Mathematics}.
\newblock Springer-Verlag, 1997.

\bibitem{thiry}
Y.~Thiry.
\newblock Les {\'e}quations de la th{\'e}orie unitaire de {K}aluza.
\newblock {\em Comptes Rendus Acad. Sci. Paris}, 226(216), 1948.

\bibitem{toller}
M.~Toller.
\newblock Classical field theory in the space of reference frames.
\newblock {\em Il Nuovo Cimento B Series 11}, 44(1):67--98, mar 1978.

\bibitem{Wise2009}
Derek~K. Wise.
\newblock Symmetric space cartan connections and gravity in three and four
  dimensions.
\newblock {\em Symmetry, Integrability and Geometry: Methods and Applications},
  August 2009.
\newblock \href{https://arxiv.org/abs/0904.1738}{arxiv.org/abs/0904.1738}.

\end{thebibliography}

\end{document}